\documentclass[preprint,12pt]{elsarticle}

\usepackage{amssymb}
\usepackage{amsmath}

\usepackage{subfigure}
\usepackage{booktabs}
\usepackage{float}
\usepackage{graphicx}
\usepackage{algorithm}
\usepackage{algorithmic}

\journal{Journal of Computational Physics}

\begin{document}

\begin{frontmatter}



\title{A Neural Network Framework for High-Dimensional Dynamic Unbalanced Optimal Transport}


\author[label1]{Wei Wan}
\ead{weiwan@ncepu.edu.cn}
\author[label2]{Jiangong Pan}
\ead{mathpjg@sina.com}
\author[label2]{Yuejin Zhang}
\ead{zhangyj19@mails.tsinghua.edu.cn}
\author[label3,label4]{Chenglong Bao}
\ead{clbao@tsinghua.edu.cn}
\author[label3,label4]{Zuoqiang Shi \corref{cor1}}
\ead{zqshi@tsinghua.edu.cn}
\cortext[cor1]{Corresponding author}

\affiliation[label1]{organization={School of mathematics and physics, North China Electric Power University},
            city={Beijing},
            postcode={102206}, 
            country={China}}           
\affiliation[label2]{organization={Department of Mathematical Science, Tsinghua University},
            city={Beijing},
            postcode={100084}, 
            country={China}}

\affiliation[label3]{organization={Yau Mathematical Sciences Center, Tsinghua University},
            city={Beijing},
            postcode={100084}, 
            country={China}}

\affiliation[label4]{organization={Yanqi Lake Beijing Institute of Mathematical Sciences and Applications},
            city={Beijing},
            postcode={101408}, 
            country={China}}

\begin{abstract}
In this paper, we introduce a neural network-based method to address the high-dimensional dynamic unbalanced optimal transport (UOT) problem.
Dynamic UOT focuses on the optimal transportation between two densities with unequal total mass, however, it introduces additional complexities compared to the traditional dynamic optimal transport (OT) problem.
To efficiently solve the dynamic UOT problem in high-dimensional space, we first relax the original problem by using the generalized Kullback-Leibler (GKL) divergence to constrain the terminal density. Next, we adopt the Lagrangian discretization to address the unbalanced continuity equation and apply the Monte Carlo method to approximate the high-dimensional spatial integrals. Moreover, a carefully designed neural network is introduced for modeling the velocity field and source function.  
Numerous experiments demonstrate that the proposed framework performs excellently in high-dimensional cases. Additionally, this method can be easily extended to more general applications, such as crowd motion problem.
\end{abstract}




\end{frontmatter}

\section{Introduction}

Optimal transport (OT) \cite{santambrogio2015optimal,villani2003topics} is a significant research field, widely applied to various practical problems in machine learning \cite{torres2021survey}, including domain adaptation \cite{seguy2018large}, transfer learning \cite{7586038,redko2019optimal}, generative models \cite{2017Wasserstein,salimans2018improving,lei2019geometric}, image processing \cite{Maas2015, papadakis2015optimal}, and natural language processing \cite{fwetdd}. 
Despite its success across various domains, the numerical resolution of the optimal transport problem still presents several challenges. 
The Eulerian discretization-based method suffers from the curse of dimensionality, which means that the computational complexity grows exponentially with spatial dimension. 
Furthermore, OT typically assumes that densities have the same total mass. 
However, in real-world applications, a more common scenario arises where the total mass of each density is not equal, which corresponds to so-called unbalanced optimal transport (UOT) problems.
In this paper, we propose an extension of our previous deep learning method designed to solve the high-dimensional OT problem to address the high-dimensional UOT problem.

The OT problem originated from the ``sand moving problem" proposed by Monge \cite{monge1781memoire}, which studied the minimization of the total cost of material transportation. 
The solution to the \emph{Monge problem} is called the optimal transport mapping. 
However, directly solving the Monge problem is extremely difficult, due to its highly non-convex nature, and the existence of its solution is not guaranteed.
The \emph{Kantorovich Problem} \cite{2006On} relaxes the transport mapping into transport plans, represented by joint probability distributions. 
In the discrete case, the OT problem can be relaxed into a standard linear programming problem. 
The Sinkhorn method \cite{cuturi2013sinkhorn}  introduces an entropy regularization term to the original linear programming problem, which can be efficiently solved using the Sinkhorn algorithm \cite{1964Sinkhorn}, significantly accelerating the computation of the OT problem. 
However, as the dimensionality increases, maintaining accuracy requires an exponentially growing number of grid points, leading to higher computational costs and memory demands.
An alternative formulation of the OT problem from the fluid mechanics perspective was first introduced by Benamou and Brenier \cite{benamou2000computational}. Its computation requires solving a PDE-constrained optimization problem.
Unlike the previous two static formulations, the \emph{dynamic OT problem} can dynamically describe changes in density by involving a time variable. 

The dynamic OT problem can be formulated as follows: Given two probability densities $\rho_{0}$ and $\rho_{1}$ on $\Omega \subset \mathbb{R}^{d}$, we want to find a time-dependent density $\rho: \Omega \times [0,1] \rightarrow \mathbb{R} $ and a time-dependent velocity field $\boldsymbol{v}:\Omega \times [0,1]  \rightarrow \mathbb{R}^{d}$ that transport
the mass of $\rho_{0}$ to the mass of $\rho_{1}$ at minimal transport cost $\mathcal{E}(\rho,\boldsymbol{v})$, i.e., 
\begin{align}\label{OT}
\min_{(\rho, \boldsymbol{v}) \in \mathcal{C}(\rho_0,\rho_1)} \int_0^1 \int_{\Omega} ||\boldsymbol{v}(\boldsymbol{x},t)||^2\rho(\boldsymbol{x},t)\mathrm{d}\boldsymbol{x} \mathrm{d}t :=\mathcal{E}(\rho, \boldsymbol{v}),
\end{align}
where $(\rho, \boldsymbol{v}) \in \mathcal{C}\left(\rho_{0}, \rho_{1}\right)$ satisfies the dynamical constraint, i.e. the continuity equation, the initial and terminal densities are given by $\rho_0$ and $\rho_1$, respectively:
\begin{equation}
\begin{aligned}
\label{massconserv}
    \mathcal{C}(\rho_0,\rho_1):=\{(\rho, \boldsymbol{v}): \partial_t \rho(\boldsymbol{\boldsymbol{x}},t)+\nabla\cdot (\rho(\boldsymbol{\boldsymbol{x}},t)\boldsymbol{v}(\boldsymbol{x},t))= 0,\\ \rho(\cdot,0)=\rho_0,\rho(\cdot,1)=\rho_1\},
\end{aligned}
\end{equation}
where $\mathcal{E}(\rho,\boldsymbol{v})$ denotes the generalized kinetic energy and the optimal $\mathcal{E}(\rho,\boldsymbol{v})$ gives the square of the $L^2$ Wasserstein distance.
This fluid mechanics formulation provides a natural time interpolation, which is highly beneficial for a wide range of applications, including continuous normalizing flow \cite{2020TrajectoryNet} and mean field games \cite{2007Mean}, therefore making the numerical solution for dynamic OT more attractive.

The UOT problem describes the optimal transportation between two densities with different total masses. 
Various versions of UOT have been developed based on different formulations of OT \cite{Sejourne2023}.
This work mainly focuses on solving the dynamic formulation of the UOT problem, a generalization of the dynamic OT problem proposed by Benamou and Brenier, which introduces a source function to control the creation or destruction of mass during transportation.
Several dynamic UOT models have been proposed to account for mass change by introducing different source functions $f$ into the cost functional and the continuity equation.
Given two densities $\rho_{0}$ and $\rho_{1}$, the general dynamic UOT problem can be formulated as:
\begin{align}\label{dUOT}
\min_{(\rho, \boldsymbol{v},f) \in \mathcal{C}(\rho_0,\rho_1)} \int_0^1 \int_{\Omega} ||\boldsymbol{v}(\boldsymbol{x},t)||^2\rho(\boldsymbol{x},t)\mathrm{d}\boldsymbol{x} \mathrm{d}t + \frac{1}{\alpha}\mathcal{R}(\rho,f).
\end{align}
Here, $\alpha$ is a positive constant, and $(\rho, \boldsymbol{v}, f) \in \mathcal{C}\left(\rho_{0}, \rho_{1}\right)$ satisfies the unbalanced continuity equation with initial and terminal densities $\rho_0$ and $\rho_1$, respectively:
\begin{equation}\label{unbalancedPDE}
\begin{aligned}
    \mathcal{C}(\rho_0,\rho_1):=\{(\rho, \boldsymbol{v},f): \partial_t \rho(\boldsymbol{x},t)+\nabla\cdot (\rho(\boldsymbol{x},t)\boldsymbol{v}(\boldsymbol{x},t))= \mathcal{S}(\rho,f),\\ \rho(\cdot,0)=\rho_0,\rho(\cdot,1)=\rho_1\}.
\end{aligned}
\end{equation}
When $f \equiv 0$, the dynamic UOT problem simplifies to the dynamic OT problem.
Several authors have considered density-independent source terms \cite{Piccoli2014, Piccoli2016, GANGBO2019, LEE2019}.
In \cite{GANGBO2019}, Gangbo et al. 
consider $\mathcal{R}(f)$ to be density-independent and the source function $f(t)$ to be space-independent:
\begin{align*}
\mathcal{R}(f)=\int_0^1 |f(t)|^p \mathrm{d}t\cdot |\Omega|,~~
\mathcal{S}(f)=f(t).
\end{align*}
They derive the corresponding unnormalized Monge problem, Monge-Amp$\acute{e}$re equation,
and Kantorovich formulation, and introduce a primal-dual algorithm for the numerical solution of this
dynamic UOT problem. 
In \cite{LEE2019}, Lee et al. extend this UOT model in \cite{GANGBO2019} by using a spatially dependent source function $f(\boldsymbol{x},t)$:
\begin{align*}
\mathcal{R}(f)=\int_0^1 \int_{\Omega}|f(\boldsymbol{x},t)|^p \mathrm{d}x\mathrm{d}t,~~
\mathcal{S}(f)=f(\boldsymbol{x},t).
\end{align*}
They design fast algorithms for $L^1$ and $L^2$ generalized unnormalized optimal transport using a primal-dual algorithm and a Nesterov accelerated gradient descent method, respectively.
In this work, we mainly study the dynamic UOT problem as introduced by Chizat et al. \cite{CHIZAT2016, CHIZAT2018}.
In this model, $\mathcal{R}(\rho,f)$ is density-dependent and  $f(\boldsymbol{x},t)$ is space-dependent:
\begin{align*}
\mathcal{R}(\rho,f)=\int_0^1 \int_{\Omega}f(\boldsymbol{x},t)^2\rho(\boldsymbol{x},t) \mathrm{d}x\mathrm{d}t,~~
\mathcal{S}(\rho,f)=f(\boldsymbol{x},t)\rho(\boldsymbol{x},t).
\end{align*}
This UOT model is defined on the Wasserstein-Fisher-Rao metric, also known as Hellinger-Kantorovich, initially and independently proposed by \cite{CHIZAT2016, Kondratyev2016}, providing an interpolation between the Wasserstein distance and the Fisher-Rao metric for measuring unequal masses.
In \cite{CHIZAT2016}, Chizat et al. theoretically establish the existence of solutions to this variational problem.
They employ staggered grid discretization techniques and solve this problem using the first-order proximal splitting method, demonstrating its application to image interpolation. In \cite{CHIZAT2018}, the authors further generalize this model and propose a family of dynamic UOT problems. In addition, they formulate a corresponding static Kantorovich formulation and establish the equivalence between dynamic and static formulations.
Building on this UOT model, Chen et al. \cite{Chen2024UOT} propose the addition of a diffusion term, $\sigma \Delta \rho$, into the source term $\mathcal{S}(\rho,f)$. This incorporates both advection and diffusion into the transport process. In their model, $f(\boldsymbol{x},t)$ is a relative source function with a specified indicator function. They have developed a detailed numerical algorithm to implement this model and applied it to modeling the fluid flow in the brain.
In \cite{LOMBARDI2015}, Lombardi et al. propose a genetic source term and consider three distinct types of source terms specifically designed for interpolating growing tumor images. 
Additionally, they further prove the existence and uniqueness of solutions for two of these terms, providing numerical evidence of their convergence and ability to handle mass growth in various test cases.
The Eulerian discretization-based methods mentioned above are effective in low-dimensional cases but cannot be readily extended to high dimensions due to the exponential growth in the number of samples required to achieve a given tolerance as the dimensionality increases. 
In recent years, deep learning has demonstrated remarkable success on a variety of computational problems.
In \cite{Jing2022UOT}, Jing et al. develop a deep learning framework to compute the geodesics under the spherical Wasserstein-Fisher-Rao metric, where $f$ is restricted to the family of zero mean functions. The learned geodesics can be adopted to generate weighted samples from some target distribution.
In \cite{Wan2023}, we propose a novel deep learning approach for addressing the high-dimensional dynamic OT problem. 
We adopt the Lagrangian discretization to solve the continuity equation and use Monte Carlo sampling to approximate the high-dimensional integrals rather than Eulerian discretization-based methods. 
Moreover, we carefully design neural networks to parametrize the velocity field, achieving more accurate results in high-dimensional cases with excellent scalability.

\textbf{Contributions:} In this paper, our first contribution is the introduction of a novel neural network approach to address the high-dimensional dynamic UOT problem used with the Wasserstein-Fisher-Rao metric.
In contrast to our previous method \cite{Wan2023}, which was designed for high-dimensional dynamic OT problems, the UOT problem introduces additional complexities by incorporating a source term into both the objective function and the continuity equation. 
There are three main differences compared to the OT problem: 
(i) To alleviate the dynamic formulation, we employ the generalized Kullback-Leibler (GKL) divergence to constrain the terminal density, rather than the standard Kullback-Leibler (KL) divergence, which is ineffective for the UOT problem.
(ii) To overcome the curse of dimensionality, we utilize Lagrangian discretization to solve the unbalanced PDE constraint and employ Monte Carlo sampling to approximate the high-dimensional integral. However, this introduces the need to solve an additional ODE compared to the OT problem. 
(iii) In addition to parameterizing the velocity field, the source function also requires parameterization which increase the number of parameters.
Our second contribution is the good scalability of our proposed framework in high dimensions with constant number of samples. Furthermore, it can be easily extended to more general applications, such as Gaussian mixtures and crowd motion problems with different total masses. 

\textbf{Organization:} The structure of this paper is as follows:  In Section \ref{se:pre}, we introduce the dynamic UOT problem and the GKL divergence.
Section \ref{lag} presents the Lagrangian discretization approach.
The neural network architecture is described in Section \ref{NNArc}. Section \ref{Experimentresult} provides the experimental results along with implementation details. Finally, we conclude this work in Section \ref{Conclusion}.


\section{Review of dynamic UOT and generalized Killback-Leibler divergence}\label{se:pre}

In this work, we mainly study the dynamic UOT problem as introduced by Chizat et al. \cite{CHIZAT2016, CHIZAT2018}. Given two densities $\rho_{0}$ and $\rho_{1}$, the fundamental problem can be formulated as minimizing the Wasserstein-Fisher-Rao metric:
\begin{equation}
\begin{aligned}\label{dOT}
\min_{(\rho, \boldsymbol{v},f) \in \mathcal{C}(\rho_0,\rho_1)} & \int_0^1 \int_{\Omega} ||\boldsymbol{v}(\boldsymbol{x},t)||^2\rho(\boldsymbol{x},t)\mathrm{d}\boldsymbol{x} \mathrm{d}t + \frac{1}{\alpha}\int_0^1\int_{\Omega}f(\boldsymbol{x},t)^2\rho(\boldsymbol{x},t)\mathrm{d}\boldsymbol{x}\mathrm{d}t\\
&:=\mathcal{E}(\rho, \boldsymbol{v}) + \frac{1}{\alpha} \mathcal{R}(\rho,f),
\end{aligned}
\end{equation}
where $\rho: \Omega \times [0,1] \rightarrow \mathbb{R}$ is a time-dependent density, $\boldsymbol{v}:\Omega \times [0,1]  \rightarrow \mathbb{R}^{d}$ is a velocity field that describes the movement of mass and $f: \Omega \times [0,1] \rightarrow \mathbb{R}$ represents a spatially dependent source function that captures local growth and destruction of mass. In addition, $\alpha$ is a positive constant, and $(\rho, \boldsymbol{v}, f) \in \mathcal{C}\left(\rho_{0}, \rho_{1}\right)$ satisfies the unbalanced continuity equation, with the initial and terminal densities denoted as $\rho_0$ and $\rho_1$ respectively:
\begin{equation}
\begin{aligned}\label{UOTcons}
    \mathcal{C}(\rho_0,\rho_1):=\{(\rho, \boldsymbol{v},f): \partial_t \rho(\boldsymbol{x},t)+\nabla\cdot (\rho(\boldsymbol{x},t)\boldsymbol{v}(\boldsymbol{x},t))= f(\boldsymbol{x},t)\rho(\boldsymbol{x},t), \\\rho(\cdot,0)=\rho_0,\rho(\cdot,1)=\rho_1\}.
\end{aligned}
\end{equation}
In particular, when $f(\boldsymbol{x},t) \equiv 0$, the dynamic UOT problem reduces to the dynamic OT problem.

To efficiently solve the dynamic UOT problem, we first relax the constraint on the terminal density $\rho(\cdot,1)=\rho_1$ by transforming it into an implicit condition. 
Unlike the standard OT problem, where the total mass is conserved and $\rho(\cdot,1)$ represents a normalized probability distribution, the UOT problem allows for the creation or destruction of mass. This means that during the transportation process, the total mass can change, leading to a terminal density $\rho(\cdot,1)$ that may not by a normalized probability distribution.
To address this issue, instead of using the original KL divergence to enforce the terminal density constraint in the OT problem, we apply the generalized Kullback-Leibler (GKL) divergence for the dynamic UOT problem.

Let $p(\boldsymbol{x})$ and $q(\boldsymbol{x})$ be unnormalized probability distributions. The GKL divergence \cite{Miller2023GKL} between $p(\boldsymbol{x})$ and $q(\boldsymbol{x})$ is defined by:
\begin{align*}
&\mathbf{GKL}(p(\boldsymbol{x})||q(\boldsymbol{x})):=\int_{\Omega}\phi \Bigg (\frac{q(\boldsymbol{x})}{p(\boldsymbol{x})} \Bigg )p(\boldsymbol{x})\mathrm{d}\boldsymbol{x},\\
&\phi(r):=-\log r +r -1.
\end{align*}

The key properties of the GKL divergence are as follows: $(a)$ Gibbs' inequality holds: $\mathbf{GKL}(p(\boldsymbol{x})||q(\boldsymbol{x})) \geq 0$, and $(b)~\mathbf{GKL}(p(\boldsymbol{x})||q(\boldsymbol{x})) = 0$ if and only if $p(\boldsymbol{x})=q(\boldsymbol{x})$ holds almost everywhere for all $\boldsymbol{x}$. The GKL divergence is general because for normalized probability distributions $p(\boldsymbol{x})$ and $q(\boldsymbol{x})$, it reduces to the original KL divergence:
\begin{align*}
\mathbf{GKL}(p(\boldsymbol{x})||q(\boldsymbol{x}))=\int_{\Omega} p(\boldsymbol{x}) \log \frac{p(\boldsymbol{x})}{q(\boldsymbol{x})} \mathrm{d}\boldsymbol{x} :=\mathbf{KL}(p(\boldsymbol{x})||q(\boldsymbol{x})).
\end{align*}

Therefore, by relaxing the terminal density constraint $\rho(\cdot,1)=\rho_1$ as an implicit condition through the introduction of the GKL divergence term:
\begin{align*}
\mathcal{P}(\rho):&=\mathbf{GKL}(\rho(\boldsymbol{x},1)||\rho_1(\boldsymbol{x}))\\
&=\int_{\Omega}\rho(\boldsymbol{x},1)\log \frac{\rho(\boldsymbol{x},1)}{\rho_1(\boldsymbol{x})}
-\rho(\boldsymbol{x},1)+\rho_1(\boldsymbol{x})\mathrm{d}\boldsymbol{x},
\end{align*}
we obtain the following optimal control formulation:
\begin{align}\label{objective}
\min_{(\rho,\boldsymbol{v},f) \in \mathcal{C}(\rho_0)}\mathcal{F}(\rho,\boldsymbol{v},f) = \mathcal{E}(\rho, \boldsymbol{v})  + \frac{1}{\alpha}\mathcal{R}(\rho,f) + \lambda\mathcal{P}(\rho),
\end{align}
where $\mathcal{C}(\rho_0):=\{(\rho, \boldsymbol{v},f): \partial_t \rho(\boldsymbol{x},t)+\nabla\cdot (\rho(\boldsymbol{x},t) \boldsymbol{v}(\boldsymbol{x},t))=f(\boldsymbol{x},t)\rho(\boldsymbol{x},t), \rho(\cdot,0)=\rho_0\}$. 
Additionally, $\lambda$ serves as the penalty parameter.

\section{Discretization  of Optimization Problem}\label{lag}

\subsection{Lagrangian Method}

To solve UOT in high dimensional space, we solve the transport equation 
$$\partial_t \rho(\boldsymbol{x},t)+\nabla\cdot (\rho(\boldsymbol{x},t)\boldsymbol{v}(\boldsymbol{x},t))= f(\boldsymbol{x},t)\rho(\boldsymbol{x},t)$$
along the characteristics. The characteristics can be obtained by solving an ODE system:
\begin{align}
&\frac{\mathrm{d}}{\mathrm{d}t}\mathbf{z}(\boldsymbol{x},t)=\boldsymbol{v}(\mathbf{z}(\boldsymbol{x},t),t),\mathbf{z}(\boldsymbol{x},0)=\boldsymbol{x},  \label{char}
\end{align}
Then the original PDE can be transformed to an ODE system along the characteristics:
\begin{equation}\label{rhos}
\begin{aligned}
\frac{\mathrm{d}}{\mathrm{d}t}\rho(\mathbf{z}(\boldsymbol{x},t),t)=&-\rho(\mathbf{z}(\boldsymbol{x},t),t)\nabla\cdot \boldsymbol{v}(\mathbf{z}(\boldsymbol{x},t),t) + f(\mathbf{z}(\boldsymbol{x},t),t)\rho(\mathbf{z}(\boldsymbol{x},t),t),\\
\rho(\mathbf{z}(\boldsymbol{x},0),0)=&\rho_0(\boldsymbol{x}). 
\end{aligned}
\end{equation}
Diving $\rho$ on both sides, above ODE can be simplified further. 
\begin{align}
\label{lnrhos_f}
\frac{\mathrm{d}}{\mathrm{d}t}\ln\rho(\mathbf{z}(\boldsymbol{x},t),t)=&-\nabla\cdot \boldsymbol{v}(\mathbf{z}(\boldsymbol{x},t),t)+f(\mathbf{z}(\boldsymbol{x},t),t),\\
\ln\rho(\mathbf{z}(\boldsymbol{x},0),0)=&\ln\rho_0(\boldsymbol{x}). \nonumber
\end{align}
Note that \eqref{rhos} and \eqref{lnrhos_f} are equivalent, however \eqref{lnrhos_f} offers some advantages in numerical computation. In high-dimensional cases, since $\rho_0$ is a positive number very close to zero, solving equation \eqref{rhos} can lead to numerical overflow issues. Therefore, in practice, we solve equation \eqref{lnrhos_f} instead.

The other difficult is to evaluate the objective function \eqref{dOT} in high dimensional space. To address this, we use the Monte Carlo method to calculate the integrals in high dimensional space. 
For the standard OT problem, the density along the characteristics satisfies the equation:
\begin{align}
\rho(\mathbf{z}(\boldsymbol{x},t),t)\det(\nabla \mathbf{z}(\boldsymbol{x},t))= \rho_0(\boldsymbol{x}),
\end{align}
for all $t \in [0,1]$. However, this relationship does not hold for the dynamic UOT problem.
As a result, implementing the Monte Carlo method requires solving an additional ODE system:
\begin{align}
     &\frac{\mathrm{d}}{\mathrm{d}t}\ln\mu(\mathbf{z}(\boldsymbol{x},t),t)=-\nabla \cdot \boldsymbol{v}(\mathbf{z}(\boldsymbol{x},t),t),~\ln\mu(\mathbf{z}(\boldsymbol{x},0),0)=\ln\mu_0(\boldsymbol{x}).\label{lnmu}
\end{align}
Here, $\mathbf{z}(\boldsymbol{x},t)$ represents the characteristics solved in \eqref{char}. The distribution $\mu_0$ can be another distribution, such as uniform distribution, or it can be chosen as $\mu_0=\rho_0$.

From \eqref{char} and \eqref{lnmu}, we can solve that
\begin{align}\label{var_subst}
\mu(\mathbf{z}(\boldsymbol{x},t),t)\det(\nabla \mathbf{z}(\boldsymbol{x},t))= \mu_0(\boldsymbol{x}),
\end{align}
for all $t \in [0,1]$. 

\subsection{Monte Carlo Method}

Now, we can give the Monte Carlo scheme to compute the objective function.
By substituting equation \eqref{var_subst} into the objective function \eqref{objective} and let $\boldsymbol{x}_1,\boldsymbol{x}_2,...,\boldsymbol{x}_r \in \mathbb{R}^{d}$ be random samples independently drawn from $\mu_0$, we can then have
\begin{align*}
 \mathcal{E}(\rho,\boldsymbol{v})
 &= \int_0^1 \int_{\Omega} ||\boldsymbol{v}(\boldsymbol{x},t)||^2\rho(\boldsymbol{x},t)\mathrm{d}\boldsymbol{x} \mathrm{d}t\\
 &=\int_0^1 \int_{\Omega}||\boldsymbol{v}(\boldsymbol{x},t)||^2\frac{\rho(\boldsymbol{x},t)}{\mu(\boldsymbol{x},t)}\mu(\boldsymbol{x},t)\mathrm{d}\boldsymbol{x} \mathrm{d}t\\
 &=\int_0^1 \int_{\Omega} ||\boldsymbol{v}(\mathbf{z}(\boldsymbol{x},t),t)||^2 \frac{\rho(\mathbf{z}(\boldsymbol{x},t),t)}{\mu(\mathbf{z}(\boldsymbol{x},t),t)} \mu(\mathbf{z}(\boldsymbol{x},t),t)\det(\nabla \mathbf{z}(\boldsymbol{x},t))\mathrm{d}\boldsymbol{x} \mathrm{d}t\\
  &=\int_0^1 \int_{\Omega} ||\boldsymbol{v}(\mathbf{z}(\boldsymbol{x},t),t)||^2 \frac{\rho(\mathbf{z}(\boldsymbol{x},t),t)}{\mu(\mathbf{z}(\boldsymbol{x},t),t)} \mu_0(\boldsymbol{x}) \mathrm{d}\boldsymbol{x} \mathrm{d}t\\
 & = \int_0^1 \mathbb{E}_{\mathbf{\boldsymbol{x}} \sim \mu_0}\big[ ||\boldsymbol{v}(\mathbf{z}(\mathbf{\boldsymbol{x}},t),t)||^2 \frac{\rho(\mathbf{z}(\mathbf{\boldsymbol{x}},t),t)}{\mu(\mathbf{z}(\mathbf{\boldsymbol{x}},t),t)} \big] \mathrm{d}t \\
 &\approx \int_0^1 \sum_{i=1}^r \frac{1}{r} ||\boldsymbol{v}(\mathbf{z}(\boldsymbol{x}_i,t),t)||^2 \frac{\rho(\mathbf{z}(\boldsymbol{x}_i,t),t)}{\mu(\mathbf{z}(\boldsymbol{x}_i,t),t)}\mathrm{d}t.
\end{align*}

For the other two spatial integrals in the objective function, we can similarly apply Monte Carlo discretization.

\begin{align*}
 \mathcal{R}(\rho,f)
 &=\int_0^1 \int_{\Omega} f(\boldsymbol{x},t)^2\rho(\boldsymbol{x},t)\mathrm{d}\boldsymbol{x} \mathrm{d}t\\
 &=\int_0^1 \int_{\Omega}f(\boldsymbol{x},t)^2\frac{\rho(\boldsymbol{x},t)}{\mu(\boldsymbol{x},t)}\mu(\boldsymbol{x},t)\mathrm{d}\boldsymbol{x} \mathrm{d}t\\
 &=\int_0^1 \int_{\Omega} f(\mathbf{z}(\boldsymbol{x},t),t)^2 \frac{\rho(\mathbf{z}(\boldsymbol{x},t),t)}{\mu(\mathbf{z}(\boldsymbol{x},t),t)} \mu(\mathbf{z}(\boldsymbol{x},t),t)\det(\nabla \mathbf{z}(\boldsymbol{x},t))\mathrm{d}\boldsymbol{x} \mathrm{d}t\\
  &=\int_0^1 \int_{\Omega} f(\mathbf{z}(\boldsymbol{x},t),t)^2 \frac{\rho(\mathbf{z}(\boldsymbol{x},t),t)}{\mu(\mathbf{z}(\boldsymbol{x},t),t)} \mu_0(\boldsymbol{x}) \mathrm{d}\boldsymbol{x} \mathrm{d}t\\
 & =\int_0^1 \mathbb{E}_{\mathbf{\boldsymbol{x}} \sim \mu_0}\big[ f(\mathbf{z}(\mathbf{\boldsymbol{x}},t),t)^2 \frac{\rho(\mathbf{z}(\mathbf{\boldsymbol{x}},t),t)}{\mu(\mathbf{z}(\mathbf{\boldsymbol{x}},t),t)} \big] \mathrm{d}t \\
 &\approx \int_0^1 \sum_{i=1}^r \frac{1}{r} f(\mathbf{z}(\boldsymbol{x}_i,t),t)^2 \frac{\rho(\mathbf{z}(\boldsymbol{x}_i,t),t)}{\mu(\mathbf{z}(\boldsymbol{x}_i,t),t)}\mathrm{d}t,
\end{align*}
and
\begin{align*}
\mathcal{P}(\rho)
& = \int_{\Omega}\rho(\boldsymbol{x},1)\log \frac{\rho(\boldsymbol{x},1)}{\rho_1(\boldsymbol{x})}-\rho(\boldsymbol{x},1)+\rho_1(\boldsymbol{x})
\mathrm{d}\boldsymbol{x}\\
& = \int_{\Omega}\big[\frac{\rho(\boldsymbol{x},1)}{\mu(\boldsymbol{x},1)}\log \frac{\rho(\boldsymbol{x},1)}{\rho_1(\boldsymbol{x})}-\frac{\rho(\boldsymbol{x},1)}{\mu(\boldsymbol{x},1)}+\frac{\rho_1(\boldsymbol{x})}{\mu(\boldsymbol{x},1)}\big]\mu(\boldsymbol{x},1)
\mathrm{d}\boldsymbol{x}\\
&= \int_{\Omega}\big[
\frac{\rho(\mathbf{z}(\boldsymbol{x},1),1)}{\mu(\mathbf{z}(\boldsymbol{x},1),1)}\log \frac{\rho(\mathbf{z}(\boldsymbol{x},1),1)}{\rho_1(\mathbf{z}(\boldsymbol{x},1))}-\frac{\rho(\mathbf{z}(\boldsymbol{x},1),1)}{\mu(\mathbf{z}(\boldsymbol{x},1),1)}+\frac{\rho_1(\mathbf{z}(\boldsymbol{x},1))}{\mu(\mathbf{z}(\boldsymbol{x},1),1)}\big]\\
&~~~~\mu(\mathbf{z}(\boldsymbol{x},1),1) \det(\nabla \mathbf{z}(\boldsymbol{x},1))
\mathrm{d}\boldsymbol{x}\\
&=\int_{\Omega}\big[
\frac{\rho(\mathbf{z}(\boldsymbol{x},1),1)}{\mu(\mathbf{z}(\boldsymbol{x},1),1)}\log \frac{\rho(\mathbf{z}(\boldsymbol{x},1),1)}{\rho_1(\mathbf{z}(\boldsymbol{x},1))}-\frac{\rho(\mathbf{z}(\boldsymbol{x},1),1)}{\mu(\mathbf{z}(\boldsymbol{x},1),1)}+\frac{\rho_1(\mathbf{z}(\boldsymbol{x},1))}{\mu(\mathbf{z}(\boldsymbol{x},1),1)}\big]\mu_0(\boldsymbol{x})
\mathrm{d}\boldsymbol{x}\\
&=\mathbb{E}_{\mathbf{\boldsymbol{x}}\sim \mu_0}\big[
\frac{\rho(\mathbf{z}(\mathbf{\boldsymbol{x}},1),1)}{\mu(\mathbf{z}(\mathbf{\boldsymbol{x}},1),1)}\log \frac{\rho(\mathbf{z}(\mathbf{\boldsymbol{x}},1),1)}{\rho_1(\mathbf{z}(\mathbf{\boldsymbol{x}},1))}-\frac{\rho(\mathbf{z}(\mathbf{\boldsymbol{x}},1),1)}{\mu(\mathbf{z}(\mathbf{\boldsymbol{x}},1),1)}+\frac{\rho_1(\mathbf{z}(\mathbf{\boldsymbol{x}},1))}{\mu(\mathbf{z}(\mathbf{\boldsymbol{x}},1),1)}\big] \\
&\approx \sum_{i=1}^r \frac{1}{r} \big[
\frac{\rho(\mathbf{z}(\boldsymbol{x}_i,1),1)}{\mu(\mathbf{z}(\boldsymbol{x}_i,1),1)}\log \frac{\rho(\mathbf{z}(\boldsymbol{x}_i,1),1)}{\rho_1(\mathbf{z}(\boldsymbol{x}_i,1))}-\frac{\rho(\mathbf{z}(\boldsymbol{x}_i,1),1)}{\mu(\mathbf{z}(\boldsymbol{x}_i,1),1)}+\frac{\rho_1(\mathbf{z}(\boldsymbol{x}_i,1))}{\mu(\mathbf{z}(\boldsymbol{x}_i,1),1)}\big].
\end{align*}

Therefore, coupled with the Lagrangian method \eqref{char}, \eqref{lnrhos_f} and \eqref{lnmu}, this results in an optimization problem regarding the velocity field $\boldsymbol{v}$ and the source function $f$, which will be approximated by neural networks. The structures of neural networks will be discussed in next section.
Consequently, we get the following semi-discretized version of the dynamic UOT problem:
\begin{align}\label{semi_obj} 
 \min_{\rho,\boldsymbol{v},f }\overline{\mathcal{F}} (\rho,\boldsymbol{v},f) :=  \overline{\mathcal{E}} (\rho,\boldsymbol{v})  + \frac{1}{\alpha}\overline{\mathcal{R}} (\rho,f)+ \lambda \overline{\mathcal{P}} (\rho) ,
\end{align} 
where $\overline{\mathcal{E}} (\rho,\boldsymbol{v}),  \overline{\mathcal{P}} (\rho)$ and $\overline{\mathcal{R}} (\rho,f)$ denote the discretized versions of $\mathcal{E}(\rho,\boldsymbol{v}), \mathcal{P}(\rho)$ and $\mathcal{R}(\rho,f)$ in the spatial direction, respectively, i.e.,
\begin{align}
\overline{\mathcal{E}} (\rho,\boldsymbol{v})&:= \int_0^1 \sum_{i=1}^r \frac{1 }{r} ||\boldsymbol{v}(\mathbf{z}(\boldsymbol{x}_i,t),t)||^2 \frac{\rho(\mathbf{z}(\boldsymbol{x}_i,t),t)}{\mu(\mathbf{z}(\boldsymbol{x}_i,t),t)}\mathrm{d}t, \label{semdisE}\\
\overline{\mathcal{R}} (\rho,f)&:=\int_0^1 \sum_{i=1}^r \frac{1 }{r} f(\mathbf{z}(\boldsymbol{x}_i,t),t)^2 \frac{\rho(\mathbf{z}(\boldsymbol{x}_i,t),t)}{\mu(\mathbf{z}(\boldsymbol{x}_i,t),t)}\mathrm{d}t, \label{semdisR}\\
\overline{\mathcal{P}} (\rho)&:= \sum_{i=1}^r \frac{1}{r} \big[
\frac{\rho(\mathbf{z}(\boldsymbol{x}_i,1),1)}{\mu(\mathbf{z}(\boldsymbol{x}_i,1),1)}\log \frac{\rho(\mathbf{z}(\boldsymbol{x}_i,1),1)}{\rho_1(\mathbf{z}(\boldsymbol{x}_i,1))}-\frac{\rho(\mathbf{z}(\boldsymbol{x}_i,1),1)}{\mu(\mathbf{z}(\boldsymbol{x}_i,1),1)}+\frac{\rho_1(\mathbf{z}(\boldsymbol{x}_i,1))}{\mu(\mathbf{z}(\boldsymbol{x}_i,1),1)}\big],\label{semdisP}
\end{align}
where $\rho(\mathbf{z}(\boldsymbol{x}_i,t),t)$ and $\mu(\mathbf{z}(\boldsymbol{x}_i,t),t)$ can be obtained from the solution of \eqref{lnrhos_f} and \eqref{lnmu} along the characteristics $\mathbf{z}(\boldsymbol{x}_i,t)$. Then, \eqref{char}, \eqref{lnrhos_f} and \eqref{lnmu} can be reformulated as the following constraints:
\begin{align}
 &\frac{\mathrm{d}}{\mathrm{d}t}\mathbf{z}(\boldsymbol{x}_i,t)=\boldsymbol{v}(\mathbf{z}(\boldsymbol{x}_i,t),t),~\mathbf{z}(\boldsymbol{x}_i,0)=\boldsymbol{x}_i, ~\boldsymbol{x}_i \sim \mu_0(\boldsymbol{x}), \label{forwd1}\\
 &\frac{\mathrm{d}}{\mathrm{d}t}\ln\rho(\mathbf{z}(\boldsymbol{x}_i,t),t)=-\nabla \cdot \boldsymbol{v}(\mathbf{z}(\boldsymbol{x}_i,t),t) +f(\mathbf{z}(\boldsymbol{x}_i,t),t),~\ln\rho(\mathbf{z}(\boldsymbol{x}_i,0),0)=\ln\rho_0(\boldsymbol{x}_i),\label{forwd2}\\
  &\frac{\mathrm{d}}{\mathrm{d}t}\ln\mu(\mathbf{z}(\boldsymbol{x}_i,t),t)=-\nabla \cdot \boldsymbol{v}(\mathbf{z}(\boldsymbol{x}_i,t),t),~\ln\mu(\mathbf{z}(\boldsymbol{x}_i,0),0)=\ln\mu_0(\boldsymbol{x}_i).\label{forwd3}
 \end{align}
In our numerical implementation, we simply set $\mu_0=\rho_0$.
To get the fully discretized model, we need to solve ODEs along time and compute integrals in time. 
The above three ODE constraints, \eqref{forwd1}-\eqref{forwd3}, are solved using the classical fourth-order explicit Runge-Kutta scheme (RK4). The time integral in \eqref{semi_obj} is computed using the standard composite Simpson's formula.

\section{Neural network architectures for the velocity field and source function}
\label{NNArc}

To address high-dimensional UOT problems, we parameterize both the velocity field $\boldsymbol{v}$ and the source function $f$ using neural networks, leveraging the ability of neural networks to approximate high dimensional functions. 

The representation is split into spatial and temporal components. For the temporal component, we use simple linear finite element basis functions, while for the spatial component, we employ a standard multi-layer perceptron (MLP) network. More specifically, 
for the velocity field $\boldsymbol{v}$, our model architecture is as follows:
 \begin{align}
 \boldsymbol{v}(\boldsymbol{x},t;\theta) = \sum_{i=0}^M\boldsymbol{v}_i(\boldsymbol{x};\theta)\phi_i(t),
 \label{eq:v-struct}
 \end{align}
 where $(\boldsymbol{x},t) \in \mathbb{R}^{d+1}$ represents the input feature, and the basis functions $\phi_i(t)$ are defined as:
\begin{align*}
\phi_i(t)=
\begin{cases}
0,&t\leq t_{i-1},\\
\frac{t-t_{i-1}}{t_i-t_{i-1}}=M(t-t_{i-1}),&t_{i-1}<t\leq t_i,\\
\frac{t_{i+1}-t}{t_{i+1}-t_{i}}=M(t_{i+1}-t),&t_{i}<t\leq t_{i+1},\\
0,&t> t_{i+1}.
\end{cases}
\end{align*}
Here, the number of the basis functions is $M+1$. Once $M$ is determined, we can set
the number of time steps $N$ to $kM$, such as $k=2$, ensuring that the forward ODEs  \eqref{char}, \eqref{lnrhos_f} and \eqref{lnmu}  can be solved stably and accurately.

Additionally, $\boldsymbol{v}_i(\boldsymbol{x};\theta)$ is the sum of some two-layer MLP networks, defined by:
$$\boldsymbol{v}_i(\boldsymbol{x};\theta)=\sum_{l=1}^L[W_{2i}^{(l)}\sigma(W_{1i}^{(l)}\boldsymbol{x}+b_{1i}^{(l)})+b_{2i}^{(l)}].
$$
Here, the network parameters are denoted by $\theta = \{W_{1i}^{(l)}, W_{2i}^{(l)}, b_{1i}^{(l)}, b_{2i}^{(l)}\}$, where $W_{1i}^{(l)} \in \mathbb{R}^{H \times d}$, $W_{2i}^{(l)} \in \mathbb{R}^{d \times H}$, $b_{1i}^{(l)} \in \mathbb{R}^{H}$, $b_{2i}^{(l)} \in \mathbb{R}^{d}$. $H$ represents the number of hidden units, $L$ denotes the width of the network, and $\sigma(\boldsymbol{x})=\tanh(\boldsymbol{x})$ serves as the activation function. 

The advantage of the proposed architecture is that it enhances representation ability and reduces oscillations in the time direction by utilizing piecewise linear function as the basis functions.

Unlike the standard OT problem, in the dynamic UOT problem, the additional source function $f(\boldsymbol{x},t;\tilde{\theta})$ is also parameterized using basis functions as follows:
\begin{align*}
f(\boldsymbol{x},t;\tilde{\theta}) = \mathbf{w}^T\boldsymbol{u}(\boldsymbol{x},t;\hat{\theta}) + b ,
\end{align*}
where $\boldsymbol{u}(\boldsymbol{x},t;\hat{\theta})$ follows the same architecture as $\boldsymbol{v}$ in \eqref{eq:v-struct}. The parameters of the network $f$ are given by $\tilde{\theta}= \{\hat{\theta}, \mathbf{w}, b\}$, with $\mathbf{w} \in \mathbb{R}^d$, $b \in \mathbb{R}$, and $\hat{\theta} = \{\hat{W}_{1i}^{(l)}, \hat{W}_{2i}^{(l)}, \hat{b}_{1i}^{(l)}, \hat{b}_{2i}^{(l)}\}$.

As a result, the original optimization problem \eqref{semi_obj} concerning the velocity field $\boldsymbol{v}$ and the source function $f$ can be reformulated into an optimization problem focused on the network parameters. 
To compute the gradient of the objective function with respect to the parameters $\theta$ and $\tilde{\theta}$, we employ the standard back propagation method enabled by automatic differentiation.
Consequently, the parameters can be optimized using gradient descent and its variant. 
We summarize the proposed UOT algorithm in Algorithm 1.

\begin{algorithm}
    \caption{UOT algorithm} 
    \begin{algorithmic}
       \STATE \textbf{Input:} Initial density $\rho_0$, terminal density $\rho_1$, time steps $N$, network structure parameters $L, M, H$, hyperparameters $\lambda$ and $\alpha$.
       \STATE \textbf{Initialize:} 
        Parameters $\theta$ of the velocity field $\boldsymbol{v}$ and $\tilde{\theta}$ of the source function $f$.\\
       \FOR{$epoch=1,2,...,1000$}
          \STATE Randomly draw samples $\boldsymbol{x}_1,\boldsymbol{x}_2,...,\boldsymbol{x}_r$ from $\mu_0$.
          \STATE Set $t_j=\frac{j}{N}.$
          \FOR{$j=0$ to $N$}
             \STATE Solve $\mathbf{z}(\boldsymbol{x}_i,t_j)$ using RK4 by Eq. (\ref{forwd1}).
             \STATE Solve $\rho(\mathbf{z}(\boldsymbol{x}_i,t_j),t_j)$ using RK4 by Eq. (\ref{forwd2}).
             \STATE Solve $\mu(\mathbf{z}(\boldsymbol{x}_i,t_j),t_j)$ using RK4 by Eq. (\ref{forwd3}).
             
          \ENDFOR
       \STATE Compute $\overline{\mathcal{E}} (\rho,\boldsymbol{v})$, $\overline{\mathcal{R}} (\rho,f)$ and $\overline{\mathcal{P}} (\rho)$ using standard composite Simpson formula by Eqs. (\ref{semdisE}), (\ref{semdisR}), (\ref{semdisP}), respectively.
       \STATE Compute gradient of the objective function $\overline{\mathcal{F}} (\rho,\boldsymbol{v},f)$ with respect to parameters $\theta$ and $\tilde{\theta}$ using automatic differentiation.
       \STATE Update $\theta$ and $\tilde{\theta}$ by using the gradient descent method.
       \ENDFOR
       \STATE \textbf{Output:} Optimized parameters $\theta$ and $\tilde{\theta}$.
          
    \end{algorithmic}
\label{algorithm}
\end{algorithm}

\section{Numerical Experiments}\label{Experimentresult}

In this section, we illustrate the effectiveness of the proposed method. The numerical examples include Gaussian examples, Gaussian mixture examples and crowd motion examples. 
Note that when the means of the densities $\rho_0$ and $\rho_1$ are not identical, we first set the source function $f=0$ and employ the KL divergence to find the target mean, a process referred to as initialization.

In our experiments, unless otherwise stated, we adopt the following default settings: the time steps $N = 10$, the GKL penalty parameter $\lambda = 10000$, and the unbalanced coefficient $\alpha=0.01$.
The parameters of the network structure are as follows: the number of the basis $M = 5$, the width of neural network $L = 2$, and the number of hidden units $H = 10$. 
All models are trained for 1000 epochs (training iterations) using the Adam optimizer and a consistent batch size of 1024 training samples for all high-dimensional cases. The initial learning rate is set to 0.01, with a decay factor of 0.98 applied every 10 epochs.
All the experiments are conducted in PyTorch on an NVIDIA Tesla V100 GPU with 32GB memory.

\subsection{Examples of UOT with different weight}

In this section, we discuss the coefficient $\alpha$ governing the source term in the dynamic UOT problem \eqref{dOT}. 
This coefficient parameter regulates the creation and destruction of mass during transportation. 
When $\alpha$ approaches zero, the UOT solution converges towards the classical OT solution in \eqref{OT}.

\begin{figure}[H]
\begin{center}

\includegraphics[width=2.5cm]{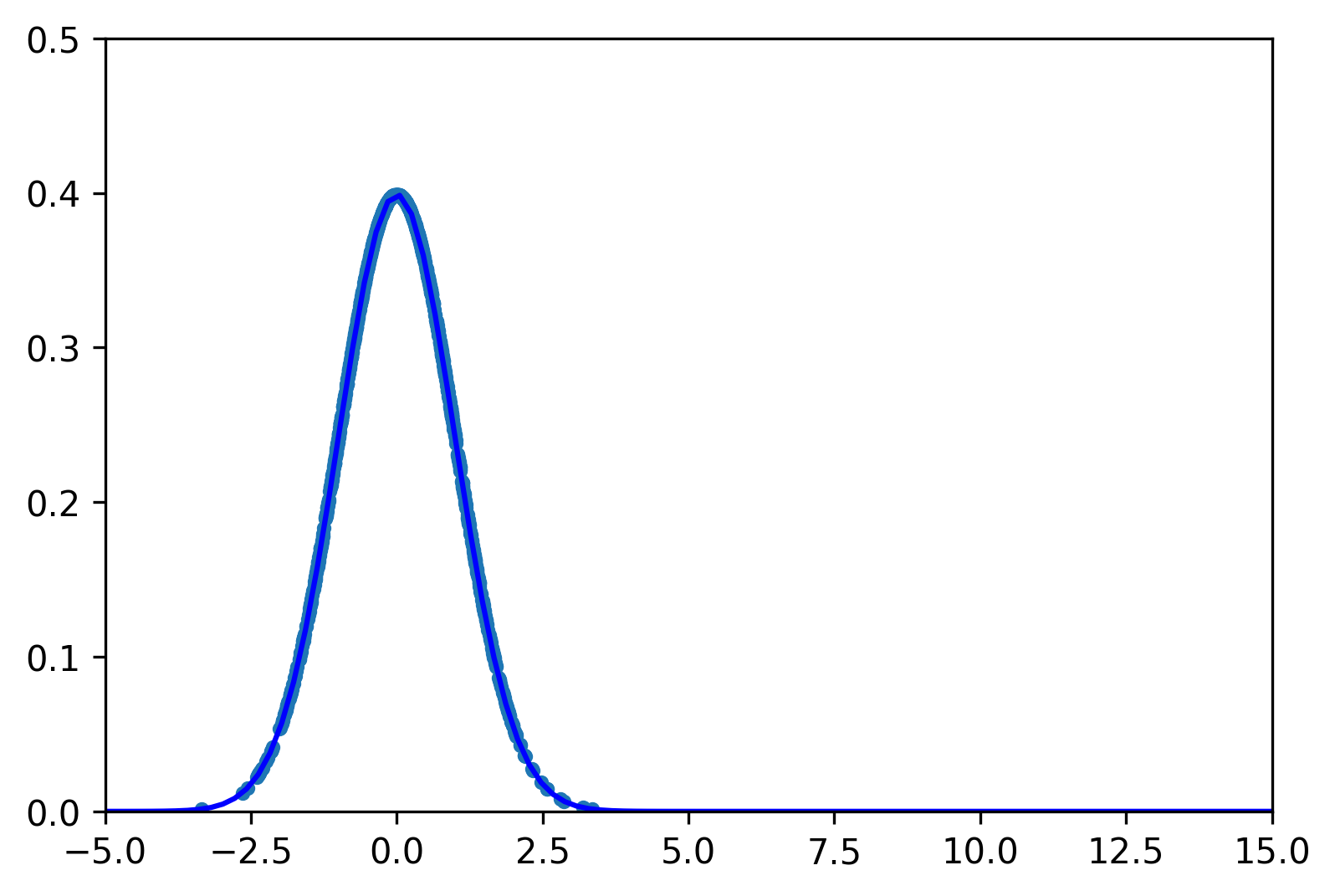}
\includegraphics[width=2.5cm]{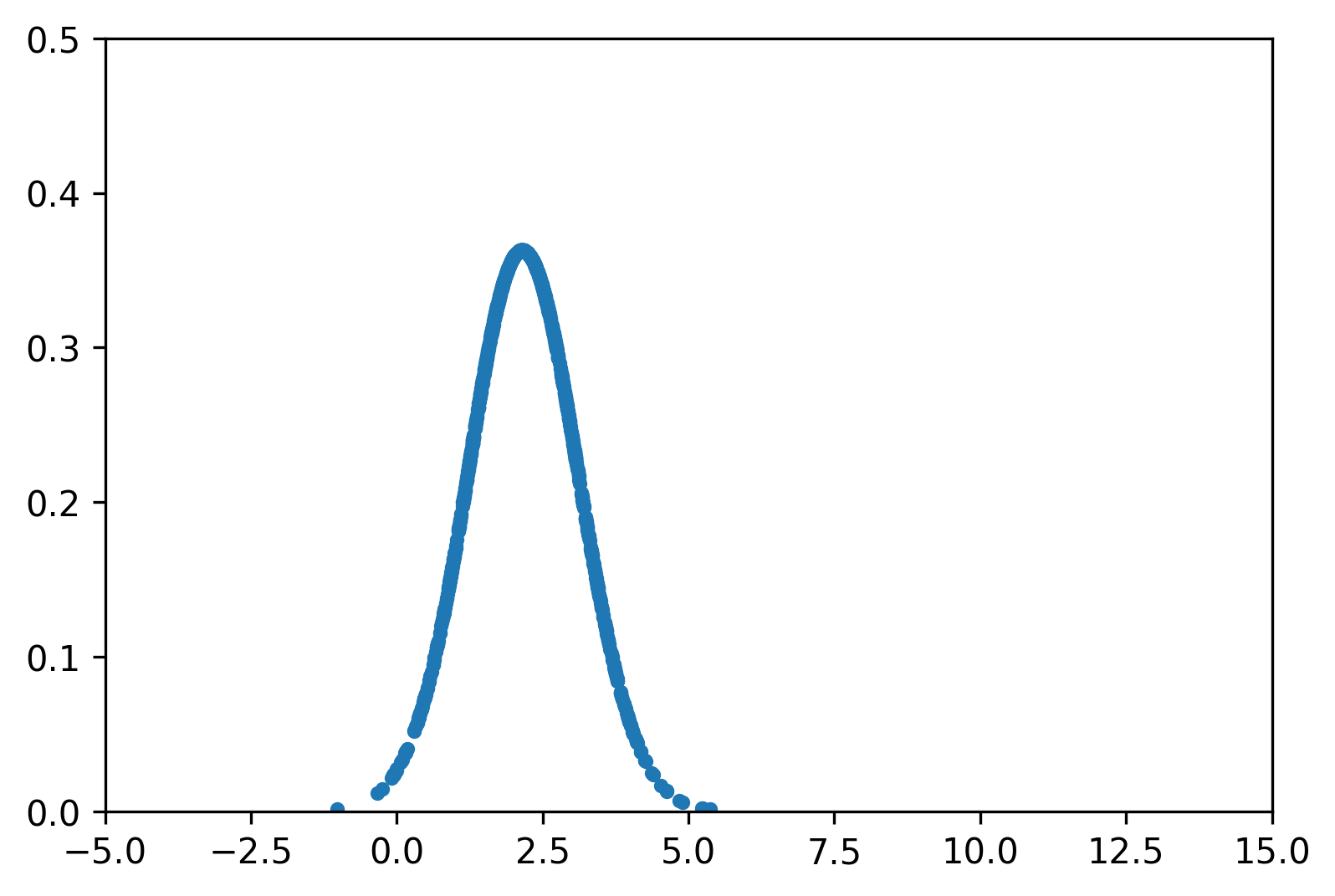}
\includegraphics[width=2.5cm]{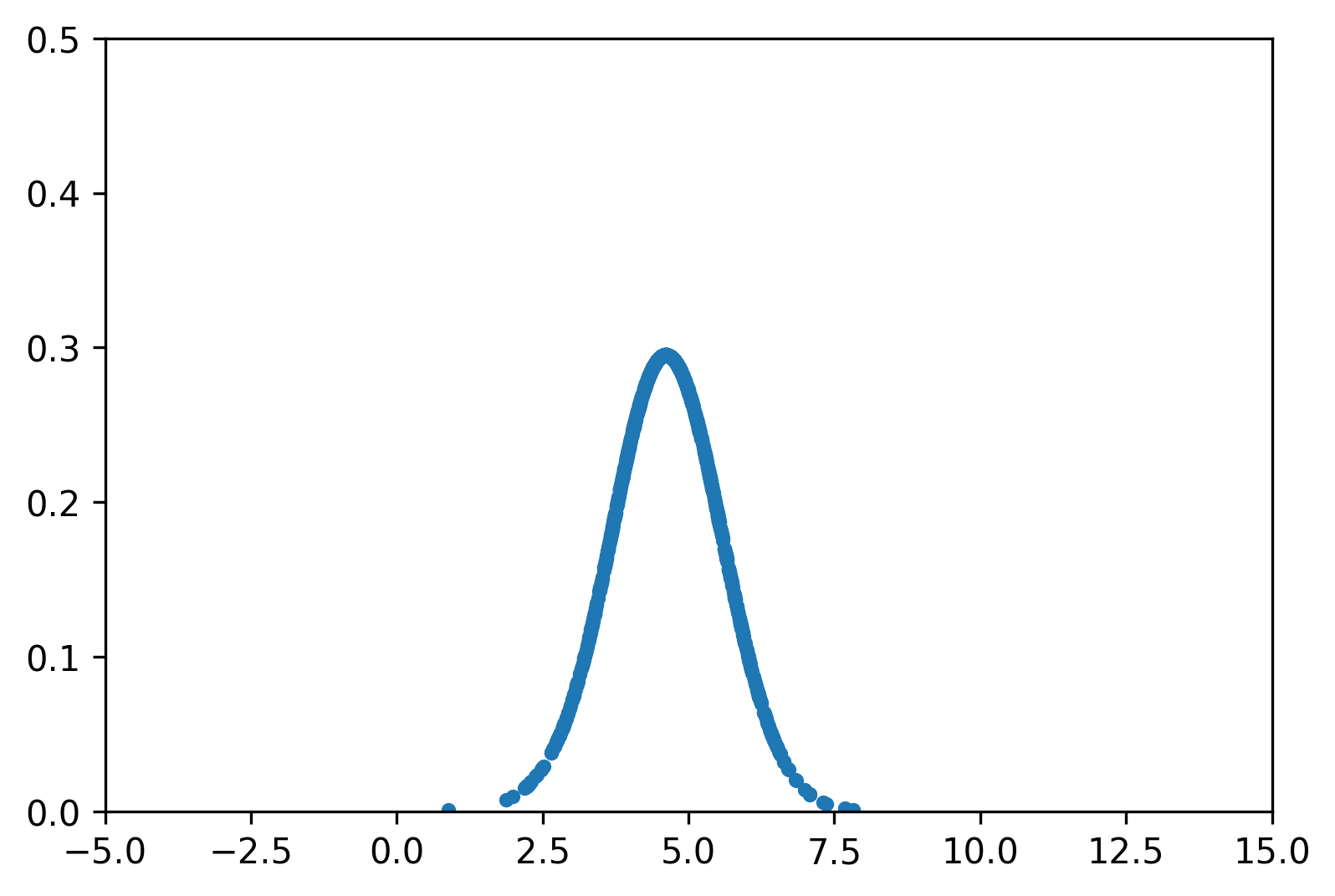}
\includegraphics[width=2.5cm]{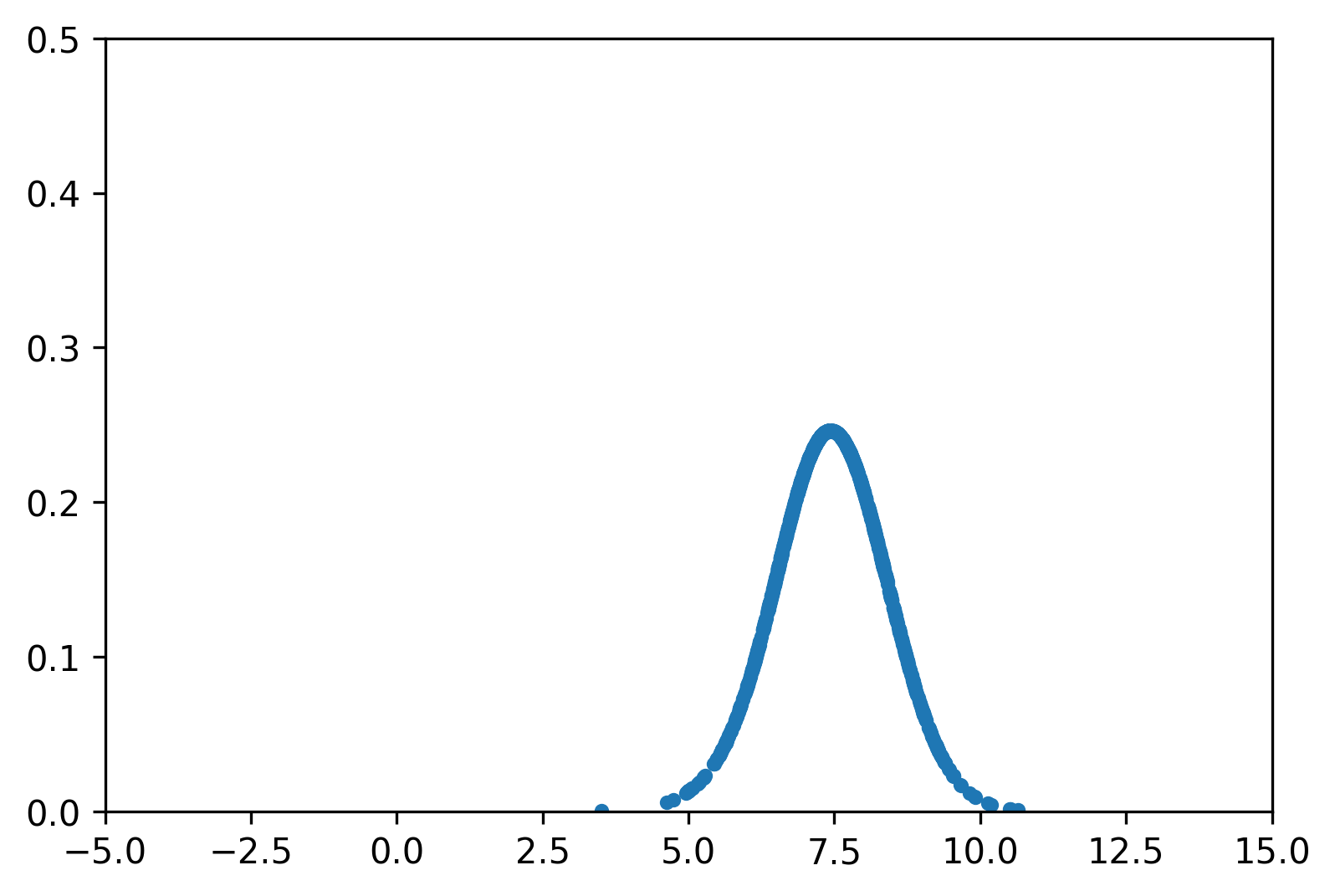}
\includegraphics[width=2.5cm]{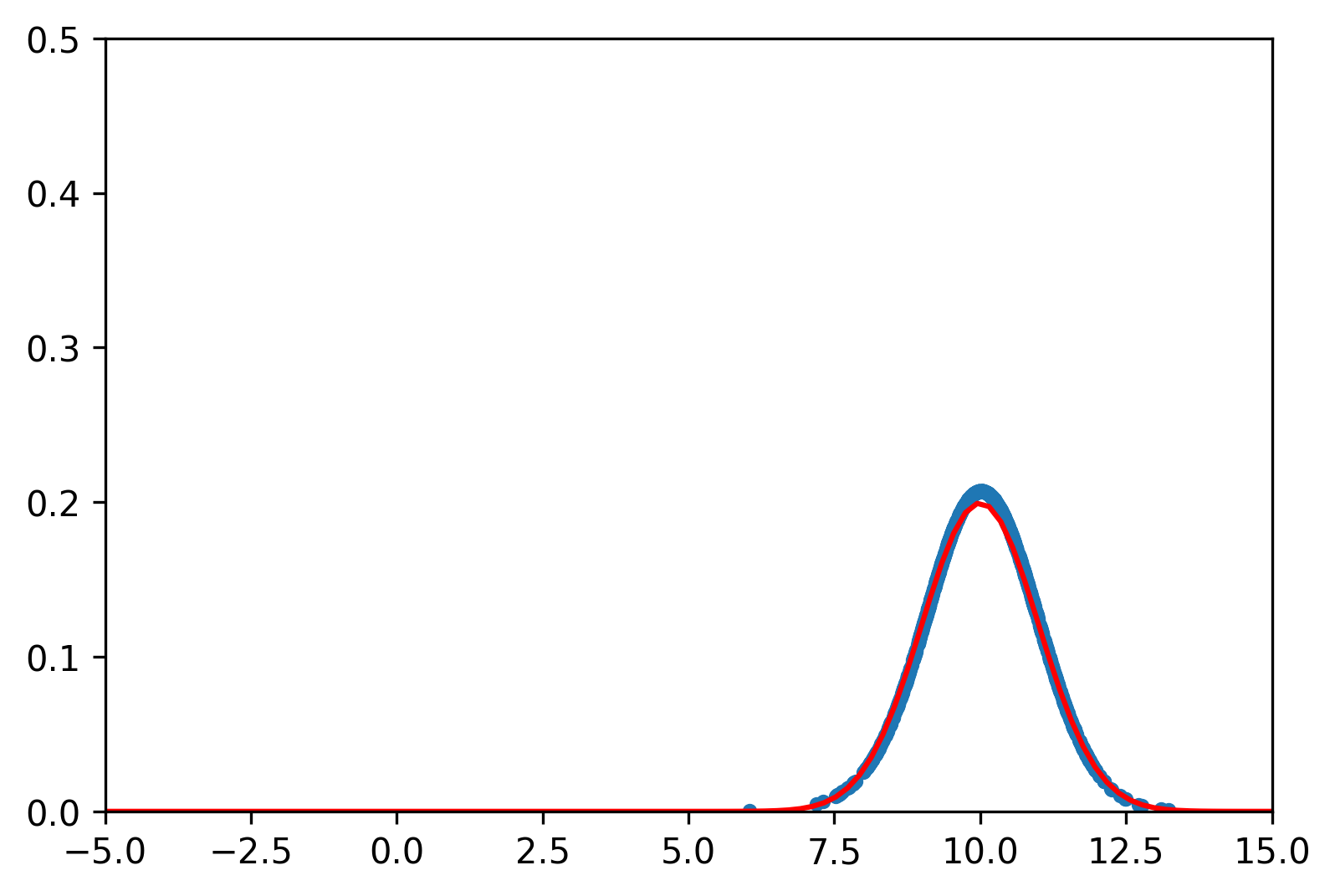}\\
\vspace{5pt}

\includegraphics[width=2.5cm]{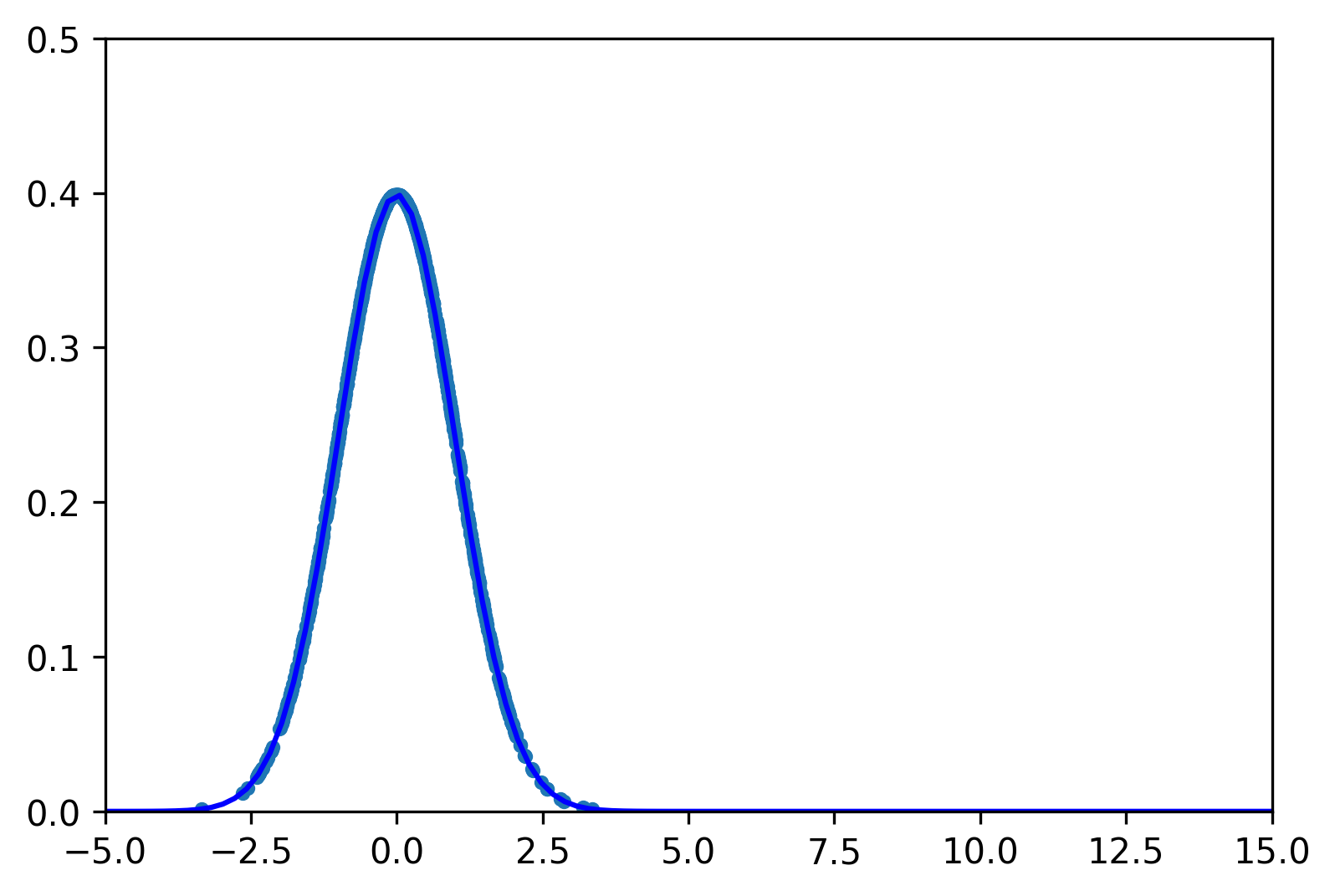}
\includegraphics[width=2.5cm]{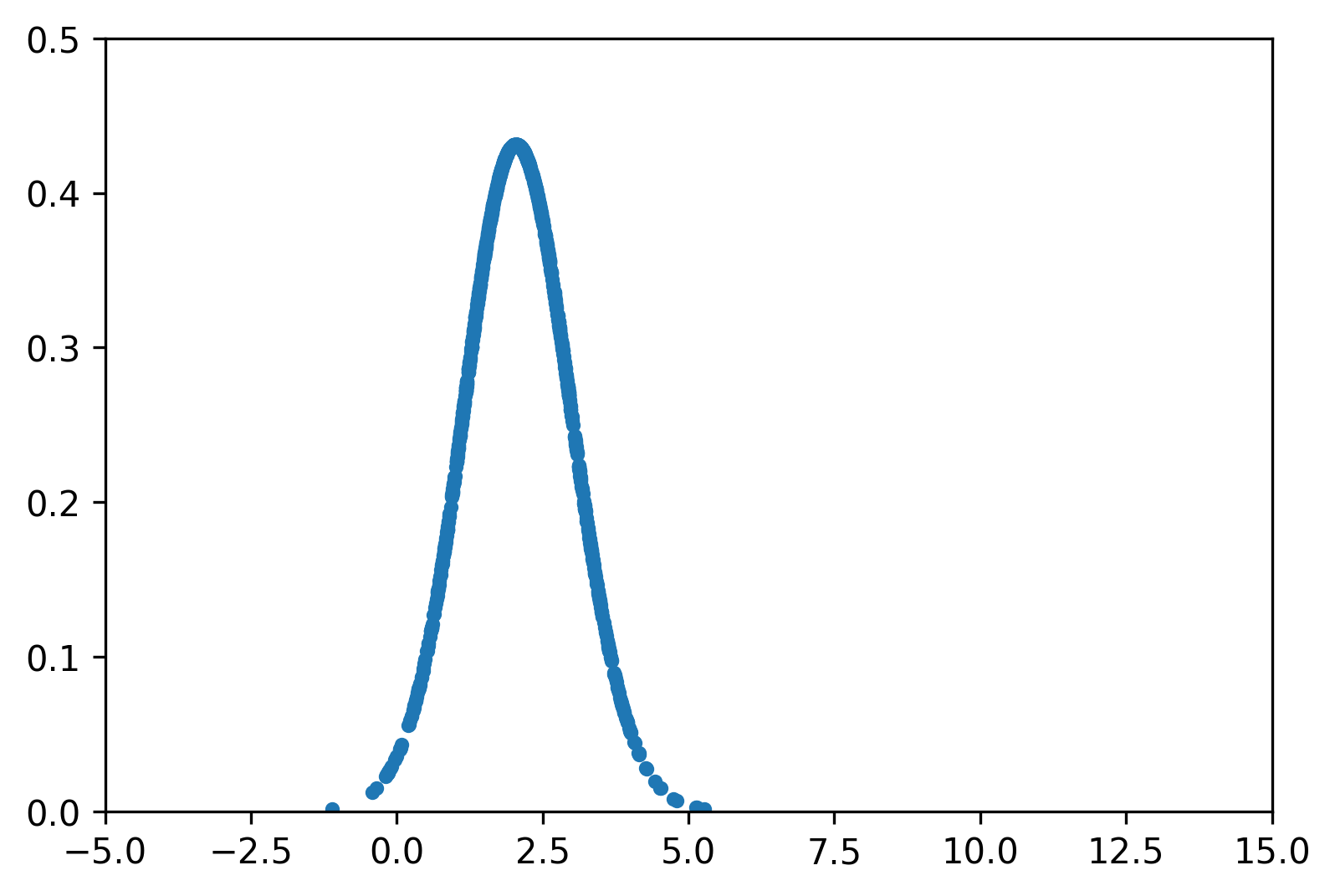}
\includegraphics[width=2.5cm]{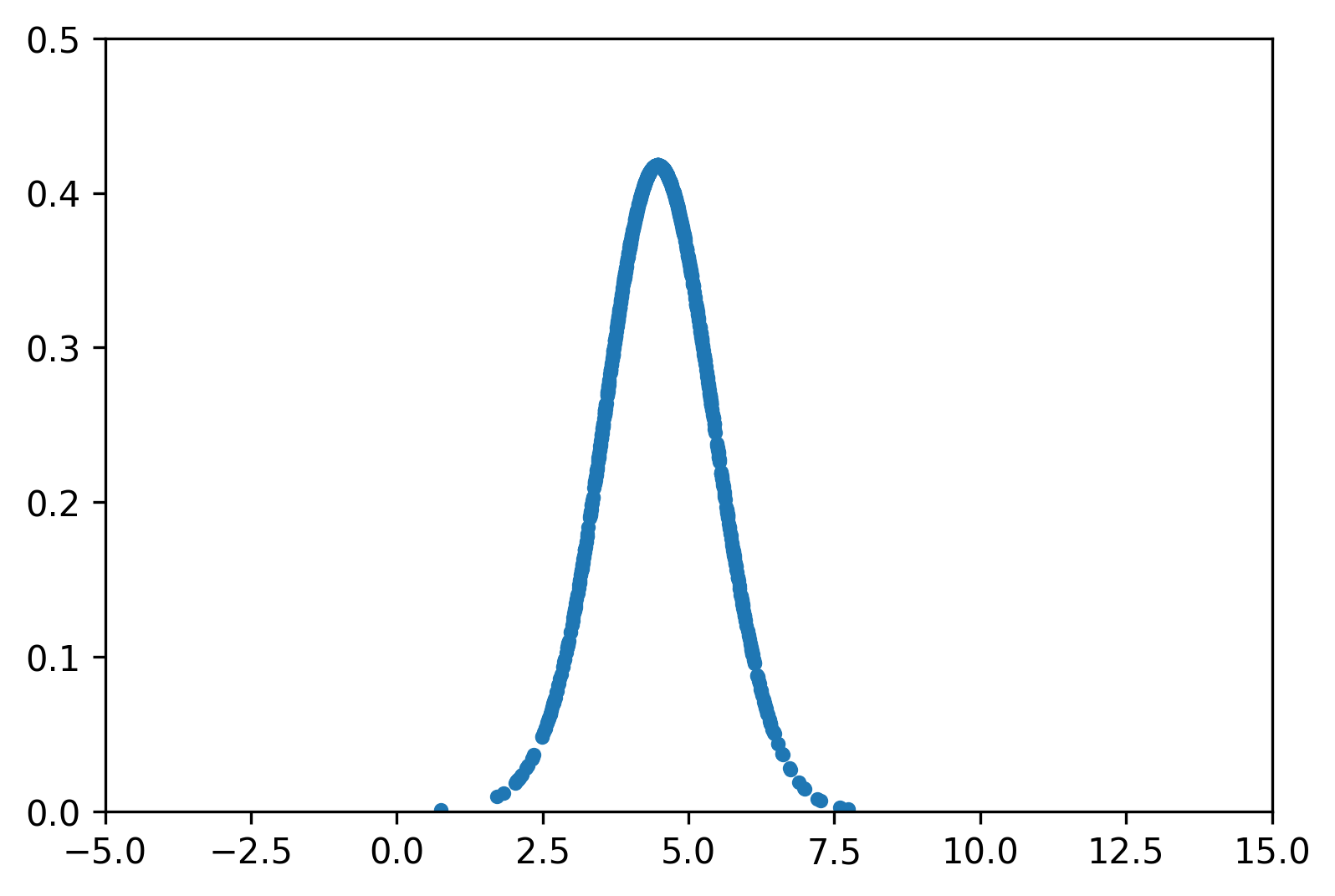}
\includegraphics[width=2.5cm]{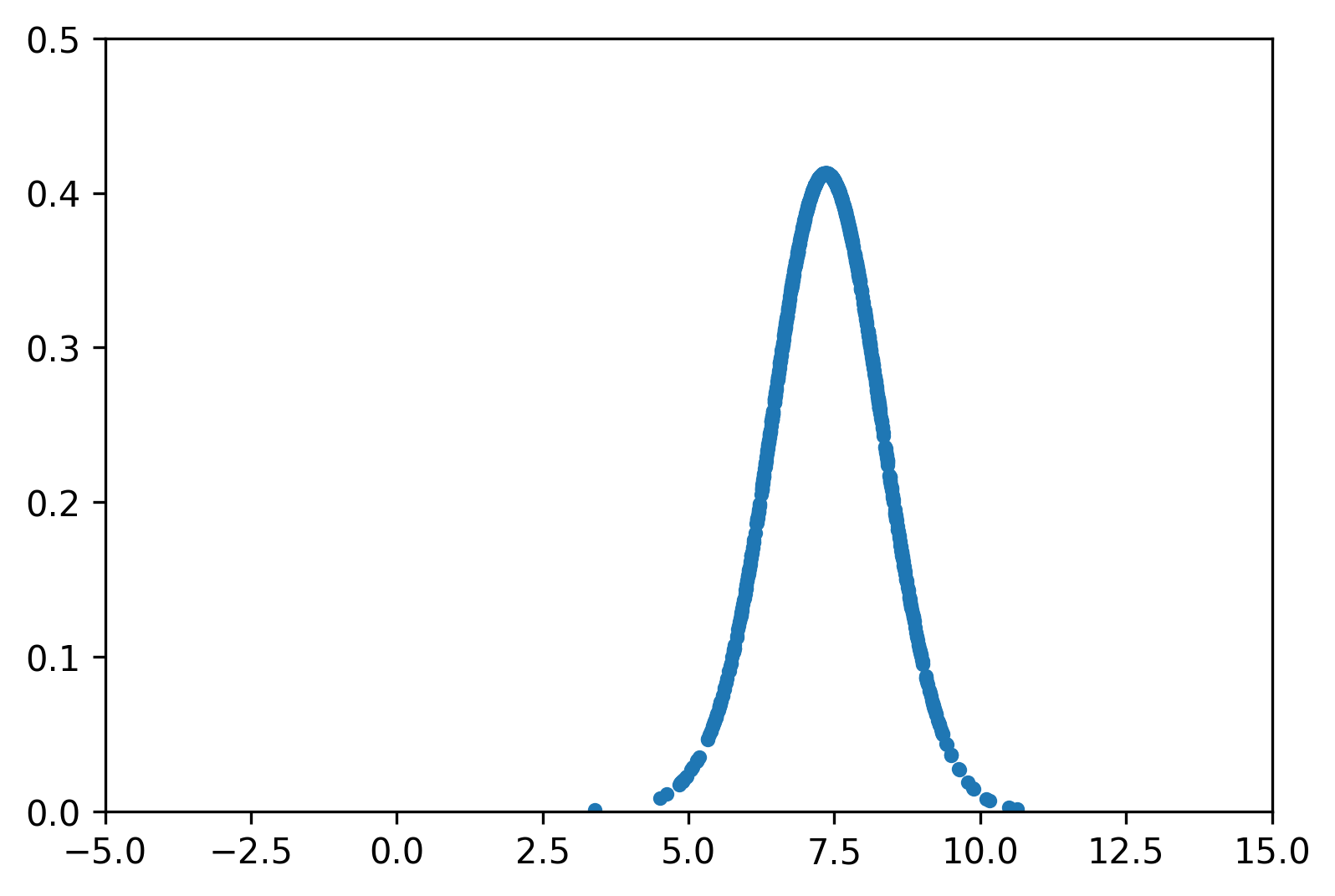}
\includegraphics[width=2.5cm]{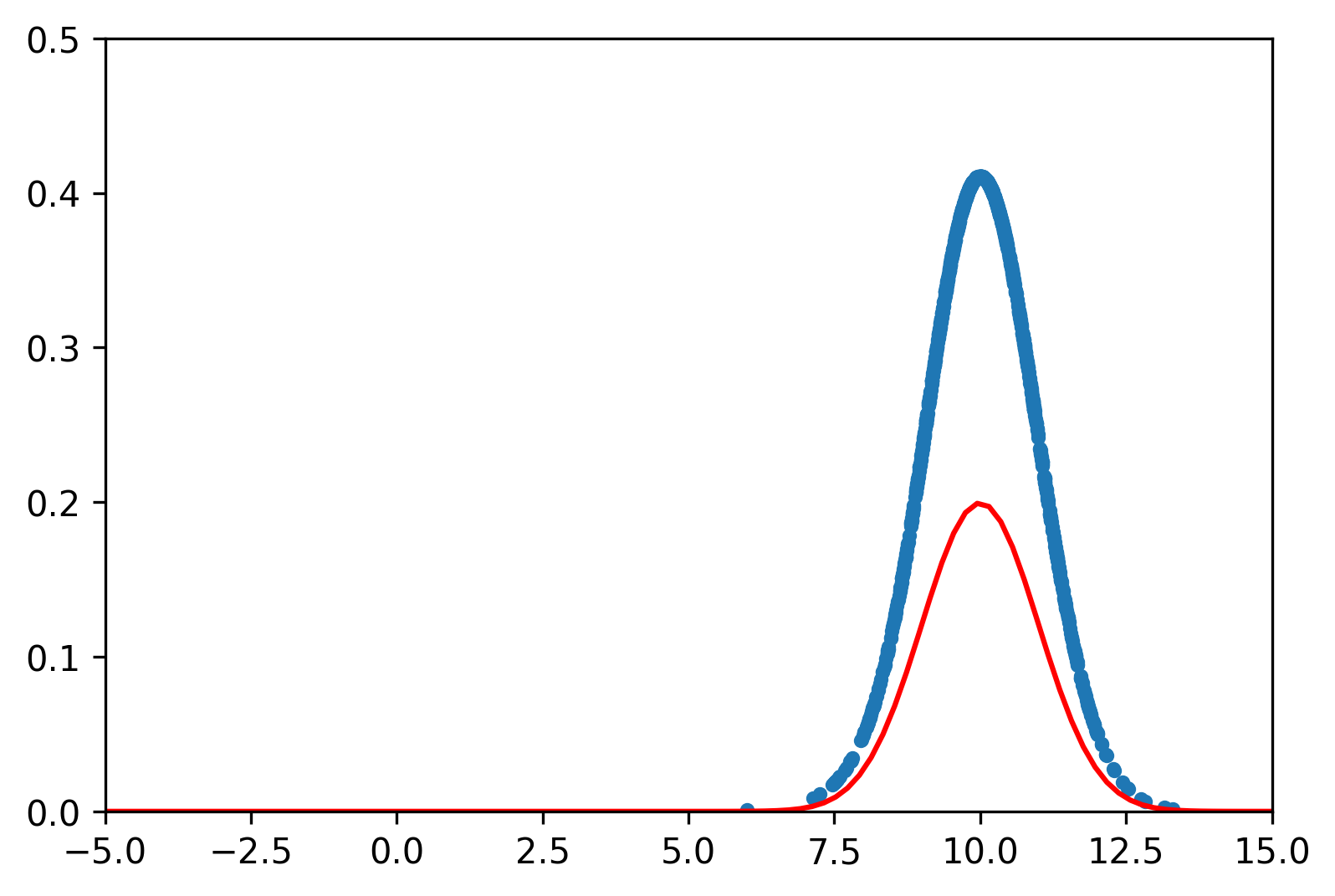}\\
\vspace{5pt}

\subfigure[$\rho_0(\boldsymbol{x}_0)$]{\includegraphics[width=2.5cm]{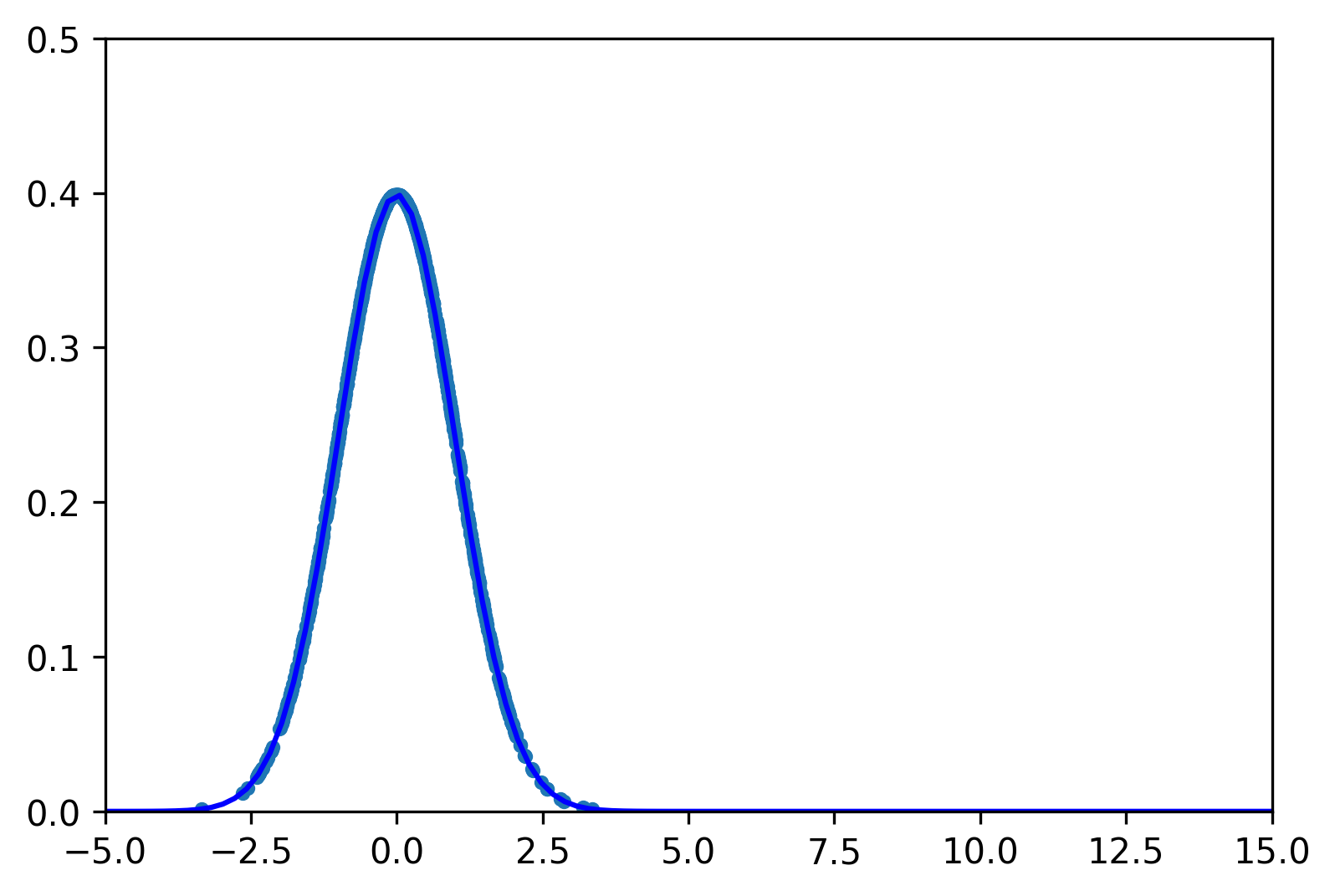}}
\subfigure[$\tilde{\rho}_{1/4}(\tilde{\boldsymbol{x}}_{1/4})$]{\includegraphics[width=2.5cm]{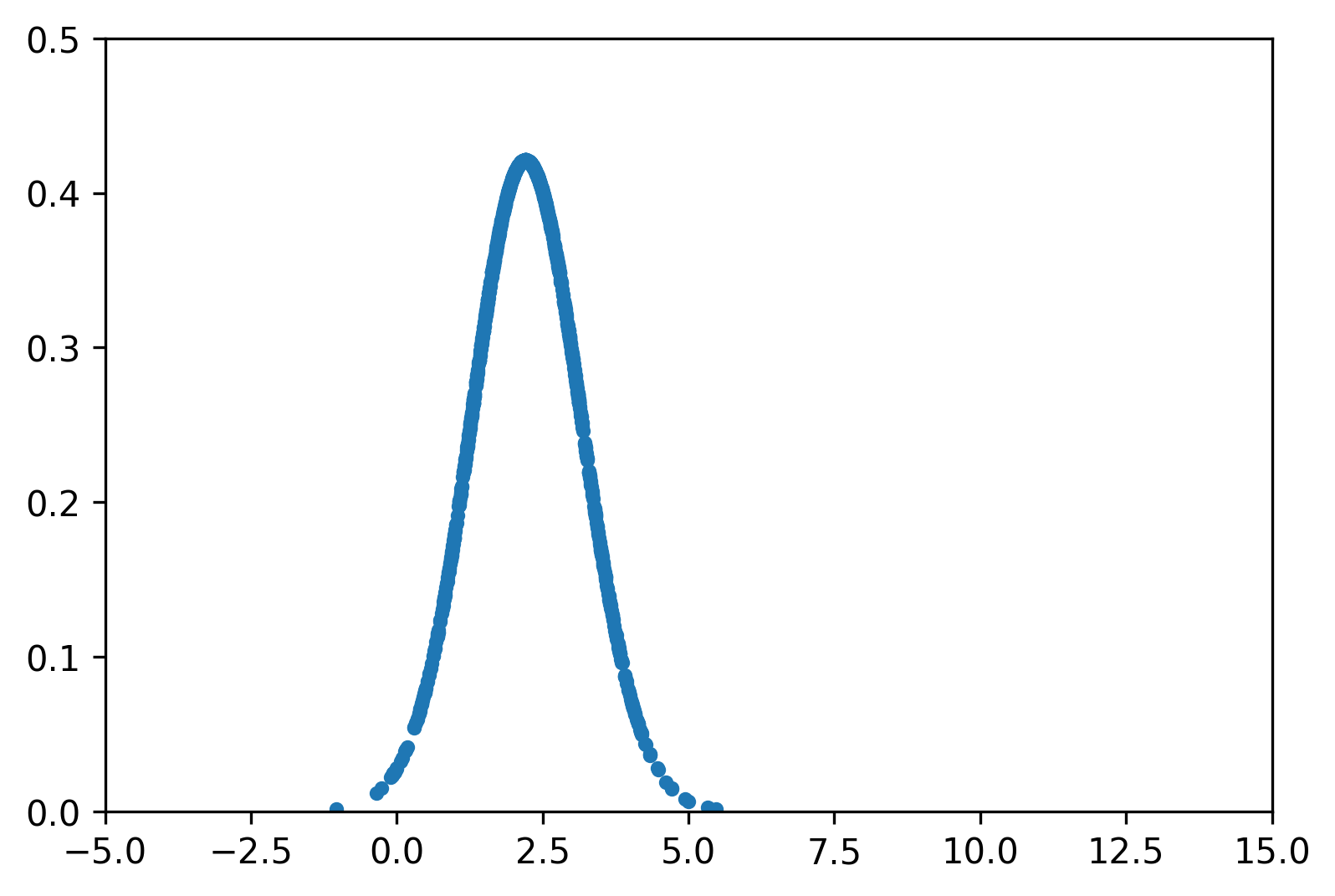}}
\subfigure[$\tilde{\rho}_{1/2}(\tilde{\boldsymbol{x}}_{1/2})$]{\includegraphics[width=2.5cm]{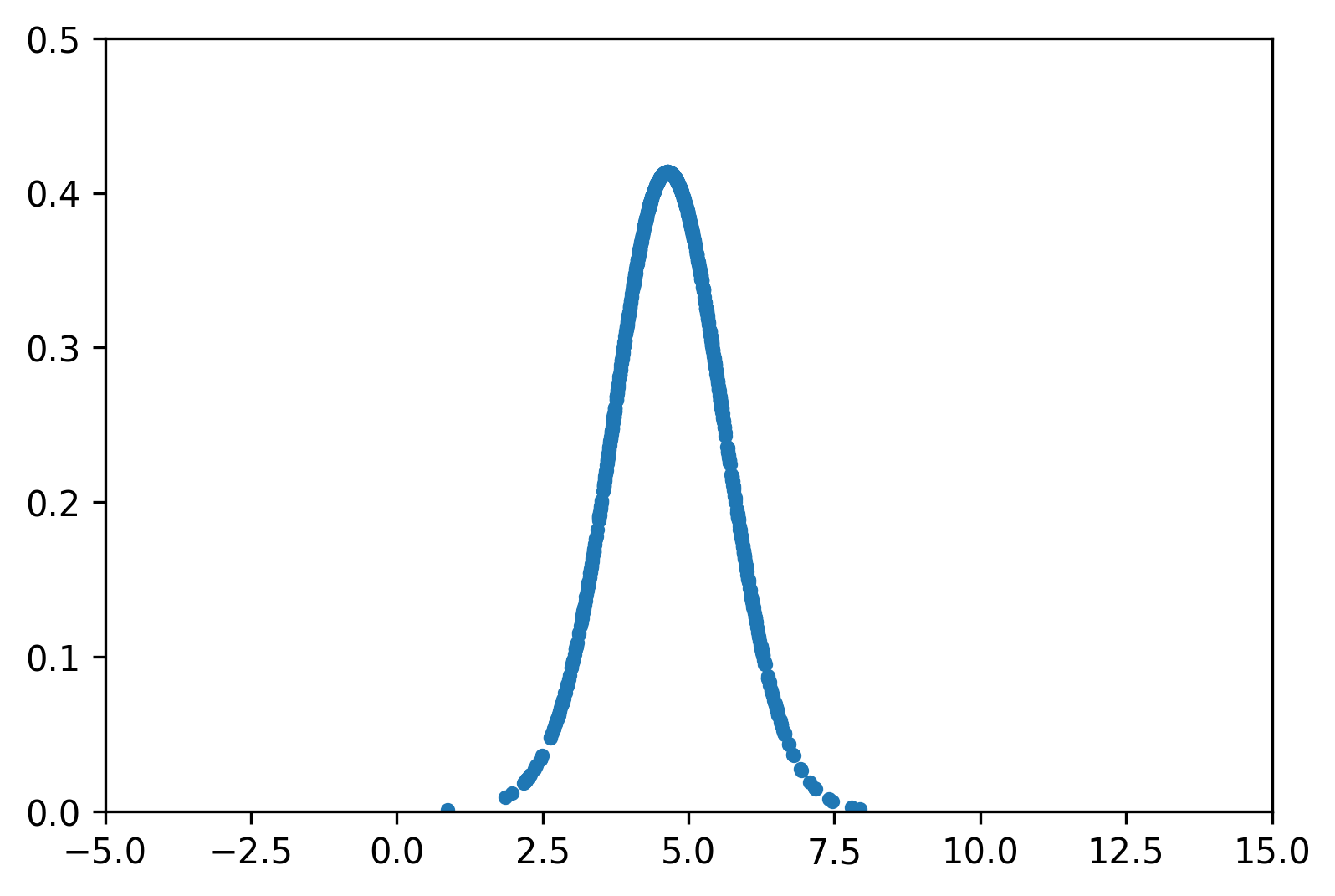}}
\subfigure[$\tilde{\rho}_{3/4}(\tilde{\boldsymbol{x}}_{3/4})$]{\includegraphics[width=2.5cm]{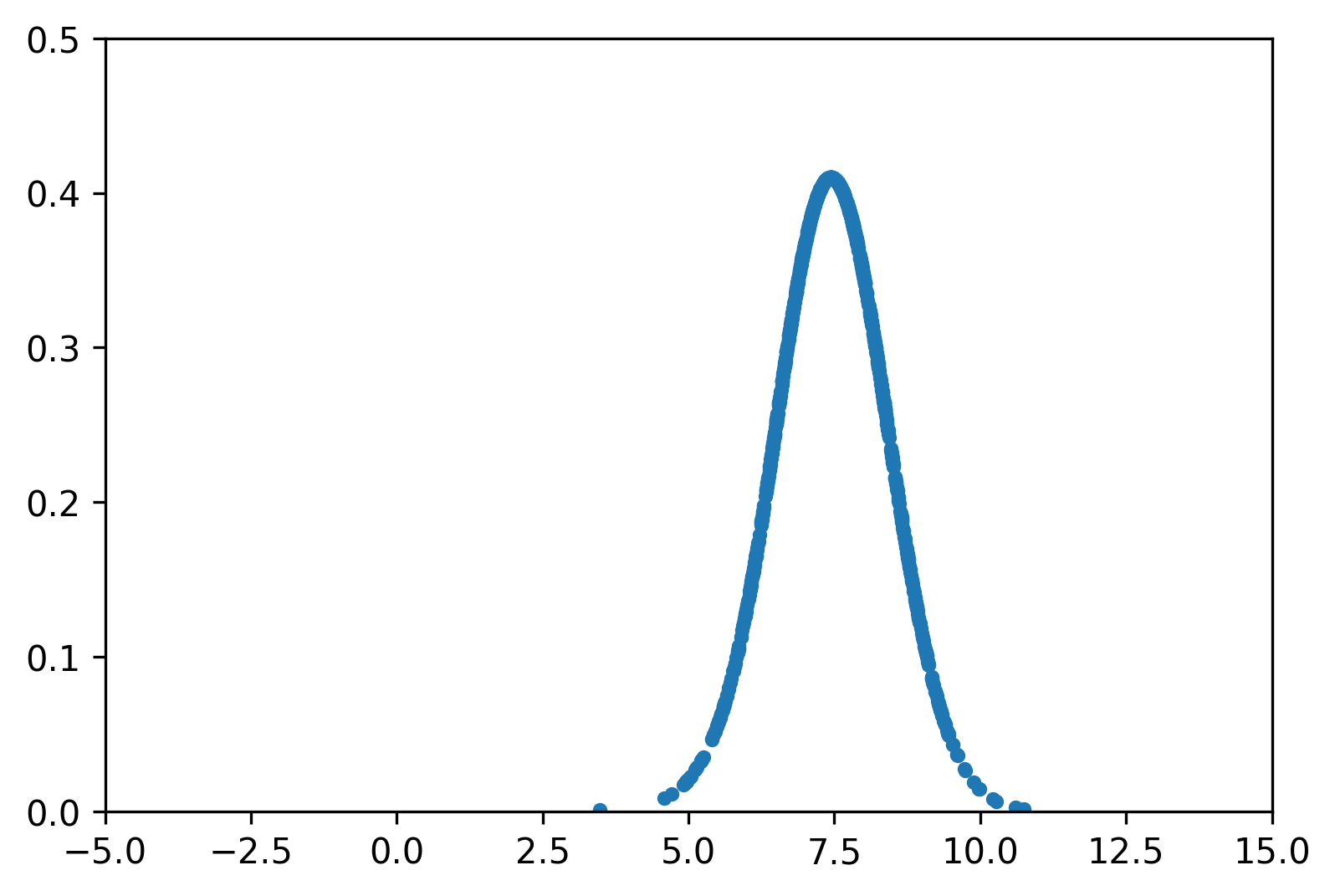}}
\subfigure[$\tilde{\rho}_{1}(\tilde{\boldsymbol{x}}_{1})$]{\includegraphics[width=2.5cm]{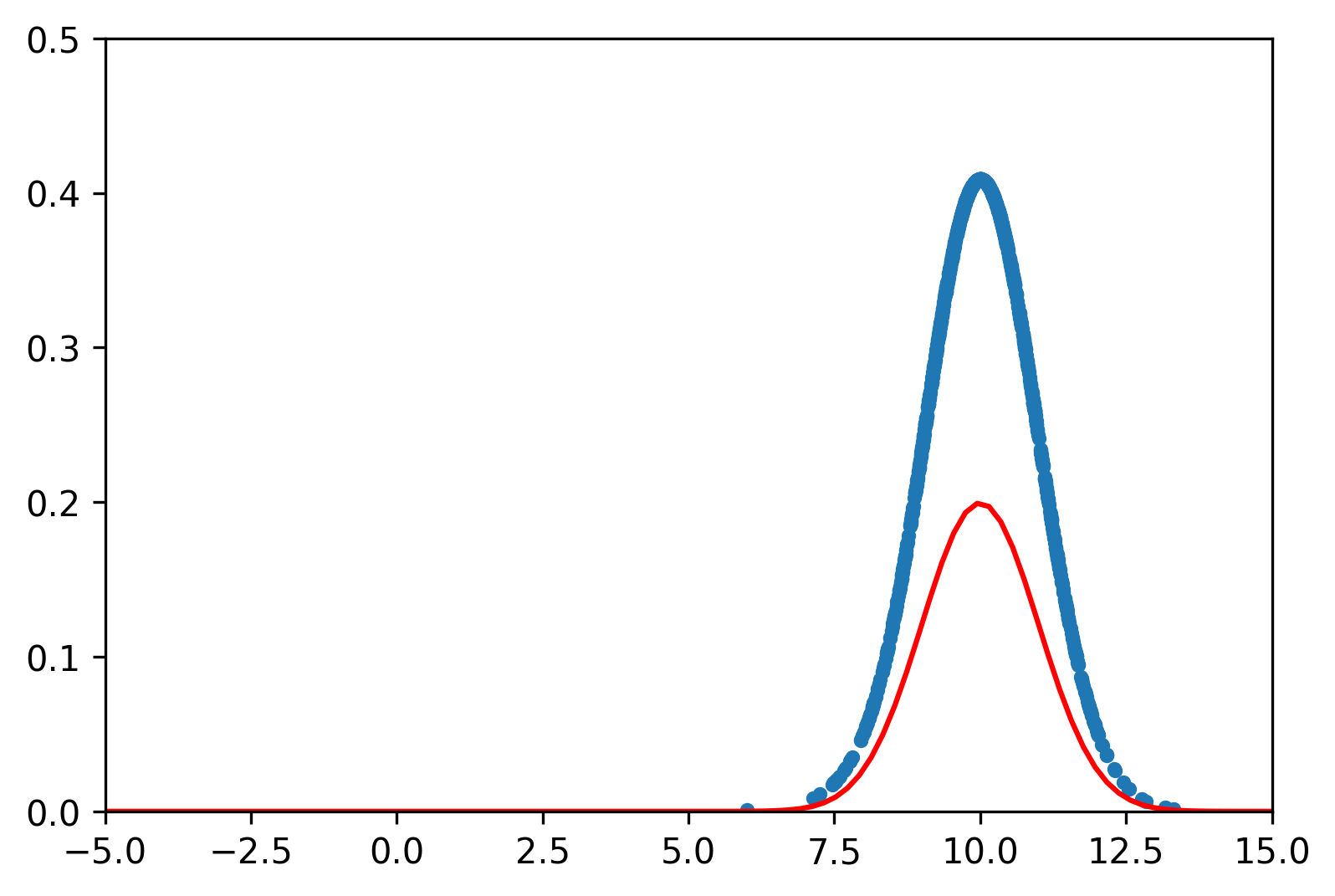}}

\caption{Top: UOT with $\alpha=10^{-2}$, middle: UOT with $\alpha=10^{-6}$, bottom: OT.}
\label{uot1}
\end{center}
\end{figure}

\begin{figure}[H]
\begin{center}
\includegraphics[width=2.5cm]{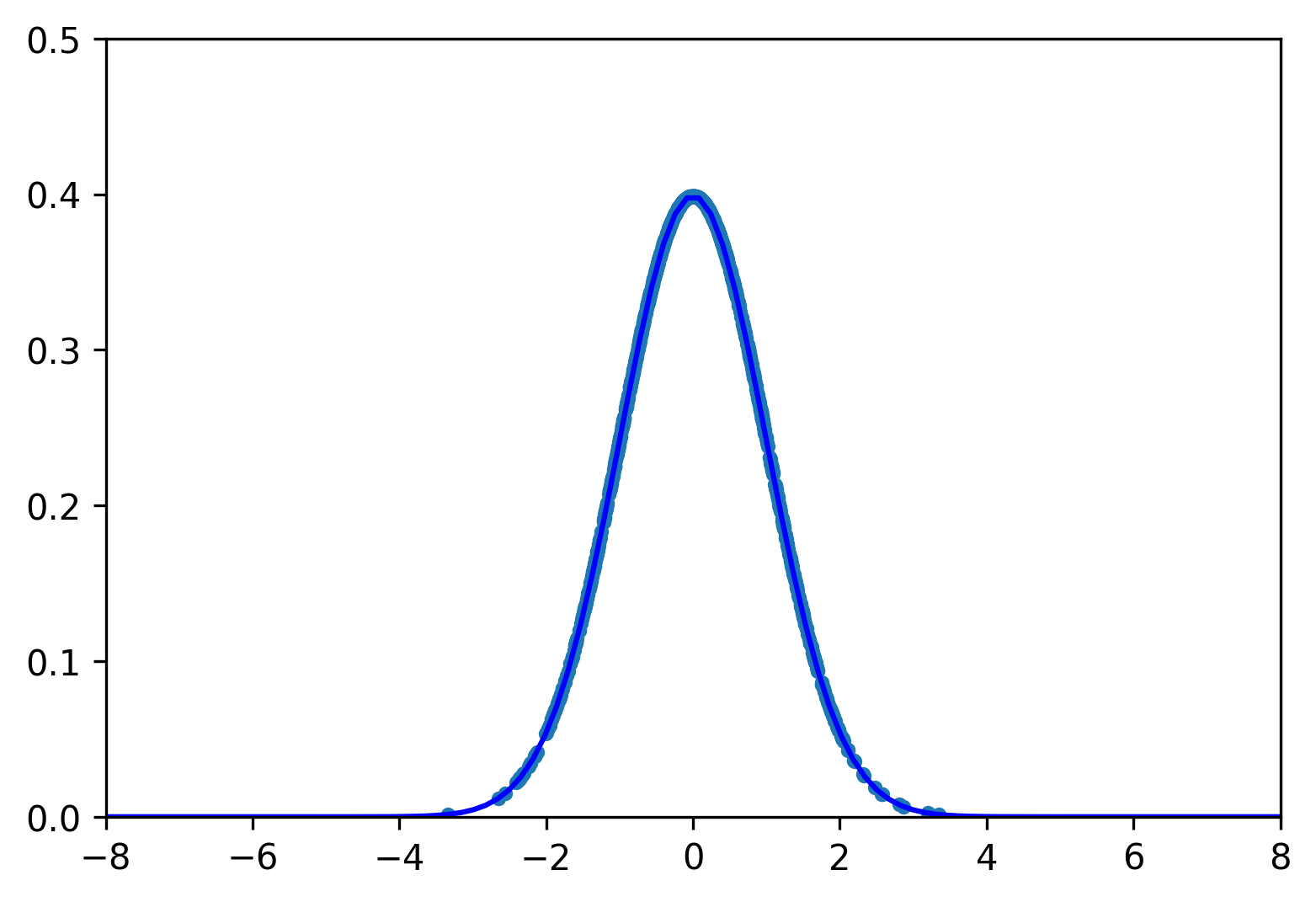}
\includegraphics[width=2.5cm]{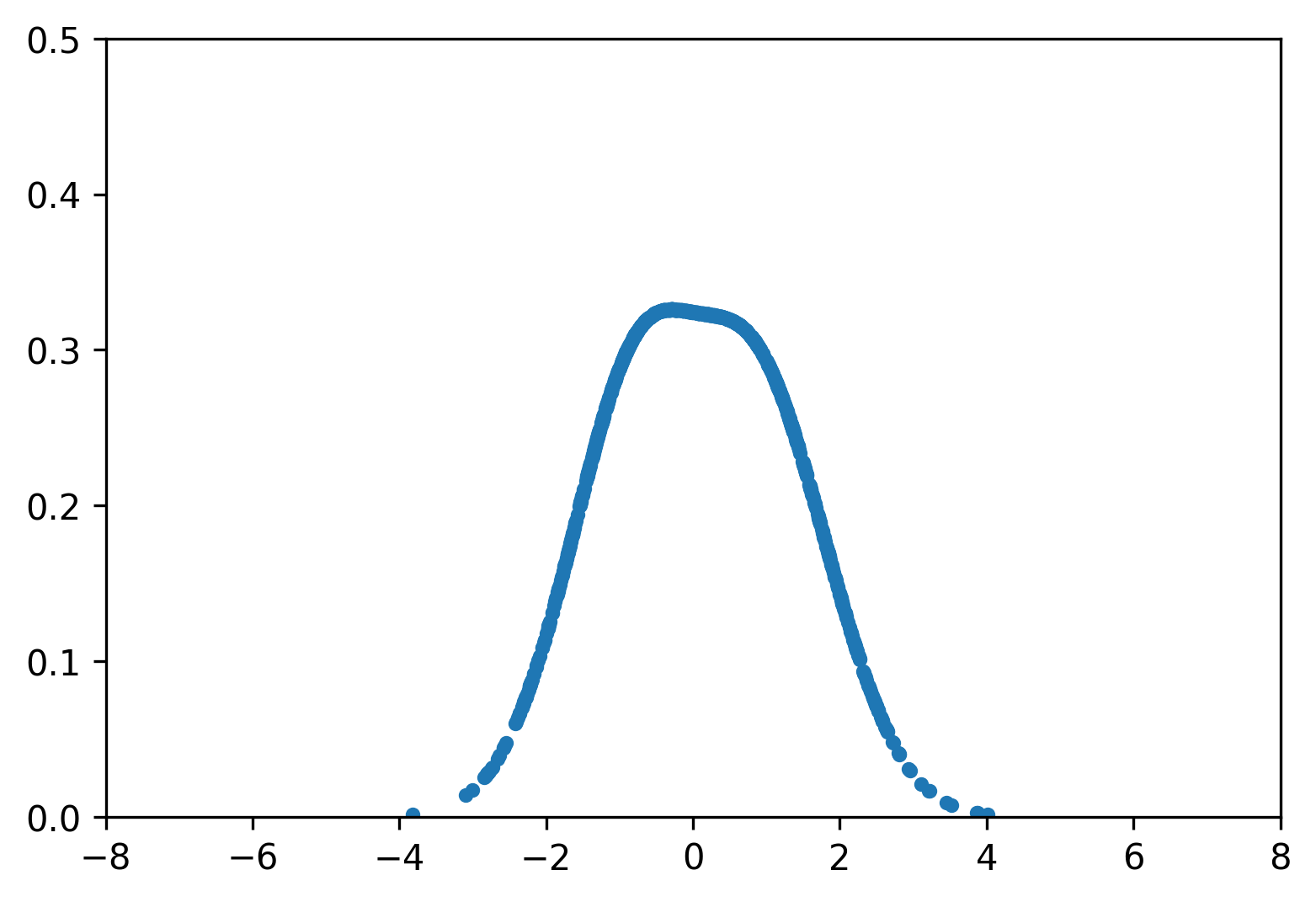}
\includegraphics[width=2.5cm]{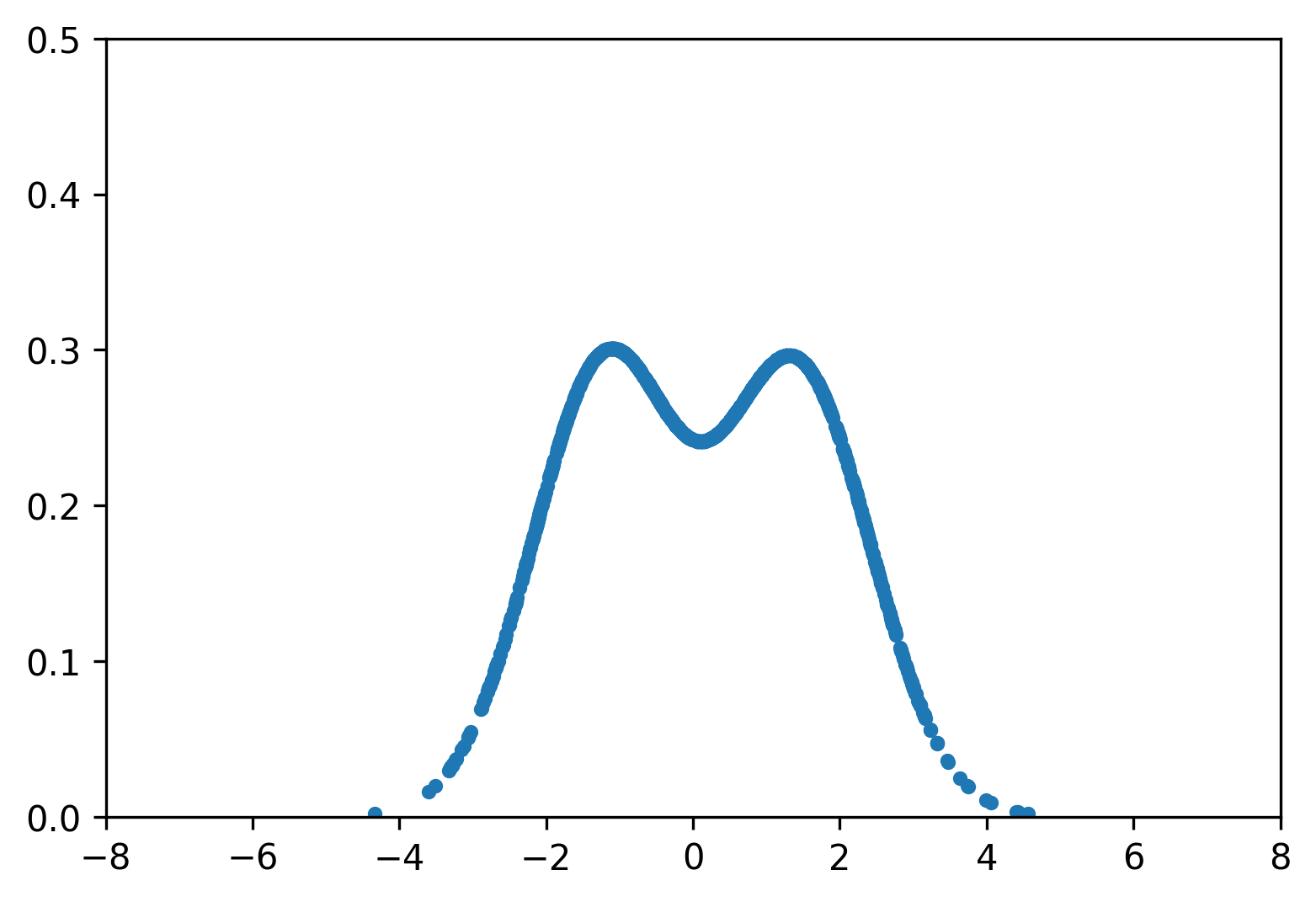}
\includegraphics[width=2.5cm]{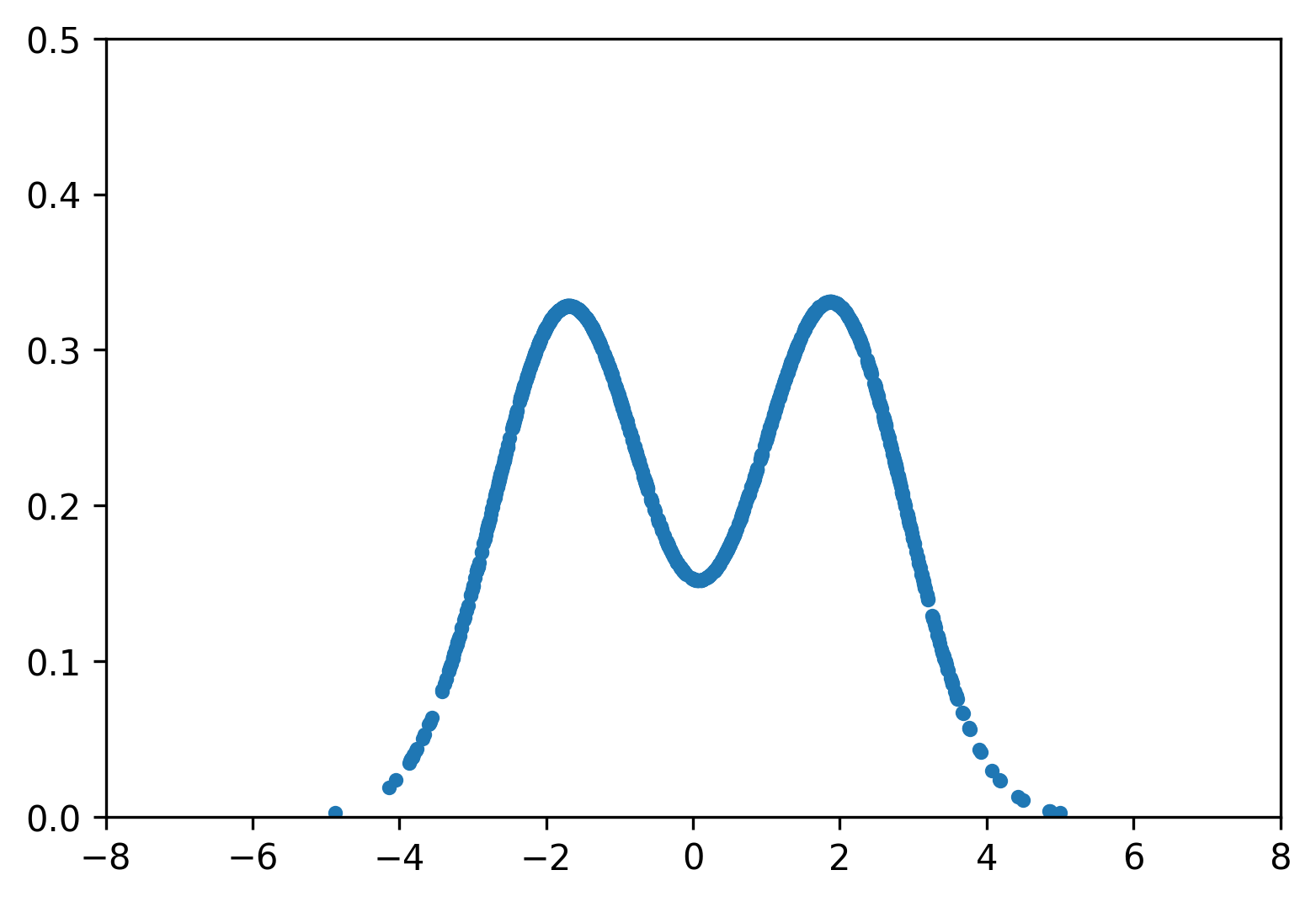}
\includegraphics[width=2.5cm]{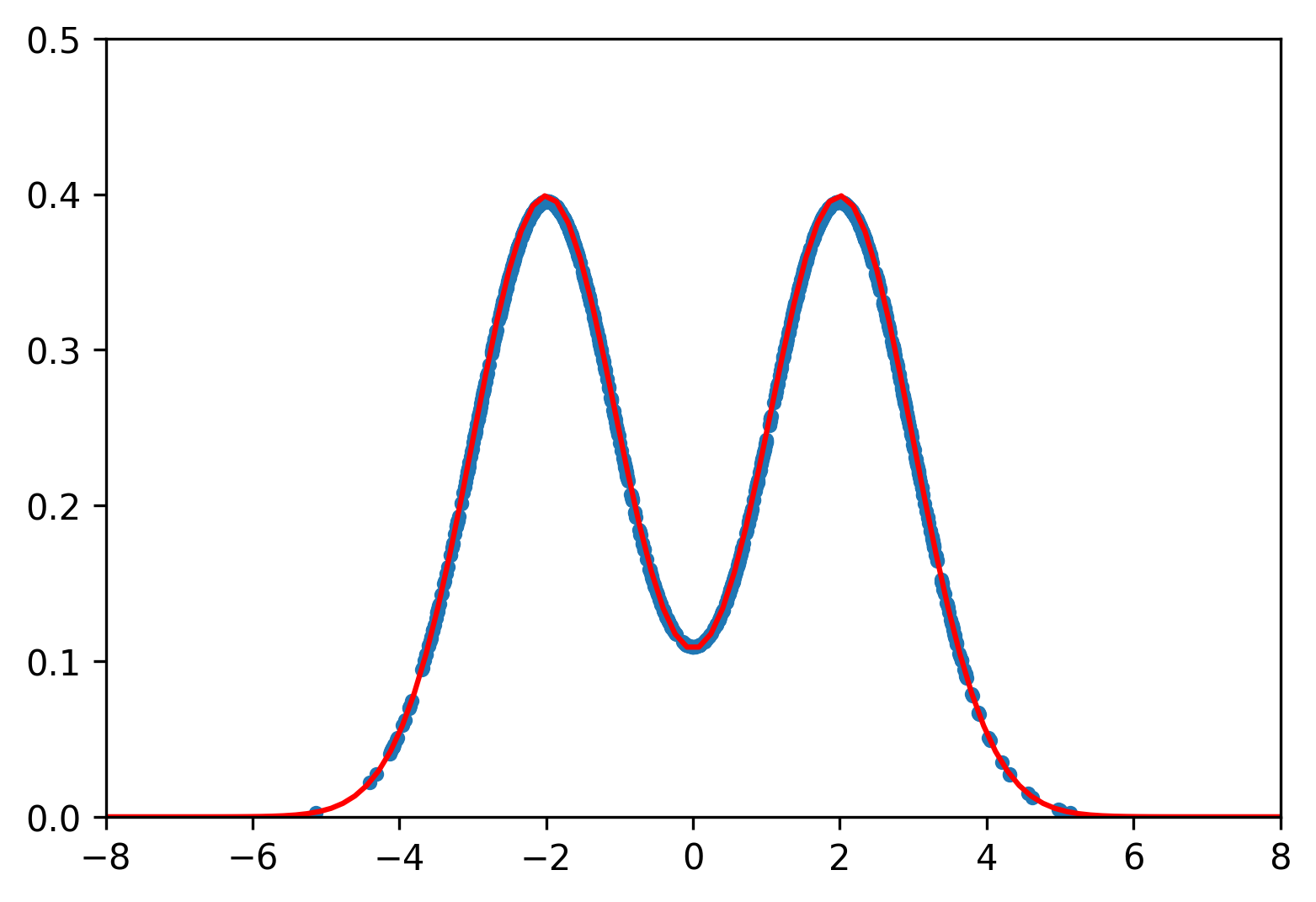}\\
\vspace{5pt}

\includegraphics[width=2.5cm]{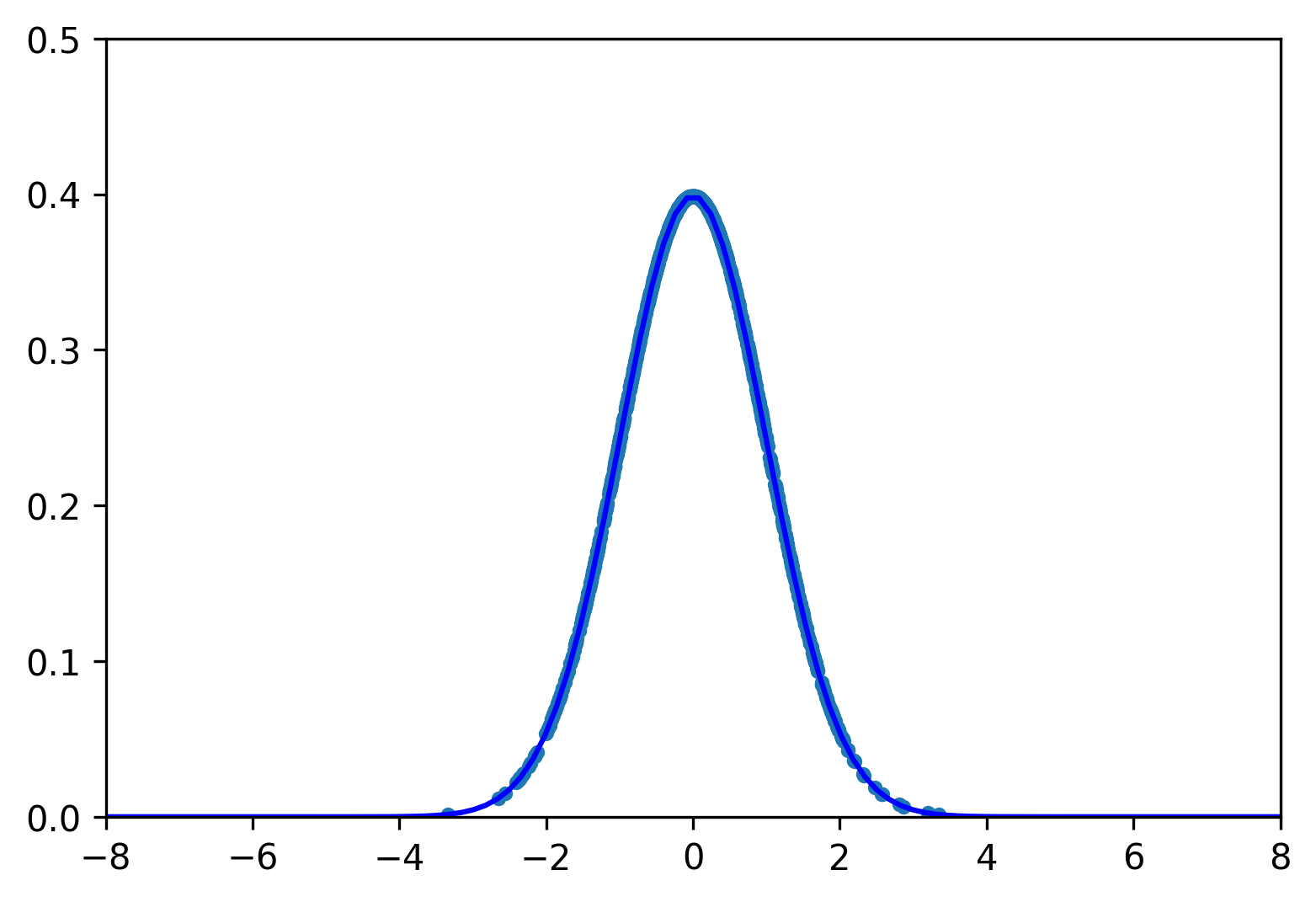}
\includegraphics[width=2.5cm]{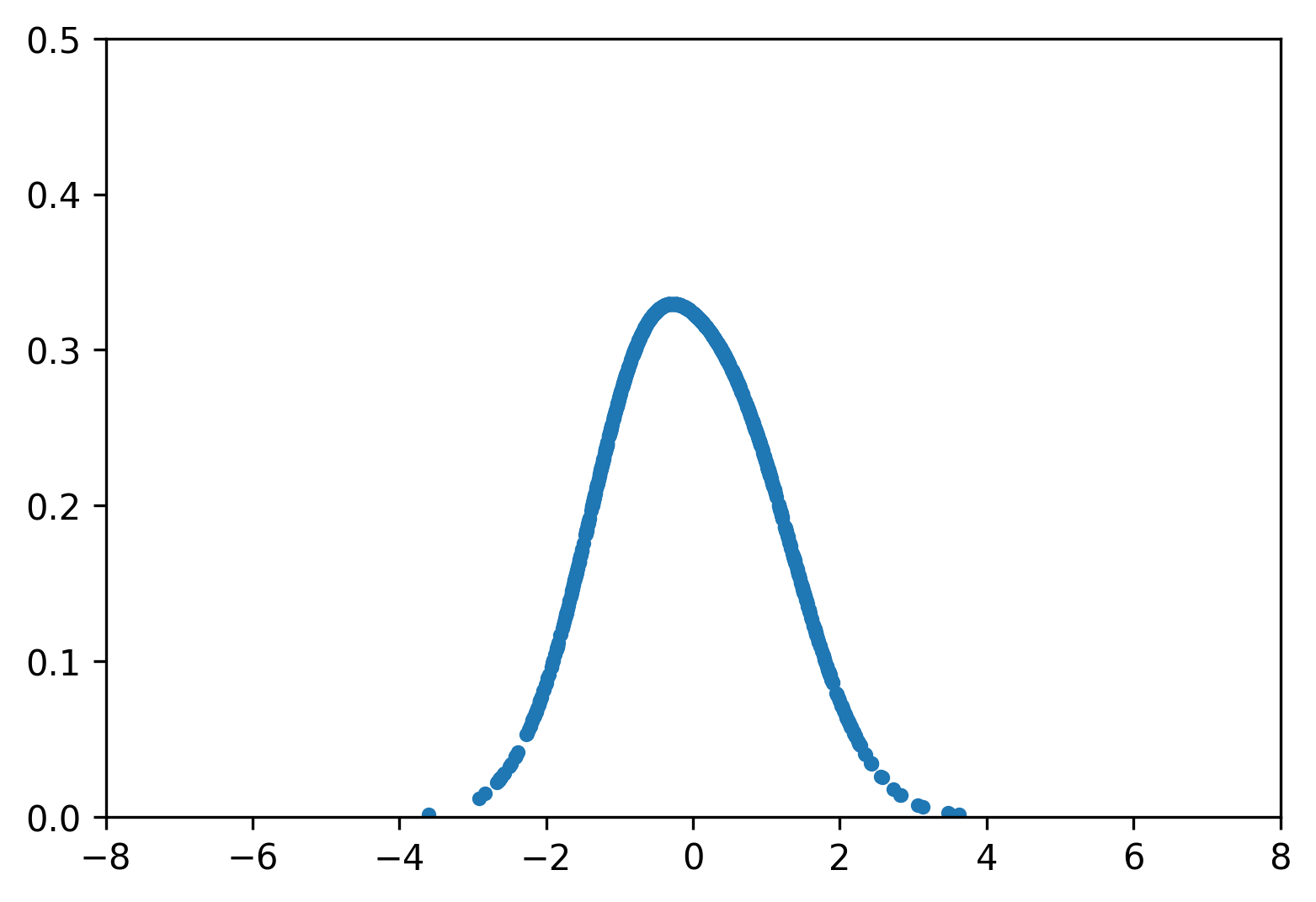}
\includegraphics[width=2.5cm]{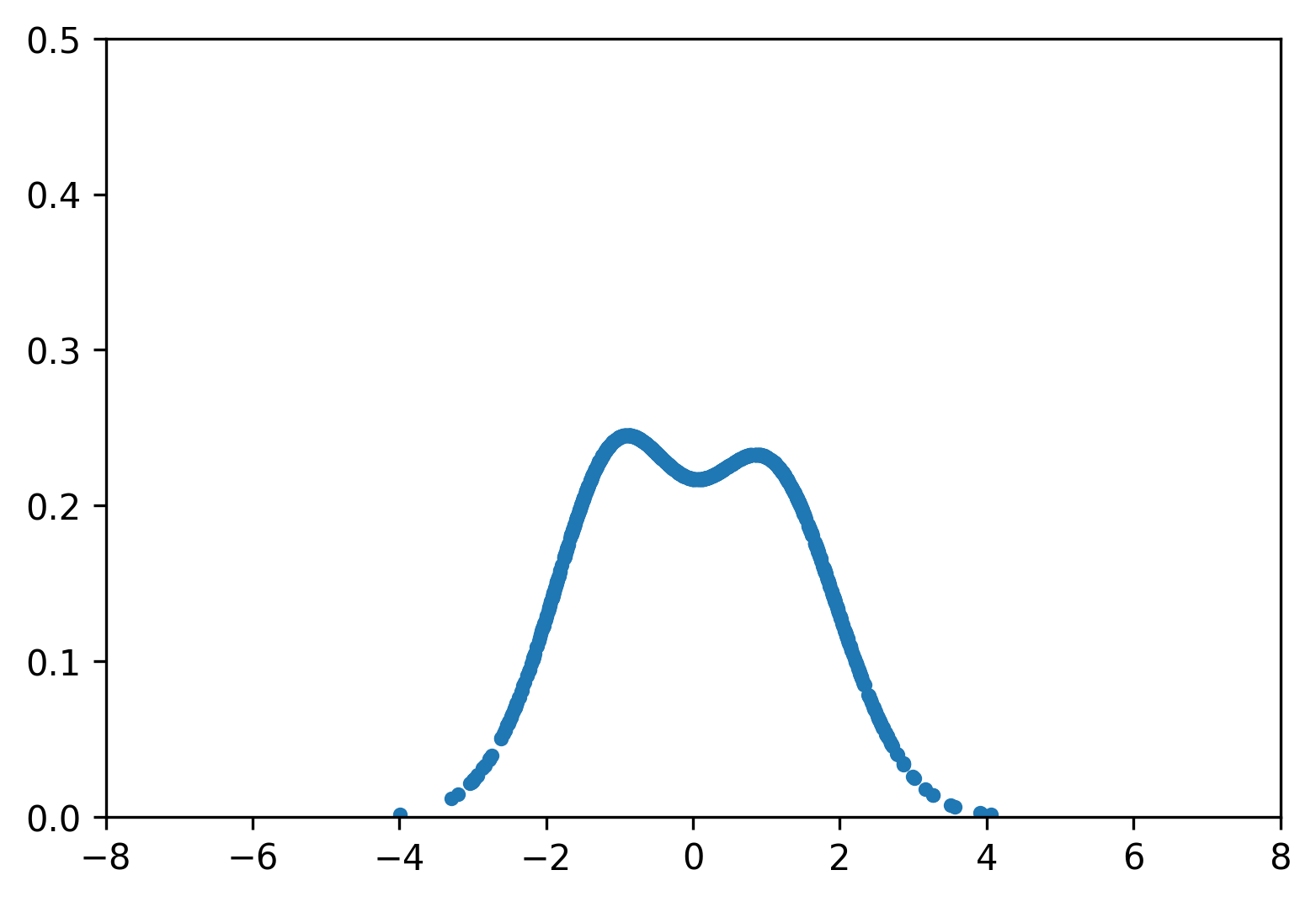}
\includegraphics[width=2.5cm]{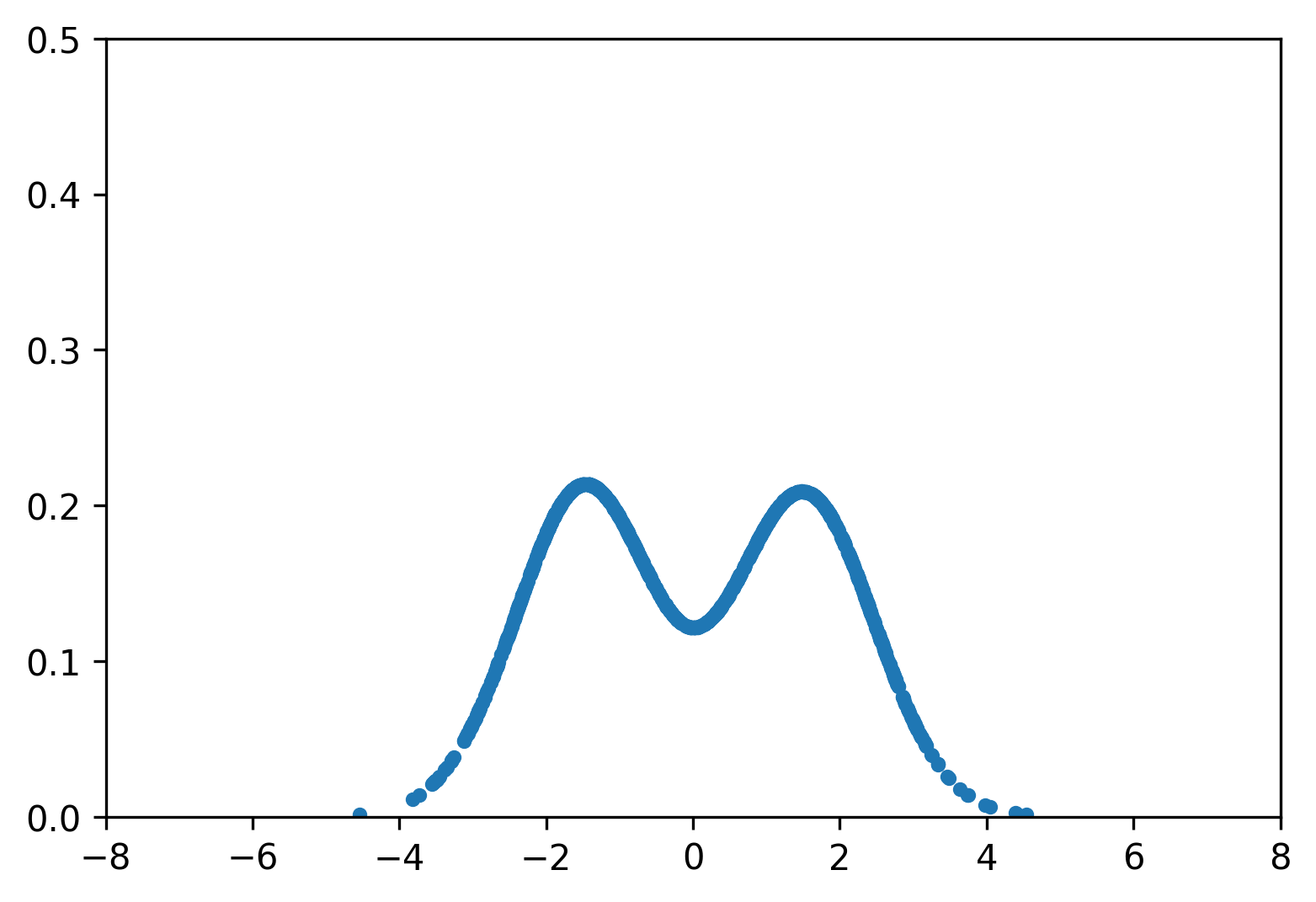}
\includegraphics[width=2.5cm]{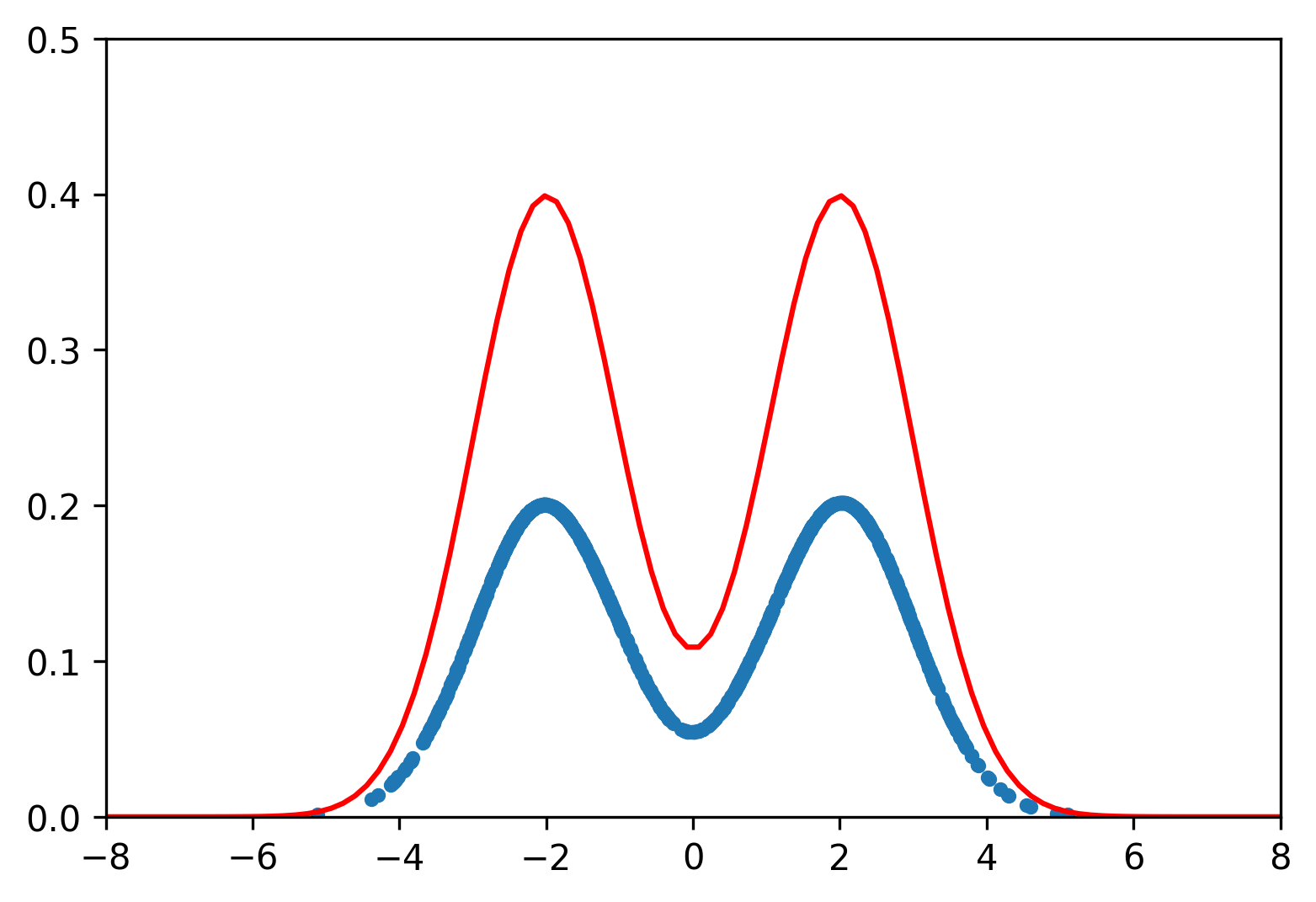}\\
\vspace{5pt}

\subfigure[$\rho_0(\boldsymbol{x}_0)$]{\includegraphics[width=2.5cm]{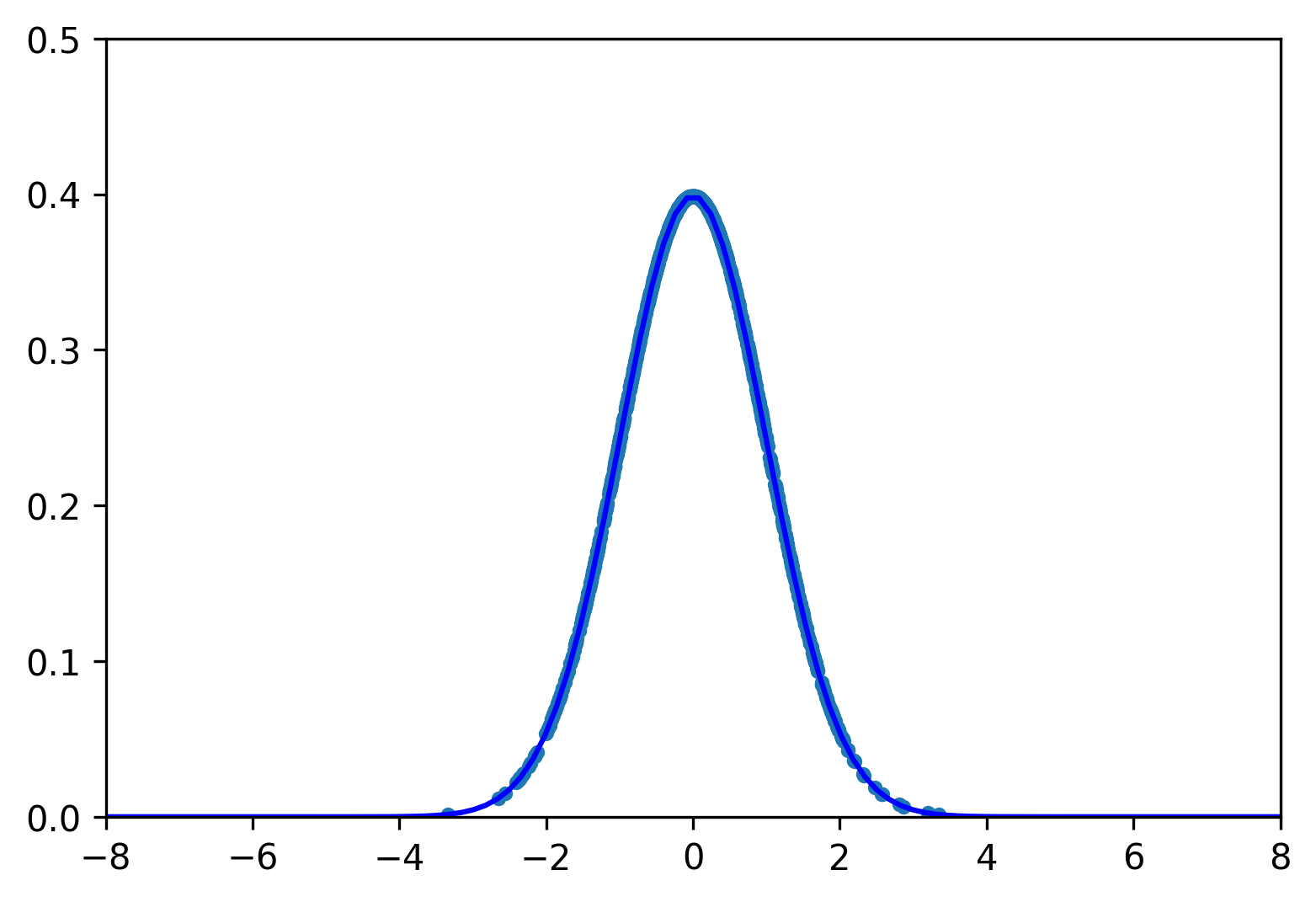}}
\subfigure[$\tilde{\rho}_{1/4}(\tilde{\boldsymbol{x}}_{1/4})$]{\includegraphics[width=2.5cm]{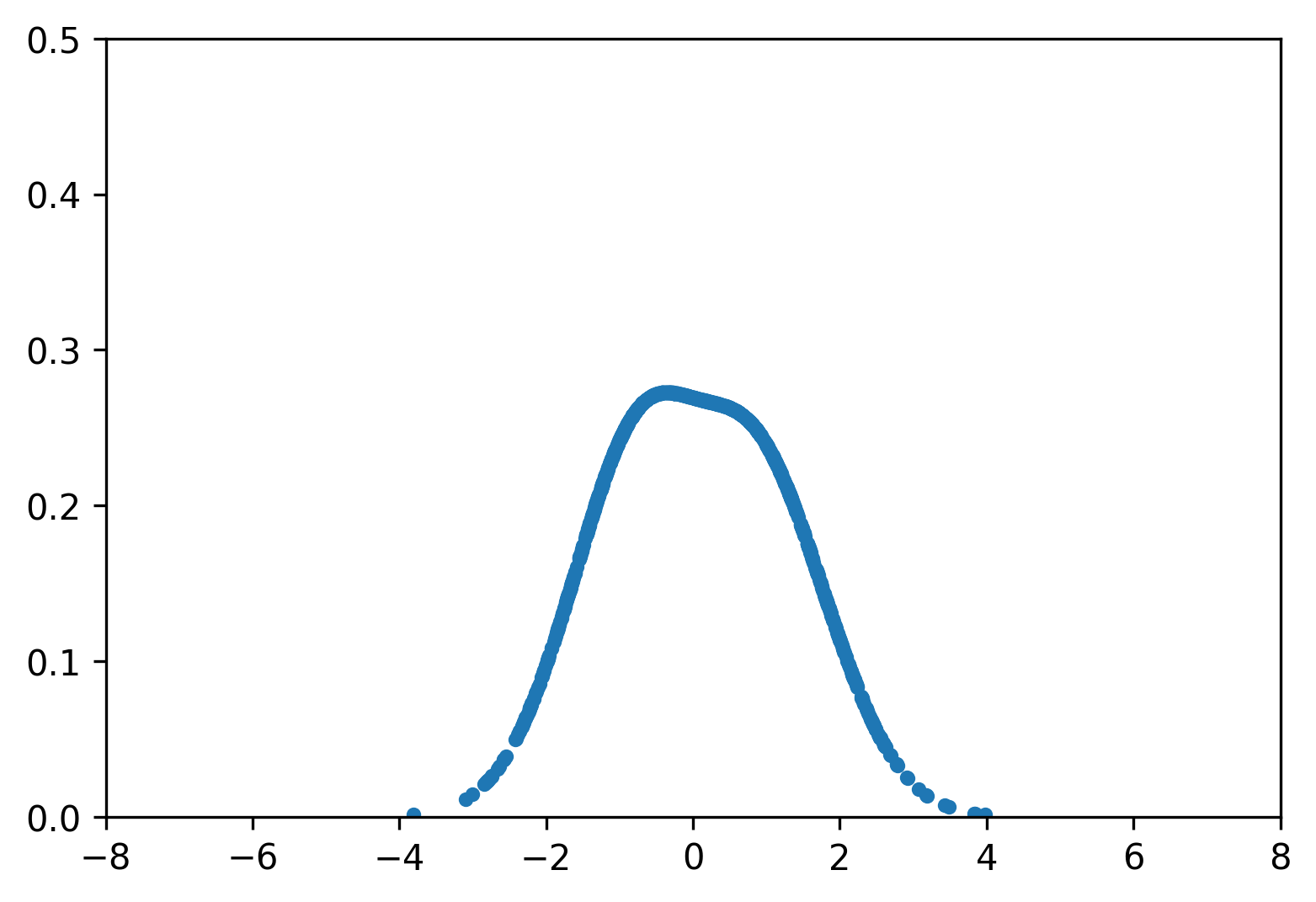}}
\subfigure[$\tilde{\rho}_{1/2}(\tilde{\boldsymbol{x}}_{1/2})$]{\includegraphics[width=2.5cm]{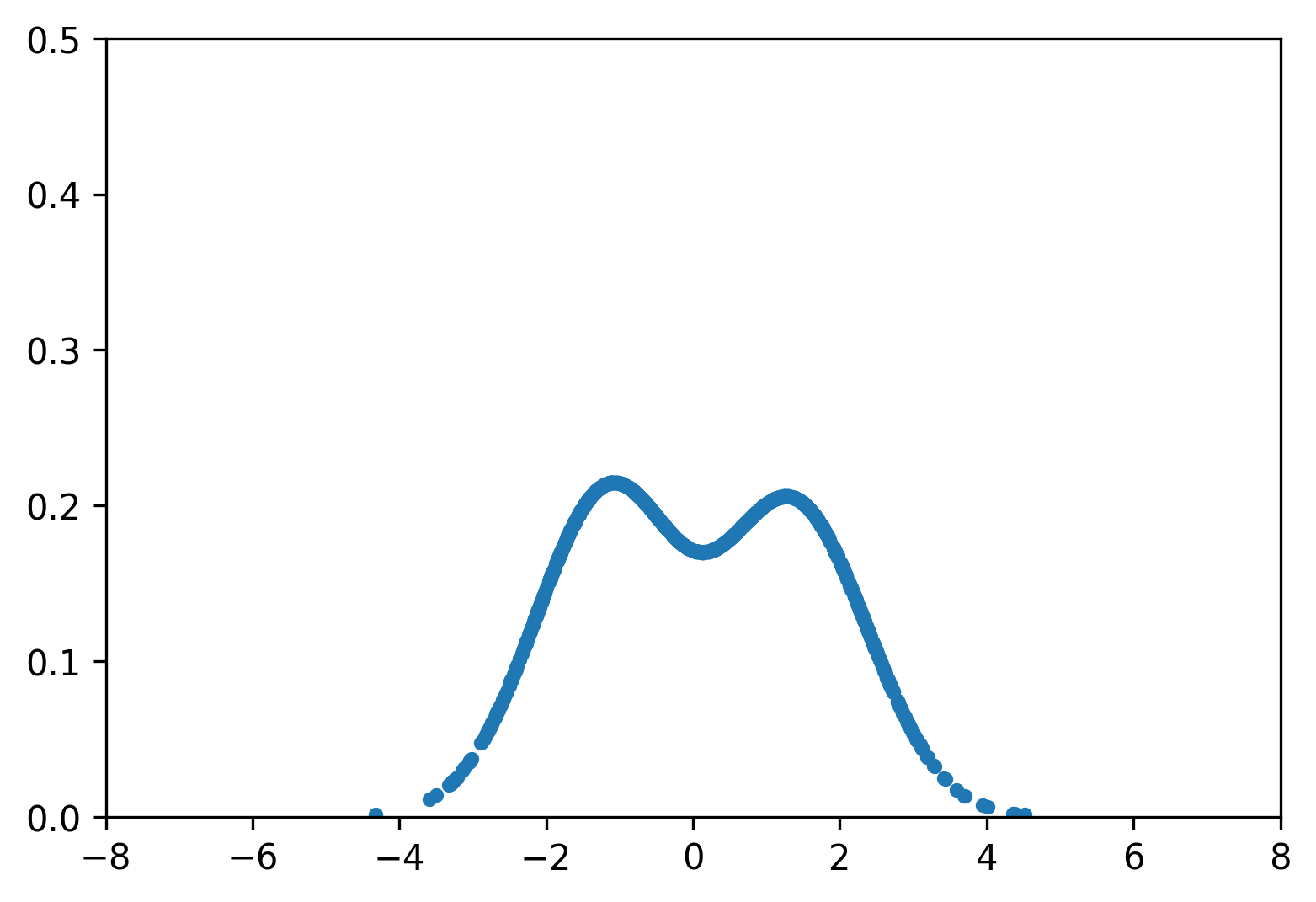}}
\subfigure[$\tilde{\rho}_{3/4}(\tilde{\boldsymbol{x}}_{3/4})$]{\includegraphics[width=2.5cm]{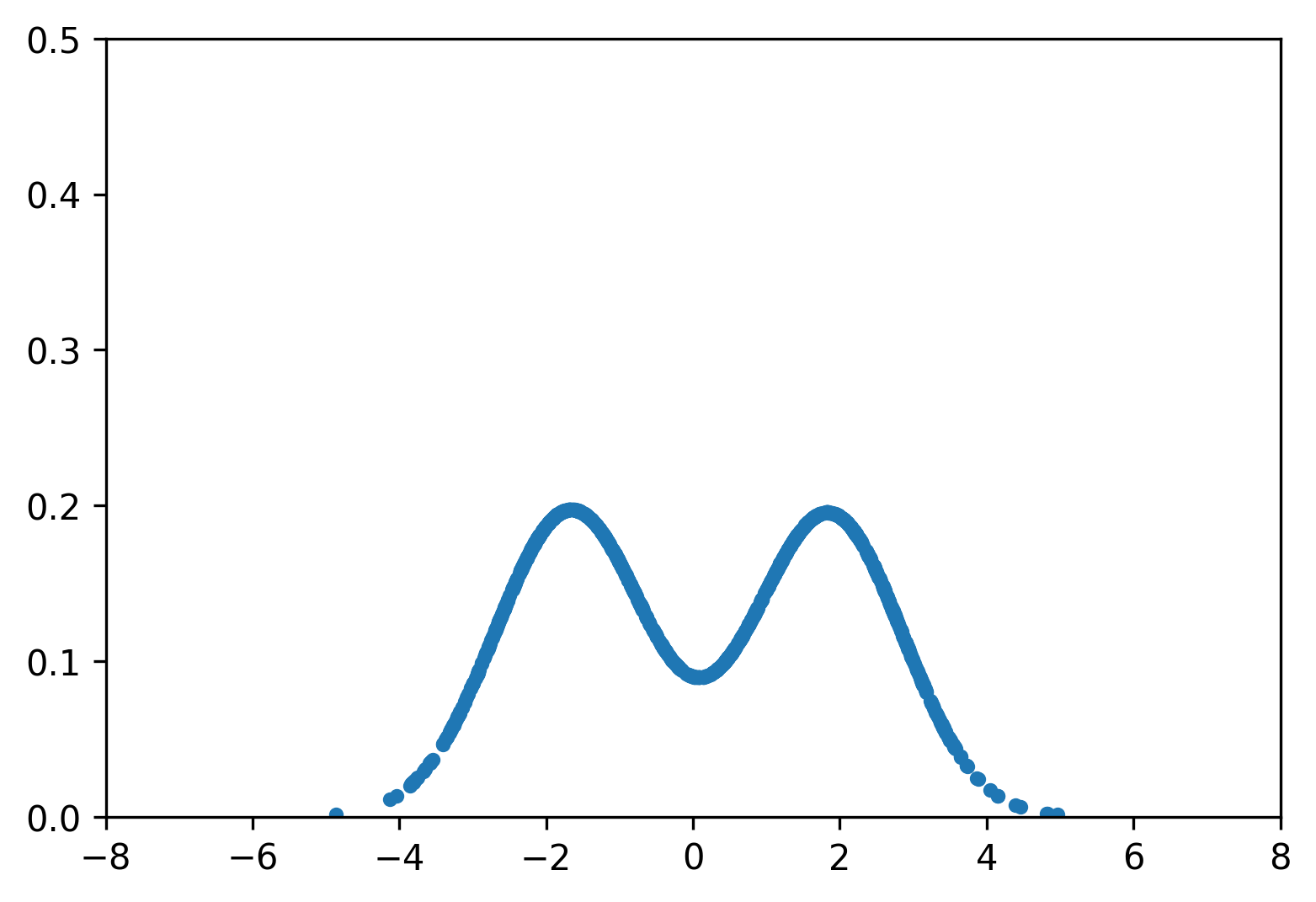}}
\subfigure[$\tilde{\rho}_{1}(\tilde{\boldsymbol{x}}_{1})$]{\includegraphics[width=2.5cm]{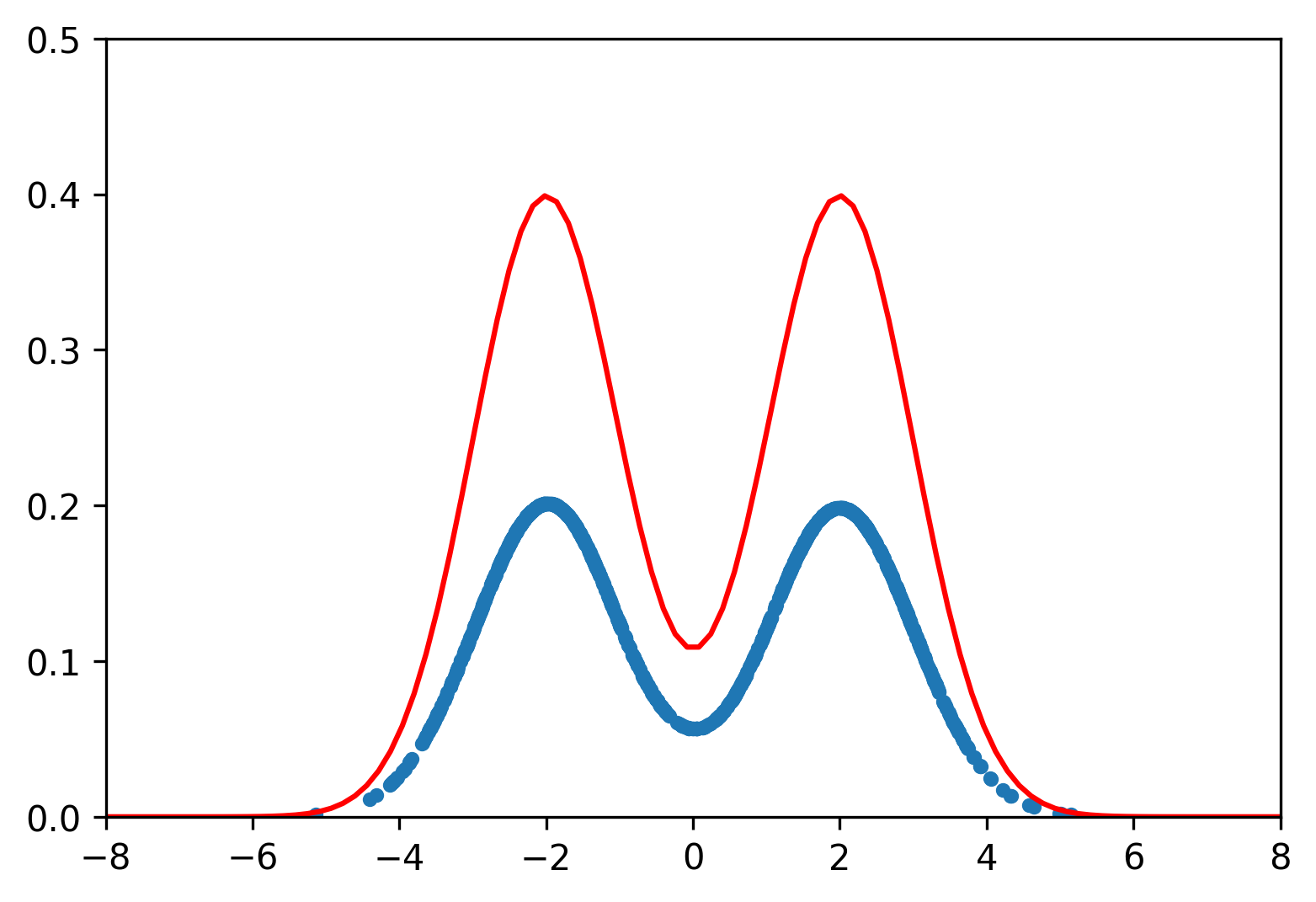}}

\caption{Top: UOT with $\alpha=10^{-2}$, middle: UOT with $\alpha=10^{-6}$, bottom: OT.}
\label{uot2}
\end{center}
\end{figure}

\begin{figure}[H]
\begin{center}
\includegraphics[width=2.5cm]{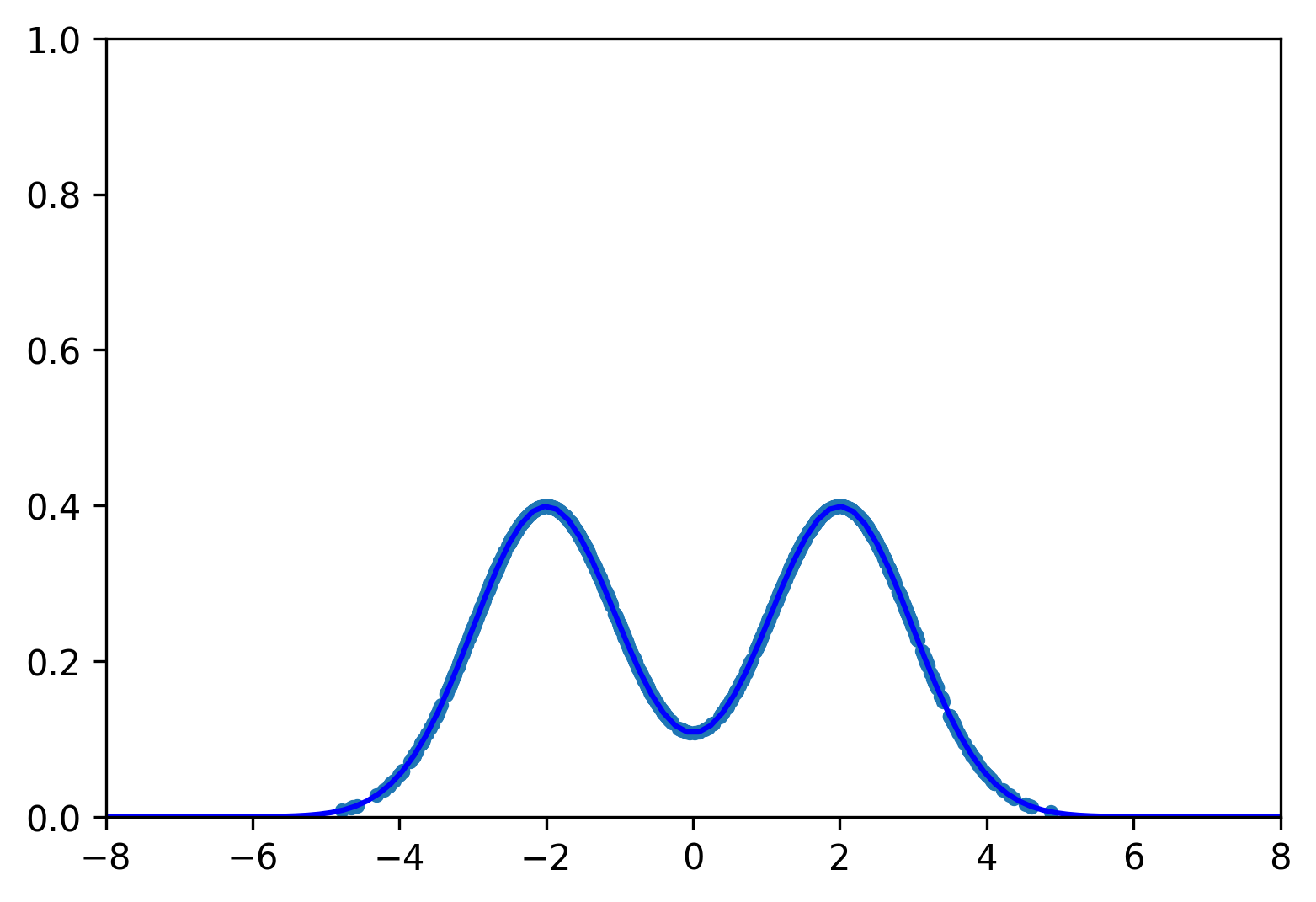}
\includegraphics[width=2.5cm]{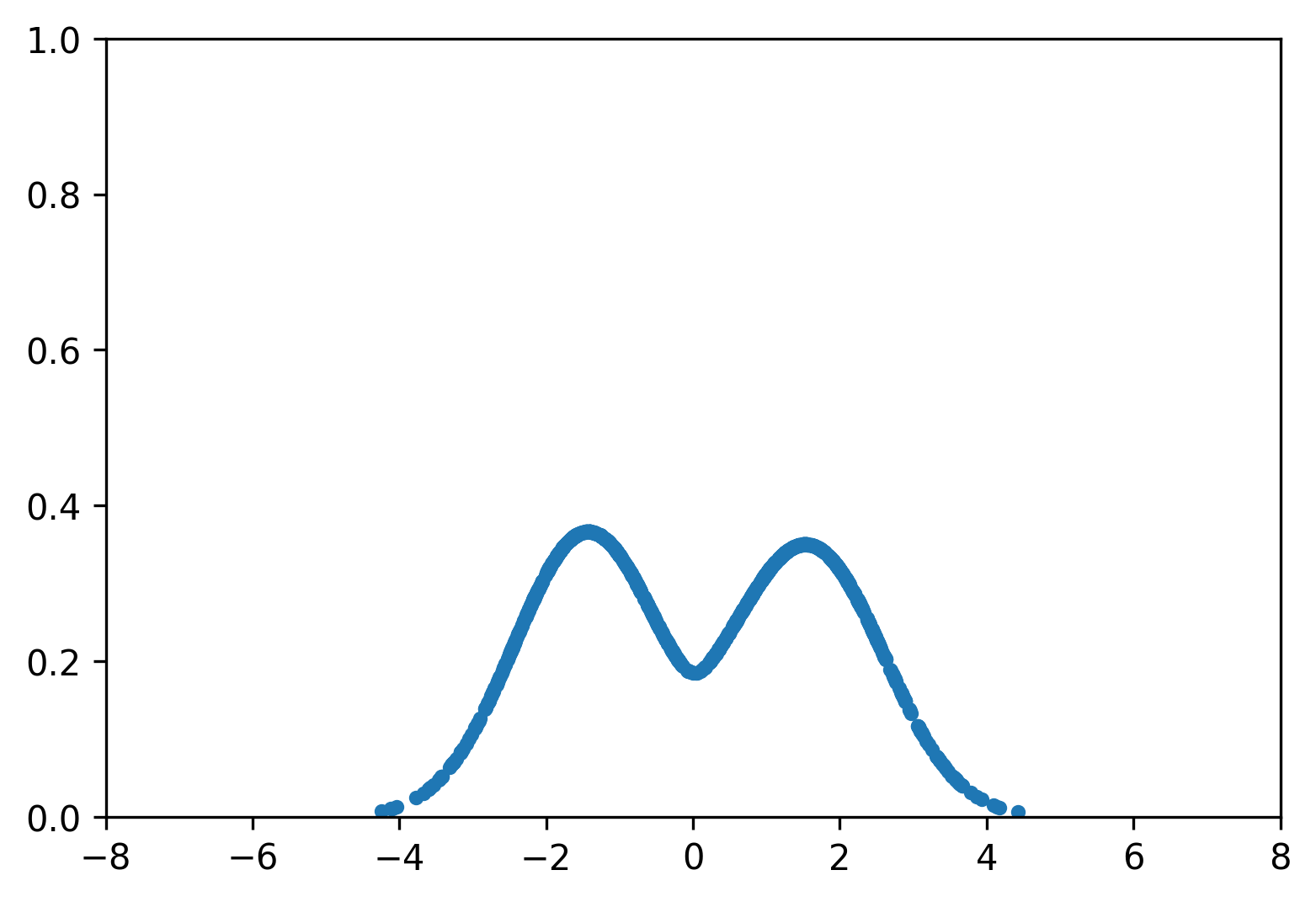}
\includegraphics[width=2.5cm]{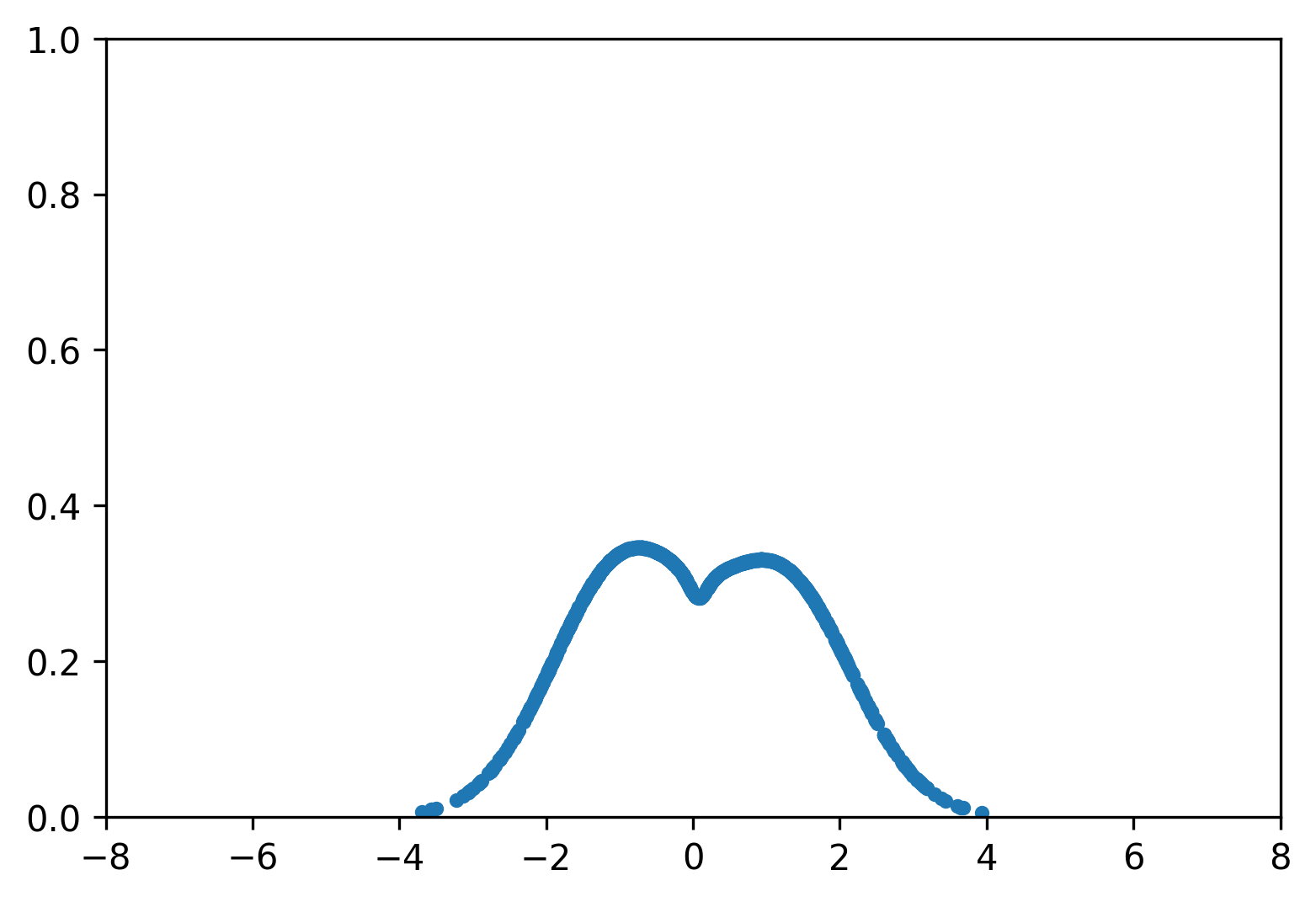}
\includegraphics[width=2.5cm]{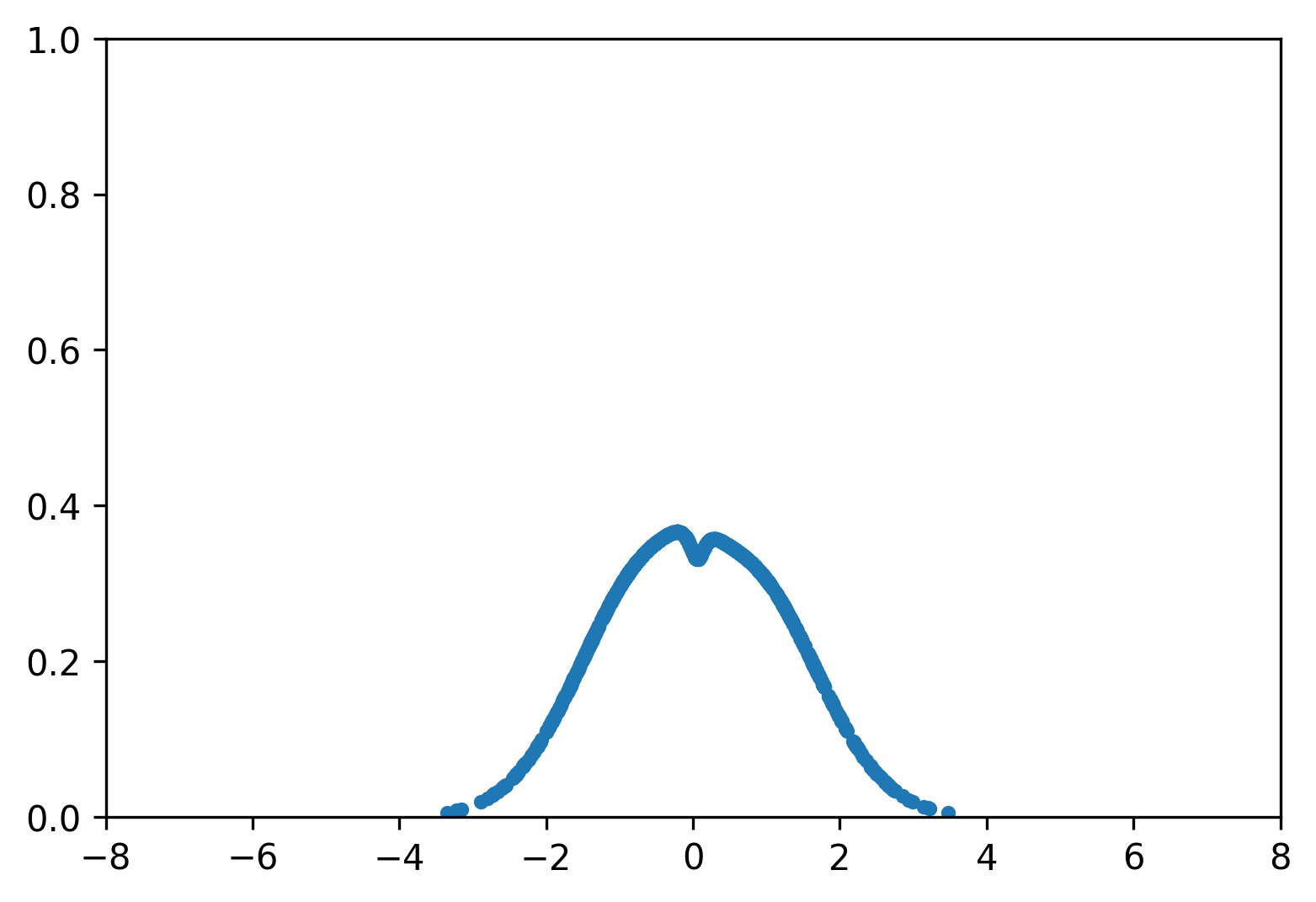}
\includegraphics[width=2.5cm]{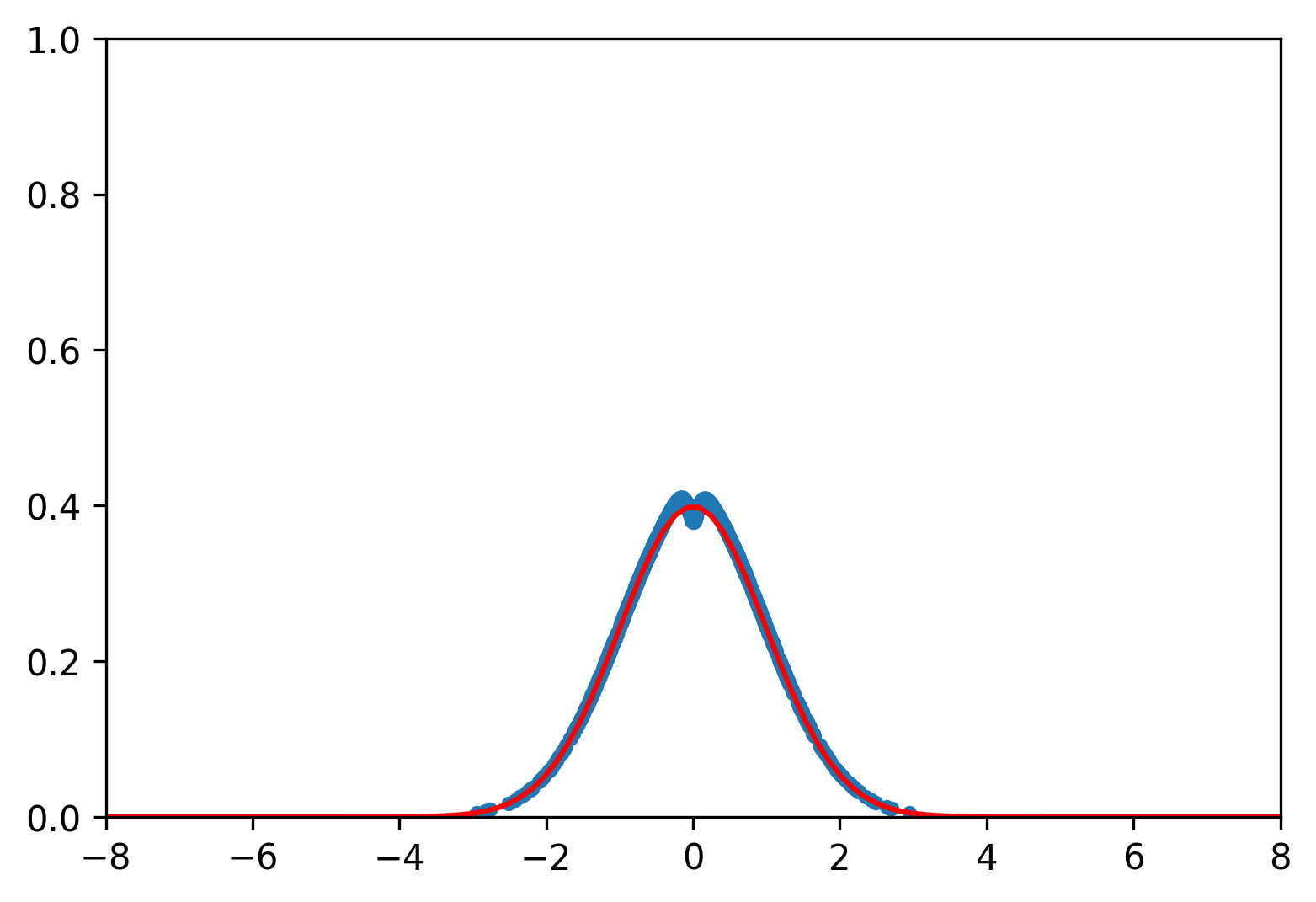}\\
\vspace{5pt}

\includegraphics[width=2.5cm]{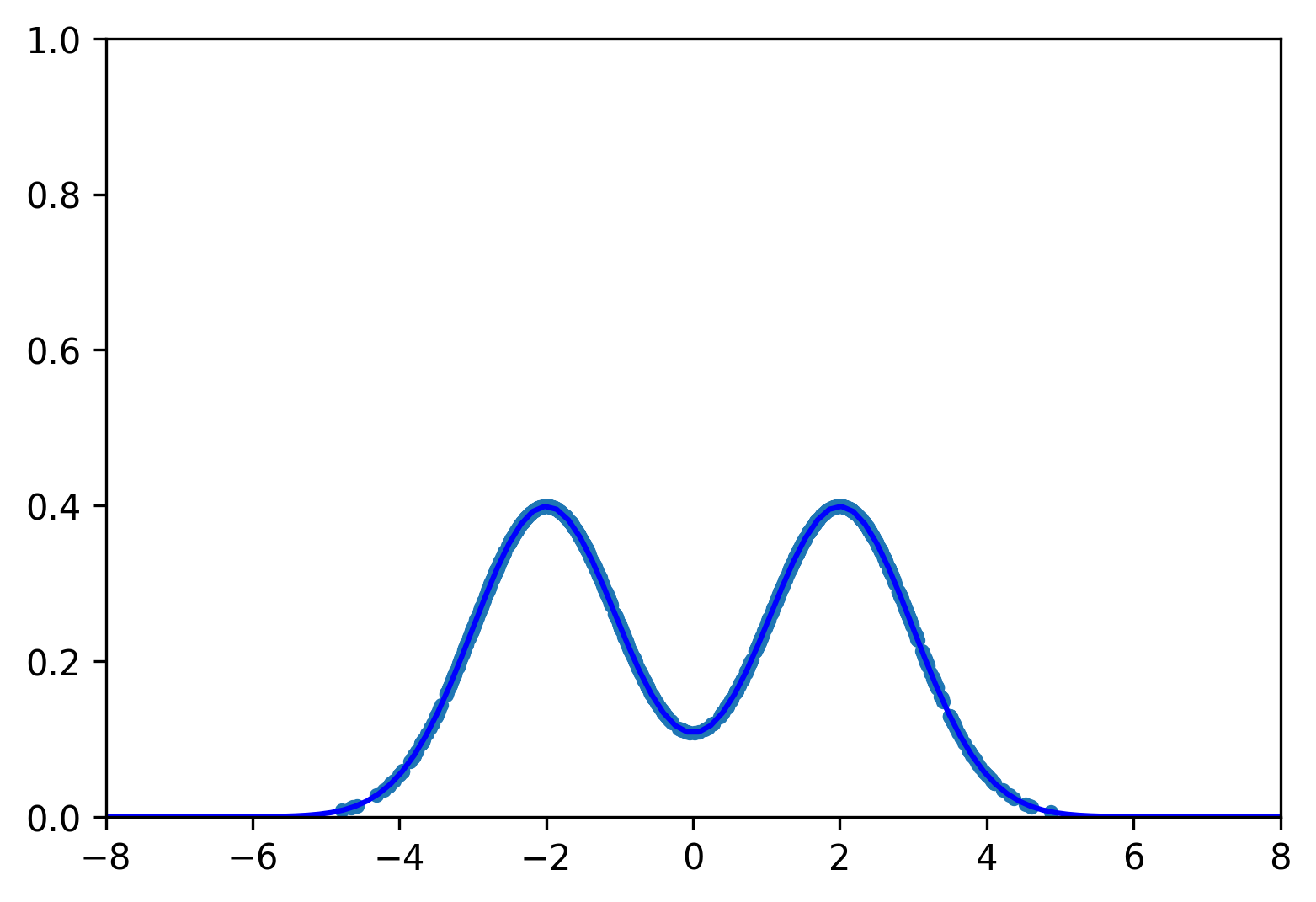}
\includegraphics[width=2.5cm]{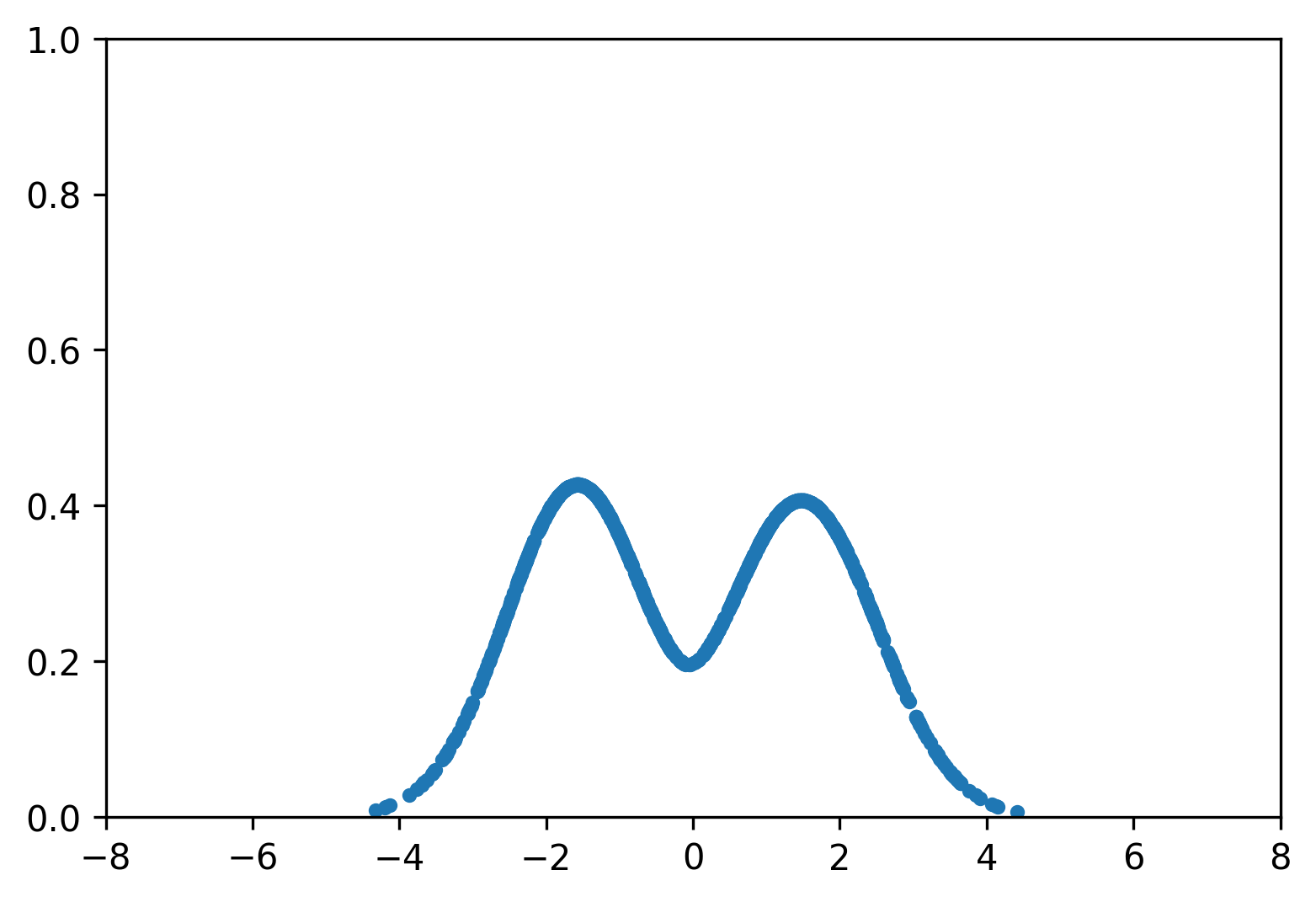}
\includegraphics[width=2.5cm]{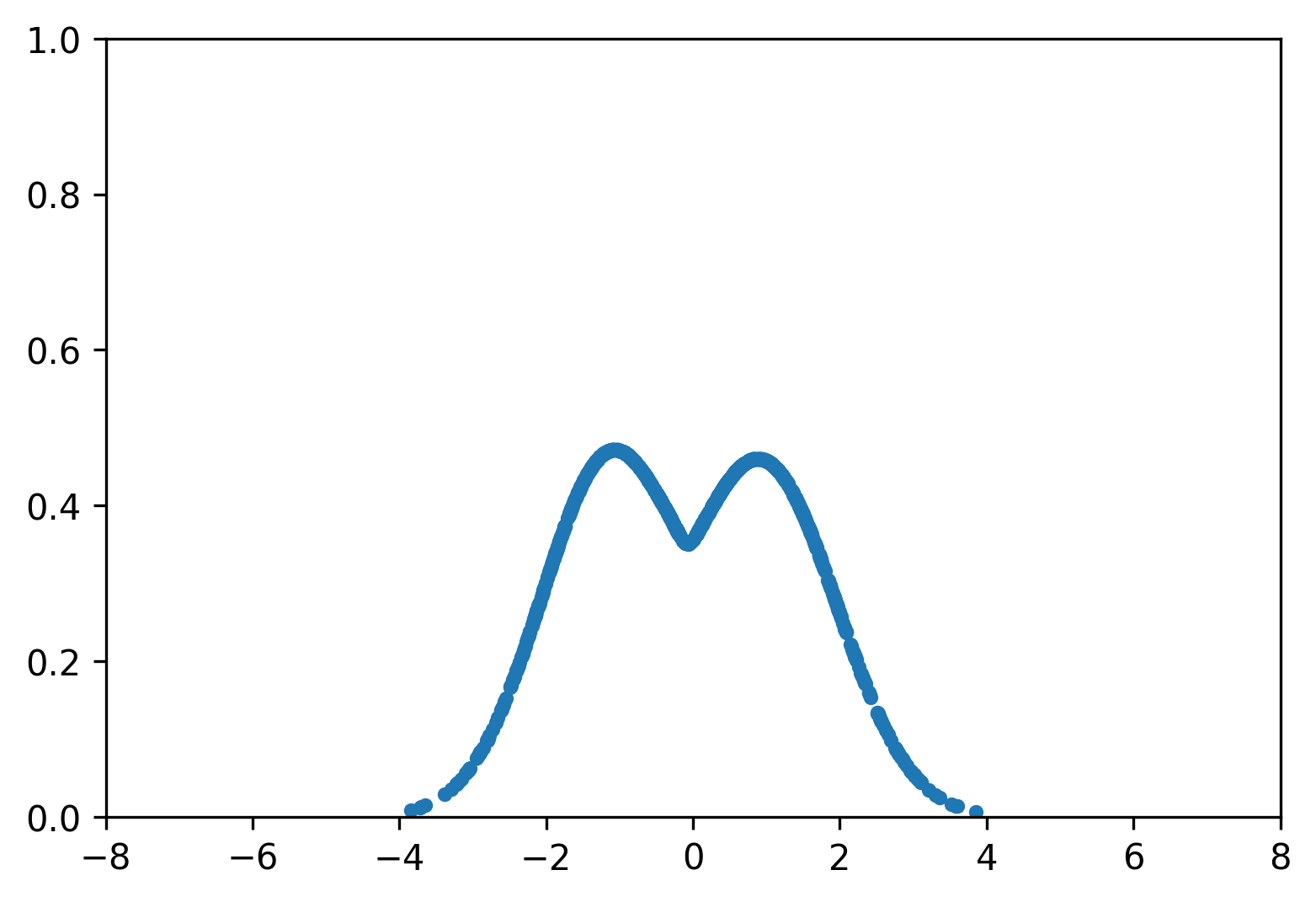}
\includegraphics[width=2.5cm]{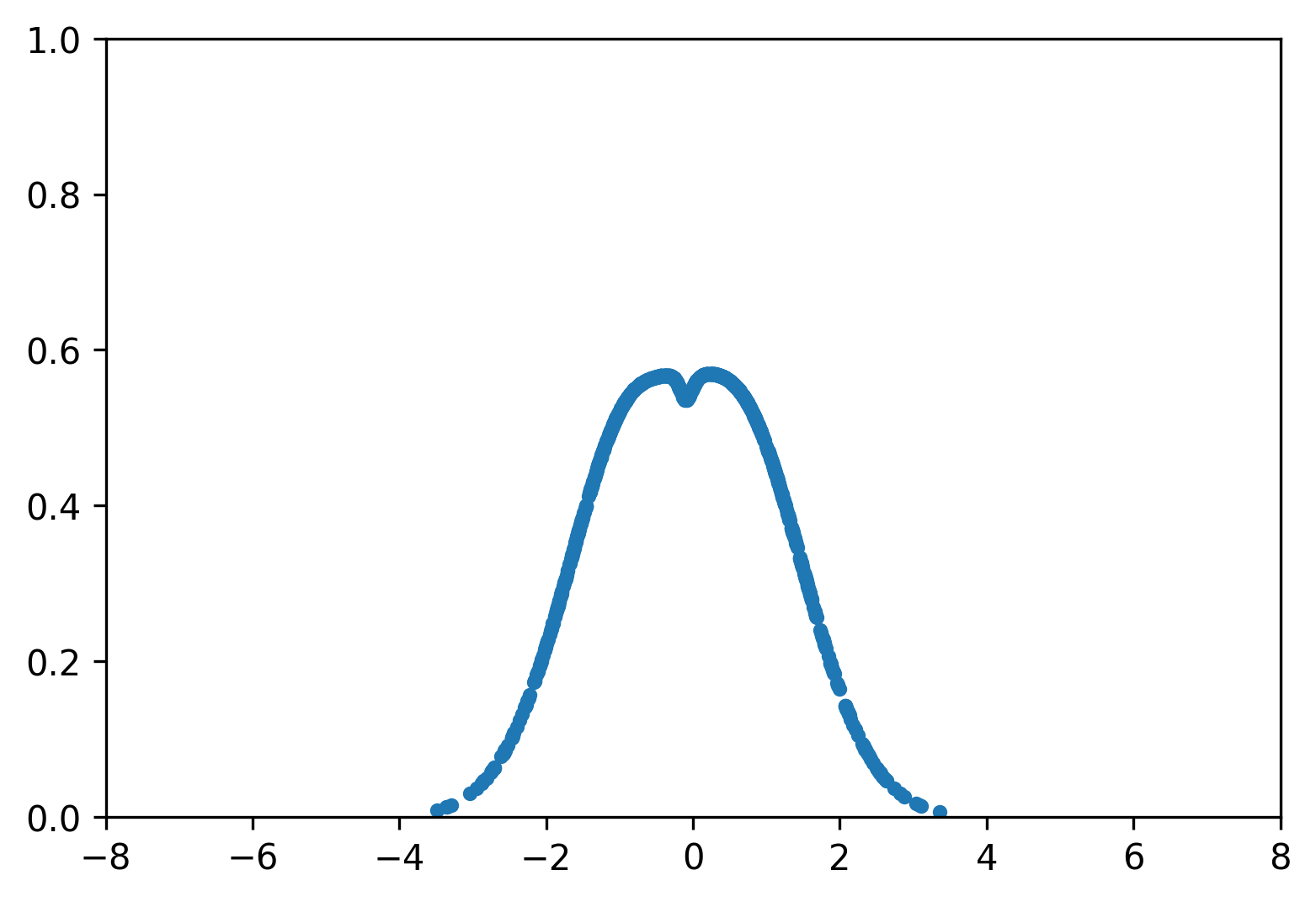}
\includegraphics[width=2.5cm]{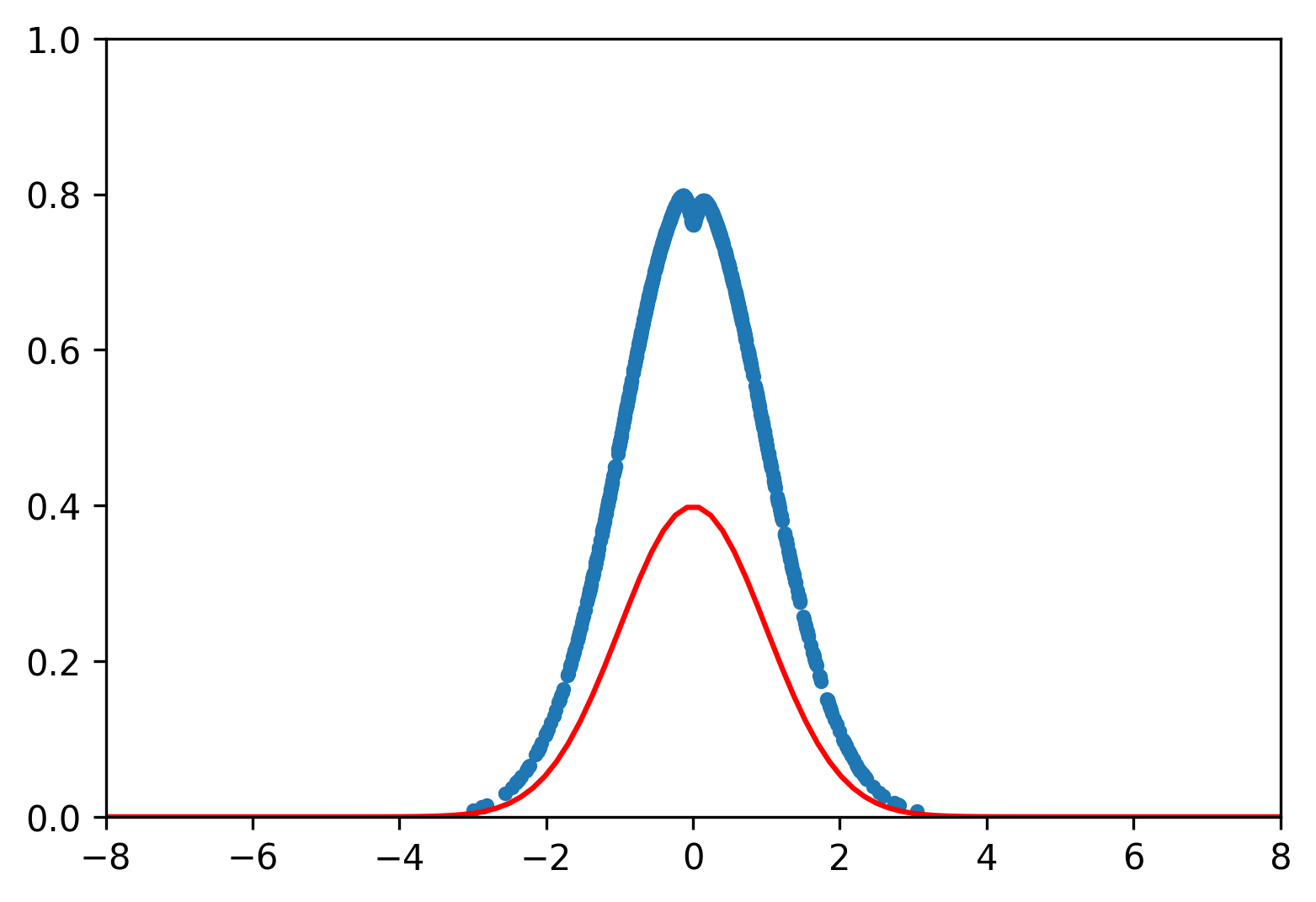}\\
\vspace{5pt}

\subfigure[$\rho_0(\boldsymbol{x}_0)$]{\includegraphics[width=2.5cm]{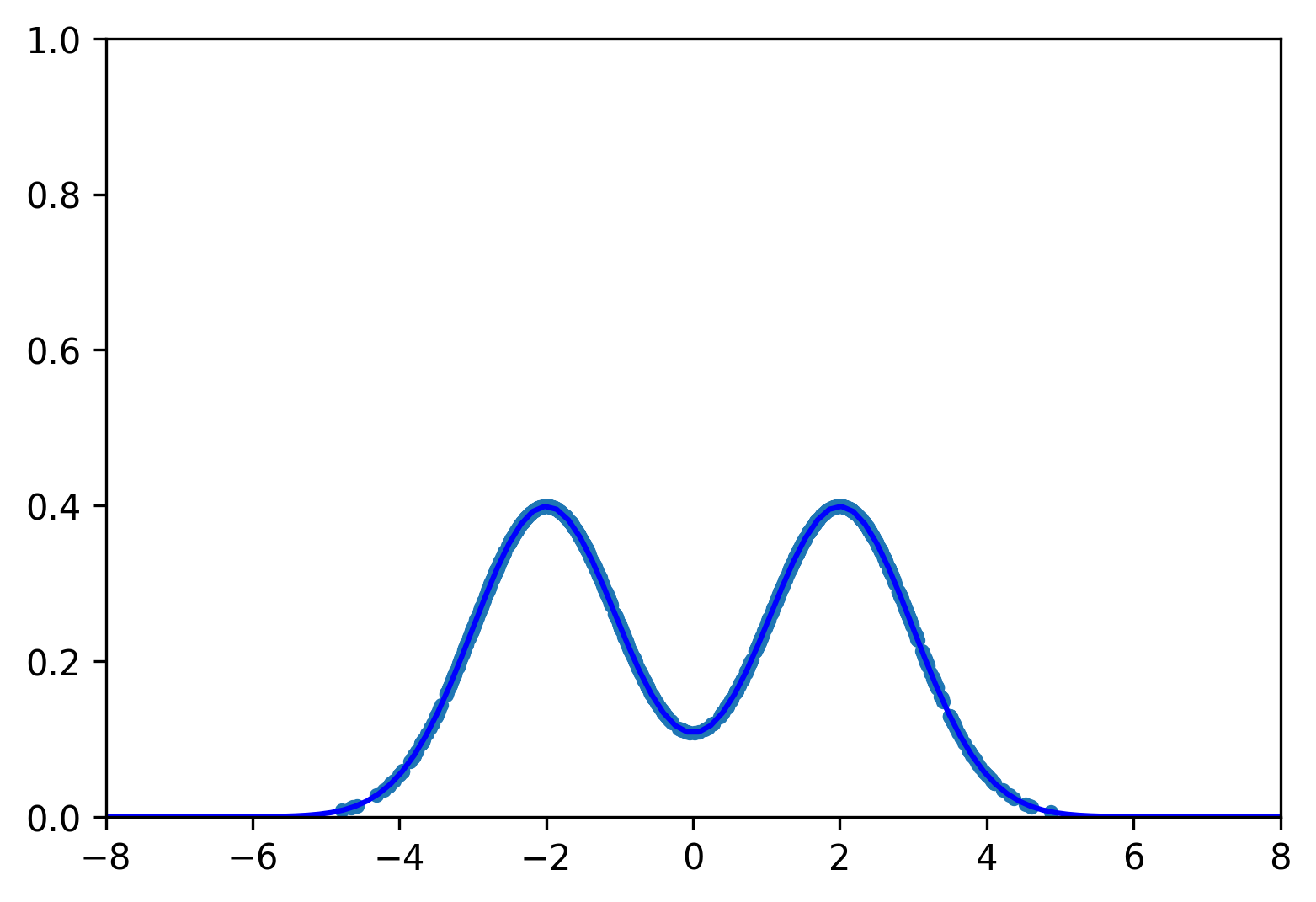}}
\subfigure[$\tilde{\rho}_{1/4}(\tilde{\boldsymbol{x}}_{1/4})$]{\includegraphics[width=2.5cm]{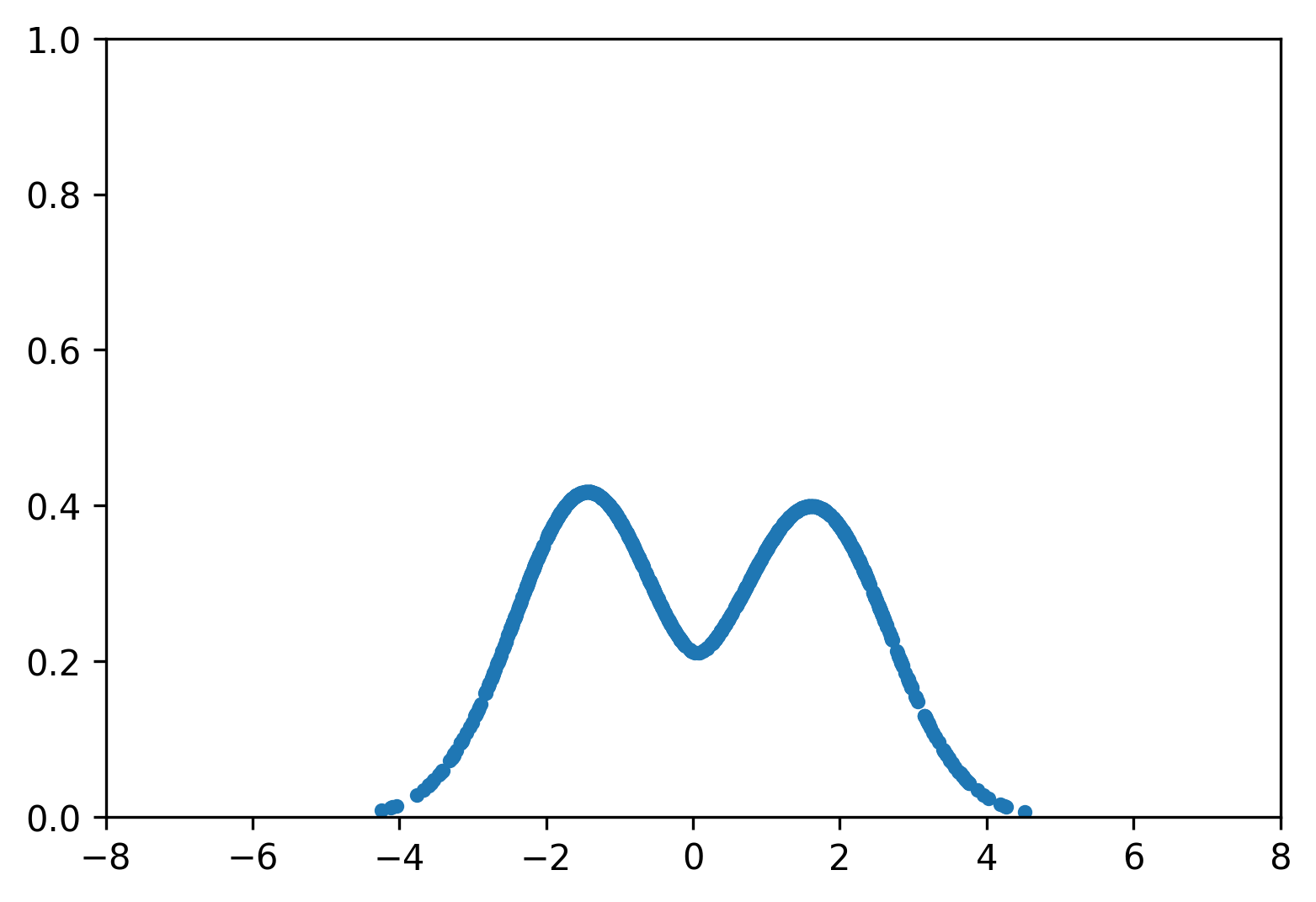}}
\subfigure[$\tilde{\rho}_{1/2}(\tilde{\boldsymbol{x}}_{1/2})$]{\includegraphics[width=2.5cm]{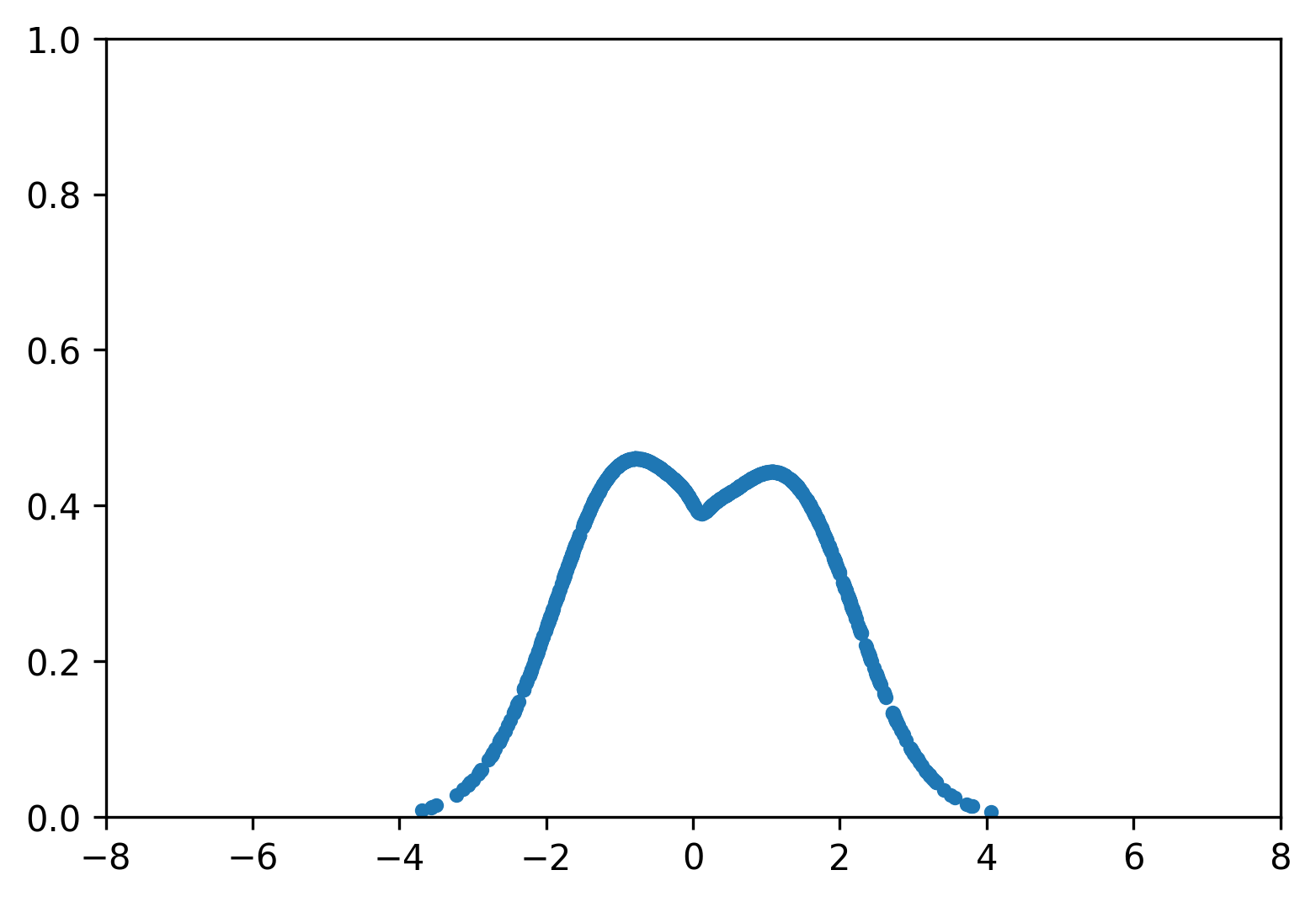}}
\subfigure[$\tilde{\rho}_{3/4}(\tilde{\boldsymbol{x}}_{3/4})$]{\includegraphics[width=2.5cm]{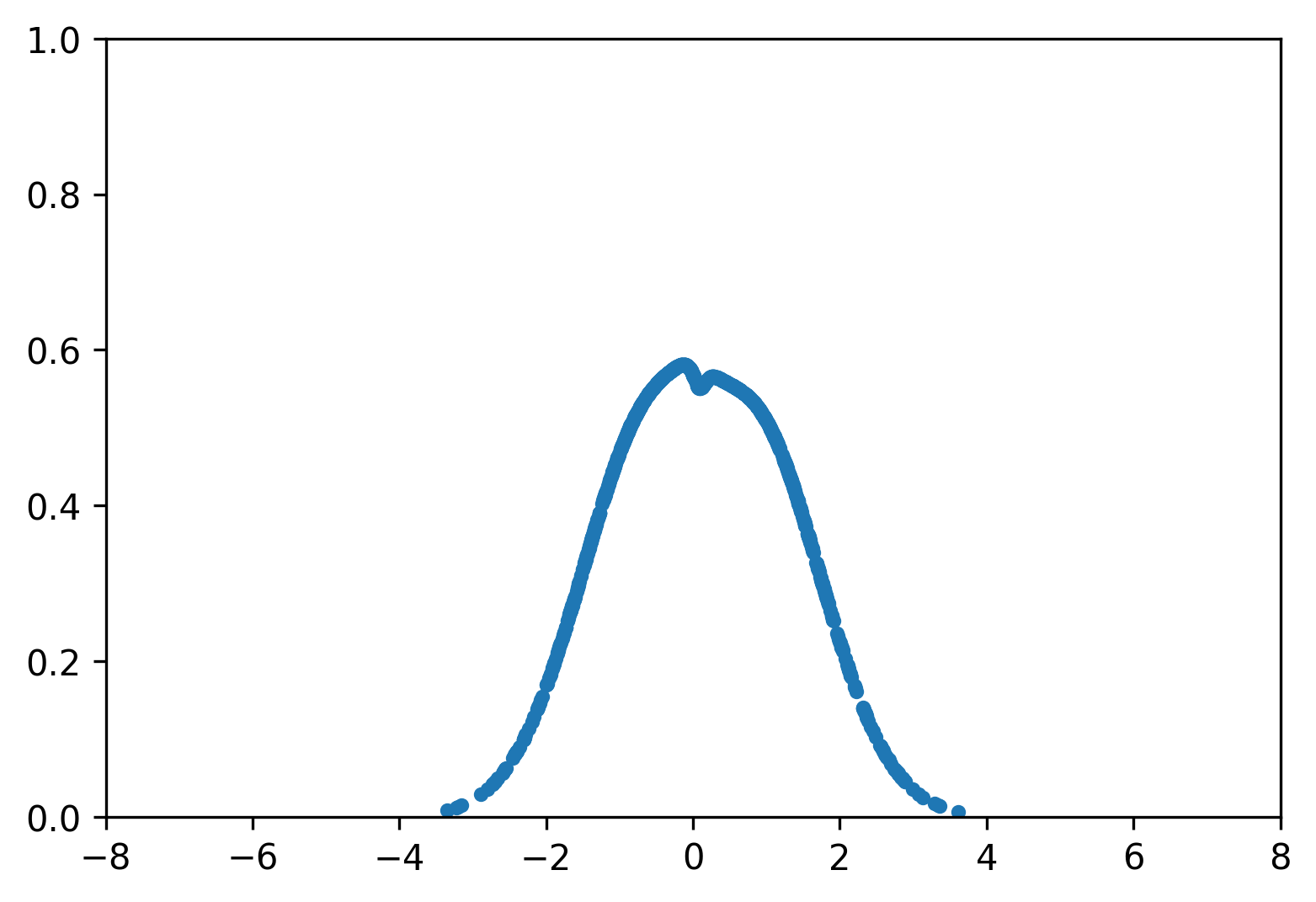}}
\subfigure[$\tilde{\rho}_{1}(\tilde{\boldsymbol{x}}_{1})$]{\includegraphics[width=2.5cm]{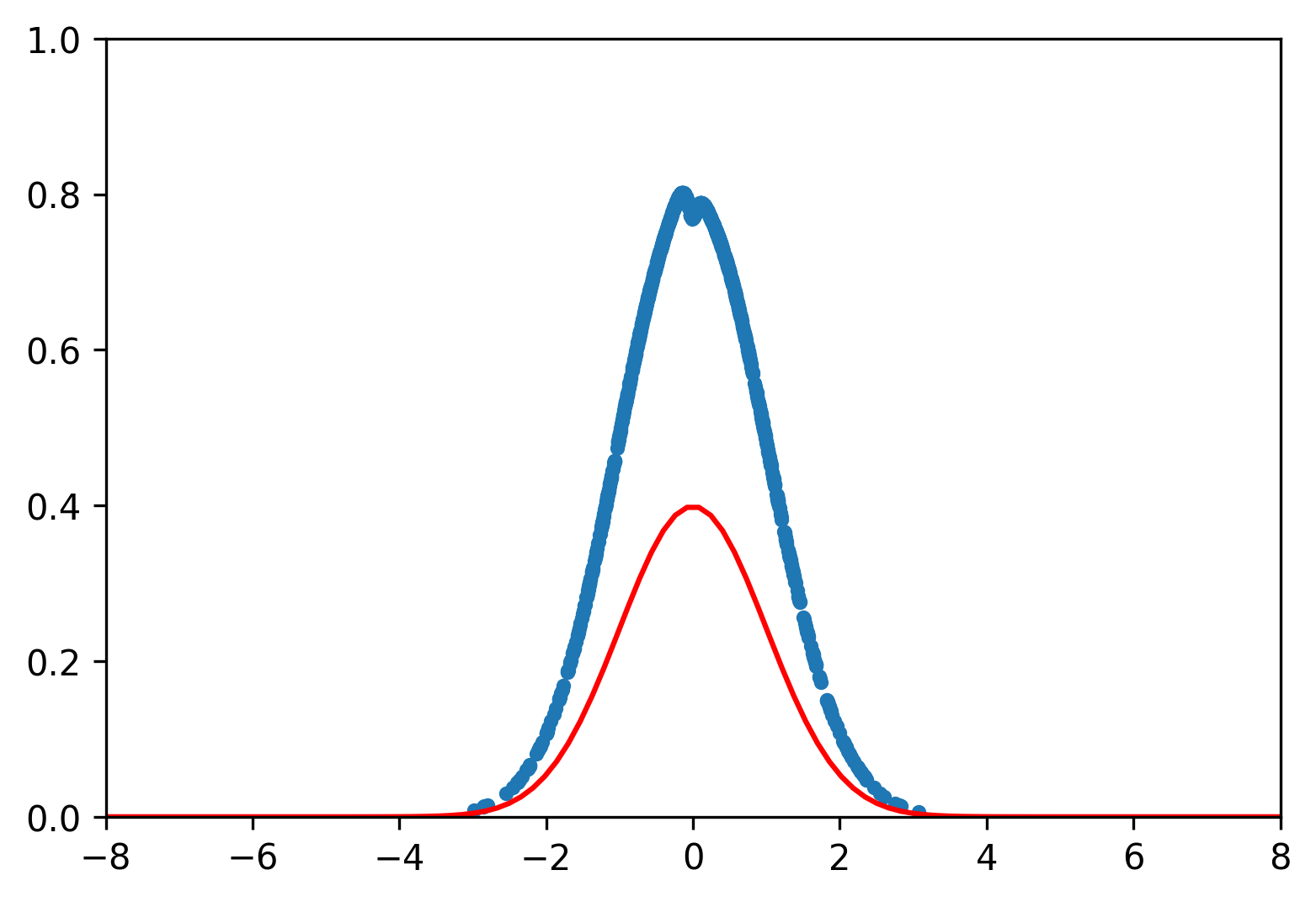}}

\caption{Top: UOT with $\alpha=10^{-2}$, middle: UOT with $\alpha=10^{-6}$, bottom: OT.}
\label{uot3}
\end{center}
\end{figure}

Figure \ref{uot1} - Figure \ref{uot3} numerically illustrate the density evolution using the proposed method, moving from the initial density (shown in blue in (a)) to the target density (shown in red in (e)) at one dimension ($d=1$). 
The first row illustrates the transportation process for UOT with $\alpha=10^{-2}$, while the second row shows the process for UOT with $\alpha=10^{-6}$.
Lastly, the third row presents the transportation process for OT.
We can clearly see that when the unbalanced coefficient $\alpha$ is sufficiently small, UOT performs similarly to OT, indicating a preference for transporting with equal mass.

\subsection{Gaussian Examples}
In this section, we consider the UOT problem on several typical Gaussian examples with different means and covariance matrices. Here, $\rho_G(\cdot, \mathbf{m}, \mathbf{\Sigma})$ represents the probability density function of a $d$-dimensional Gaussian distribution with a mean vector $\mathbf{m} \in \mathbb{R}^d$ and covariance matrix $\mathbf{\Sigma} \in \mathbb{R}^{d \times d}$, i.e., 
$$\rho_{G}(\boldsymbol{\boldsymbol{x}}, \mathbf{m}, \mathbf{\Sigma})=\frac{1}{(2\pi)^{d/2}|\mathbf{\Sigma}|^{1/2}}\mathbf{exp}\big({-\frac{1}{2}(\boldsymbol{\boldsymbol{x}}-\mathbf{m})^T\mathbf{\Sigma}^{-1}(\boldsymbol{\boldsymbol{x}}-\mathbf{m})}\big).$$

Below, we first provide specific settings for four unbalanced synthetic tests in $d=1$.

\begin{itemize}\label{list1}
    \item Test 1: $\rho_0(x)=\rho_G(x,0,1)$, $\rho_1(x)=2\rho_G(x,0,1)$;
    \item Test 2: $\rho_0(x)=\rho_G(x,0,1)$, $\rho_1(x)=\frac{1}{2}\rho_G(x,0,1)$;
    \item Test 3: $\rho_0(x)=\rho_G(x,0,1)$, $\rho_1(x)=2\rho_G(x,4,1)$;
    \item Test 4: $\rho_0(x)=\rho_G(x,0,1)$, $\rho_1(x)=\frac{1}{2}\rho_G(x,4,1)$.
\end{itemize}

\begin{figure}[H]
\begin{center}
\includegraphics[width=2.5cm]{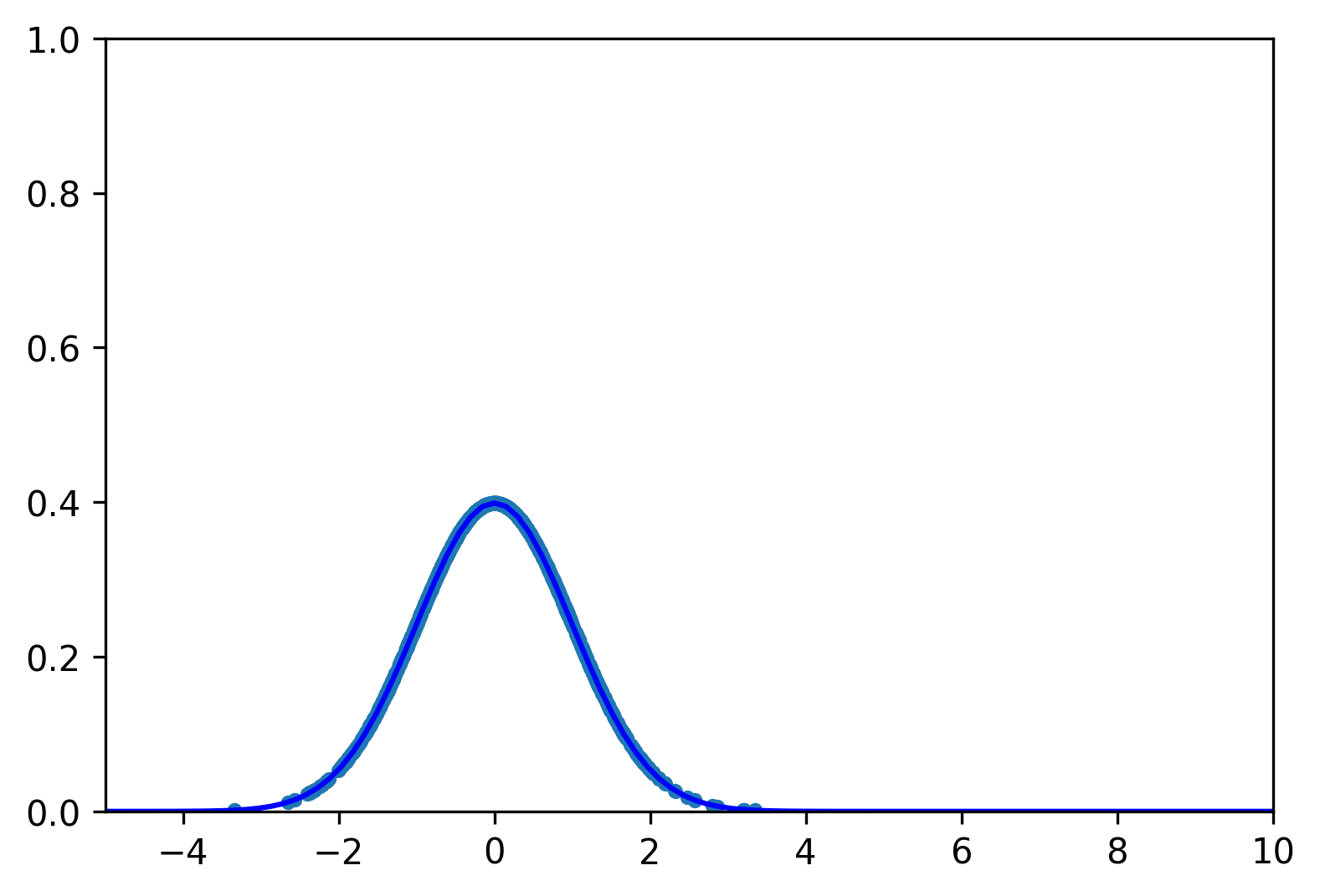}
\includegraphics[width=2.5cm]{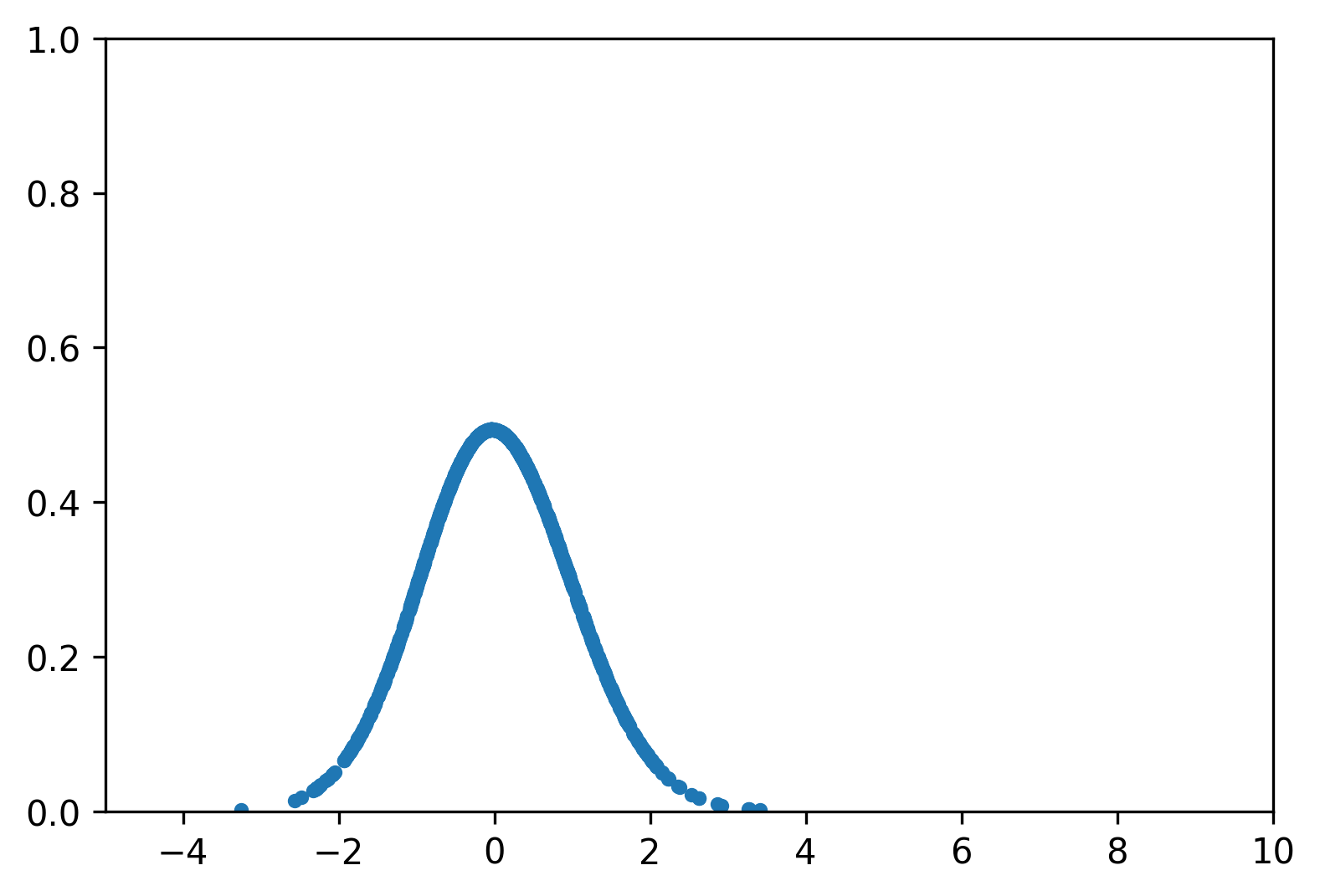}
\includegraphics[width=2.5cm]{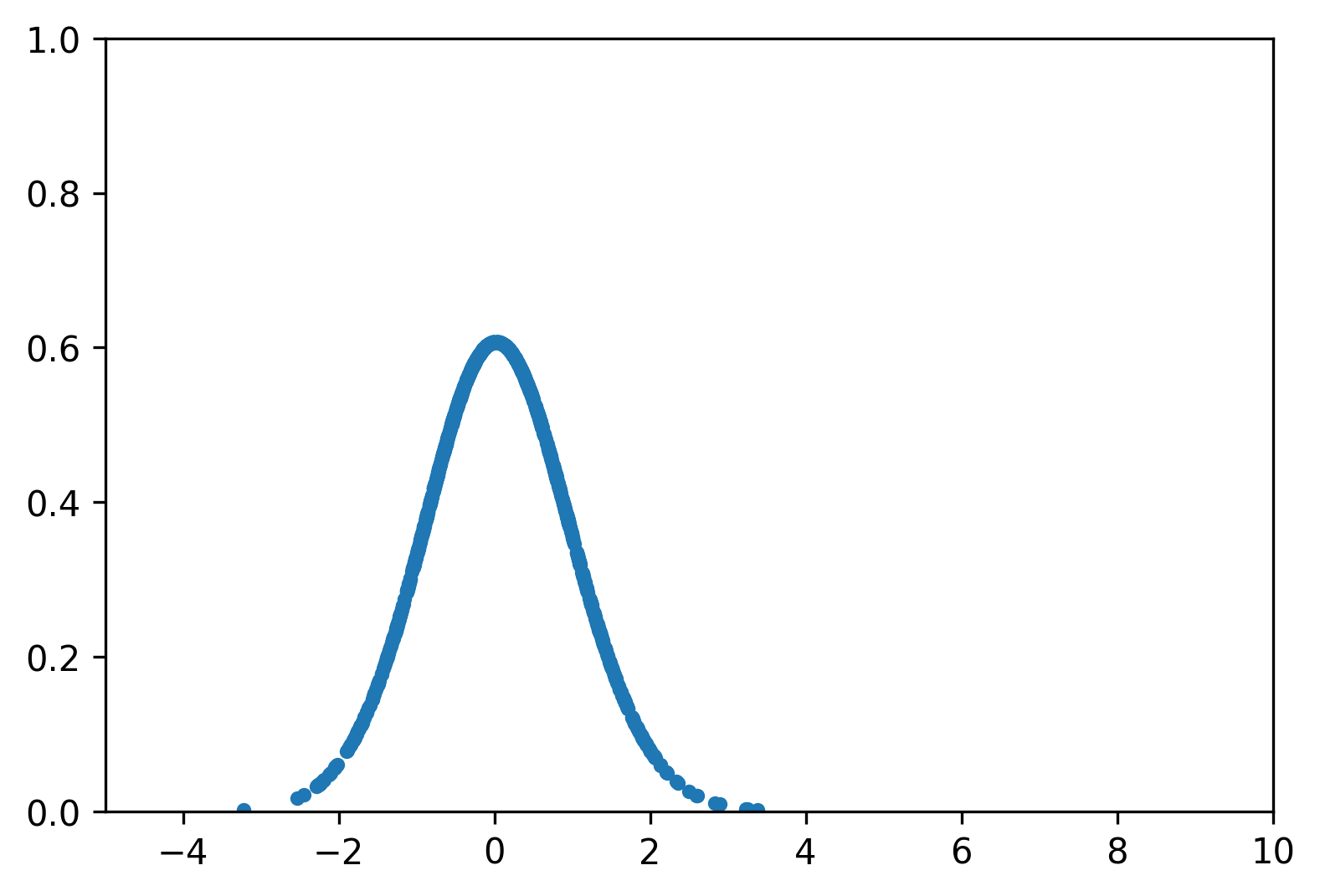}
\includegraphics[width=2.5cm]{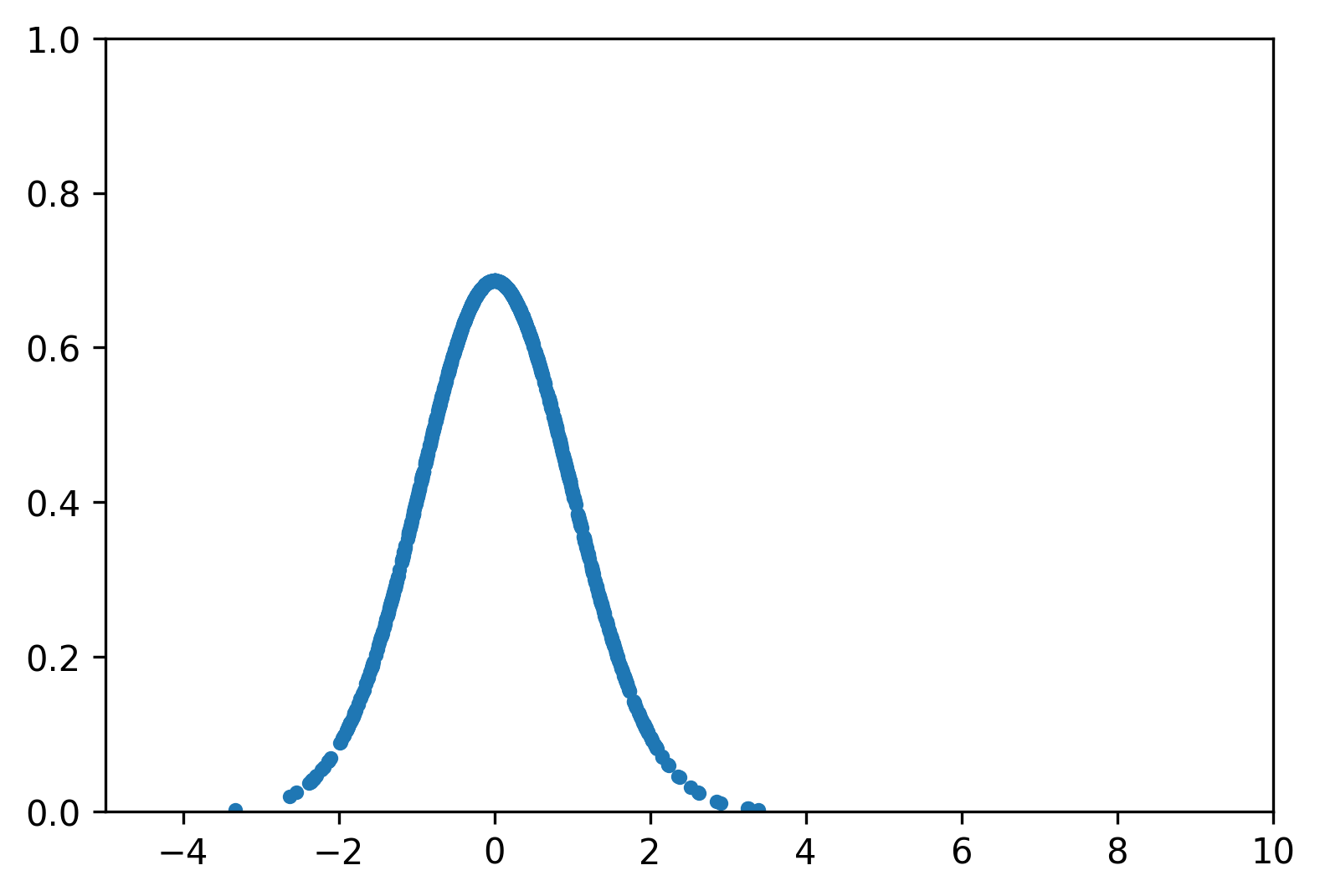}
\includegraphics[width=2.5cm]{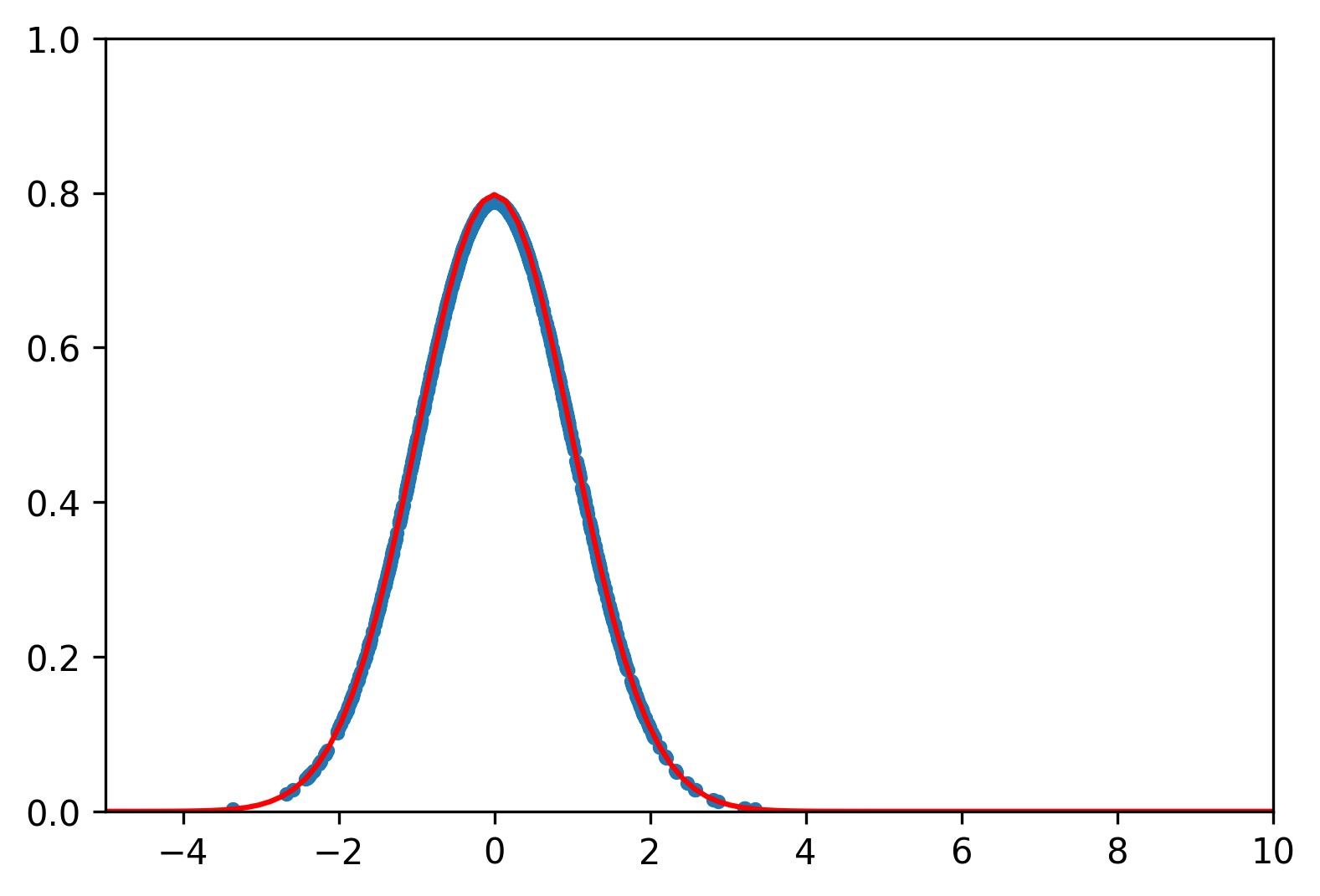}\\
\vspace{5pt}

\includegraphics[width=2.5cm]{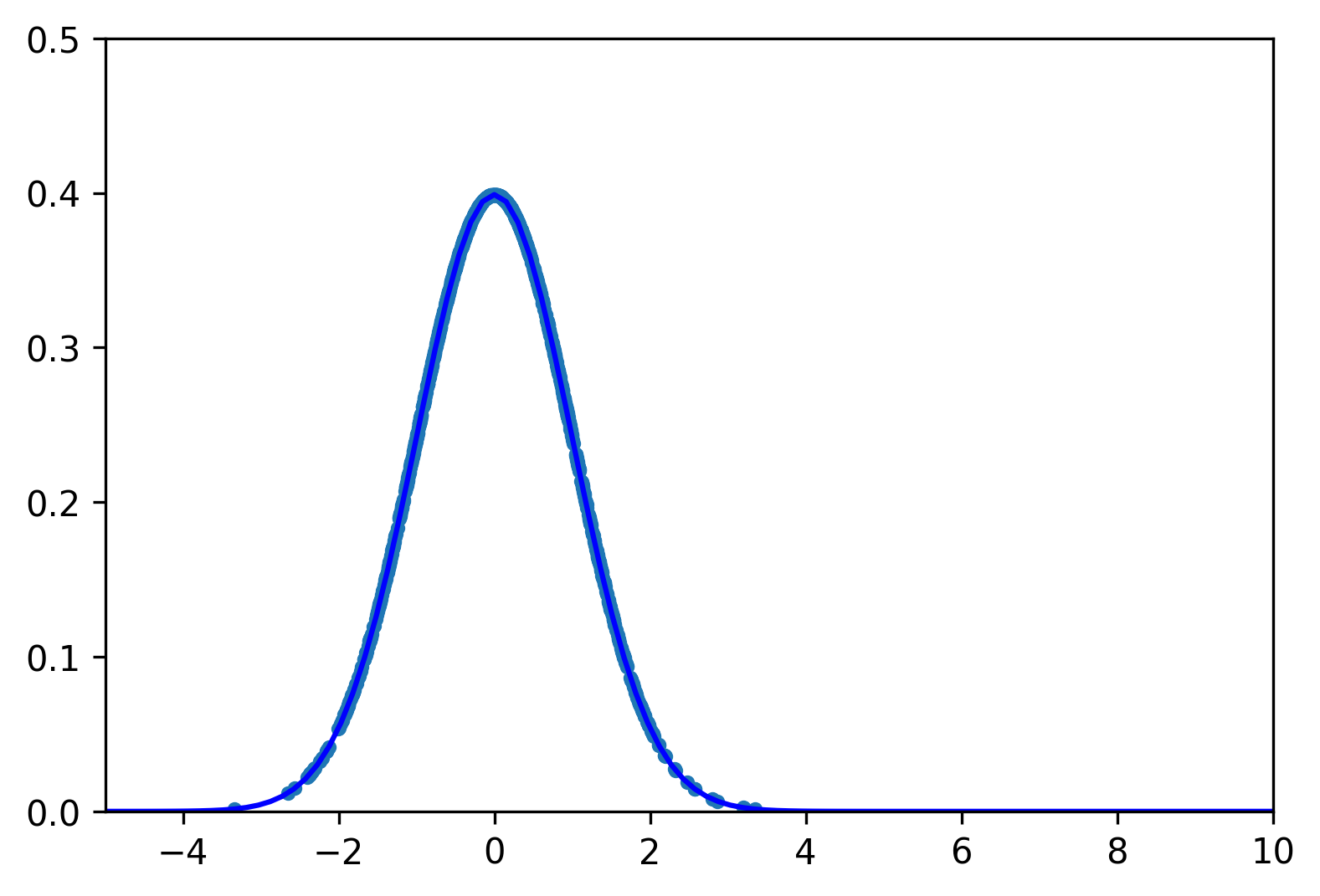}
\includegraphics[width=2.5cm]{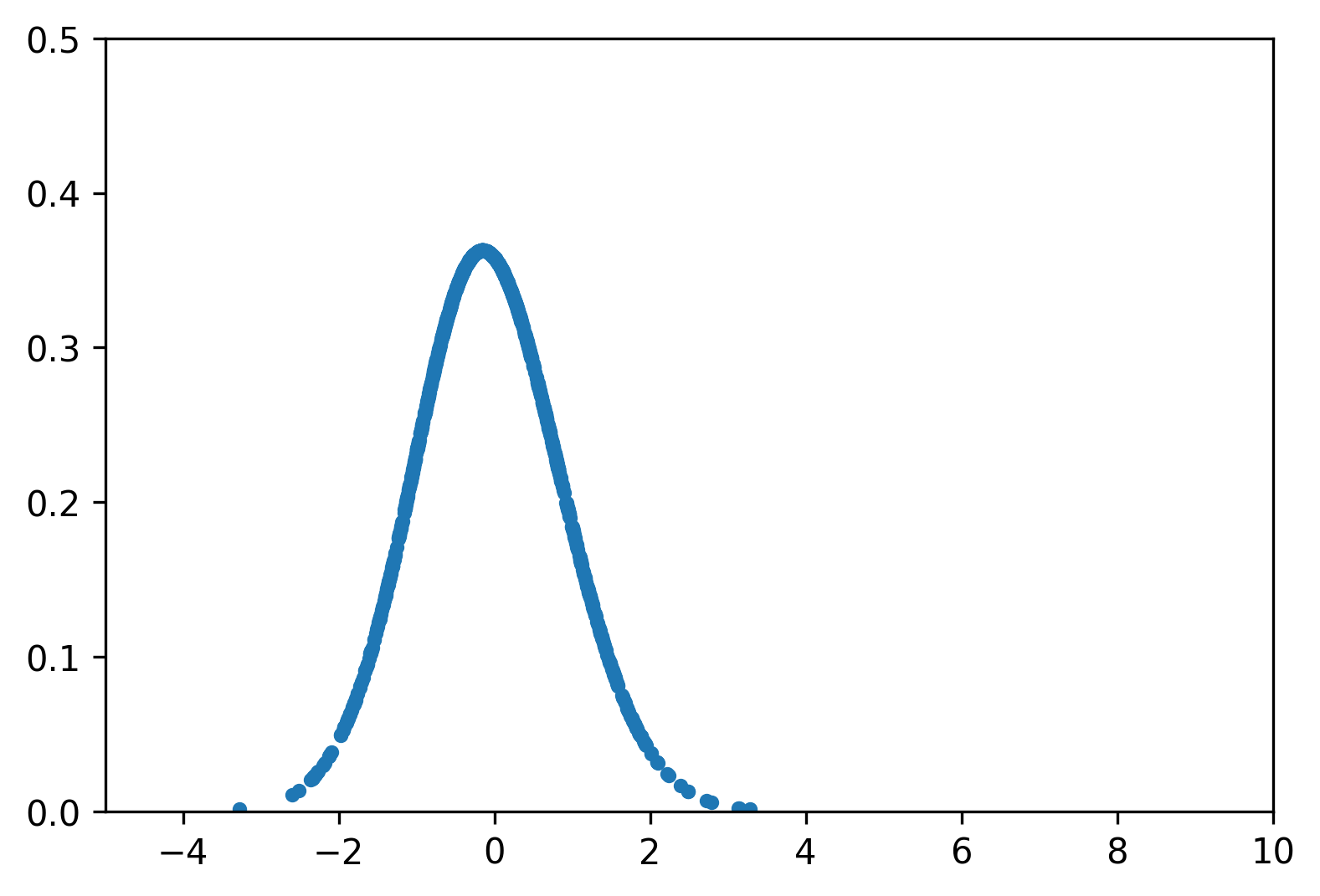}
\includegraphics[width=2.5cm]{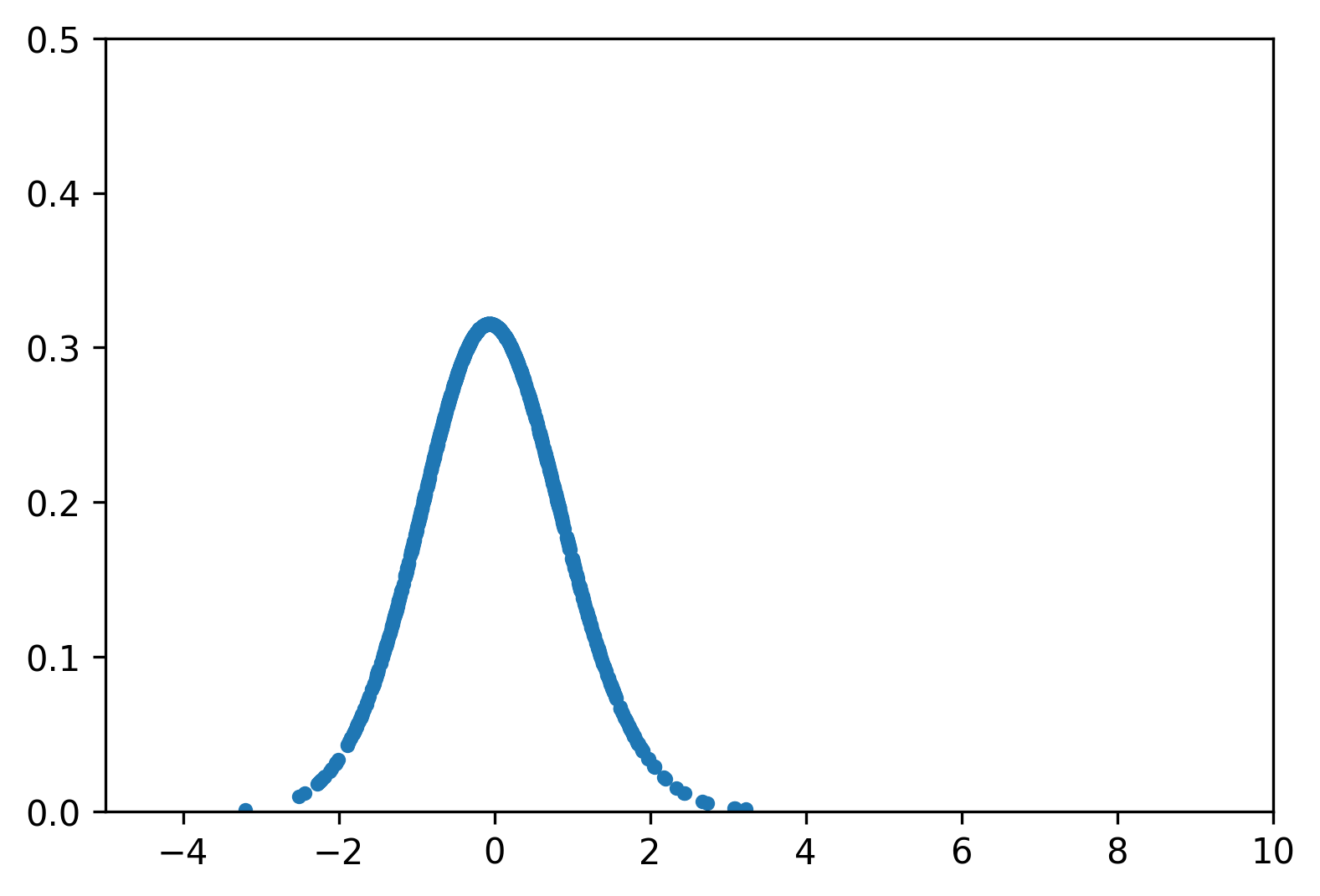}
\includegraphics[width=2.5cm]{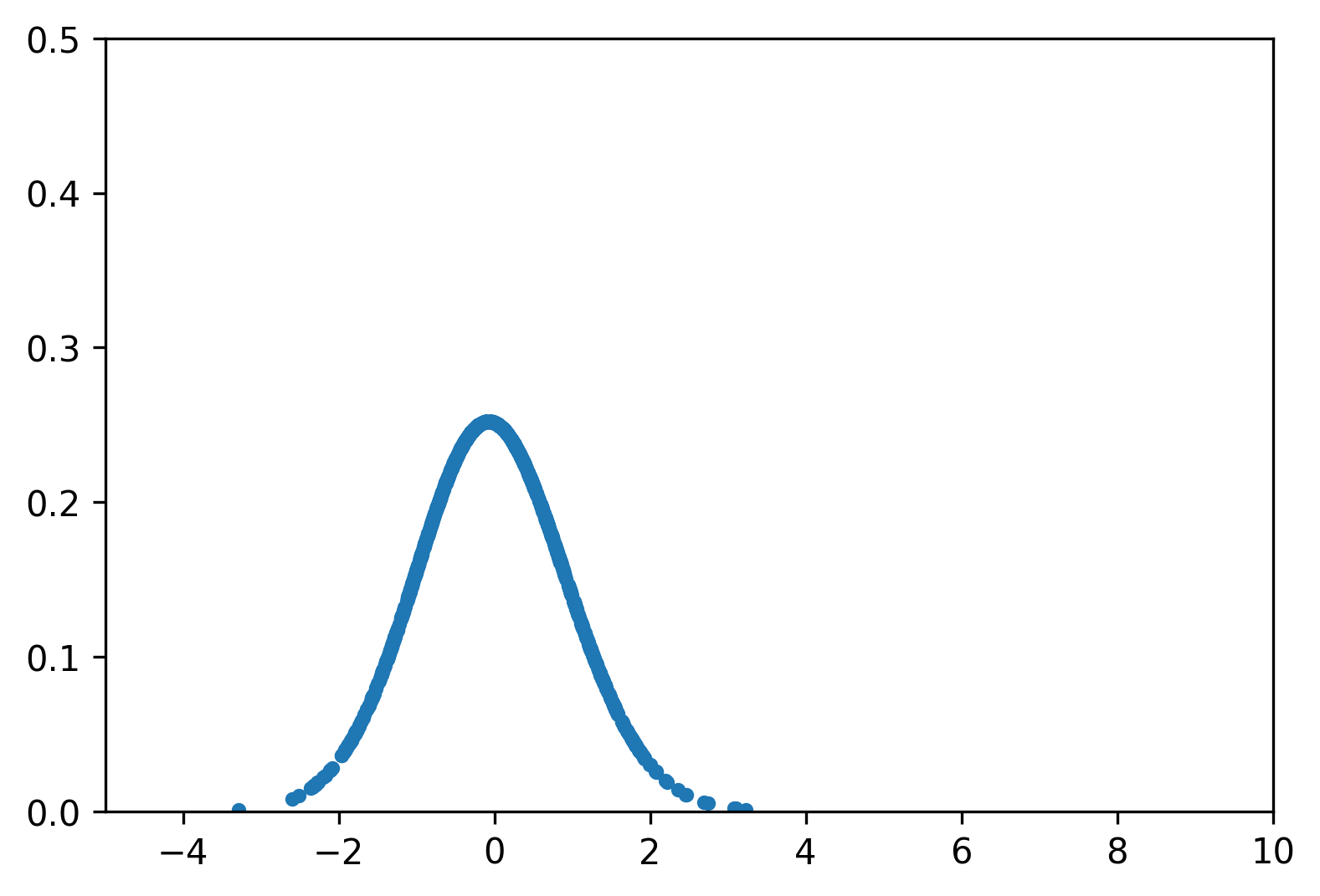}
\includegraphics[width=2.5cm]{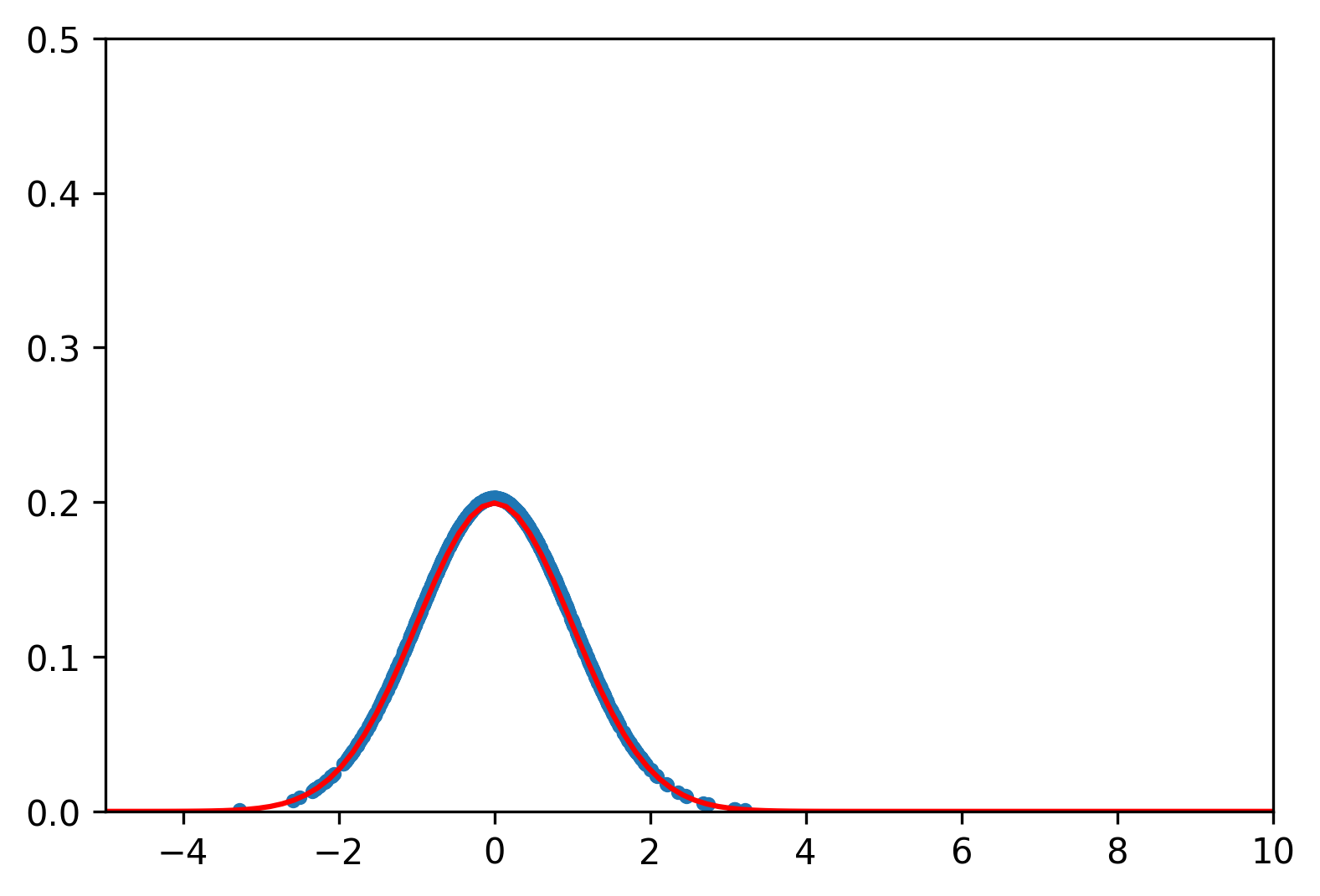}\\
\vspace{5pt}

\includegraphics[width=2.5cm]{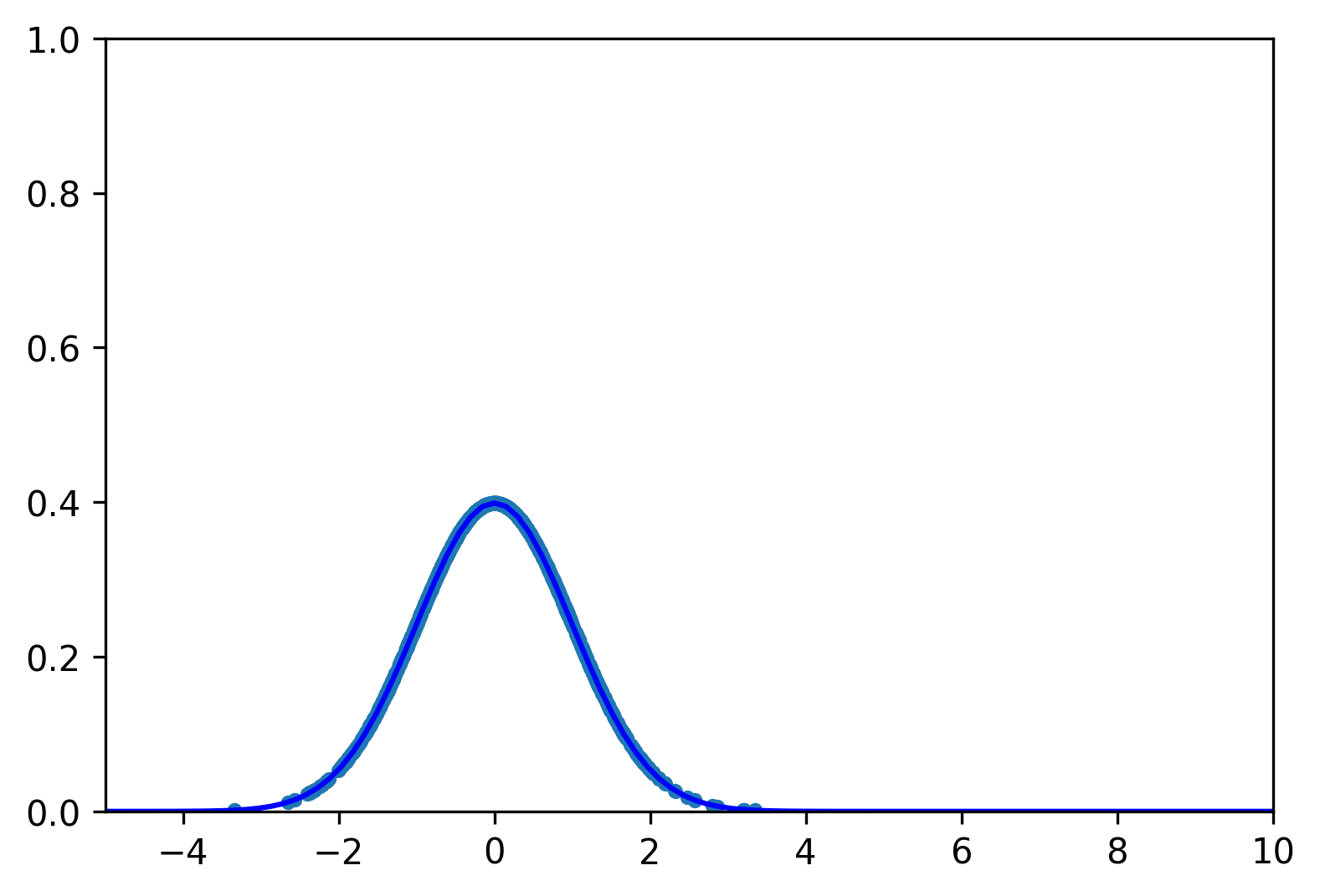}
\includegraphics[width=2.5cm]{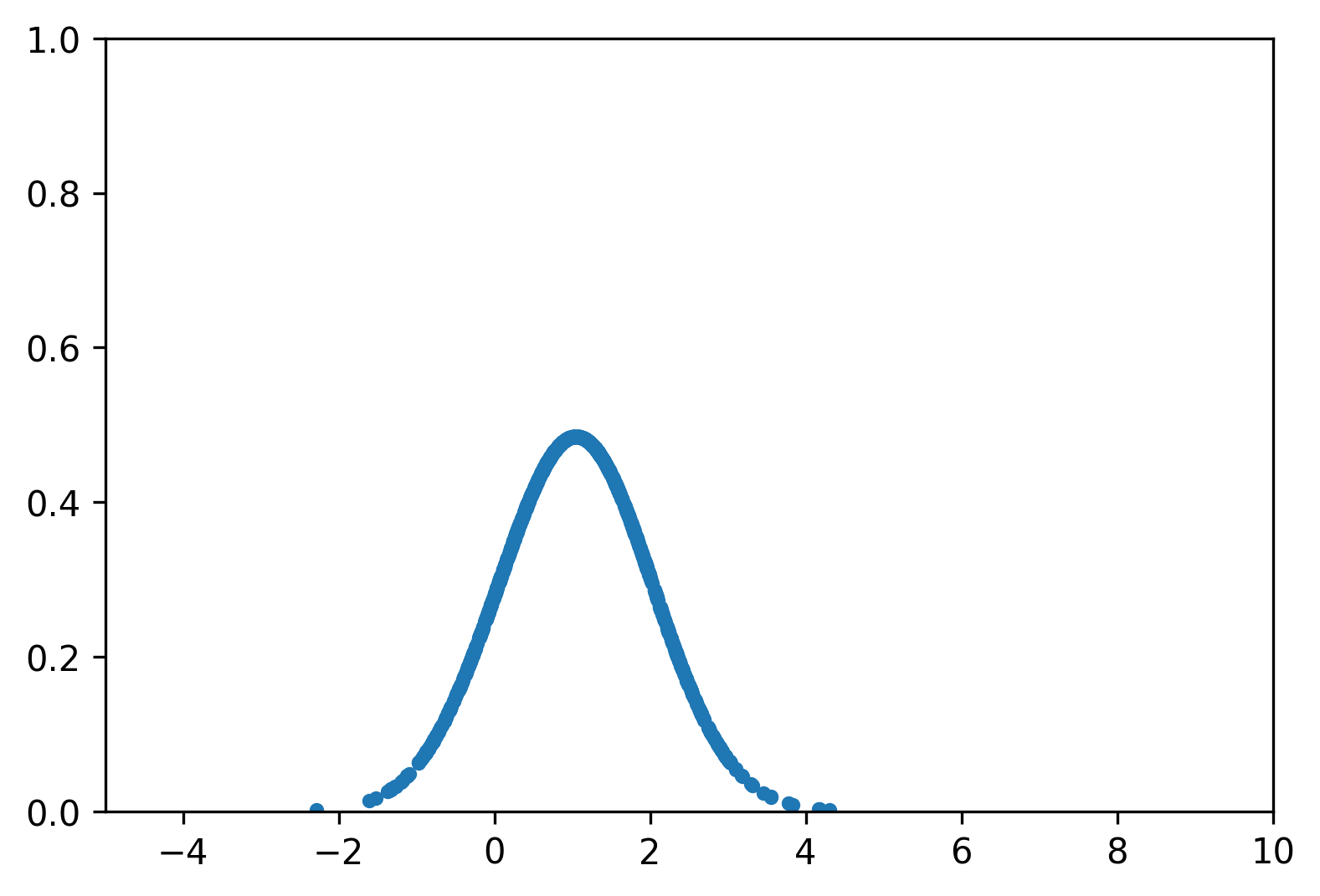}
\includegraphics[width=2.5cm]{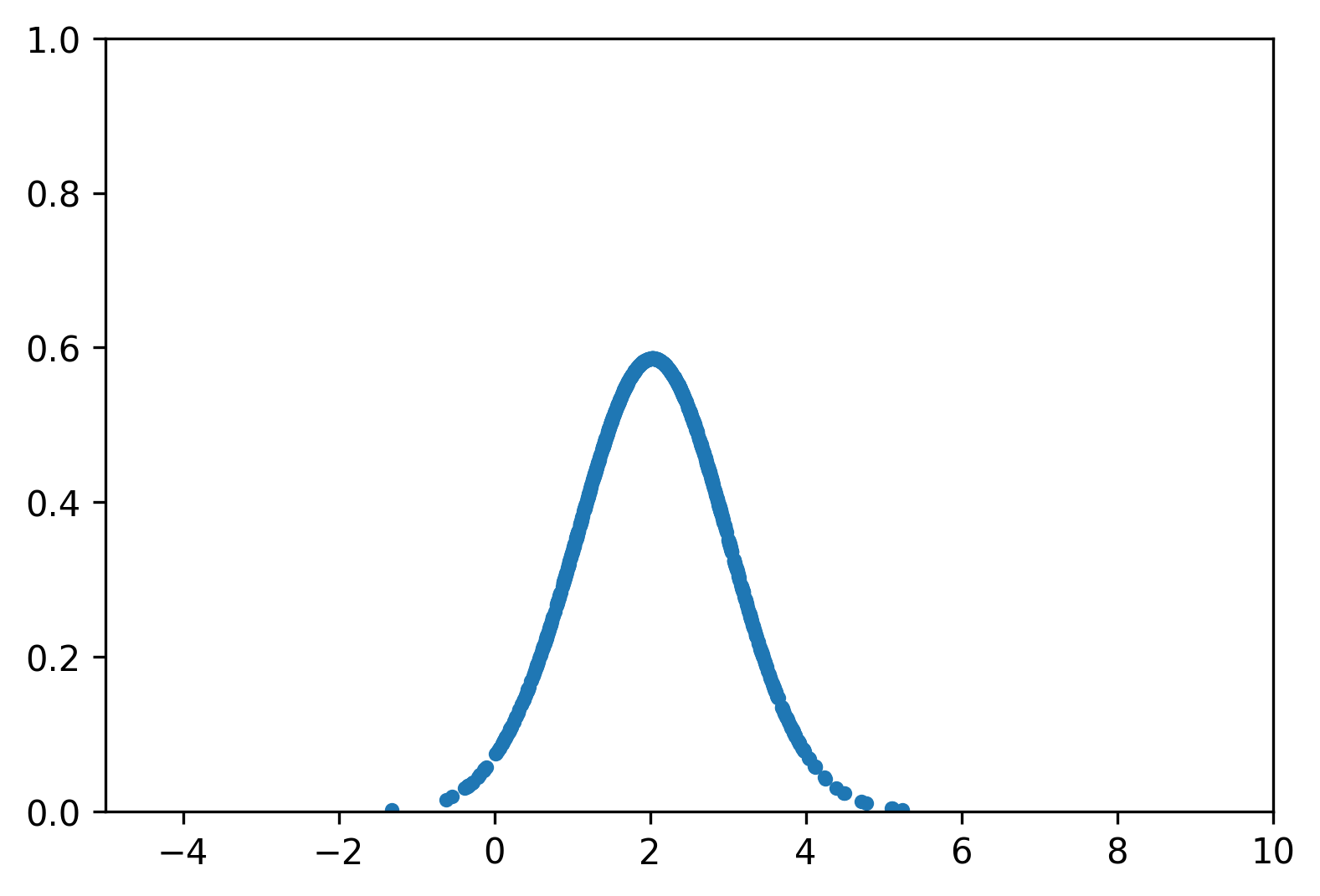}
\includegraphics[width=2.5cm]{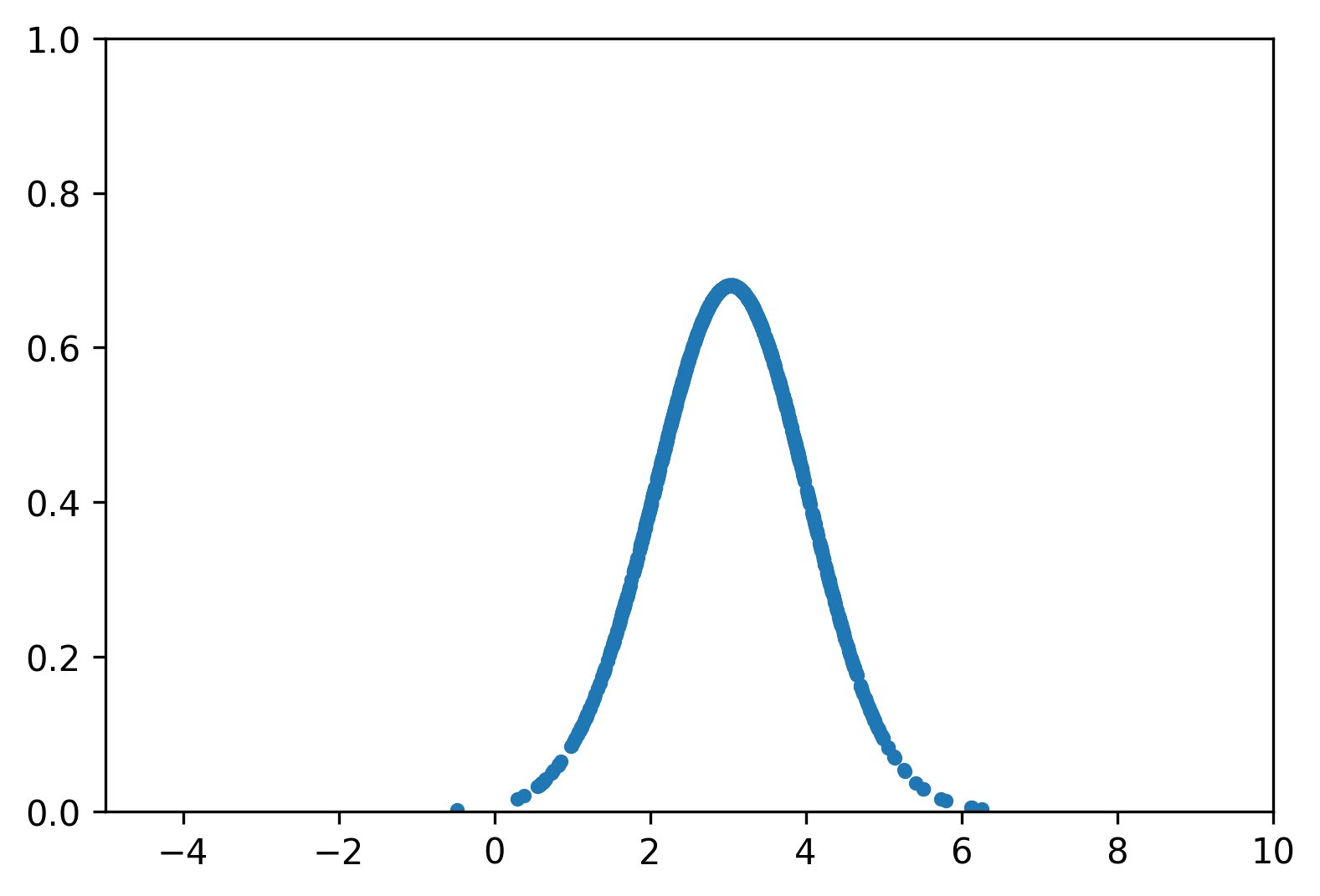}
\includegraphics[width=2.5cm]{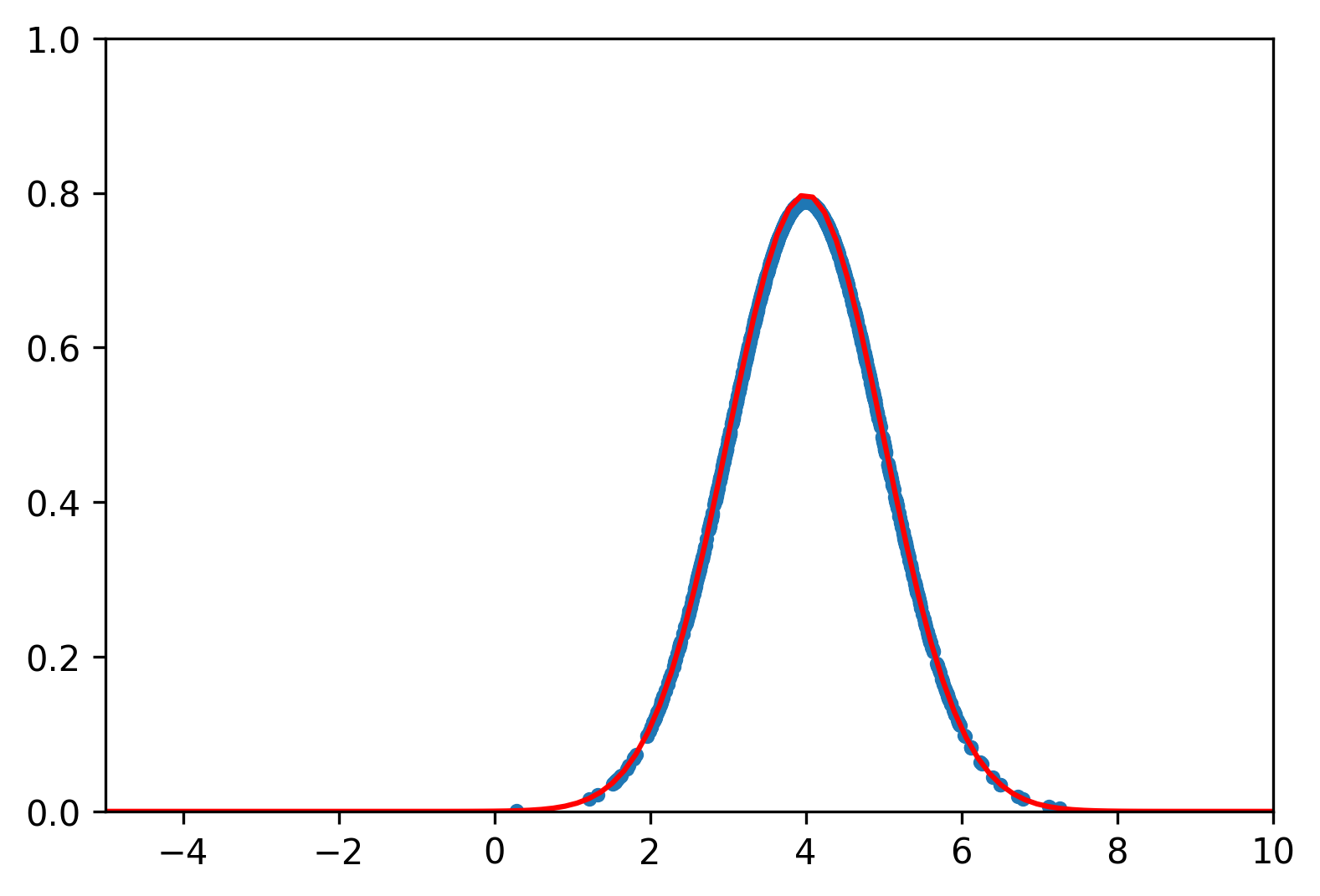}\\
\vspace{5pt}

\subfigure[$\rho_0(\boldsymbol{x}_0)$]{\includegraphics[width=2.5cm]{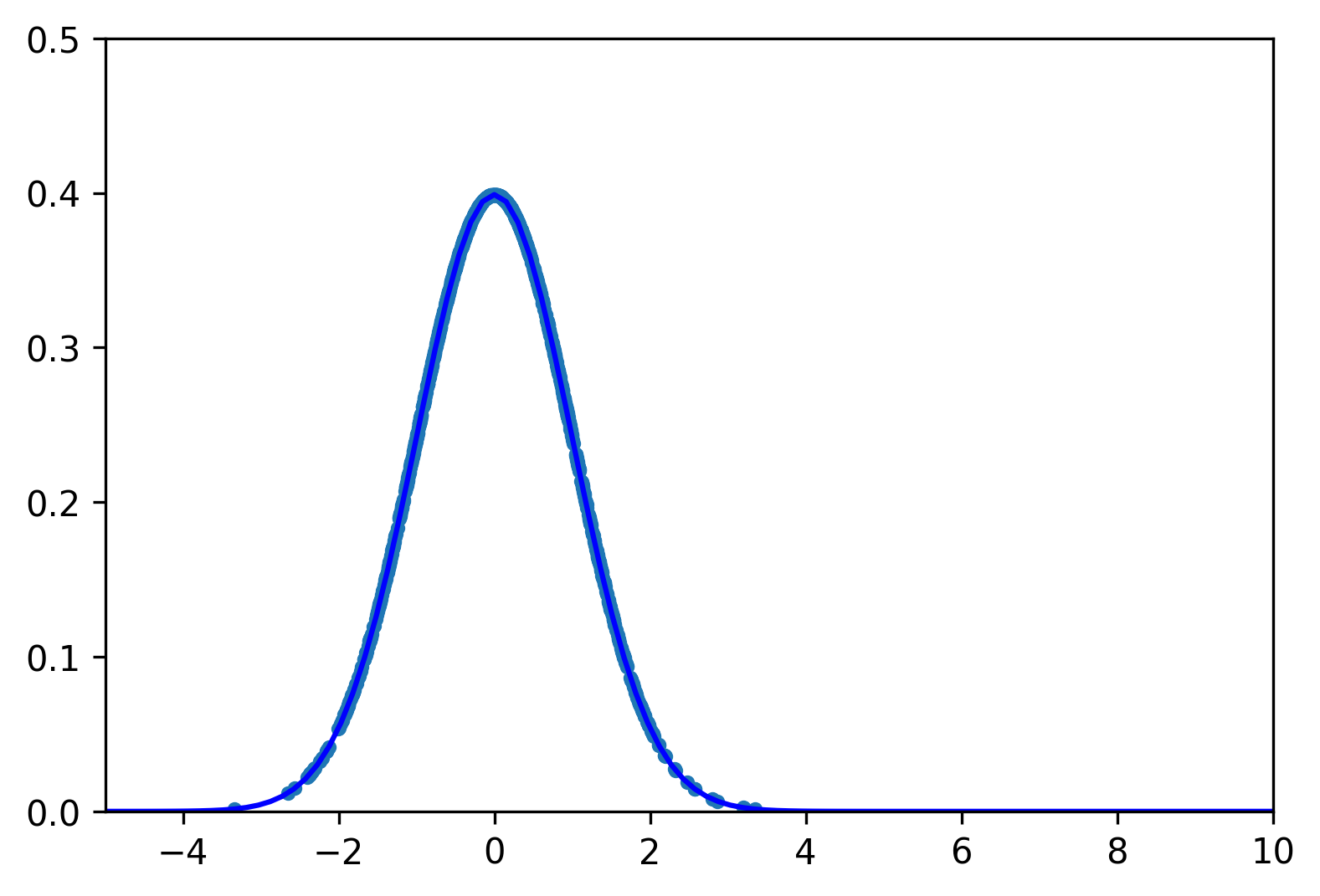}}
\subfigure[$\tilde{\rho}_{1/4}(\tilde{\boldsymbol{x}}_{1/4})$]{\includegraphics[width=2.5cm]{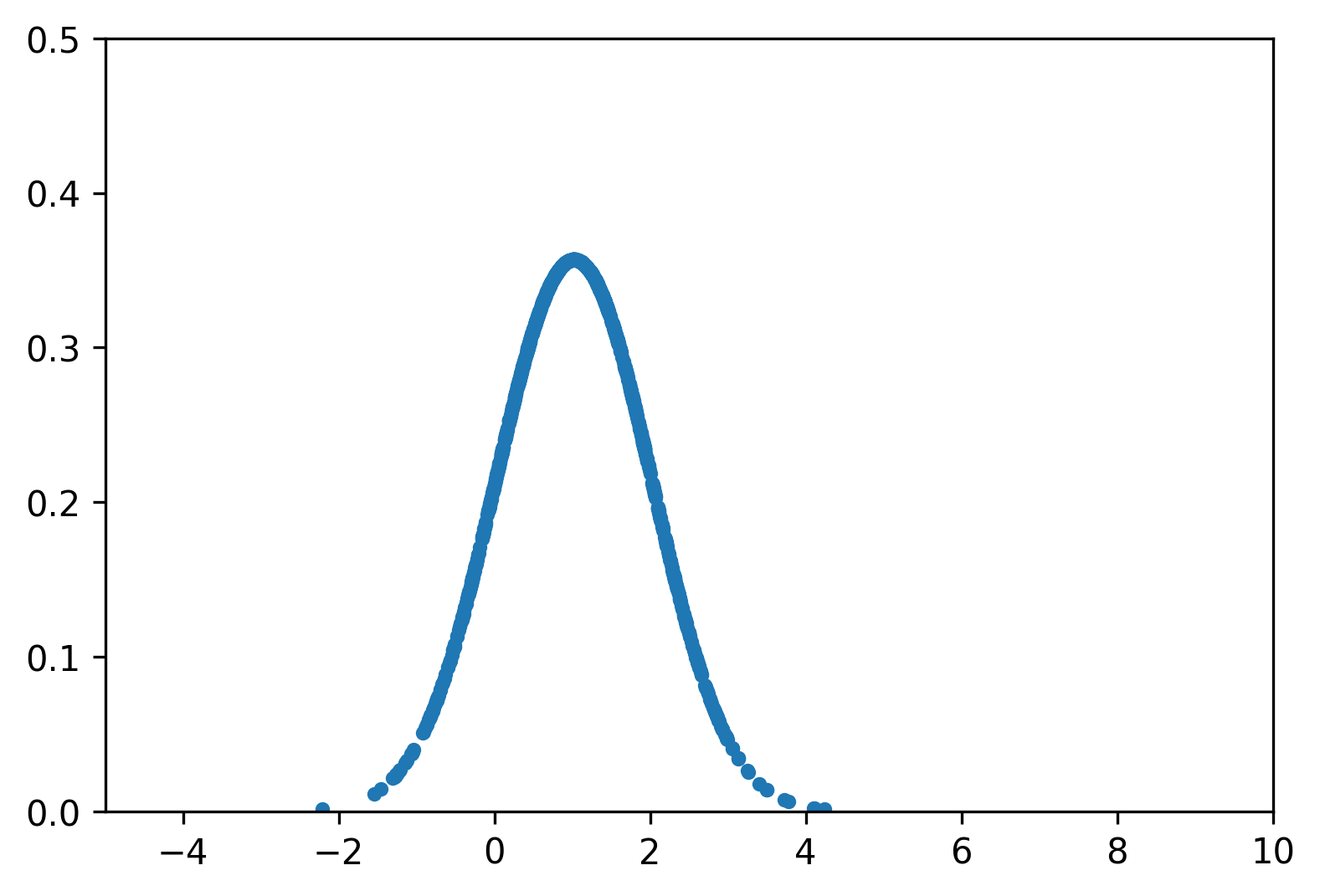}}
\subfigure[$\tilde{\rho}_{1/2}(\tilde{\boldsymbol{x}}_{1/2})$]{\includegraphics[width=2.5cm]{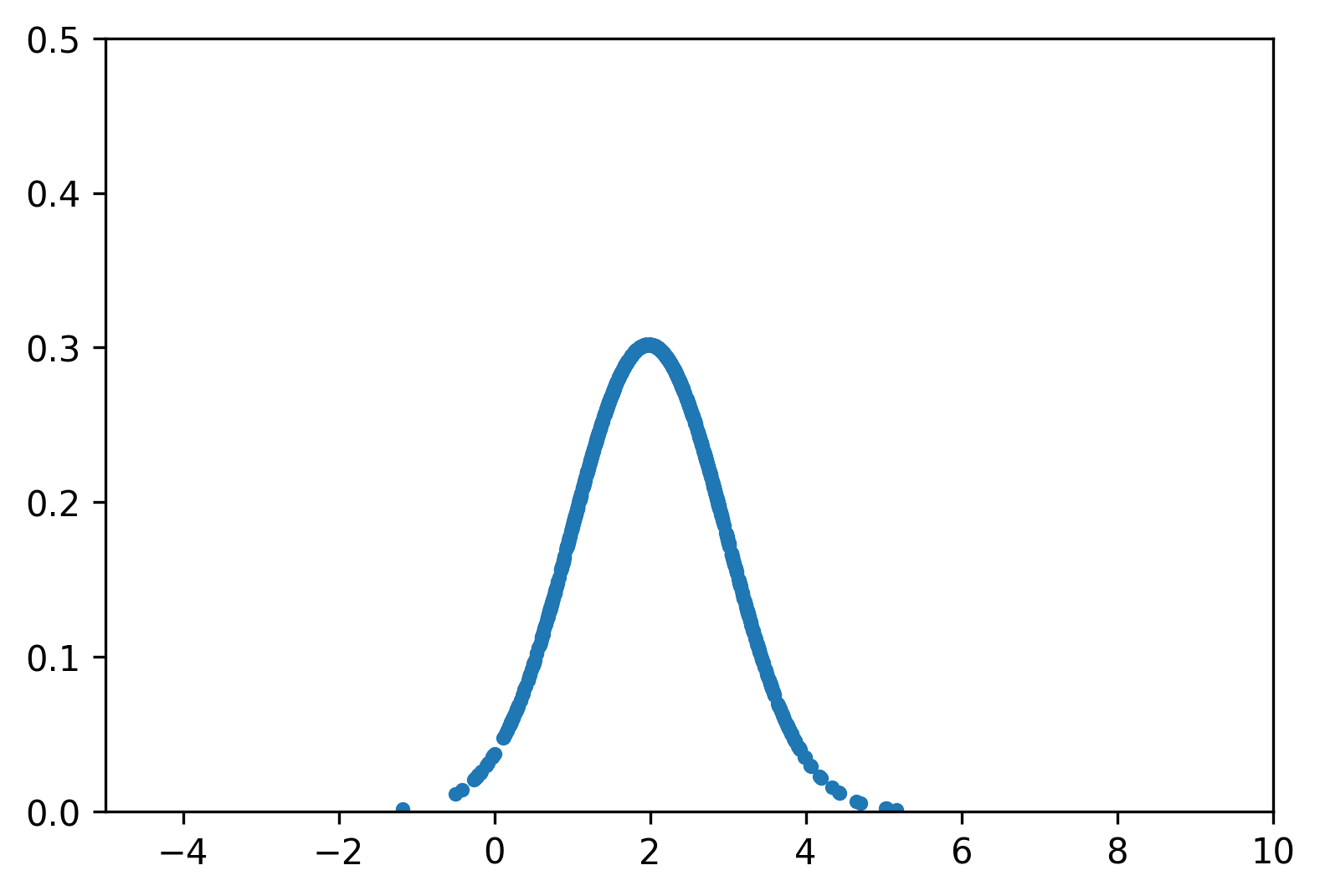}}
\subfigure[$\tilde{\rho}_{3/4}(\tilde{\boldsymbol{x}}_{3/4})$]{\includegraphics[width=2.5cm]{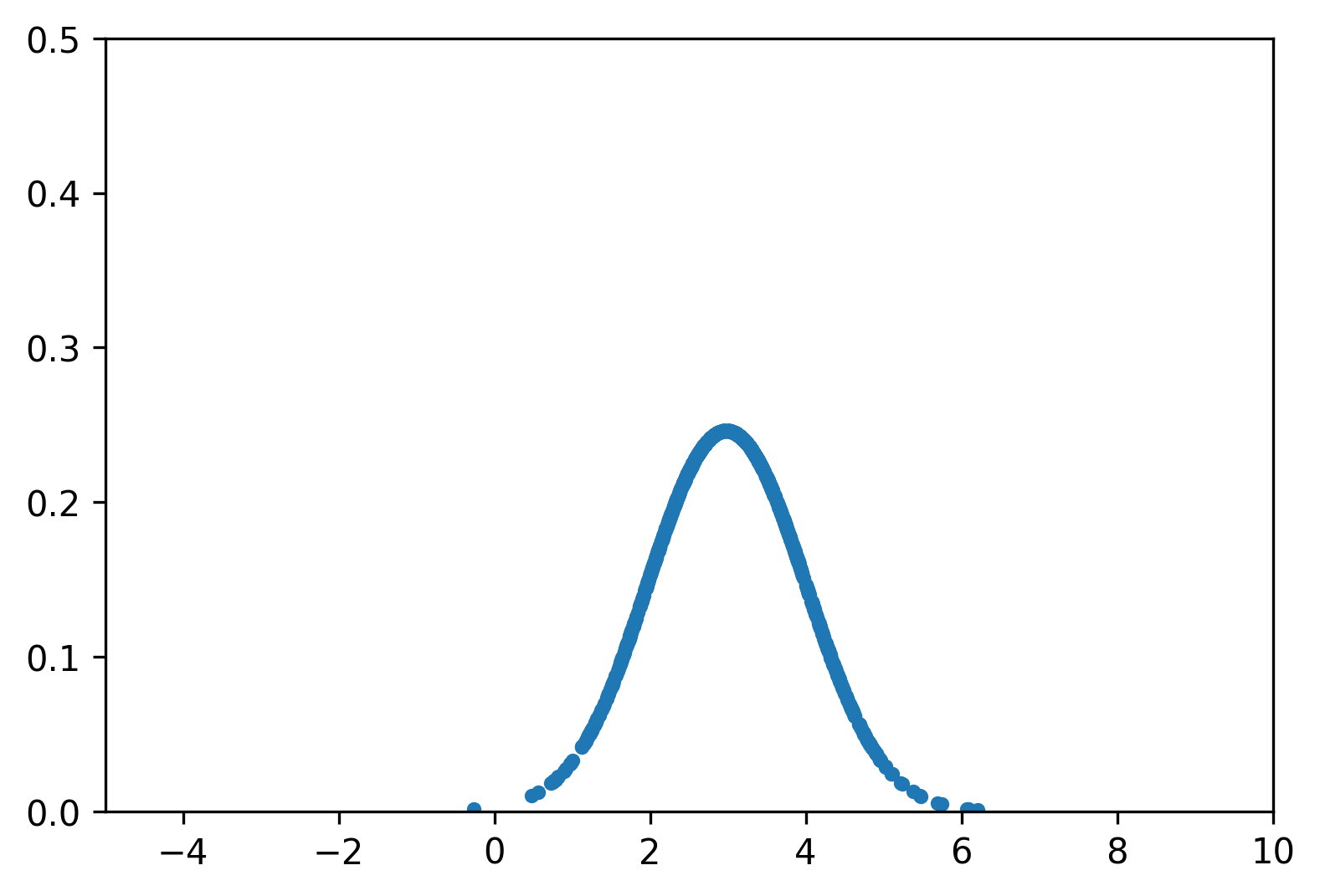}}
\subfigure[$\tilde{\rho}_{1}(\tilde{\boldsymbol{x}}_{1})$]{\includegraphics[width=2.5cm]{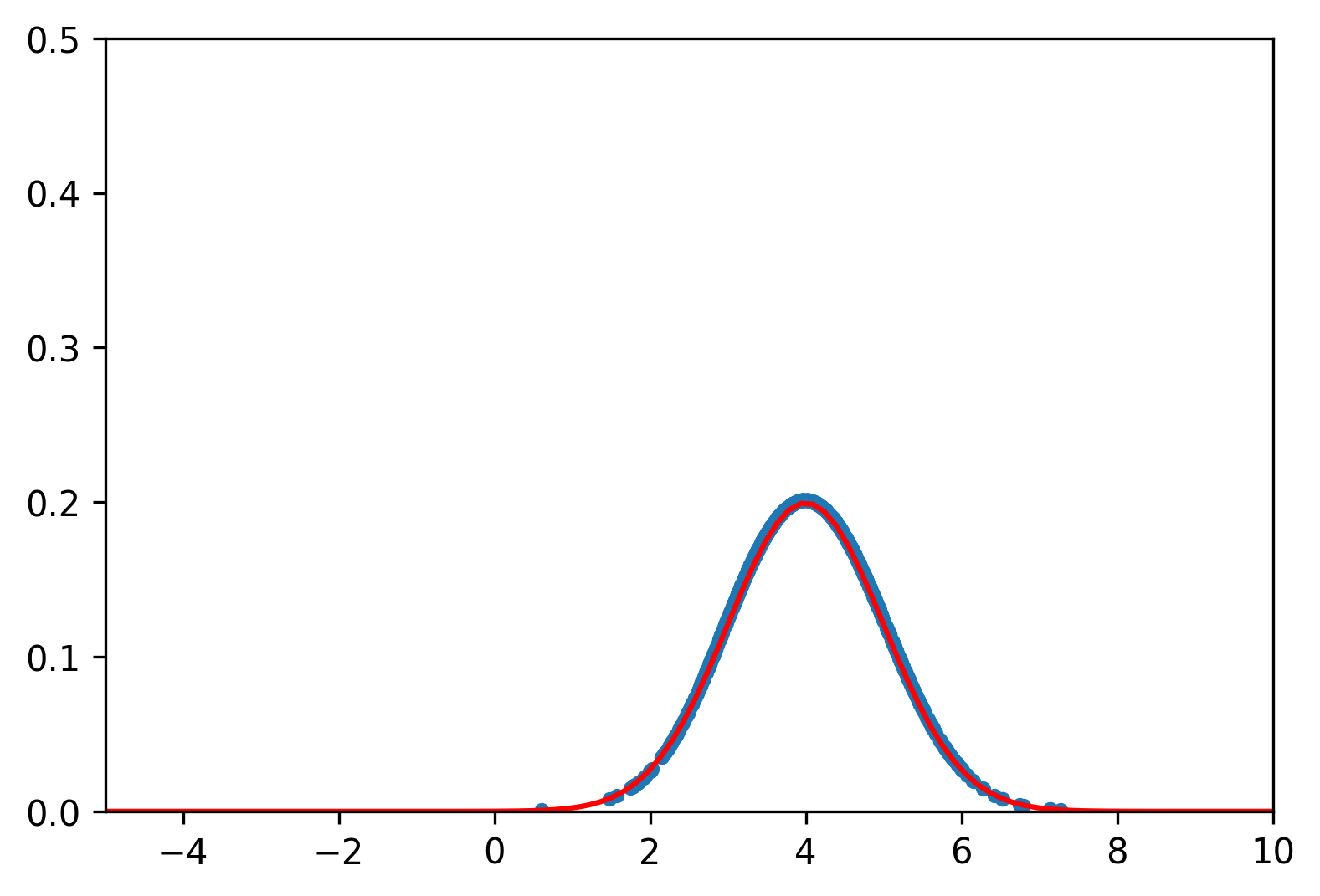}}

\caption{Illustration of these four Gaussian problems at $d=1$.}
\label{gaussian1}
\end{center}
\end{figure}

In Figure \ref{gaussian1}, we illustrate the evolution of the density transition from the initial density $\rho_0$ in (a) to the target density $\rho_1$ in (e). We observe that the unbalanced movement is reasonable, and the target density estimated by our trained network (depicted in blue in (e)) closely matches the true target density $\rho_1$ (shown in red in (e)). This result demonstrates the effectiveness of our training procedure.

Secondly, we consider the UOT problem on four typical Gaussian examples in $d \geq 2$ as listed below. 

\begin{itemize}\label{list2}
    
    \item Test 5: $\rho_0(\boldsymbol{x})=\rho_G(\boldsymbol{x},\mathbf{0},\mathbf{I})$, $\rho_1(\boldsymbol{x})=\frac{1}{2}\rho_G(\boldsymbol{x},-4 \cdot \mathbf{e_1},\mathbf{I})$;
    \item Test 6: $\rho_0(\boldsymbol{x})=\rho_G(\boldsymbol{x},\mathbf{0},\mathbf{I})$, $\rho_1(\boldsymbol{x})=2\rho_G(\boldsymbol{x}, 4 \cdot \mathbf{e_1},\mathbf{I})$;
    \item Test 7: $\rho_0(\boldsymbol{x})=\rho_G(\boldsymbol{x},\mathbf{0},0.3 \cdot\mathbf{I})$, $\rho_1(\boldsymbol{x})=2\rho_G(\boldsymbol{x},4 \cdot \mathbf{e_1},0.3\cdot\mathbf{I})$;
    \item Test 8: $\rho_0(\boldsymbol{x})=\rho_G(\boldsymbol{x},-4\cdot \mathbf{e_1}-4 \cdot \mathbf{e_2}, \mathbf{I})$, $\rho_1(\boldsymbol{x})=\frac{1}{2}\rho_G(\boldsymbol{x},4\cdot \mathbf{e_1}+4 \cdot \mathbf{e_2}, \mathbf{I})$.
  
\end{itemize}

\begin{figure}[H]
\begin{center}
\includegraphics[width=2.4cm]{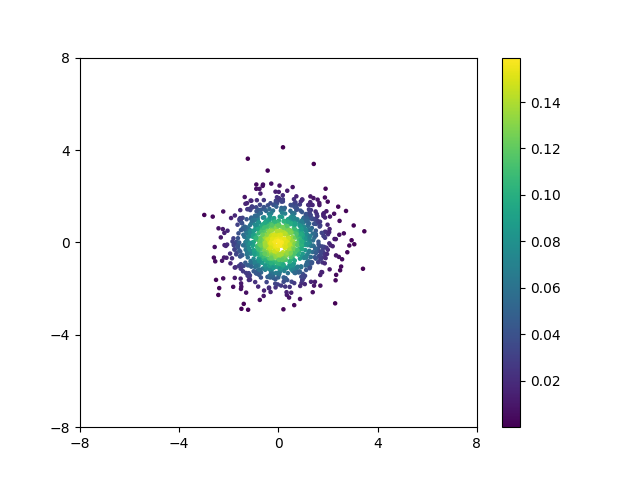}
\hspace{-5mm}
\includegraphics[width=2.4cm]{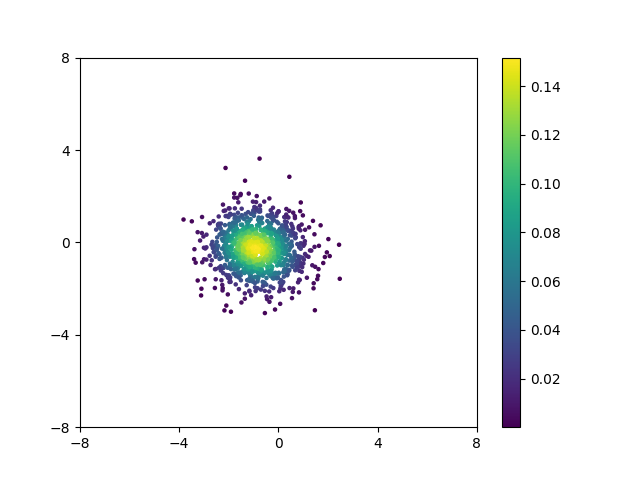}
\hspace{-5mm}
\includegraphics[width=2.4cm]{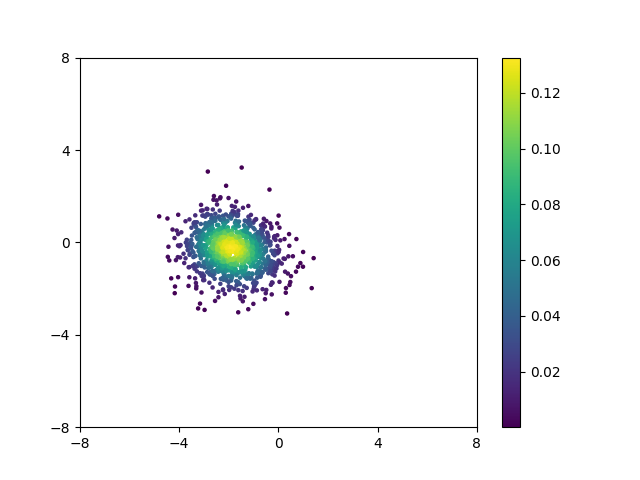}
\hspace{-5mm}
\includegraphics[width=2.4cm]{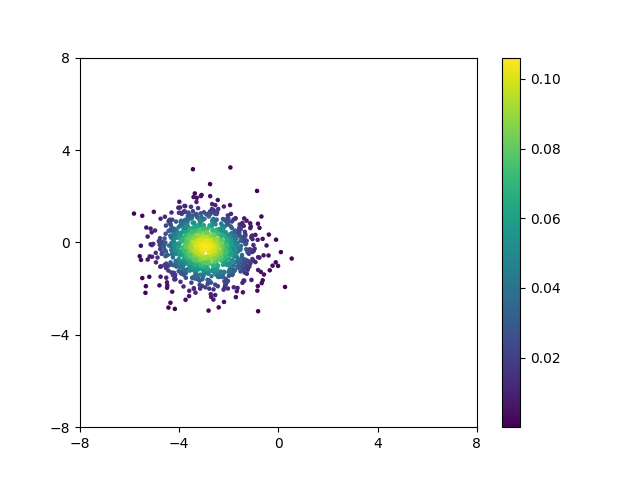}
\hspace{-5mm}
\includegraphics[width=2.4cm]{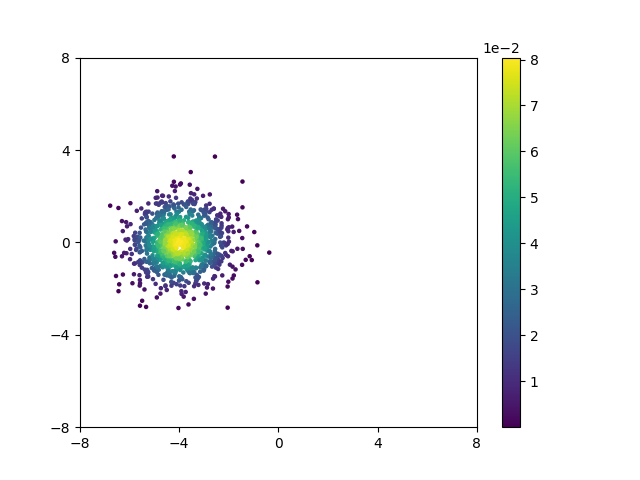}
\hspace{-5mm}
\includegraphics[width=2.4cm]{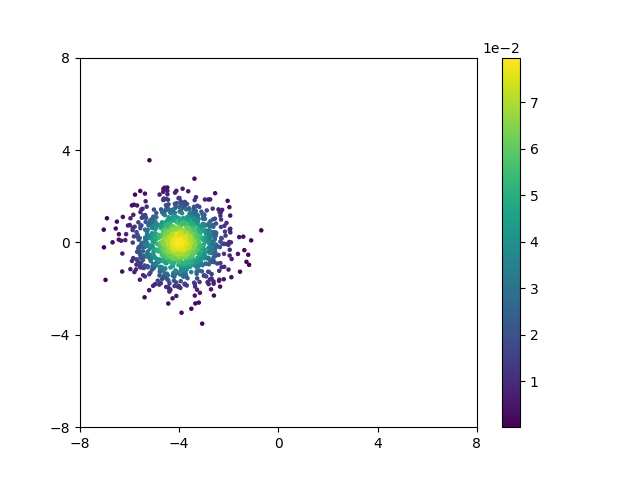}\\
\vspace{5pt}

\includegraphics[width=2.4cm]{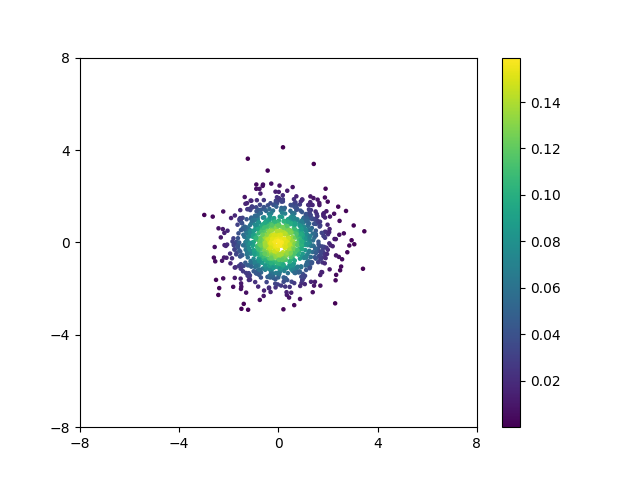}
\hspace{-5mm}
\includegraphics[width=2.4cm]{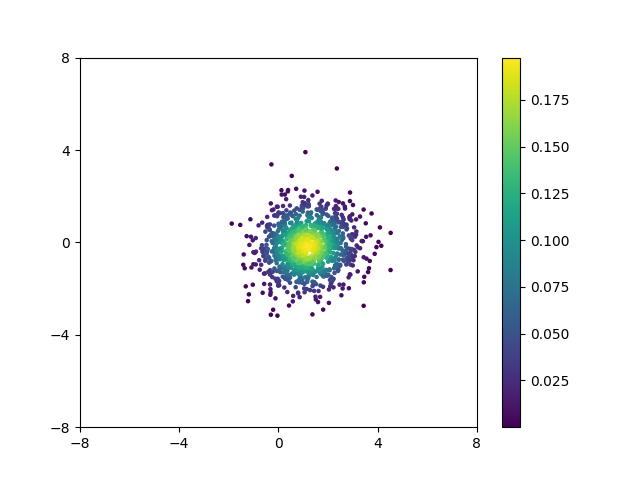}
\hspace{-5mm}
\includegraphics[width=2.4cm]{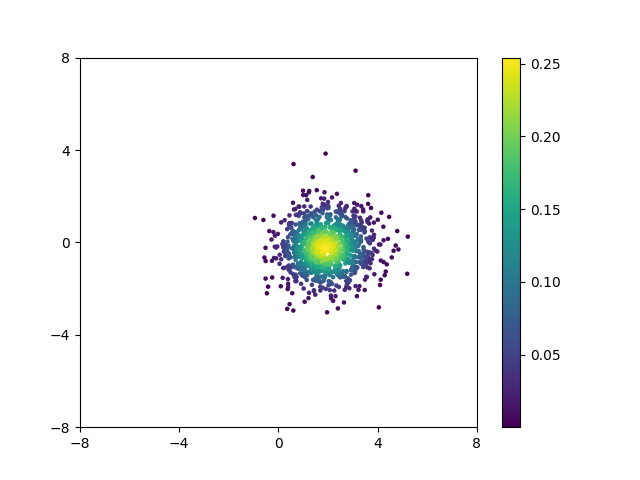}
\hspace{-5mm}
\includegraphics[width=2.4cm]{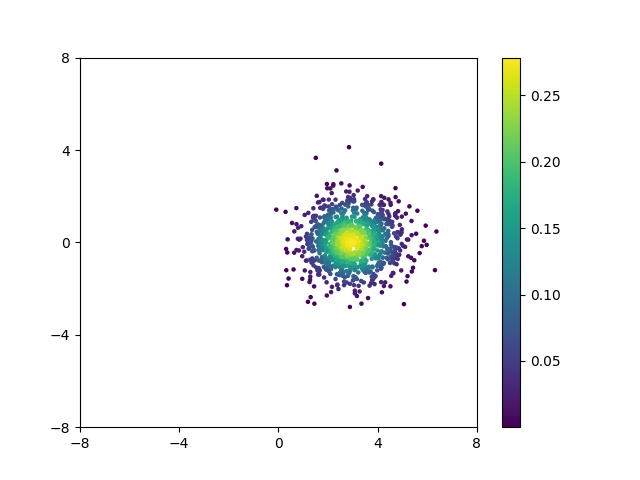}
\hspace{-5mm}
\includegraphics[width=2.4cm]{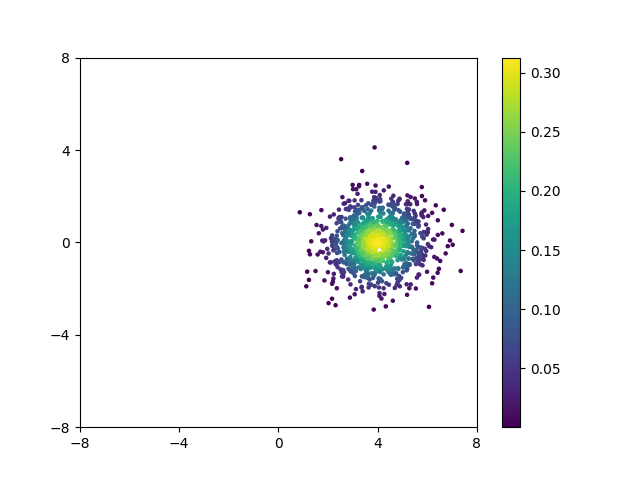}
\hspace{-5mm}
\includegraphics[width=2.4cm]{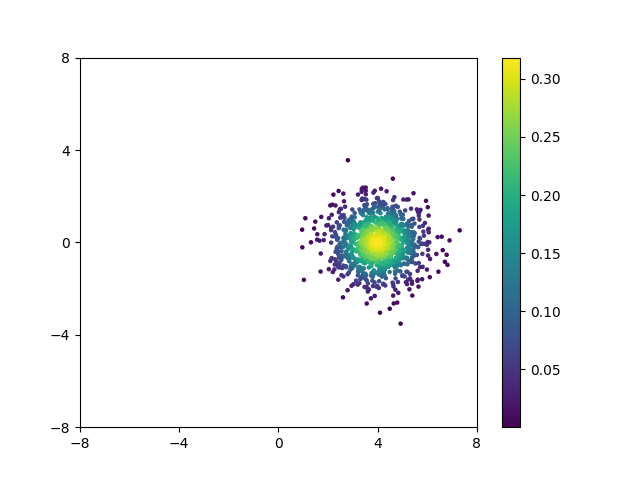}\\
\vspace{5pt}

\includegraphics[width=2.4cm]{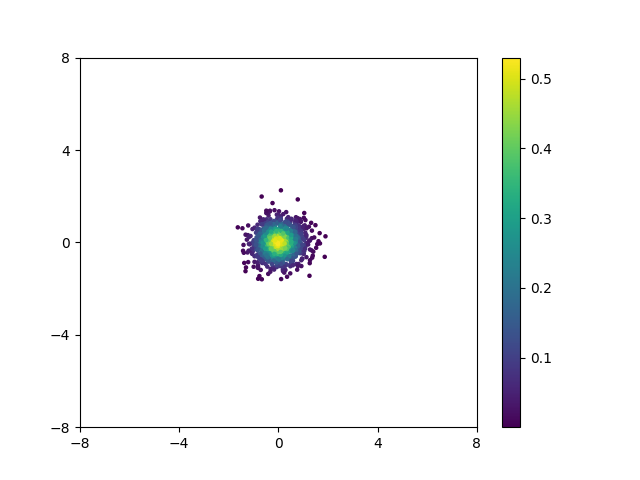}
\hspace{-5mm}
\includegraphics[width=2.4cm]{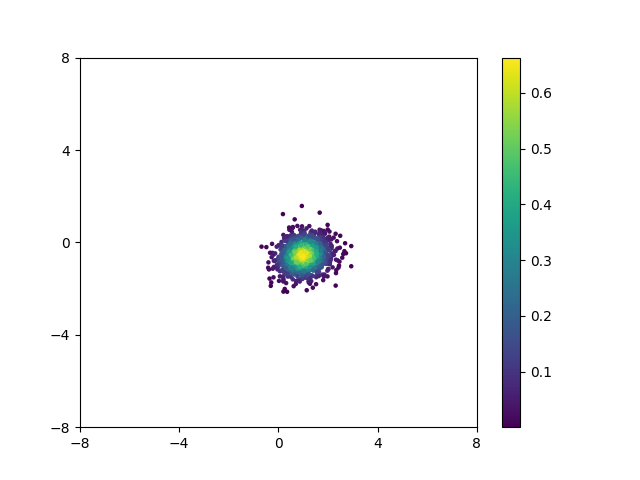}
\hspace{-5mm}
\includegraphics[width=2.4cm]{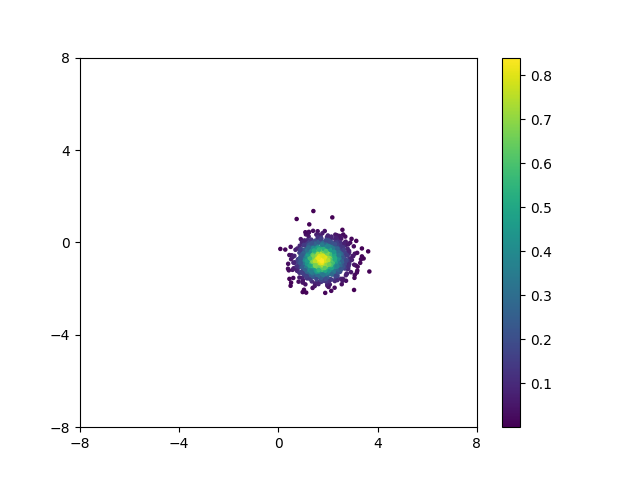}
\hspace{-5mm}
\includegraphics[width=2.4cm]{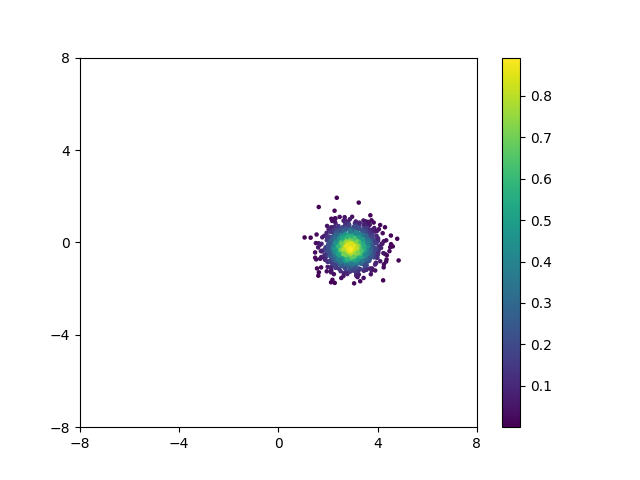}
\hspace{-5mm}
\includegraphics[width=2.4cm]{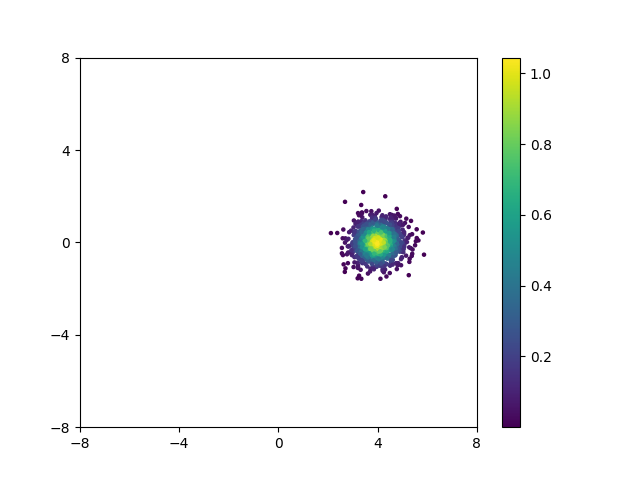}
\hspace{-5mm}
\includegraphics[width=2.4cm]{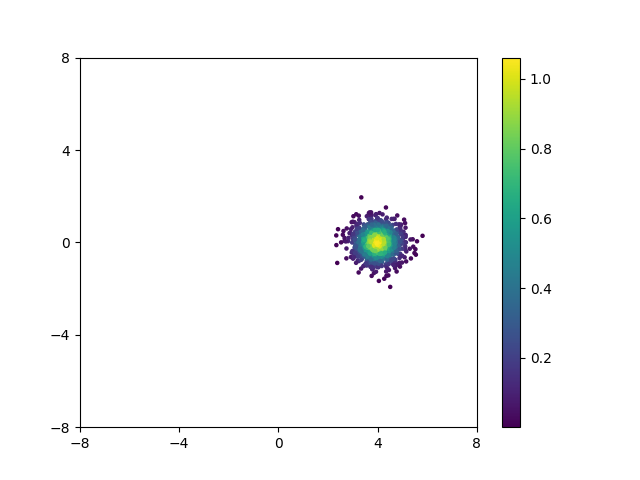}\\
\vspace{5pt}

\subfigure[$\rho_0(\boldsymbol{\boldsymbol{x}}_0)$]{\includegraphics[width=2.4cm]{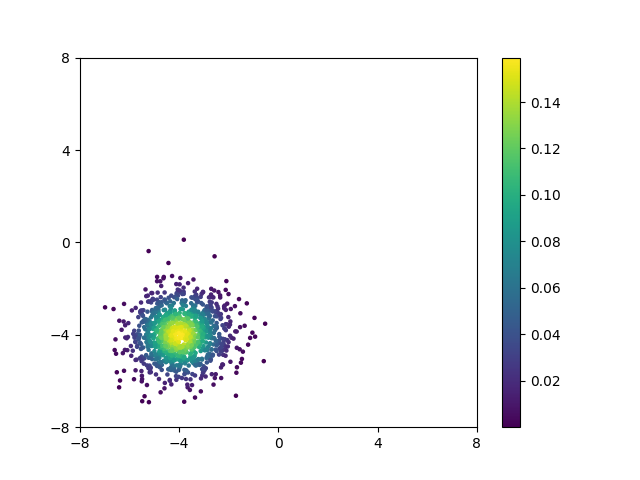}}
\hspace{-5mm}
\subfigure[$\tilde{\rho}_{1/4}(\tilde{\boldsymbol{\boldsymbol{x}}}_{1/4})$]{\includegraphics[width=2.4cm]{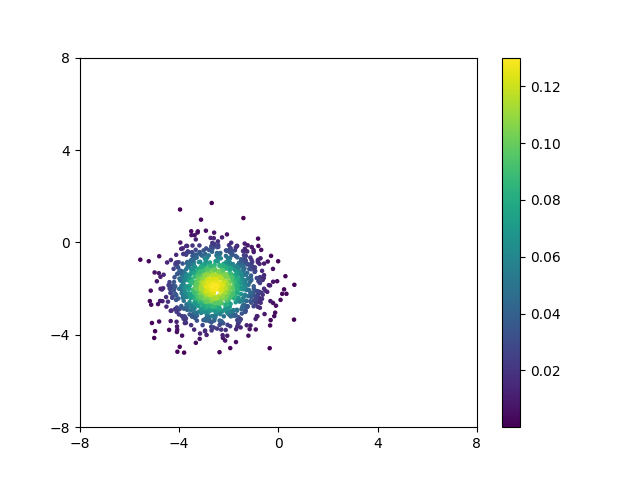}}
\hspace{-5mm}
\subfigure[$\tilde{\rho}_{1/2}(\tilde{\boldsymbol{\boldsymbol{x}}}_{1/2})$]{\includegraphics[width=2.4cm]{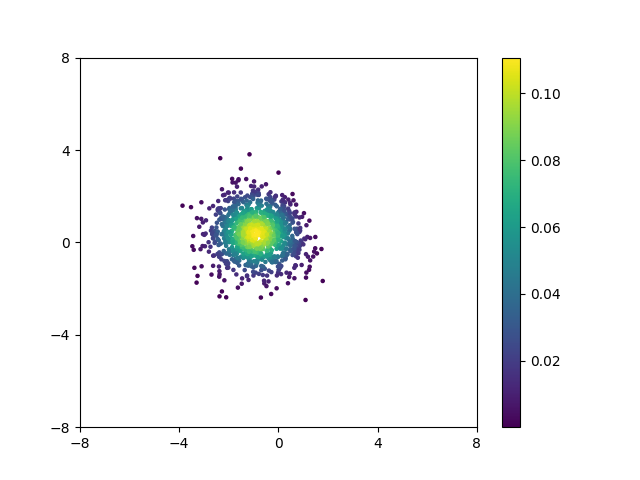}}
\hspace{-5mm}
\subfigure[$\tilde{\rho}_{3/4}(\tilde{\boldsymbol{\boldsymbol{x}}}_{3/4})$]{\includegraphics[width=2.4cm]{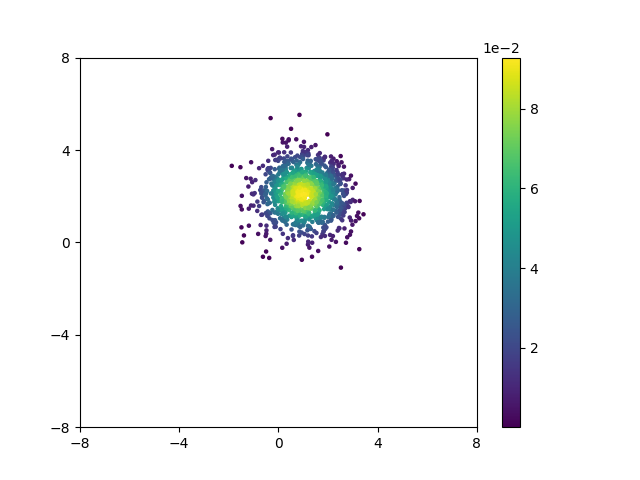}}
\hspace{-5mm}
\subfigure[$\tilde{\rho}_{1}(\tilde{\boldsymbol{\boldsymbol{x}}}_{1})$]{\includegraphics[width=2.4cm]{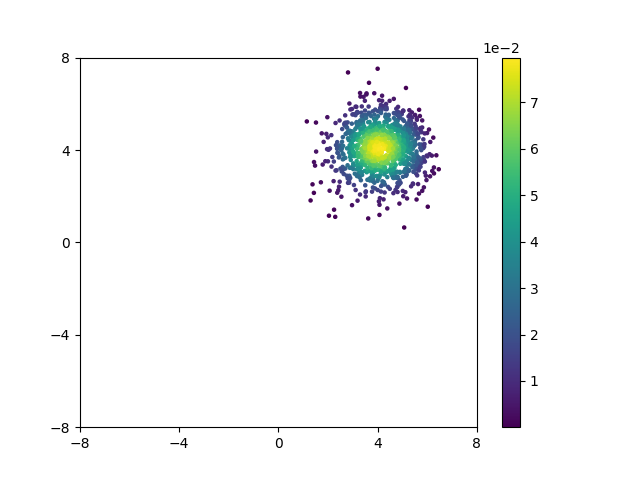}}
\hspace{-5mm}
\subfigure[$\rho_{1}(\boldsymbol{\boldsymbol{x}}_{1})$]{\includegraphics[width=2.4cm]{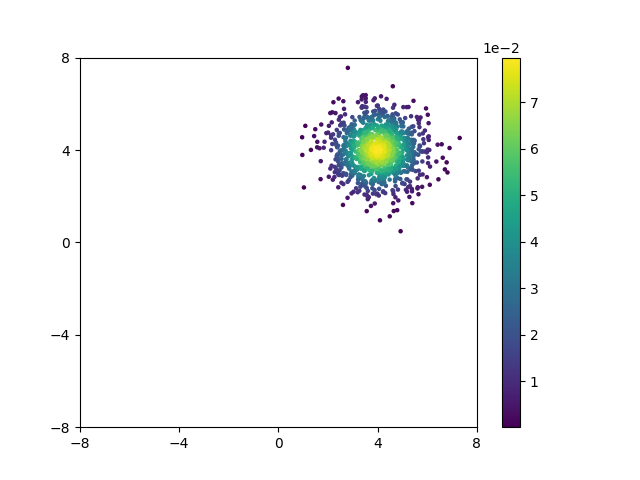}}

\caption{Illustration of these Gaussian problems at $d=2$.}
\label{compare_d2}
\end{center}
\end{figure}

\begin{figure}[H]
\begin{center}
\includegraphics[width=2.4cm]{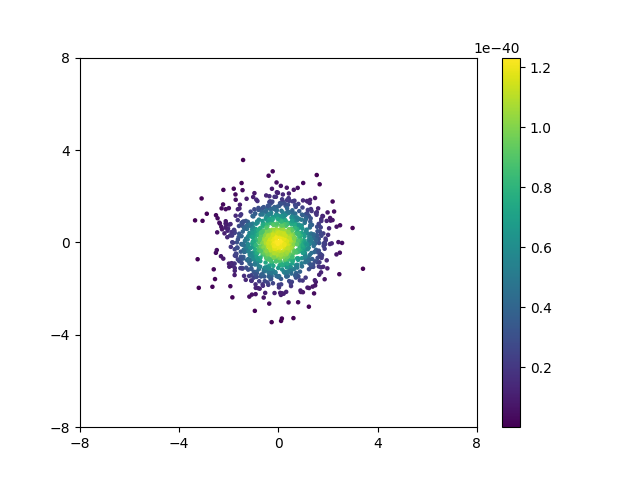}
\hspace{-5mm}
\includegraphics[width=2.4cm]{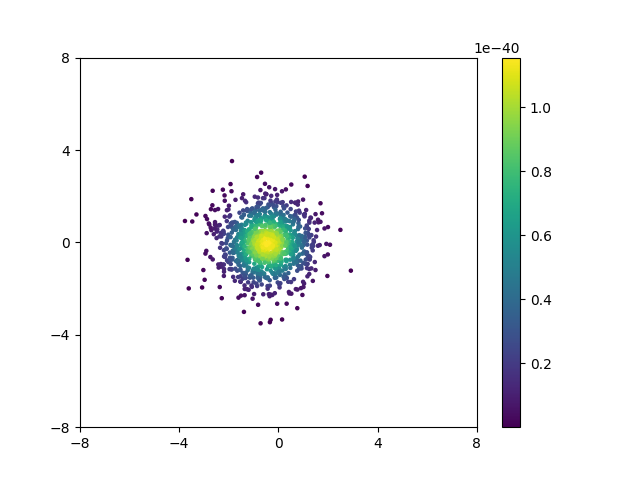}
\hspace{-5mm}
\includegraphics[width=2.4cm]{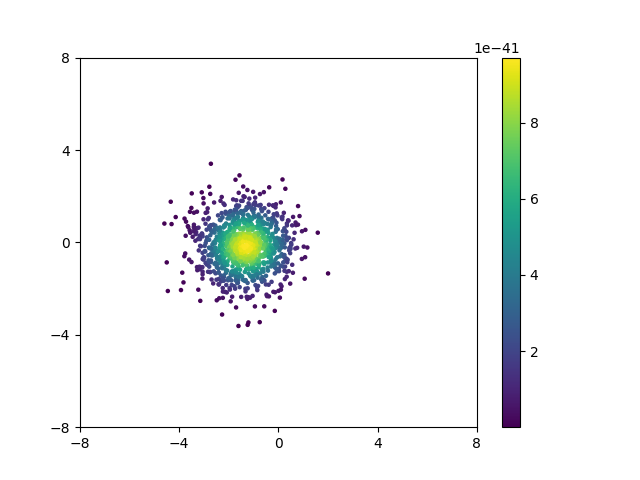}
\hspace{-5mm}
\includegraphics[width=2.4cm]{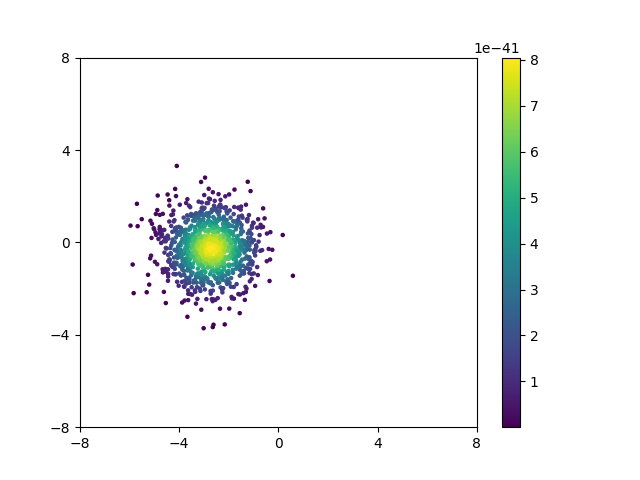}
\hspace{-5mm}
\includegraphics[width=2.4cm]{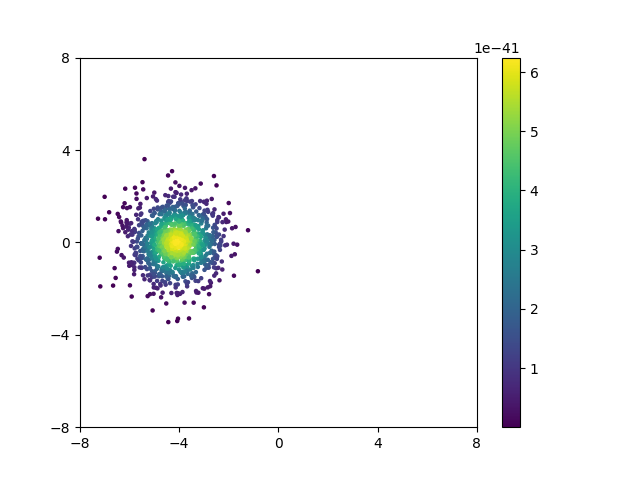}
\hspace{-5mm}
\includegraphics[width=2.4cm]{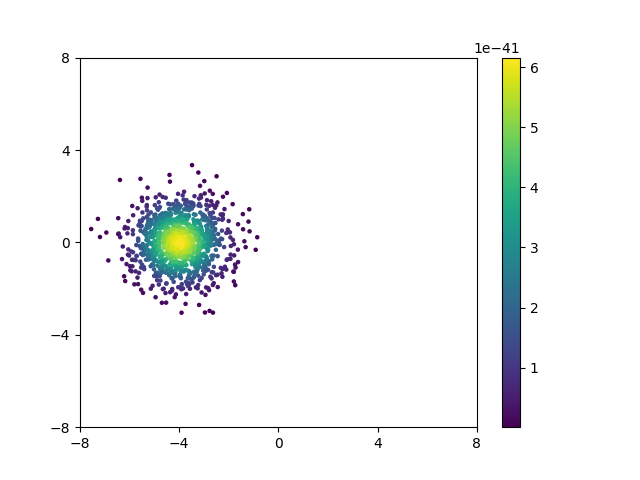}\\
\vspace{5pt}

\includegraphics[width=2.4cm]{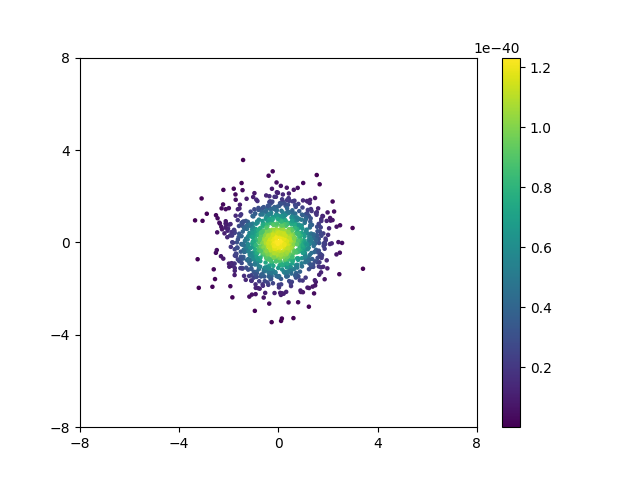}
\hspace{-5mm}
\includegraphics[width=2.4cm]{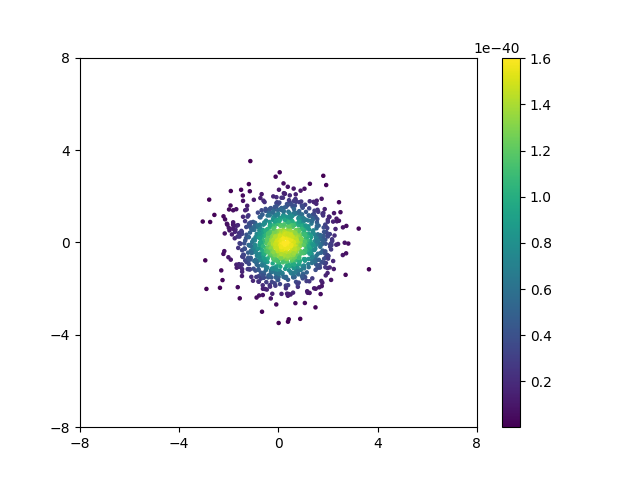}
\hspace{-5mm}
\includegraphics[width=2.4cm]{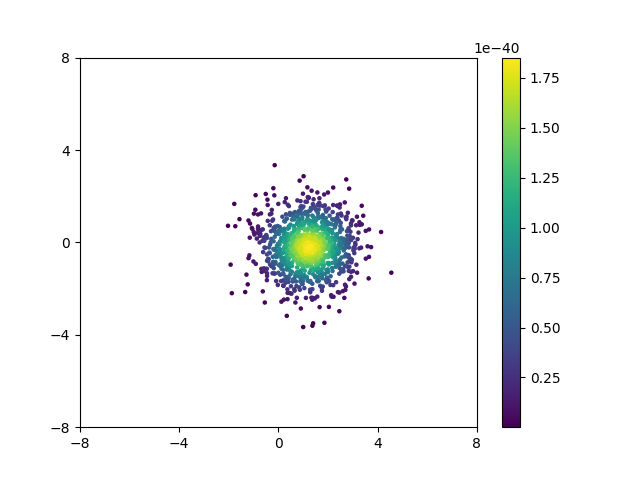}
\hspace{-5mm}
\includegraphics[width=2.4cm]{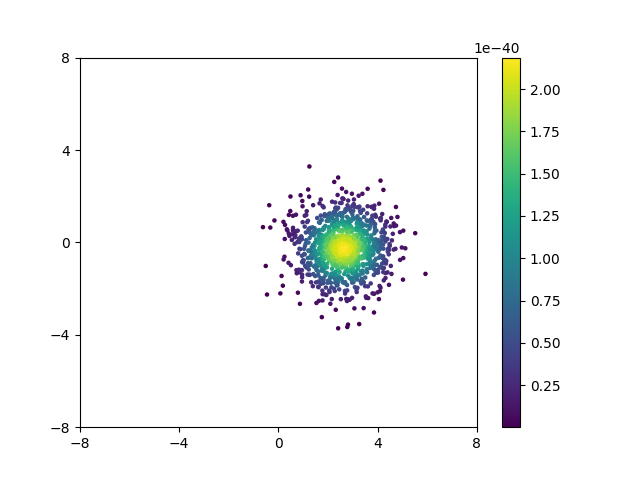}
\hspace{-5mm}
\includegraphics[width=2.4cm]{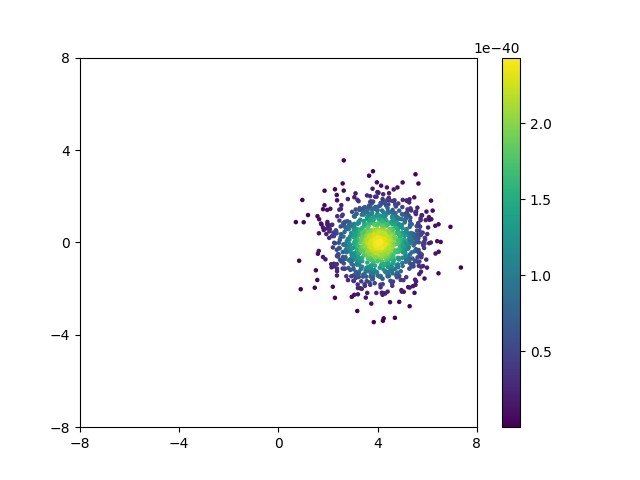}
\hspace{-5mm}
\includegraphics[width=2.4cm]{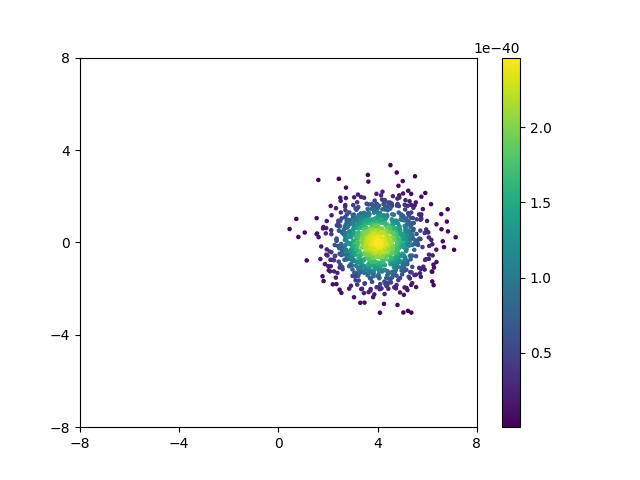}\\
\vspace{5pt}

\includegraphics[width=2.4cm]{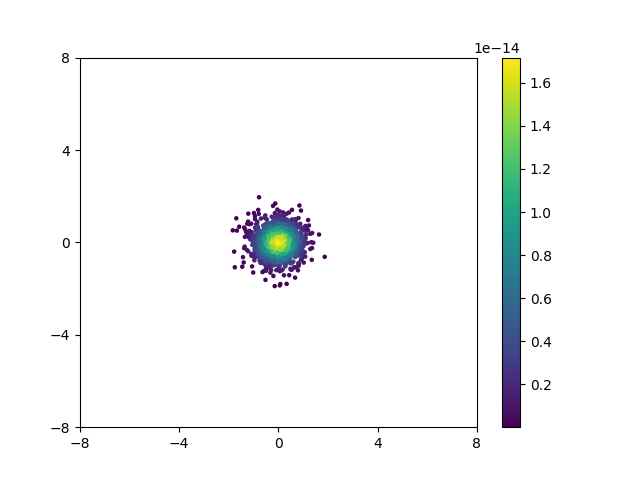}
\hspace{-5mm}
\includegraphics[width=2.4cm]{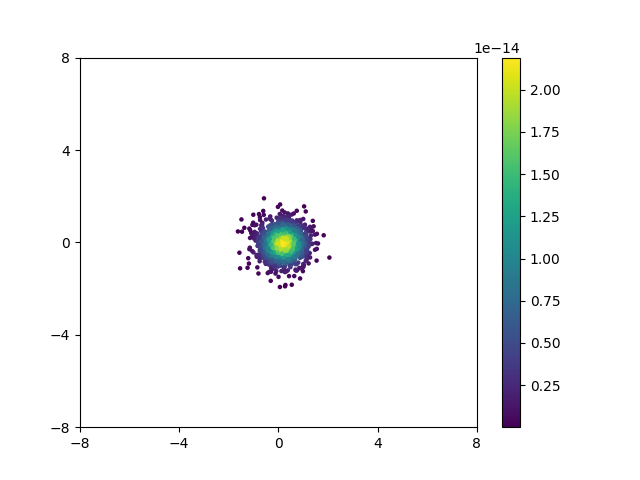}
\hspace{-5mm}
\includegraphics[width=2.4cm]{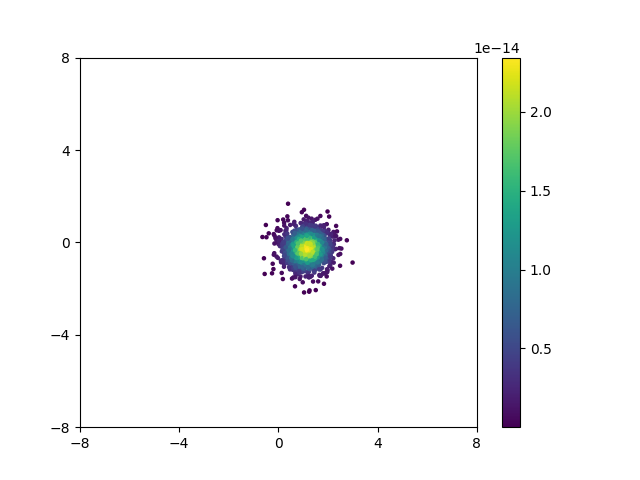}
\hspace{-5mm}
\includegraphics[width=2.4cm]{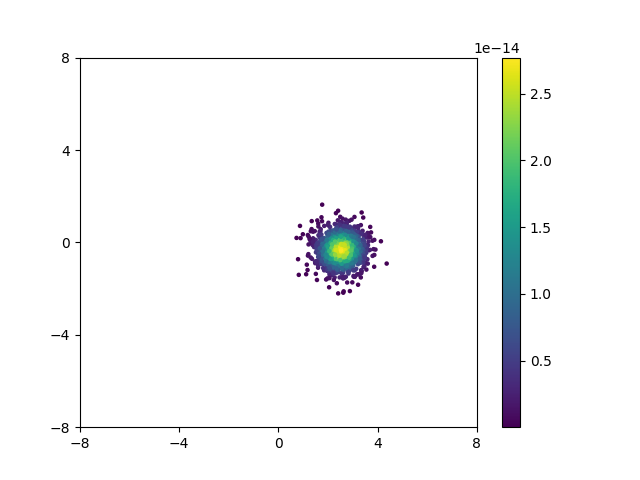}
\hspace{-5mm}
\includegraphics[width=2.4cm]{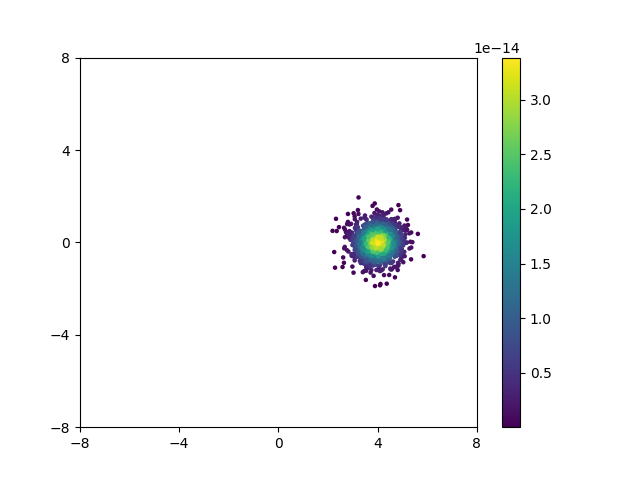}
\hspace{-5mm}
\includegraphics[width=2.4cm]{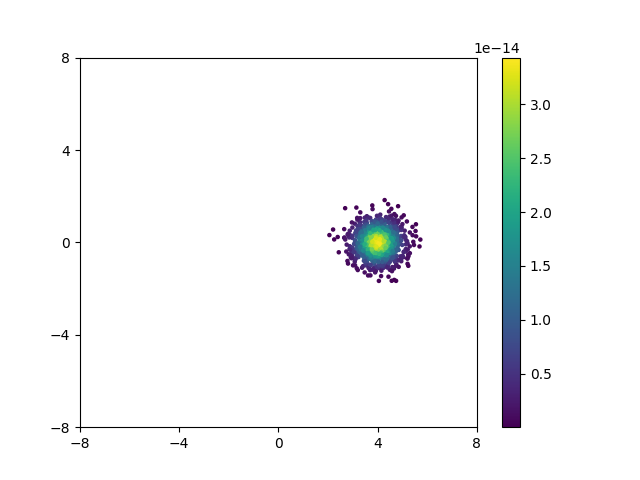}\\
\vspace{5pt}

\subfigure[$\rho_0(\boldsymbol{\boldsymbol{x}}_0)$]{\includegraphics[width=2.4cm]{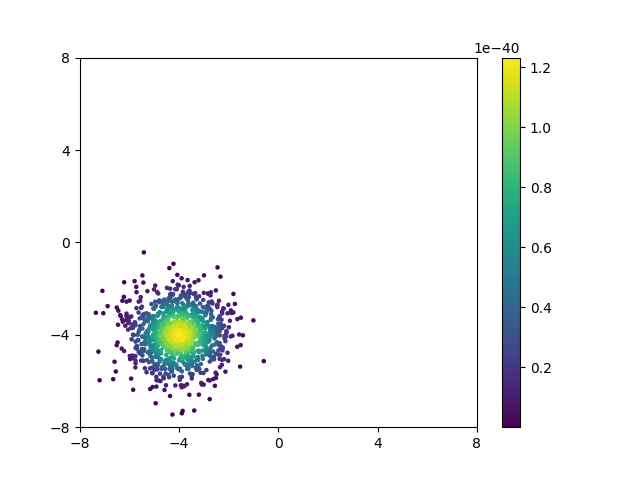}}
\hspace{-5mm}
\subfigure[$\tilde{\rho}_{1/4}(\tilde{\boldsymbol{\boldsymbol{x}}}_{1/4})$]{\includegraphics[width=2.4cm]{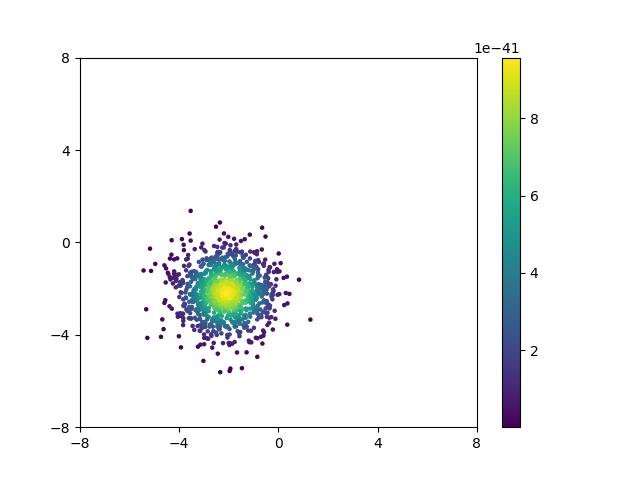}}
\hspace{-5mm}
\subfigure[$\tilde{\rho}_{1/2}(\tilde{\boldsymbol{\boldsymbol{x}}}_{1/2})$]{\includegraphics[width=2.4cm]{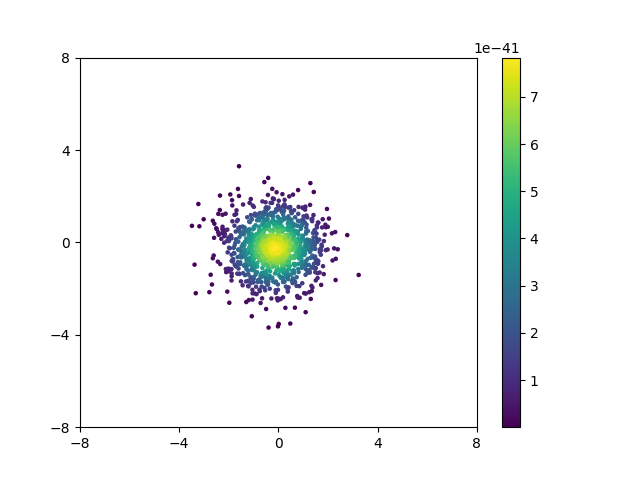}}
\hspace{-5mm}
\subfigure[$\tilde{\rho}_{3/4}(\tilde{\boldsymbol{\boldsymbol{x}}}_{3/4})$]{\includegraphics[width=2.4cm]{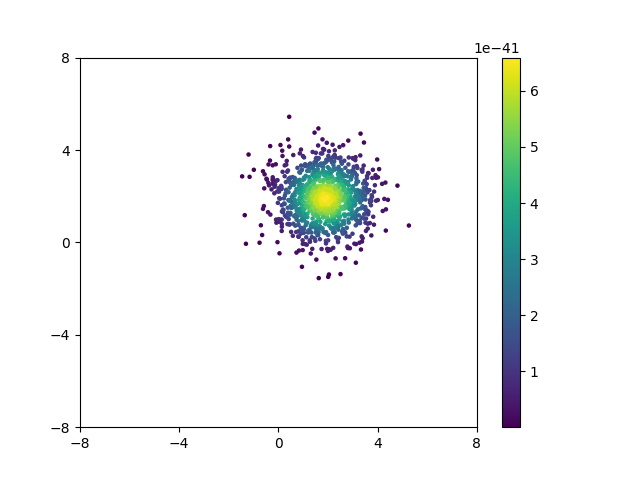}}
\hspace{-5mm}
\subfigure[$\tilde{\rho}_{1}(\tilde{\boldsymbol{\boldsymbol{x}}}_{1})$]{\includegraphics[width=2.4cm]{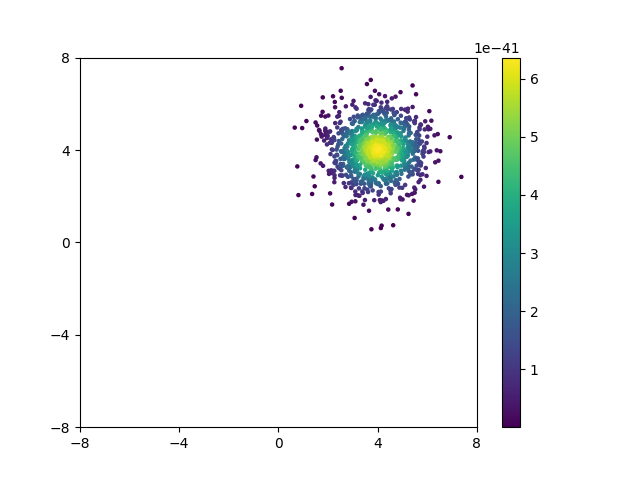}}
\hspace{-5mm}
\subfigure[$\rho_{1}(\boldsymbol{\boldsymbol{x}}_{1})$]{\includegraphics[width=2.4cm]{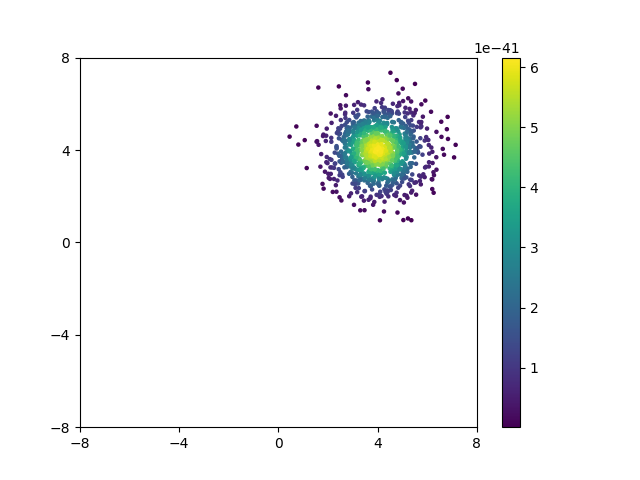}}

\caption{Illustration of these Gaussian problems at $d=100$.}
\label{compare_d100}
\end{center}
\end{figure}

In Figure \ref{compare_d2}, we present the transmission process in the two-dimensional case ($d=2$). 
The first column (a) and the last column (f) show the initial density and the target density, respectively.
The other columns correspond to the push-forward of the initial density at intermediate times $t=0.25, 0.5, 0.75$, and $1.0$, respectively.
We can clearly observe that the target density estimated by our trained network in (e) is very similar to the true target density in (f).
Additionally, the changes in mass during the transmission process can be seen from the variations in the colorbars at different time points. 
It is evident that in Test 5 and Test 8, the maximum value of the colorbar at $t=1$ is approximately half of the maximum value at $t=0$, while in Test 6 and Test 7, the maximum value of the colorbar at $t=1$ is nearly double that at $t=0$. Moreover, the observed unbalanced movement is both expected and reasonable. 
Similar results can be observed across high-dimensional cases, as shown in Figure \ref{compare_d100} for $d=100$ and in the \ref{appA} for $d=10, 30, 50, 80$. For dimensions $d > 2$, we display only slices along the first two coordinate directions.

For better illustration, the GKL loss values for dimensions ranging from 2 to 100 are recorded in Table 1. It is evident that the proposed method consistently yields promising quantitative results, even as dimensionality increases significantly. Additionally, the number of training samples is fixed at 1024 for all tested dimensions.
The running time for each epoch at different dimensions is presented in Table 2.
Since each test is trained for 1000 epochs, the total training time is 1000 times the duration of a single epoch.
The results clearly demonstrate that computational efficiency is not influenced by the increase in dimensionality.

\begin{table}[H]
\newcommand{\tabincell}[2]{\begin{tabular}{@{}#1@{}}#2\end{tabular}} 
\centering
\caption{Results of the GKL loss values for each dimension.}
\begin{tabular}{ccccccc}
\toprule
$d$ &  2 & 10 & 30 & 60 & 80 & 100\\
\midrule
Test 5 & 1.79e-04 & 9.81e-04 & 1.22e-03 & 1.23e-03 & 1.71e-03   & 1.30e-03 \\
Test 6 & 3.59e-04 & 9.94e-04 & 1.86e-03 & 3.09e-03 & 3.15e-03   & 3.19e-03 \\
Test 7 & 6.61e-04 & 1.50e-03 & 2.62e-03 & 2.92e-03 & 3.09e-03   & 2.77e-03 \\
Test 8 & 3.48e-03 & 1.92e-03 & 1.75e-03 & 2.72e-03 & 1.36e-03   & 2.23e-03 \\
\bottomrule
\end{tabular}
\end{table}

\begin{table}[H]
\newcommand{\tabincell}[2]{\begin{tabular}{@{}#1@{}}#2\end{tabular}} 
\centering
\caption{Running time (s) per epoch for different dimensions.}
\begin{tabular}{ccccccc}
\toprule
$d$ &  2 & 10 & 30 & 60 & 80 & 100\\
\midrule
Test 5 & 1.6581 &  1.7199  &  1.7099
 & 1.7073  &  1.6877 &  1.7286 \\
Test 6 & 1.7006	& 1.6982 & 1.6942 &	1.7057 & 1.7069 & 1.7395 \\
Test 7 & 1.7859	& 1.7476 & 1.7007 & 1.7379 & 1.6973 & 1.7553 \\
Test 8 & 1.7031	& 1.6661 & 1.7048 & 1.7298 & 1.7663 & 1.7985 \\
\bottomrule
\end{tabular}
\end{table}

\begin{figure}[t]
    \centering
    \includegraphics[width=0.8\textwidth]{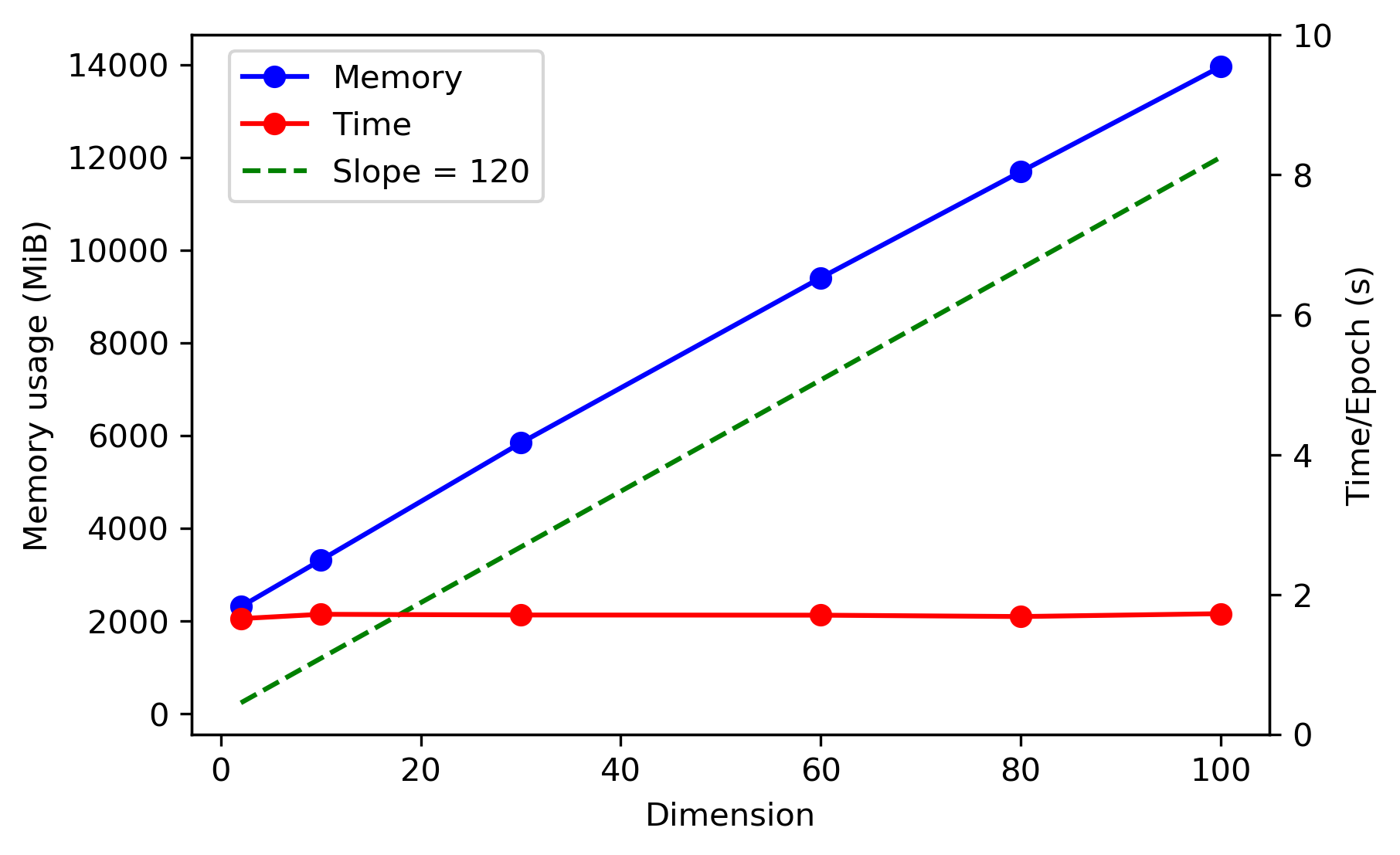} 
    \caption{Memory usage and iteration time across dimensions (2-100) with 1024 samples.}
    \label{fig:sample-image}
\end{figure}

Furthermore, we present the running time per epoch and memory usage for dimensions ranging from 2 to 100 in Figure \ref{fig:sample-image}. As shown in the figure, memory usage increases linearly with dimensionality. 
Therefore, with sufficient memory, this proposed method should remain effective for higher dimensions.
To test this, we employed the NVIDIA H800 PCIe GPU with 80GB of memory to conduct experiments at dimensions of $d=300$ and $d=400$. For both high-dimensional cases, we set $L=1$ and use 2048 samples for $d=300$ and 3096 samples for $d=400$, respectively. The corresponding results are presented in \ref{appB}.

\subsection{Gaussian mixtures}

In this section, we consider UOT problem on two Gaussian mixtures examples as listed below. 

\begin{itemize}\label{list3}

    \item Test 9: $\rho_0(\boldsymbol{x})=\rho_G(\boldsymbol{x},\mathbf{0},\mathbf{I})$, $\rho_1(\boldsymbol{x})=\rho_G(\boldsymbol{x},-2 \cdot \mathbf{e_1},\mathbf{I}) + \rho_G(\boldsymbol{x},2 \cdot \mathbf{e_1},\mathbf{I})$;
    \item Test 10: $\rho_0(\boldsymbol{x})=\rho_G(\boldsymbol{x},-2 \cdot \mathbf{e_1},\mathbf{I}) + \rho_G(\boldsymbol{x},2 \cdot \mathbf{e_1},\mathbf{I})$, $\rho_1(\boldsymbol{x})=\rho_G(\boldsymbol{x},\mathbf{0},\mathbf{I})$.
 
\end{itemize}

\begin{figure}[H]
\begin{center}
\includegraphics[width=2.4cm]{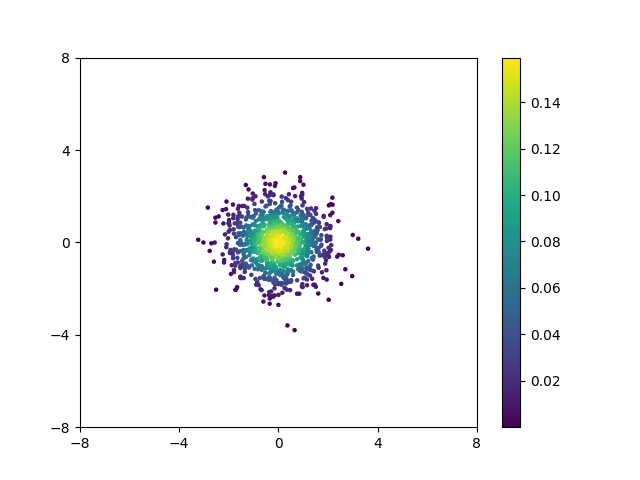}
\hspace{-5mm}
\includegraphics[width=2.4cm]{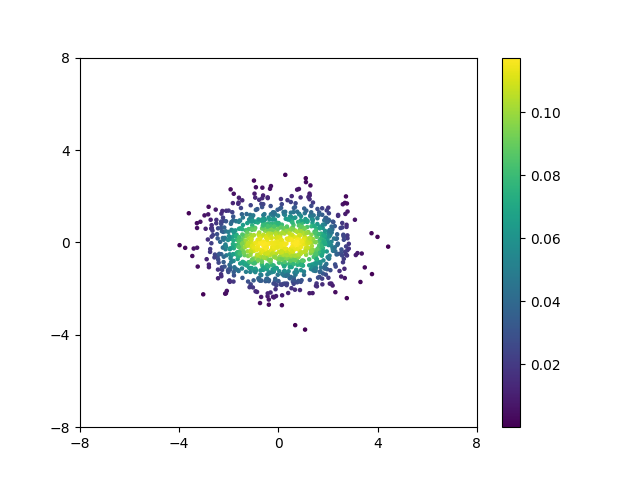}
\hspace{-5mm}
\includegraphics[width=2.4cm]{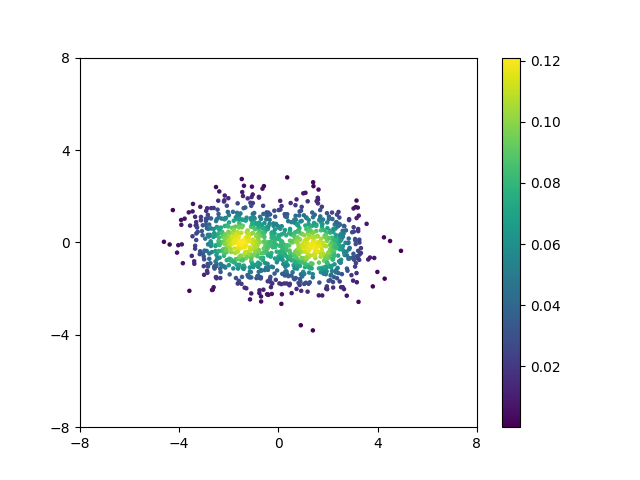}
\hspace{-5mm}
\includegraphics[width=2.4cm]{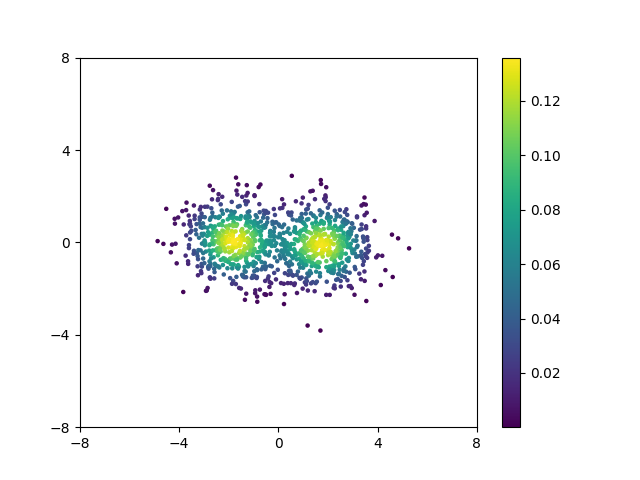}
\hspace{-5mm}
\includegraphics[width=2.4cm]{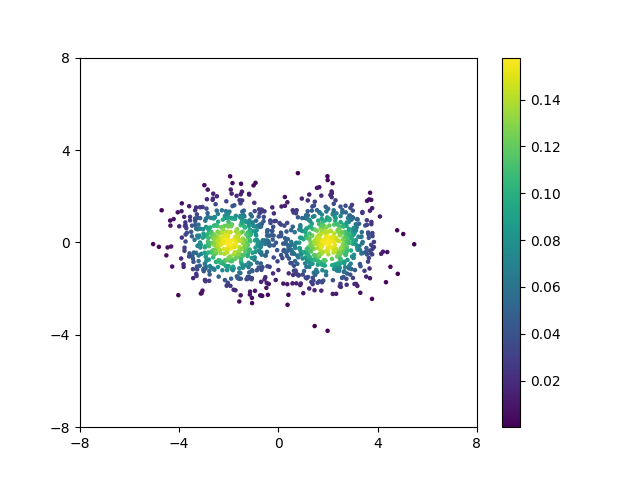}
\hspace{-5mm}
\includegraphics[width=2.4cm]{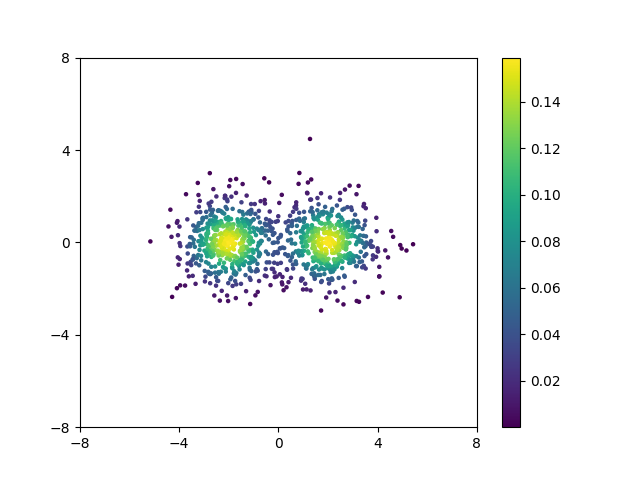}\\
\vspace{5pt}

\subfigure[$\rho_0(\boldsymbol{\boldsymbol{x}}_0)$]{\includegraphics[width=2.4cm]{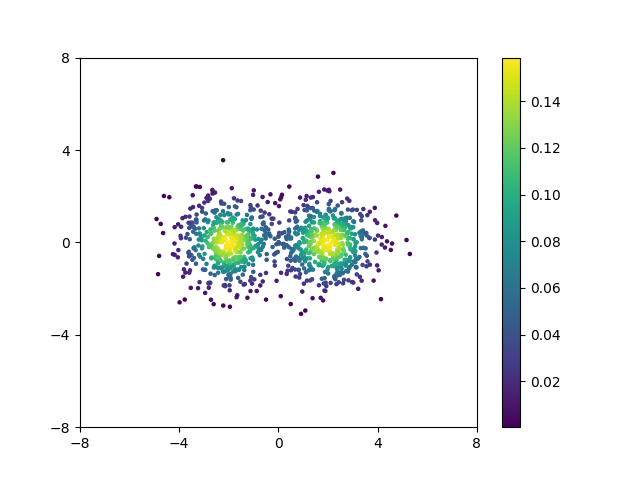}}
\hspace{-5mm}
\subfigure[$\tilde{\rho}_{1/4}(\tilde{\boldsymbol{\boldsymbol{x}}}_{1/4})$]{\includegraphics[width=2.4cm]{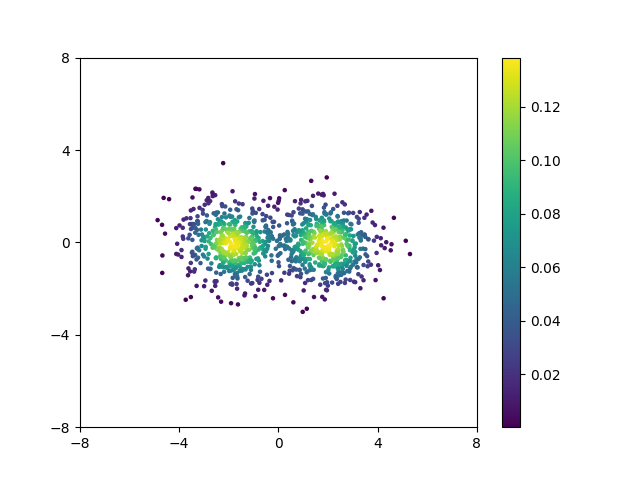}}
\hspace{-5mm}
\subfigure[$\tilde{\rho}_{1/2}(\tilde{\boldsymbol{\boldsymbol{x}}}_{1/2})$]{\includegraphics[width=2.4cm]{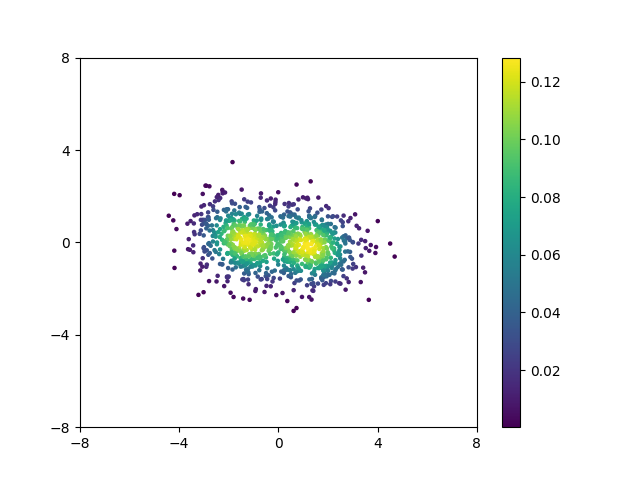}}
\hspace{-5mm}
\subfigure[$\tilde{\rho}_{3/4}(\tilde{\boldsymbol{\boldsymbol{x}}}_{3/4})$]{\includegraphics[width=2.4cm]{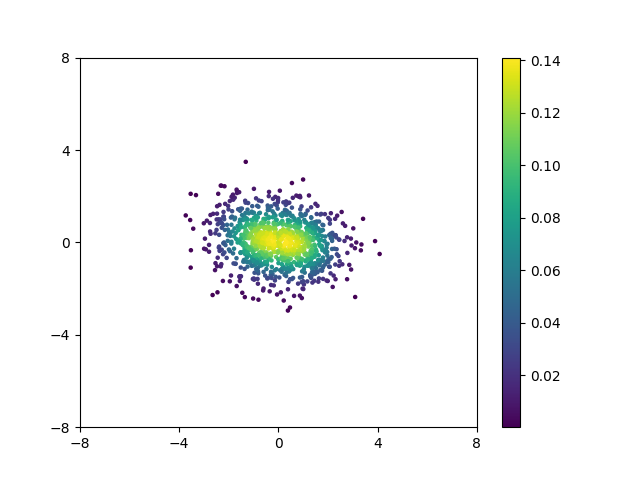}}
\hspace{-5mm}
\subfigure[$\tilde{\rho}_{1}(\tilde{\boldsymbol{\boldsymbol{x}}}_{1})$]{\includegraphics[width=2.4cm]{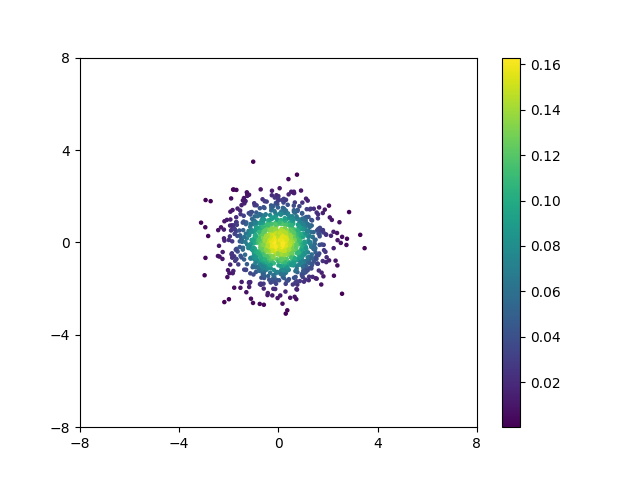}}
\hspace{-5mm}
\subfigure[$\rho_{1}(\boldsymbol{\boldsymbol{x}}_{1})$]{\includegraphics[width=2.4cm]{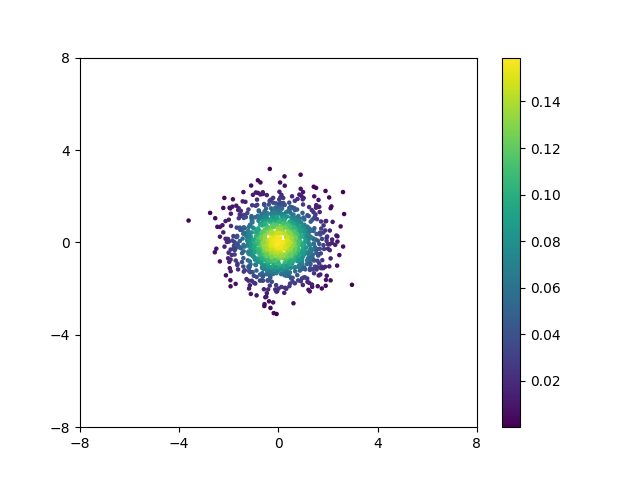}}

\caption{Illustration of these Gaussian mixtures problems at $d=2$.}
\label{gaussianmix_d2}
\end{center}
\end{figure}

\begin{figure}[H]
\begin{center}
\includegraphics[width=2.4cm]{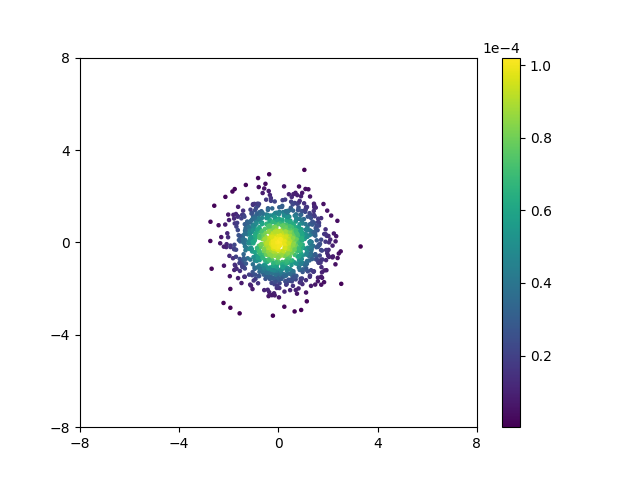}
\hspace{-5mm}
\includegraphics[width=2.4cm]{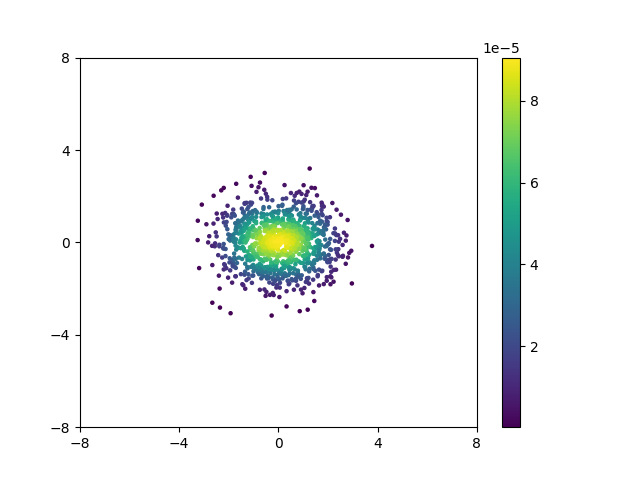}
\hspace{-5mm}
\includegraphics[width=2.4cm]{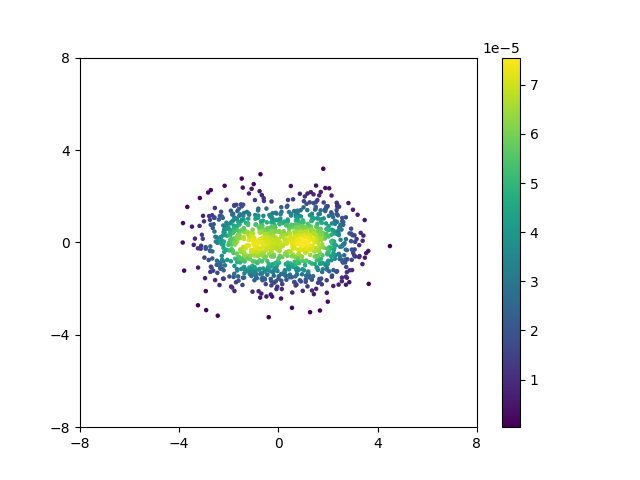}
\hspace{-5mm}
\includegraphics[width=2.4cm]{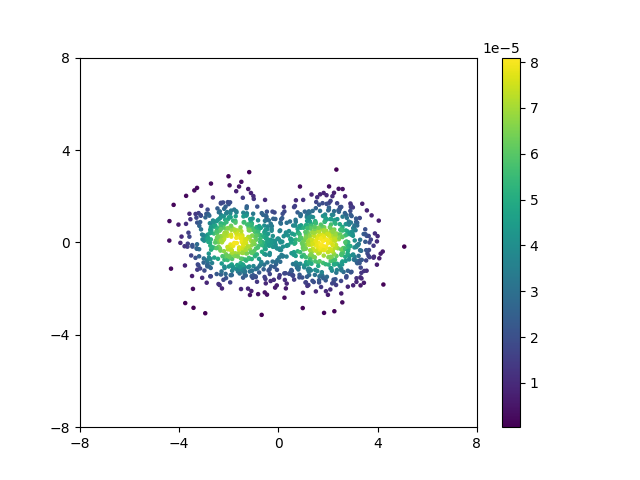}
\hspace{-5mm}
\includegraphics[width=2.4cm]{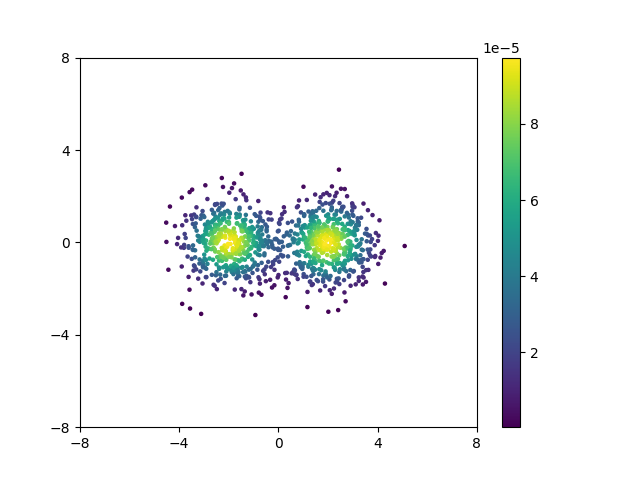}
\hspace{-5mm}
\includegraphics[width=2.4cm]{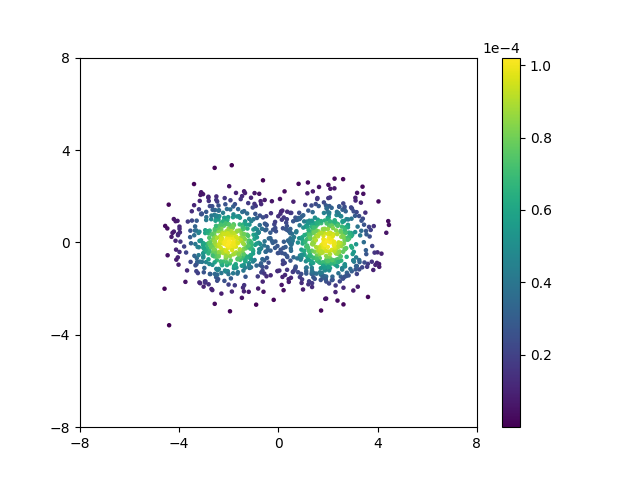}\\
\vspace{5pt}

\subfigure[$\rho_0(\boldsymbol{\boldsymbol{x}}_0)$]{\includegraphics[width=2.4cm]{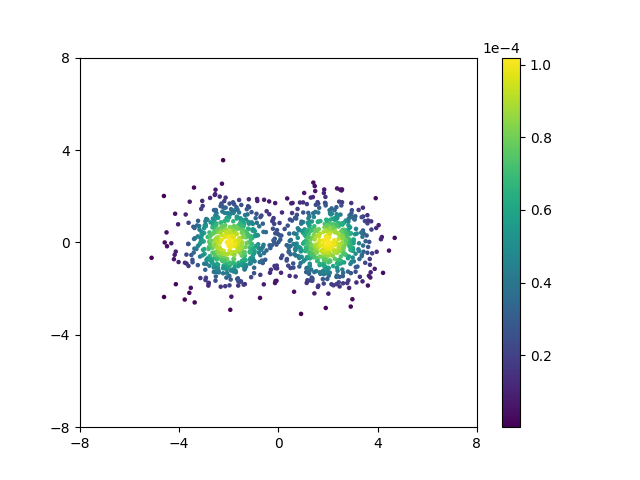}}
\hspace{-5mm}
\subfigure[$\tilde{\rho}_{1/4}(\tilde{\boldsymbol{\boldsymbol{x}}}_{1/4})$]{\includegraphics[width=2.4cm]{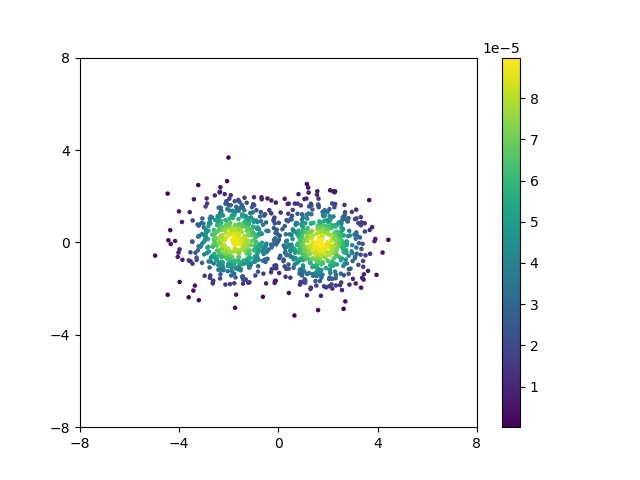}}
\hspace{-5mm}
\subfigure[$\tilde{\rho}_{1/2}(\tilde{\boldsymbol{\boldsymbol{x}}}_{1/2})$]{\includegraphics[width=2.4cm]{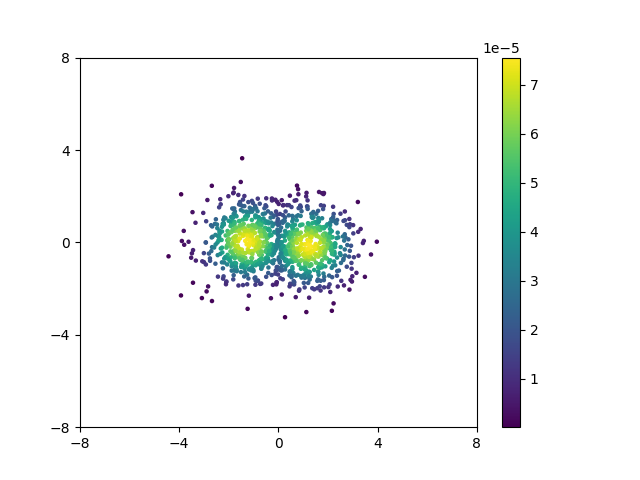}}
\hspace{-5mm}
\subfigure[$\tilde{\rho}_{3/4}(\tilde{\boldsymbol{\boldsymbol{x}}}_{3/4})$]{\includegraphics[width=2.4cm]{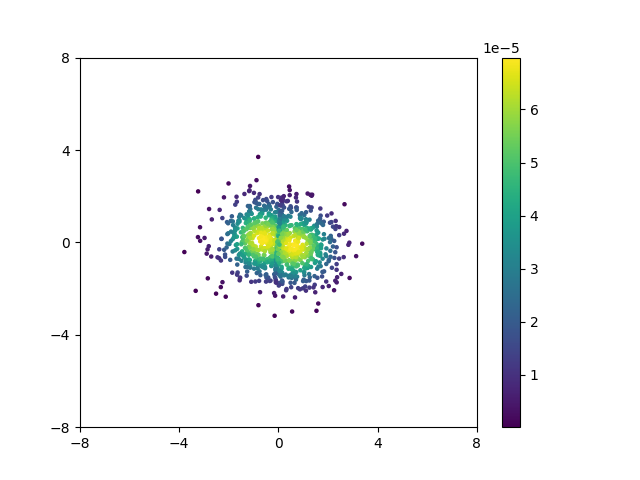}}
\hspace{-5mm}
\subfigure[$\tilde{\rho}_{1}(\tilde{\boldsymbol{\boldsymbol{x}}}_{1})$]{\includegraphics[width=2.4cm]{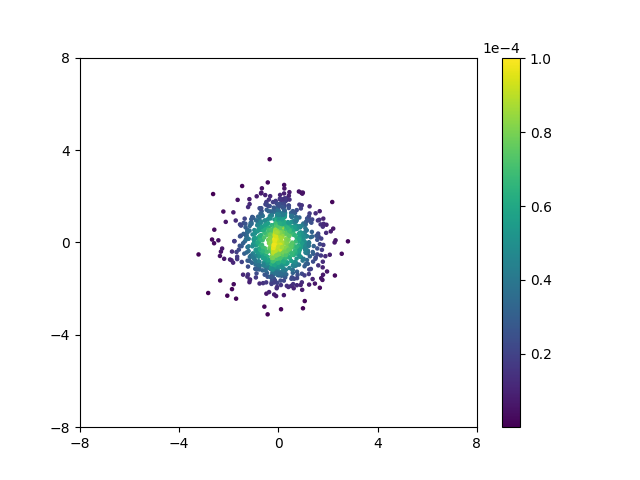}}
\hspace{-5mm}
\subfigure[$\rho_{1}(\boldsymbol{\boldsymbol{x}}_{1})$]{\includegraphics[width=2.4cm]{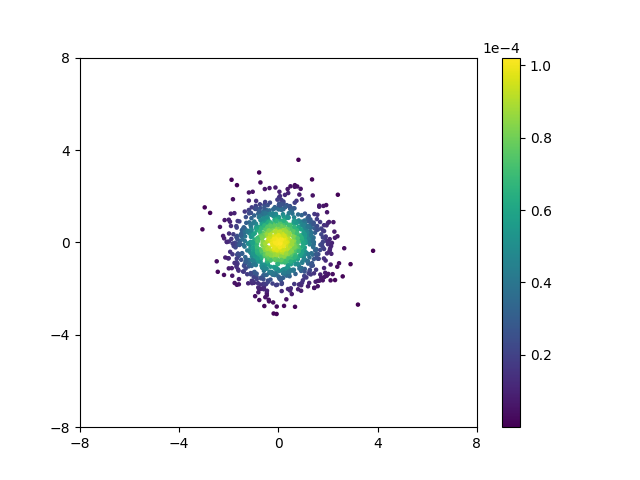}}

\caption{Illustration of these Gaussian mixtures problems at $d=10$.}
\label{gaussianmix_d10}
\end{center}
\end{figure}

In Figure \ref{gaussianmix_d2} and Figure \ref{gaussianmix_d10}, we demonstrate that the proposed method effectively achieves the splitting from one to two (Test 9) and the merging from two to one (Test 10), along with the creation and destruction in mass. 
Furthermore, we can clearly observe that the target density estimated by our trained network in (e) closely resembles the target density $\rho_1$ in (f), even in high-dimensional instances.

\subsection{Generalization on crowd motion}\label{cro}

The proposed method of UOT can be easily extended to address crowd motion problems with different total mass.
For example, in the application of traffic flow, some vehicles enter and exit the flow while traveling from one location to another, resulting in different initial and terminal vehicle counts. Similarly, in the population migration of animals, there are births of new life and deaths of weaker individuals, resulting in differences between the initial and final population counts.

In contrast to the optimal transport problems, the crowd motion problem involves not only approximating the transition from the initial density $\rho_{0}$ to the target density $\rho_{1}$ but also ensuring that obstacles are avoided during dynamic movement. 
These obstacles can be incorporated into the objective function by introducing an additional term, referred to as the preference term:
\begin{align}\label{prefernceterm}
\mathcal{Q}(\rho) = \int_0^1\int_{\Omega} Q(\boldsymbol{x})\rho(\boldsymbol{x},t)\mathrm{d}\boldsymbol{x} \mathrm{d}t,
\end{align}
where $Q:\Omega  \rightarrow \mathbb{R}$ models the spatial preferences of agents. A higher value of $Q(\boldsymbol{x})$ indicates a position closer to obstacles. 

The optimization problem can then be formulated as:
\begin{align*}
    \min_{(\rho,\boldsymbol{v},f) \in \mathcal{C}(\rho_0,\rho_1)} \mathcal{E}(\rho, \boldsymbol{v})  + \frac{1}{\alpha}\mathcal{R}(\rho,f) + \lambda_P \mathcal{Q}(\rho),
\end{align*}
where $\lambda_P$ is a penalty parameter, and $\mathcal{E}(\rho, \boldsymbol{v}), \mathcal{R}(\rho,f)$ and $\mathcal{C}(\rho_0,\rho_1)$ are defined in \eqref{dOT} and \eqref{UOTcons}.
The proposed method can naturally be applied to solve the crowd motion problem.
The semi-discrete version of $\mathcal{Q}(\rho)$ in the spatial direction is expressed as:
\begin{align*}
    \overline{\mathcal{Q}} (\rho)&:= \int_0^1 \sum_{i=1}^r \frac{1}{r}Q(\mathbf{z}(\boldsymbol{x}_i,t))\frac{\rho(\mathbf{z}(\boldsymbol{x}_i,t),t)}{\mu(\mathbf{z}(\boldsymbol{x}_i,t),t)} \mathrm{d}t.
\end{align*}

Our crowd motion experiment involves two maze examples, where the initial density $\rho_0$ and target density $\rho_1$ have different total masses.
\begin{itemize}\label{list4}

    \item Test 11: $\rho_0(\boldsymbol{x})=\rho_G(\boldsymbol{x},8\cdot \mathbf{e_1}+4 \cdot \mathbf{e_2}, \mathbf{I})$, $\rho_1(\boldsymbol{x})=\frac{1}{2}\rho_G(\boldsymbol{x},8\cdot \mathbf{e_1}+12 \cdot \mathbf{e_2}, 0.3 \cdot \mathbf{I})$;
    \item Test 12: $\rho_0(\boldsymbol{x})=\rho_G(\boldsymbol{x},-4\cdot \mathbf{e_1}-4 \cdot \mathbf{e_2}, \mathbf{I})$, $\rho_1(\boldsymbol{x})=\frac{1}{2}\rho_G(\boldsymbol{x},4\cdot \mathbf{e_1}+4 \cdot \mathbf{e_2}, 0.3 \cdot \mathbf{I})$.
 
\end{itemize}

\begin{figure}[H]
\begin{center}
\includegraphics[width=1.9cm]{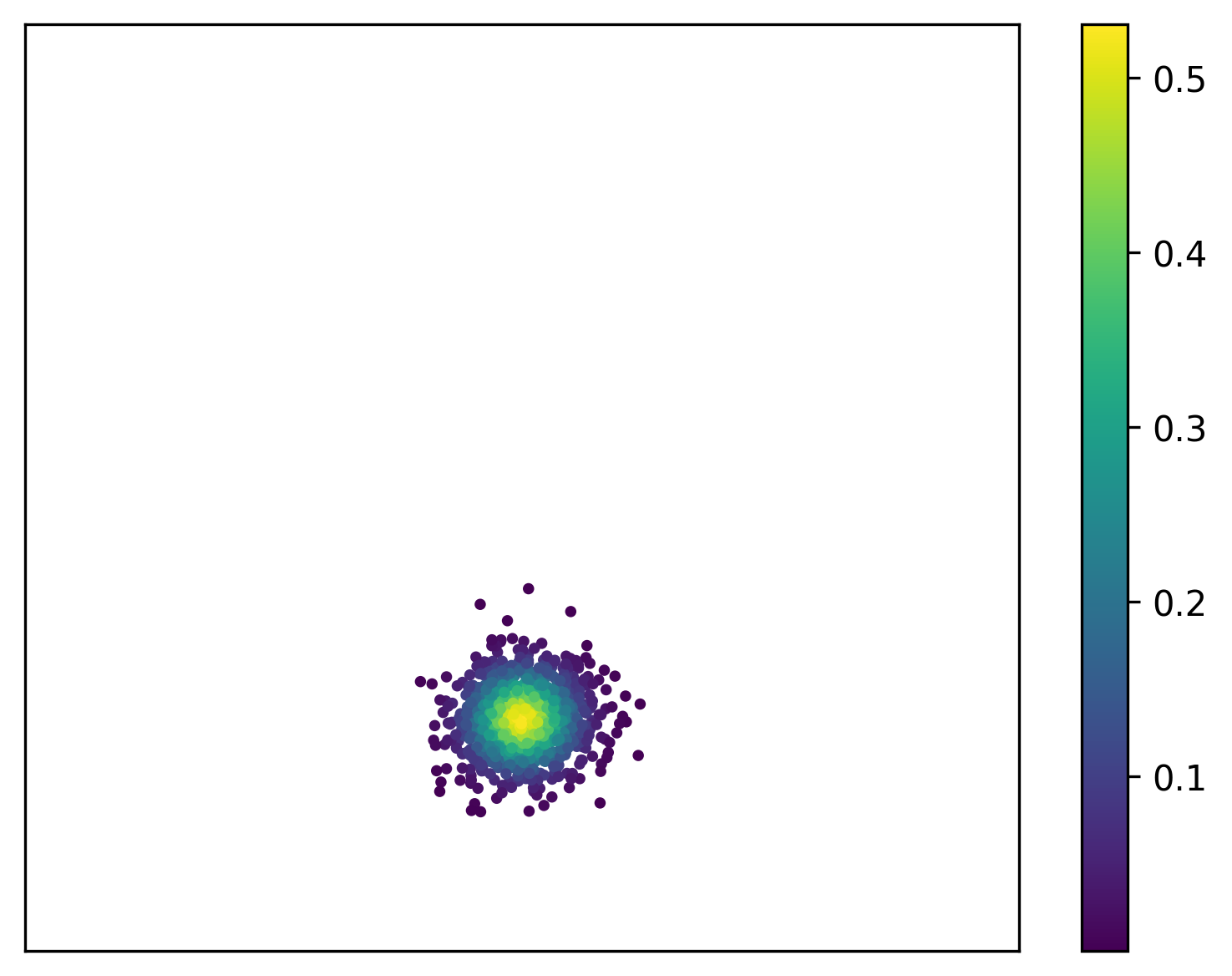}
\includegraphics[width=1.9cm]{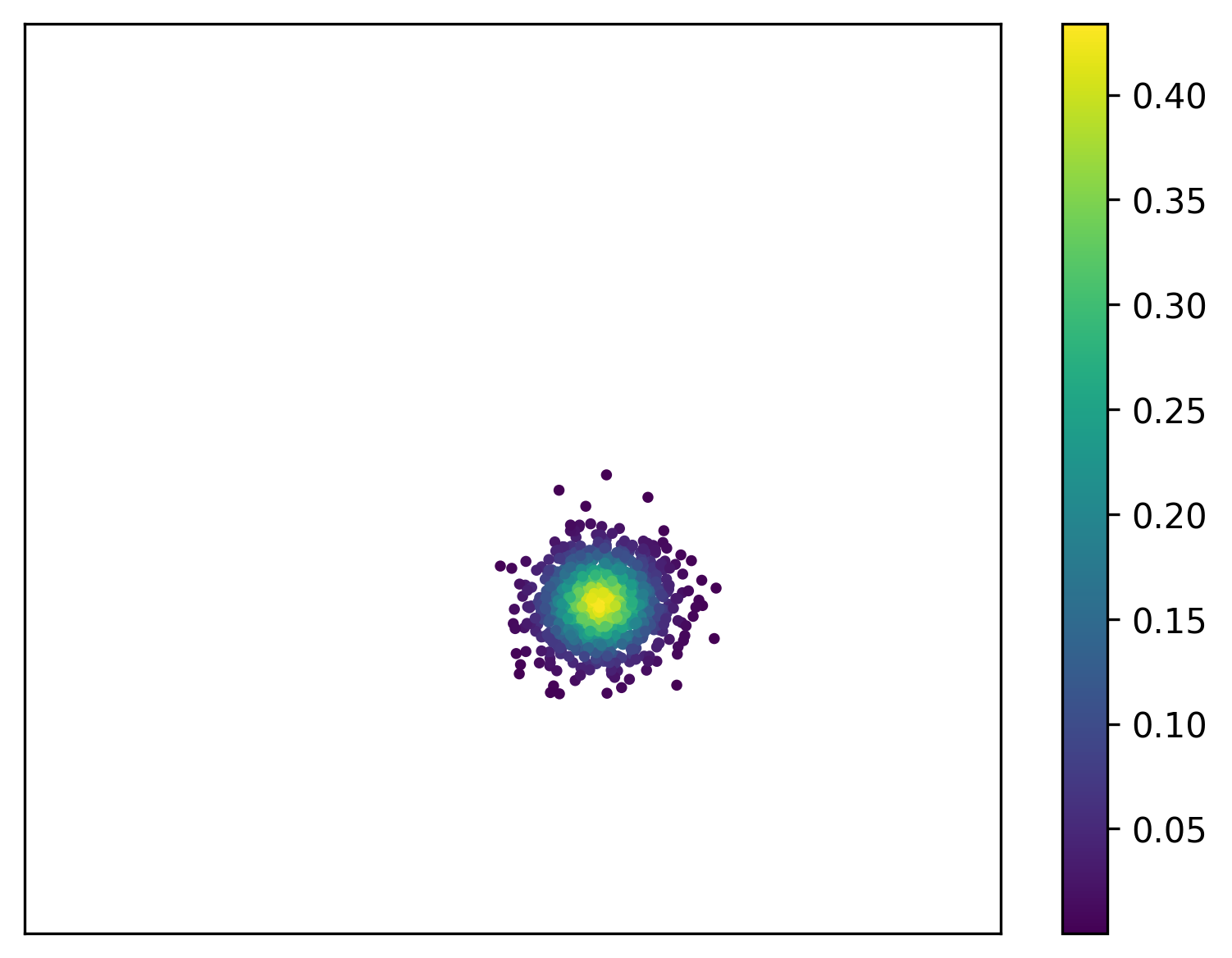}
\includegraphics[width=1.9cm]{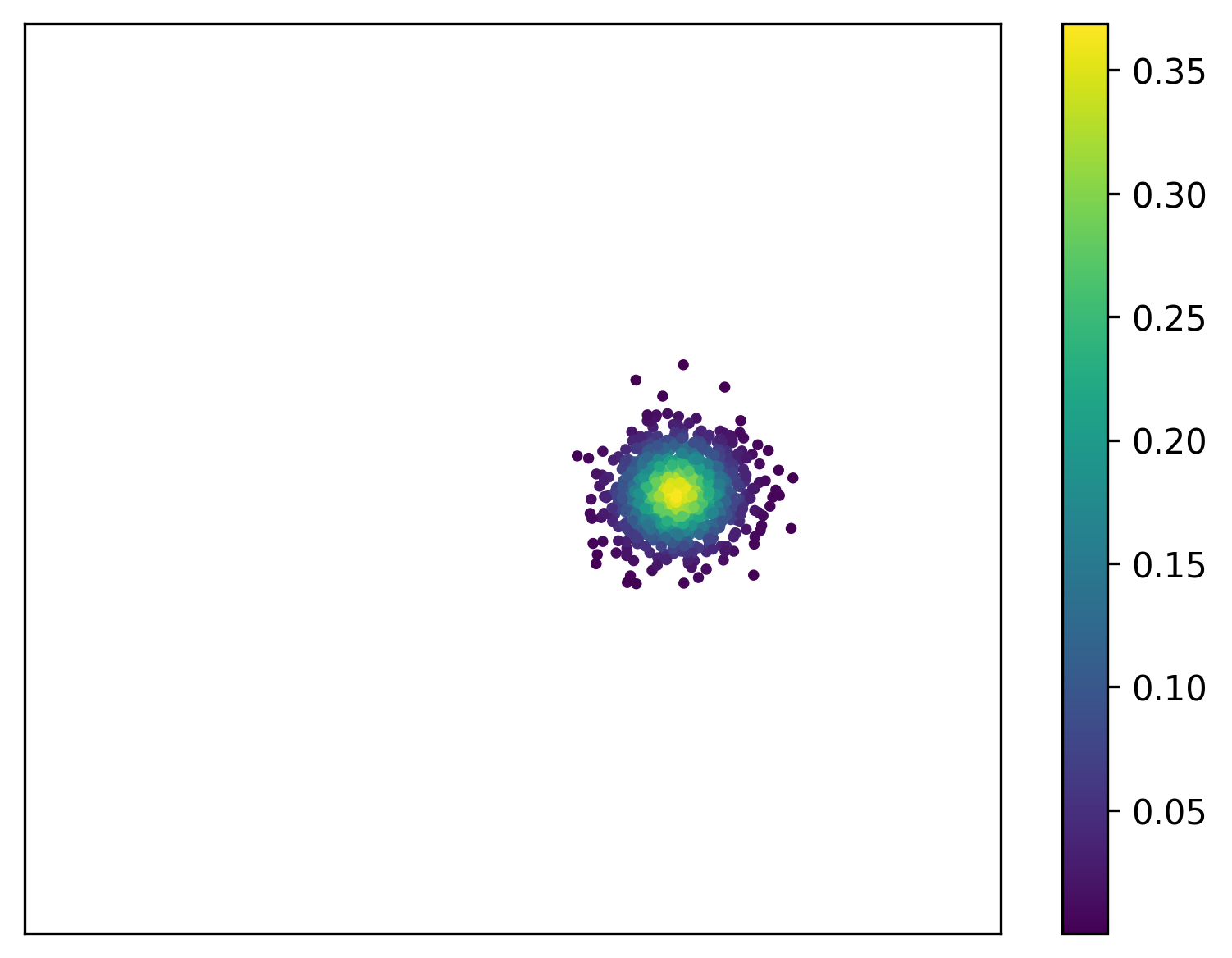}
\includegraphics[width=1.9cm]{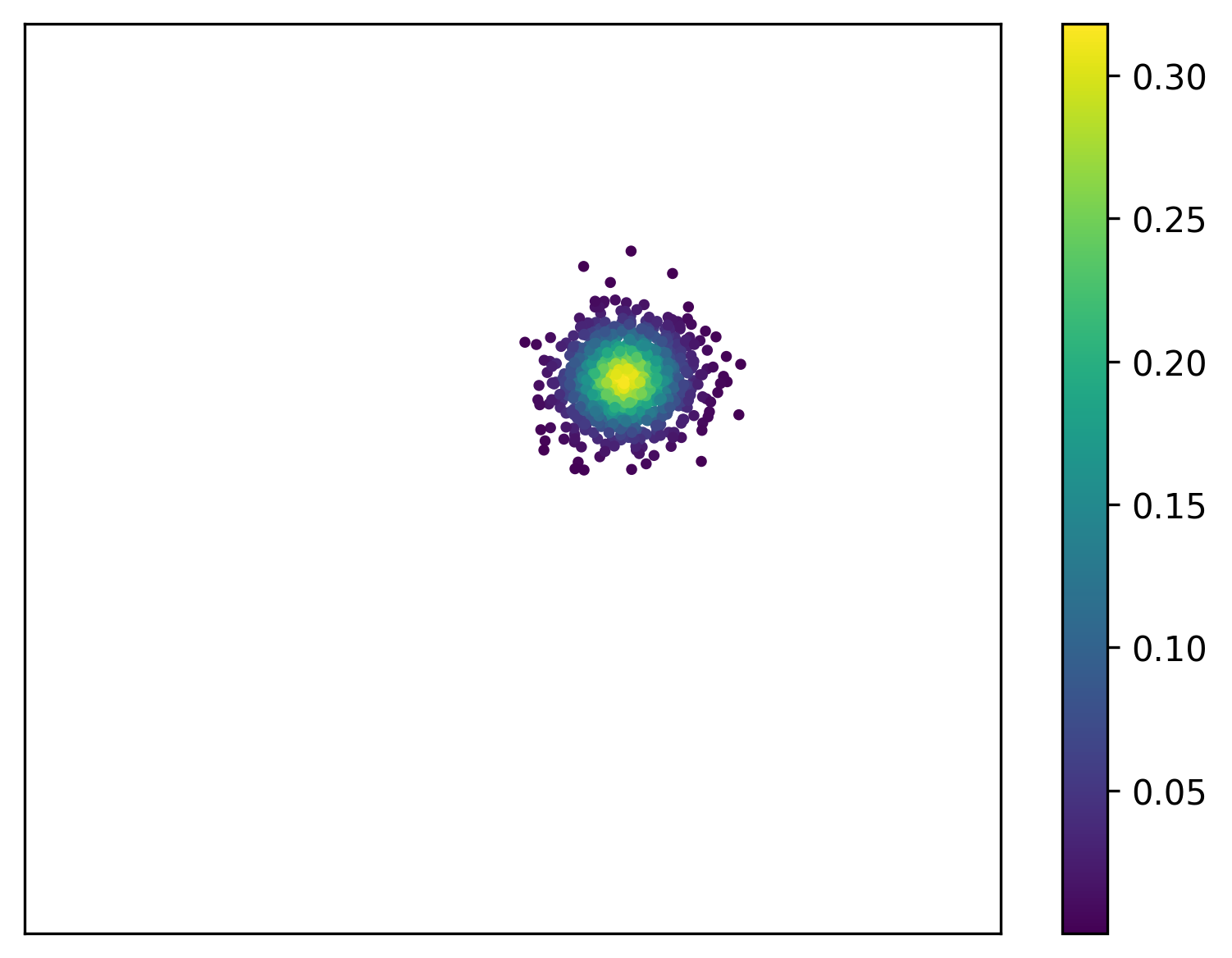}
\includegraphics[width=1.9cm]{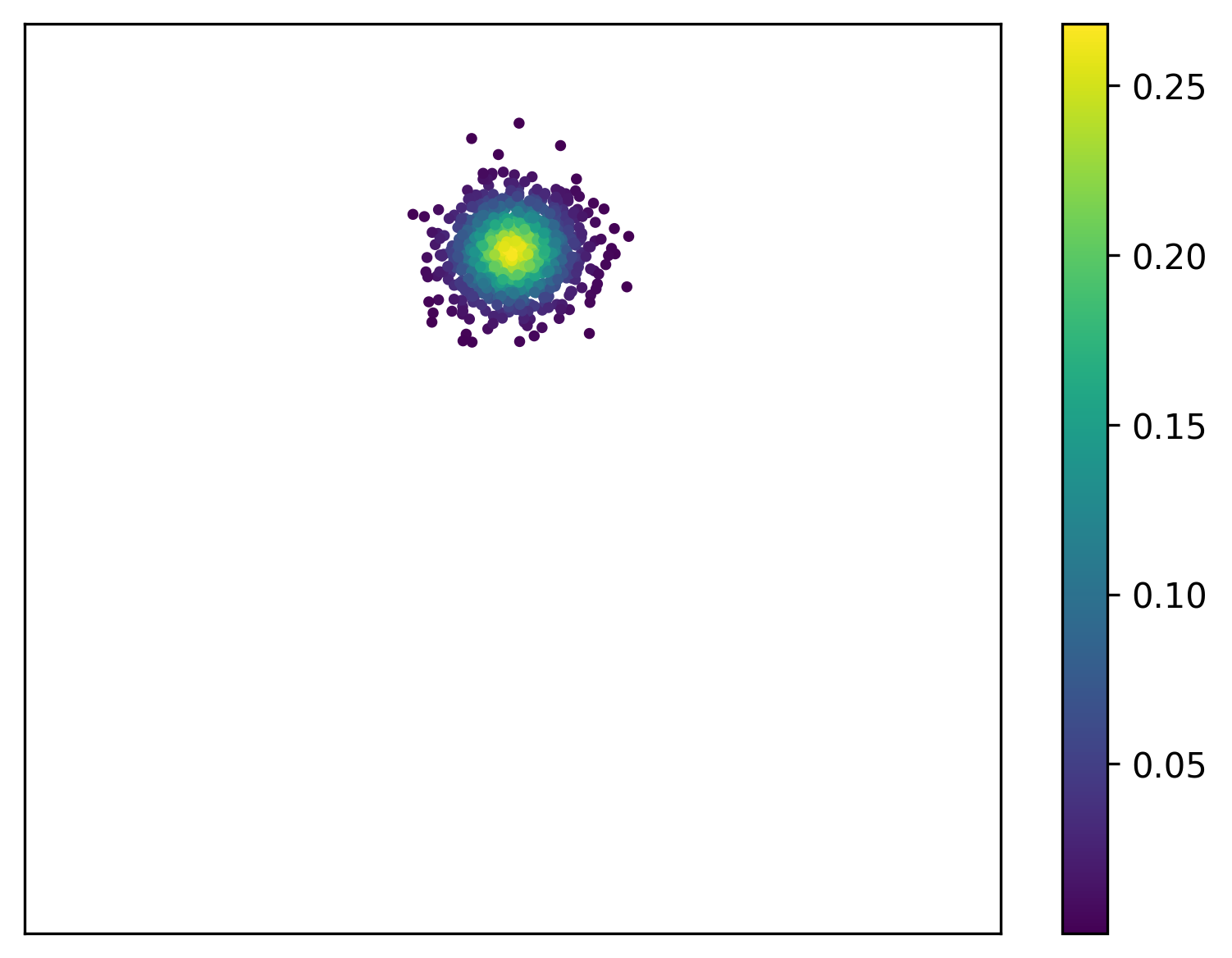}
\includegraphics[width=1.9cm]{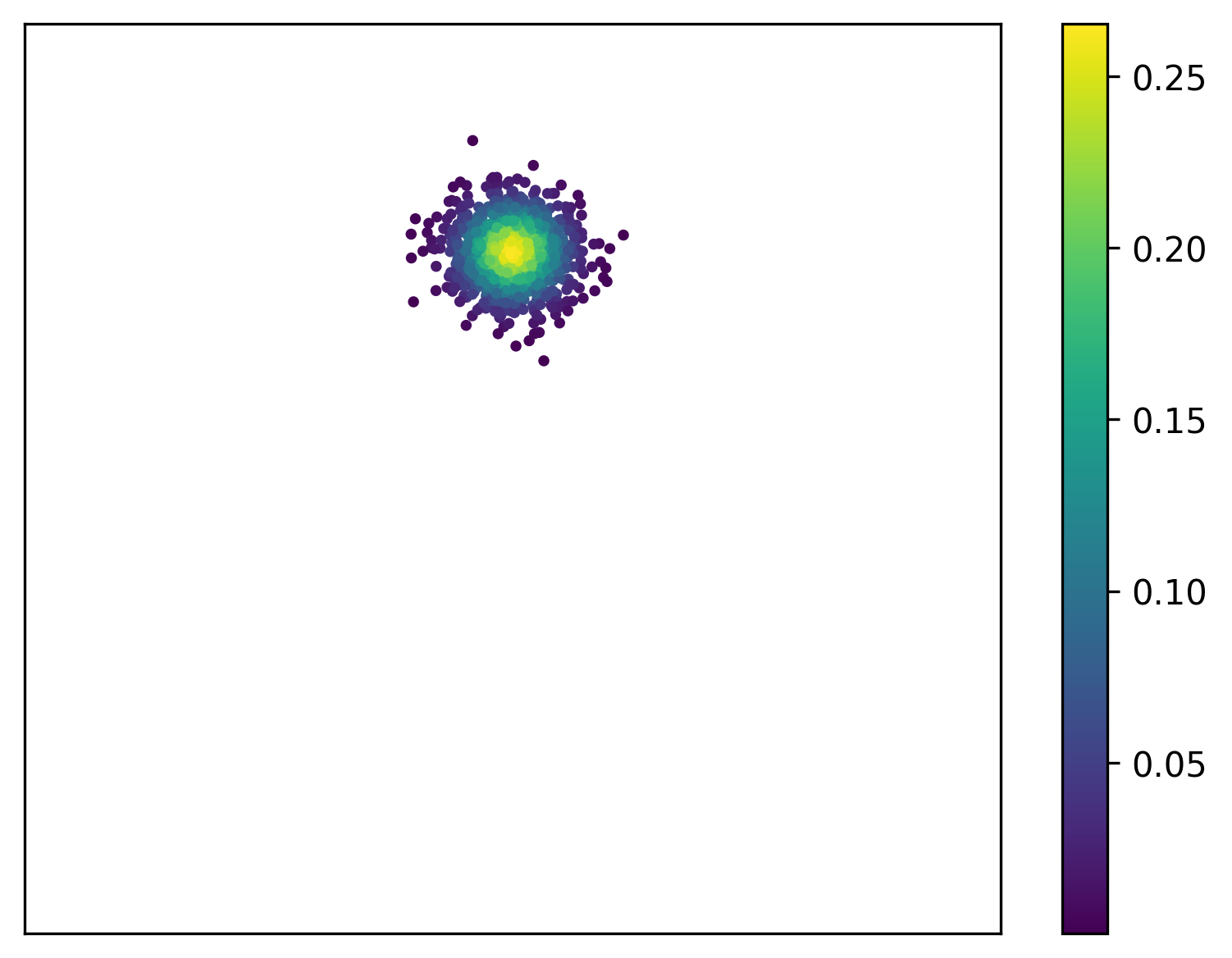}
\includegraphics[width=1.5cm]{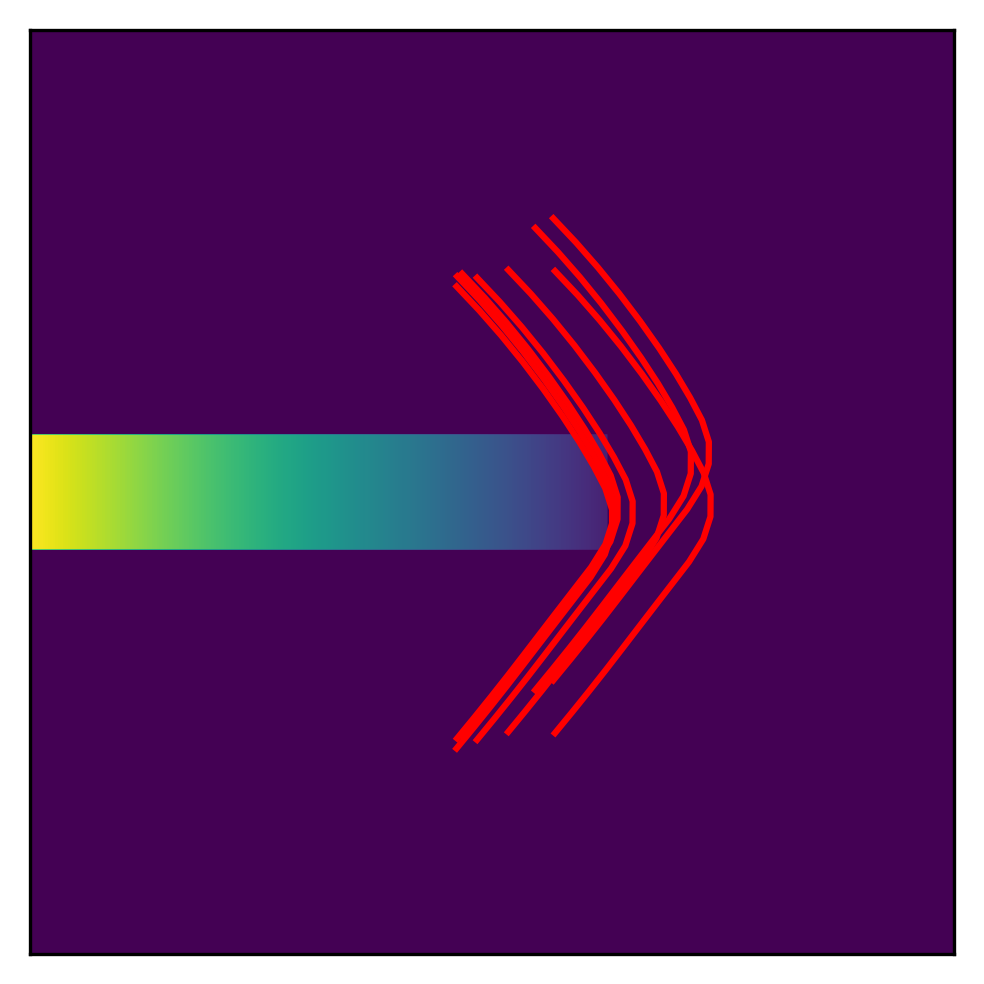}
\vspace{5pt}

\subfigure[]{\includegraphics[width=1.9cm]{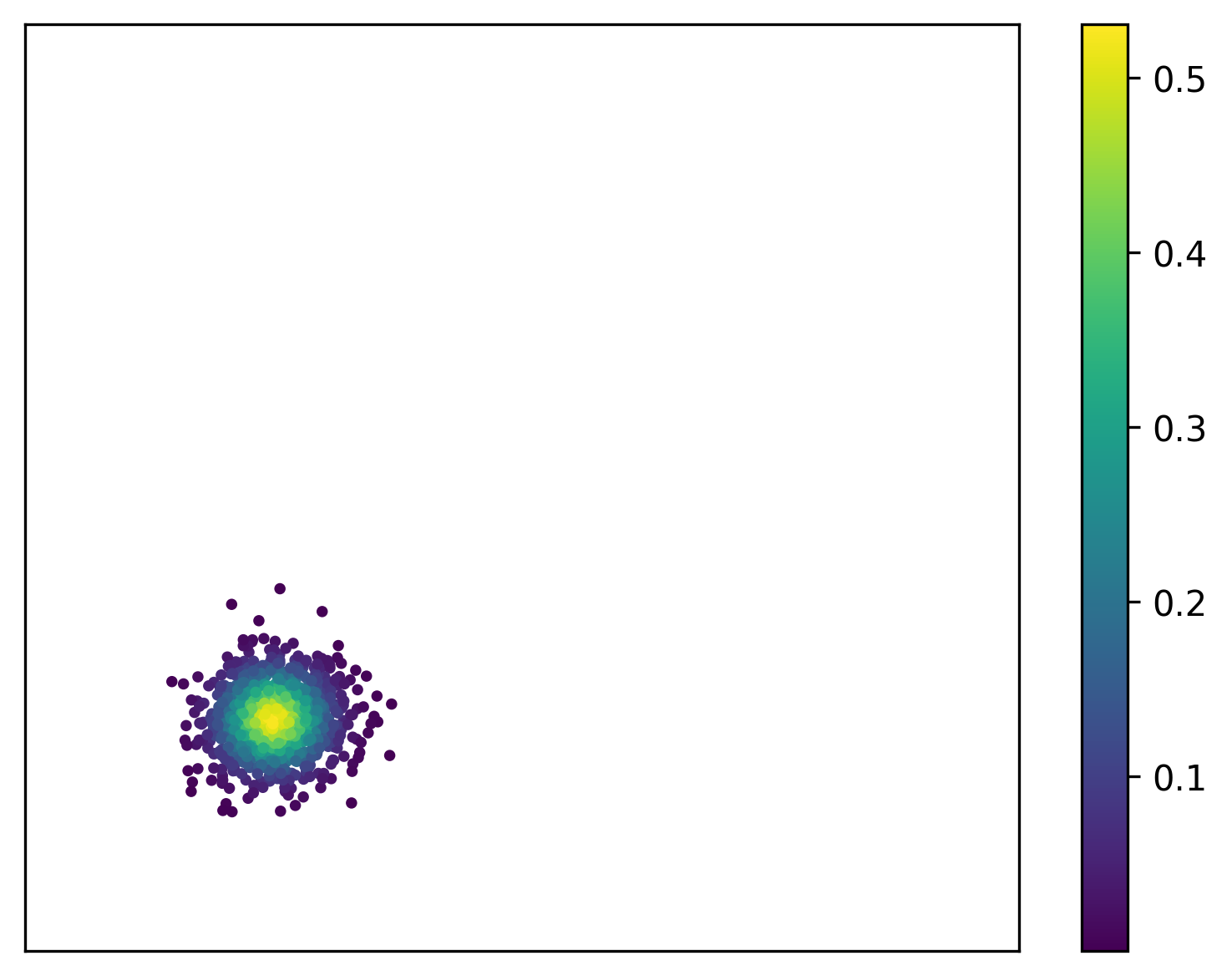}}
\subfigure[]{\includegraphics[width=1.9cm]{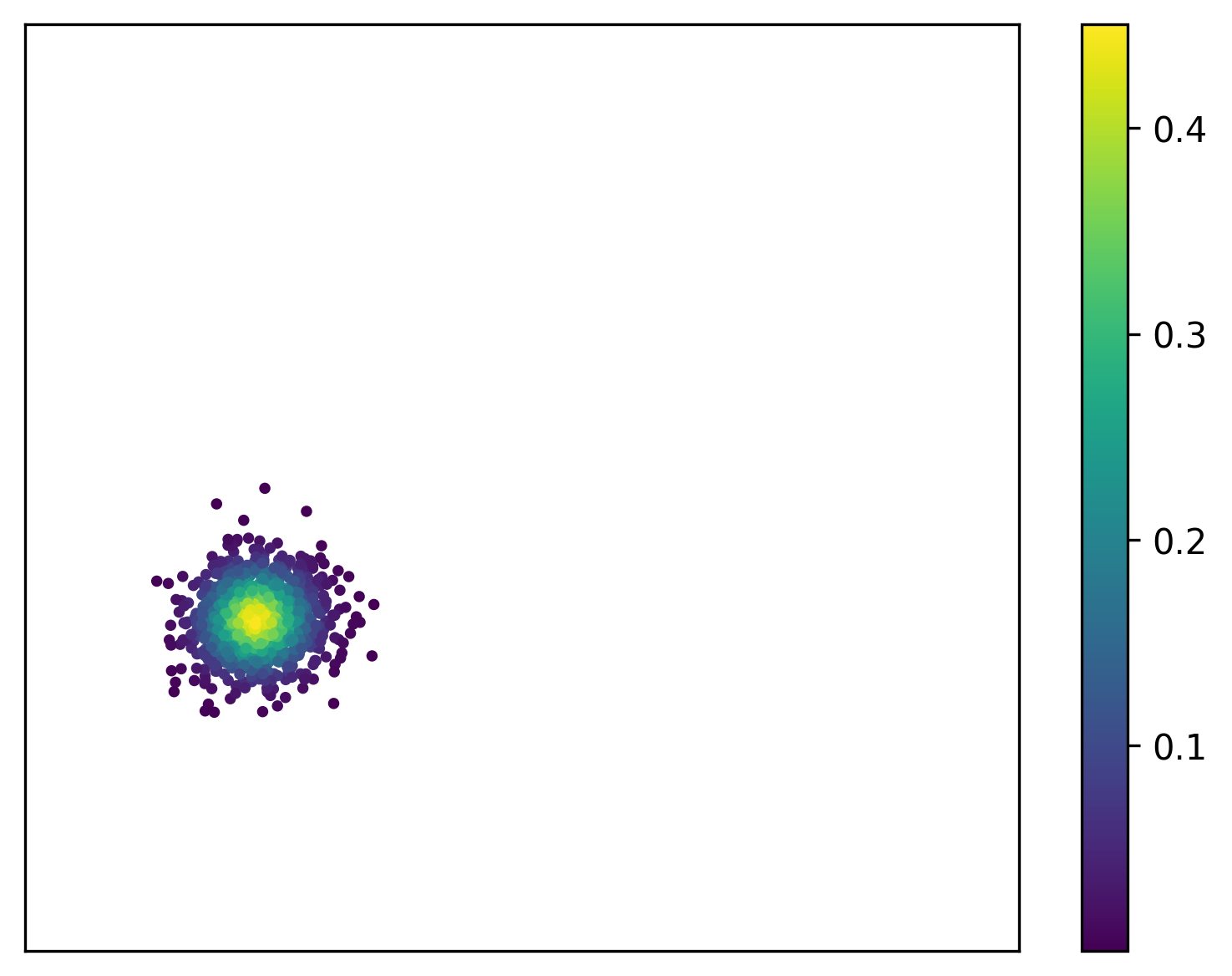}}
\subfigure[]{\includegraphics[width=1.9cm]{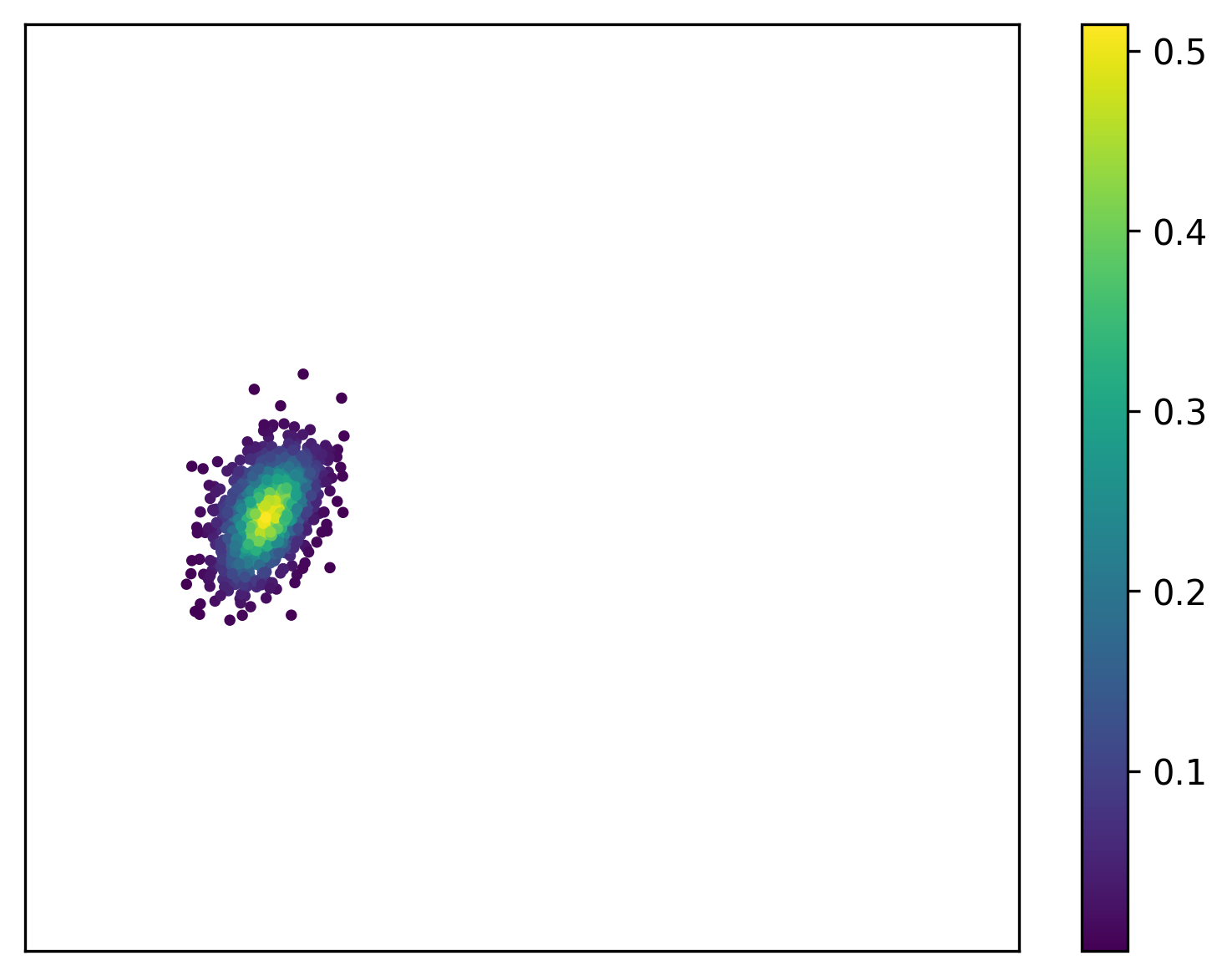}}
\subfigure[]{\includegraphics[width=1.9cm]{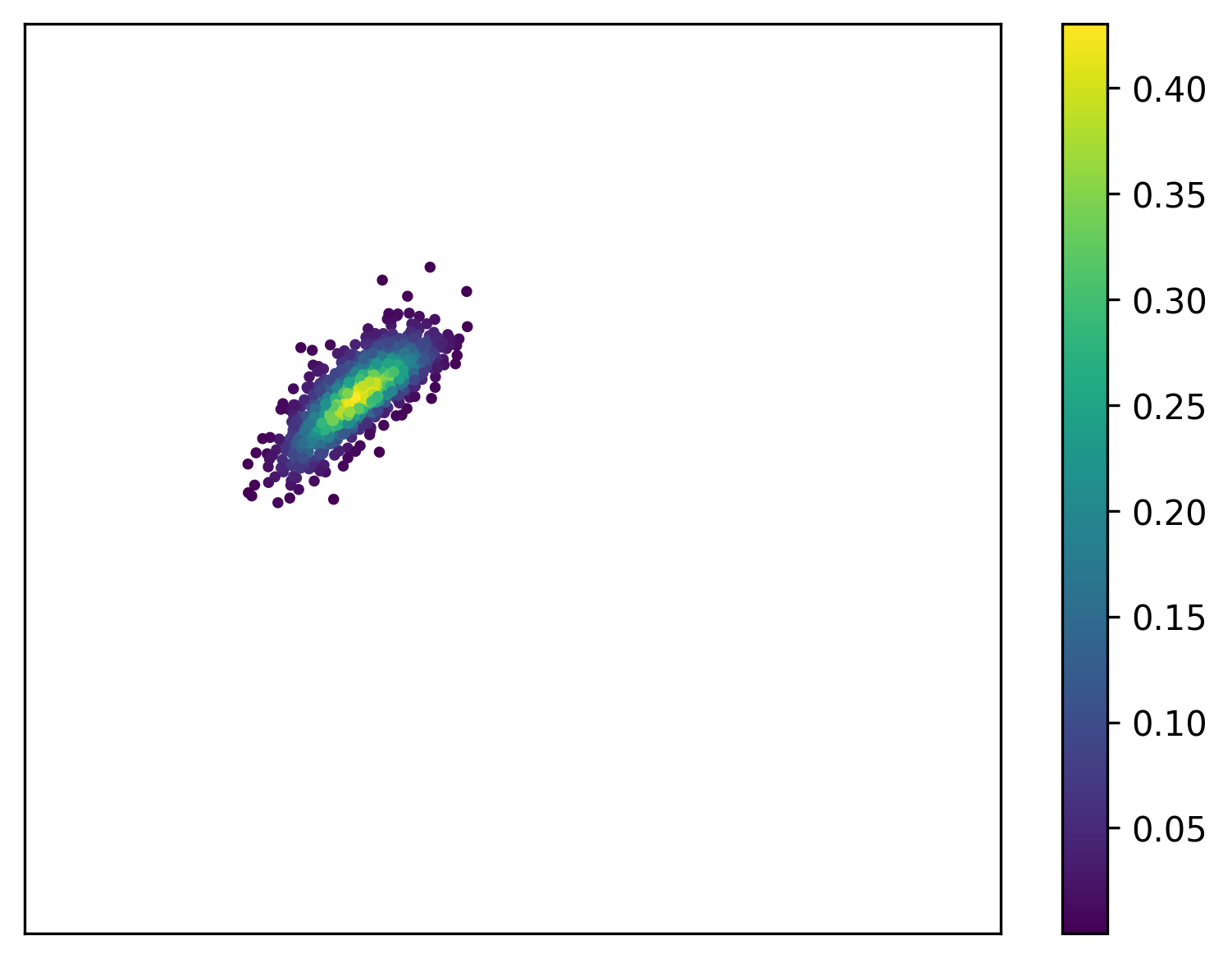}}
\subfigure[]{\includegraphics[width=1.9cm]{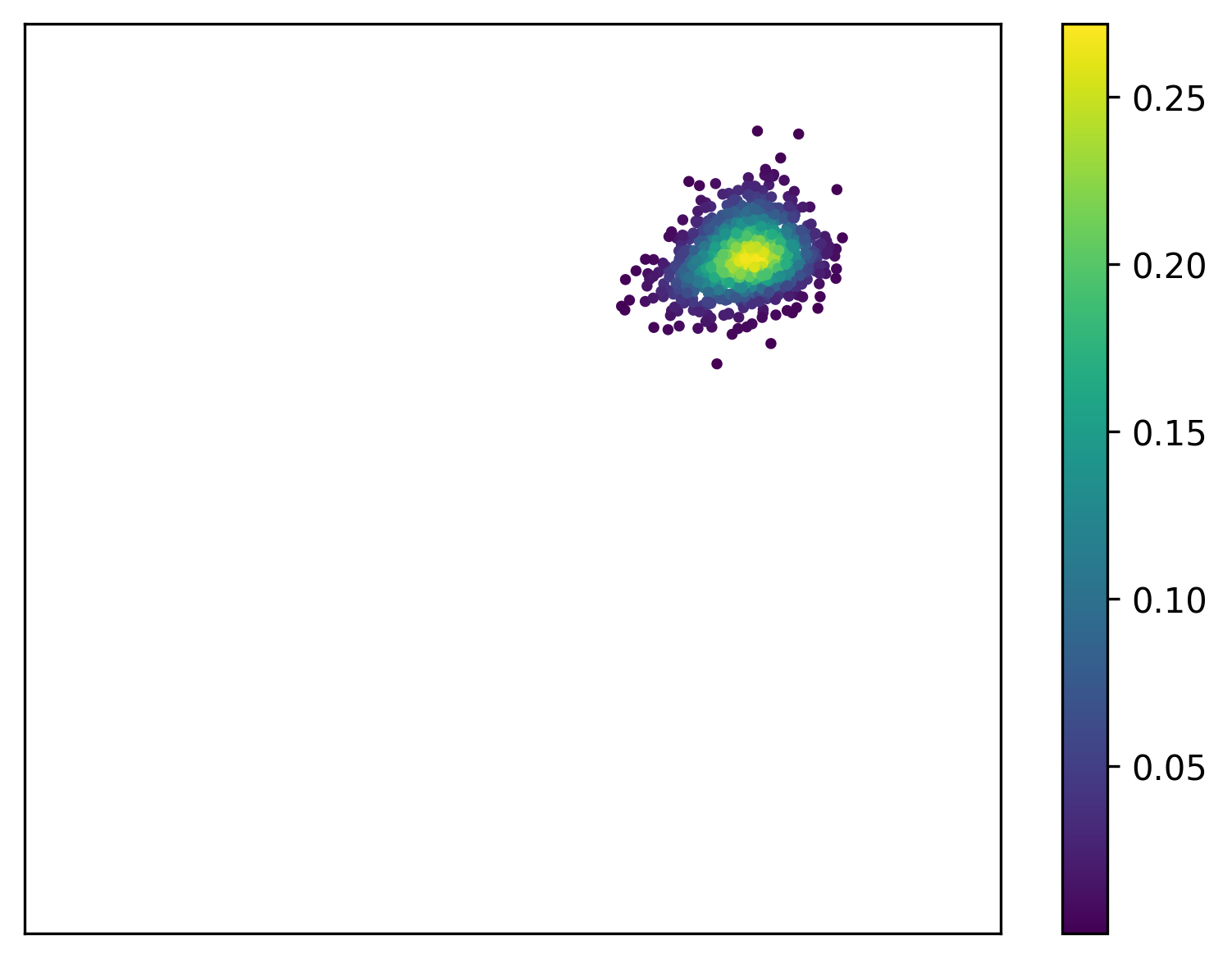}}
\subfigure[]{\includegraphics[width=1.9cm]{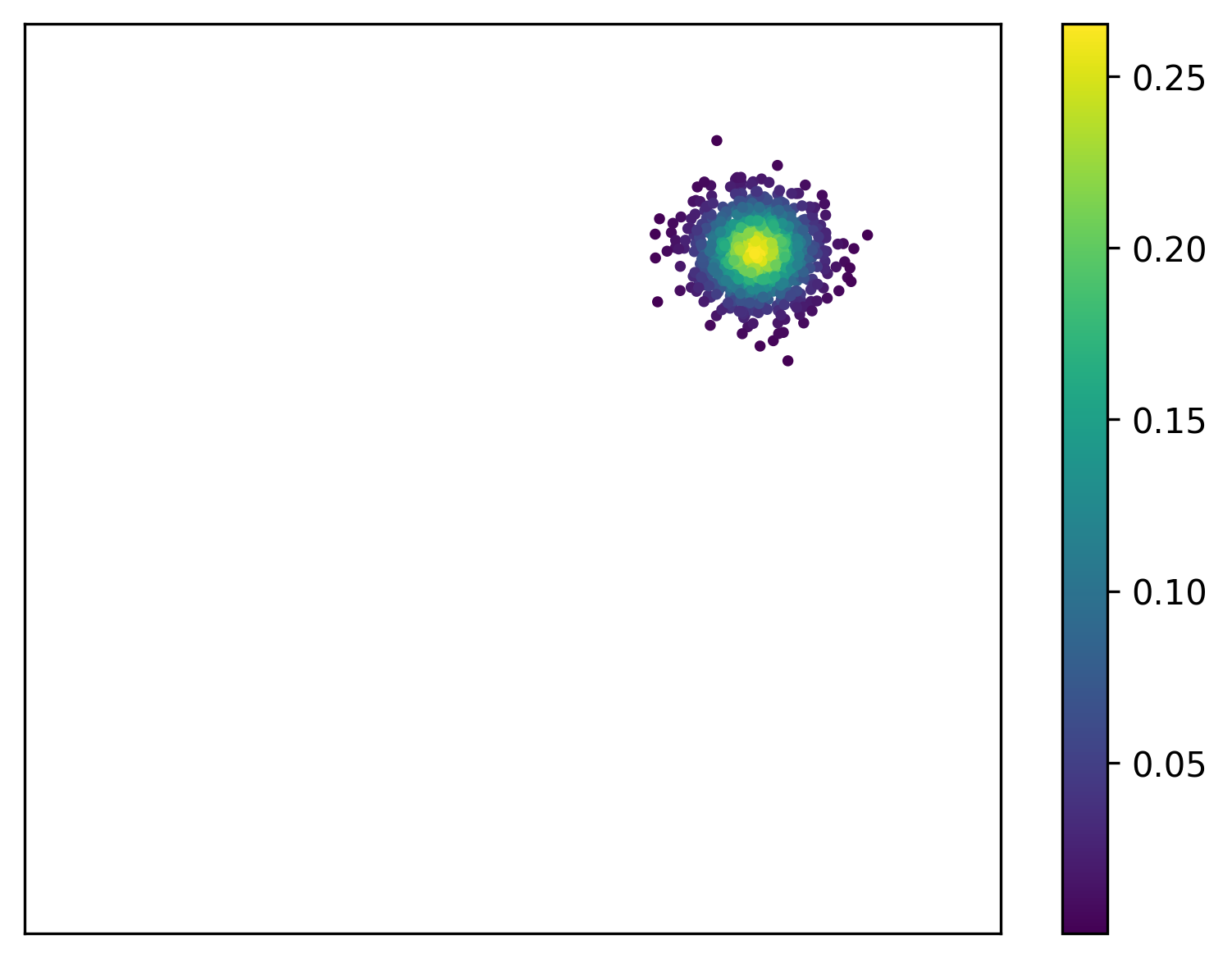}}
\subfigure[]{\includegraphics[width=1.5cm]{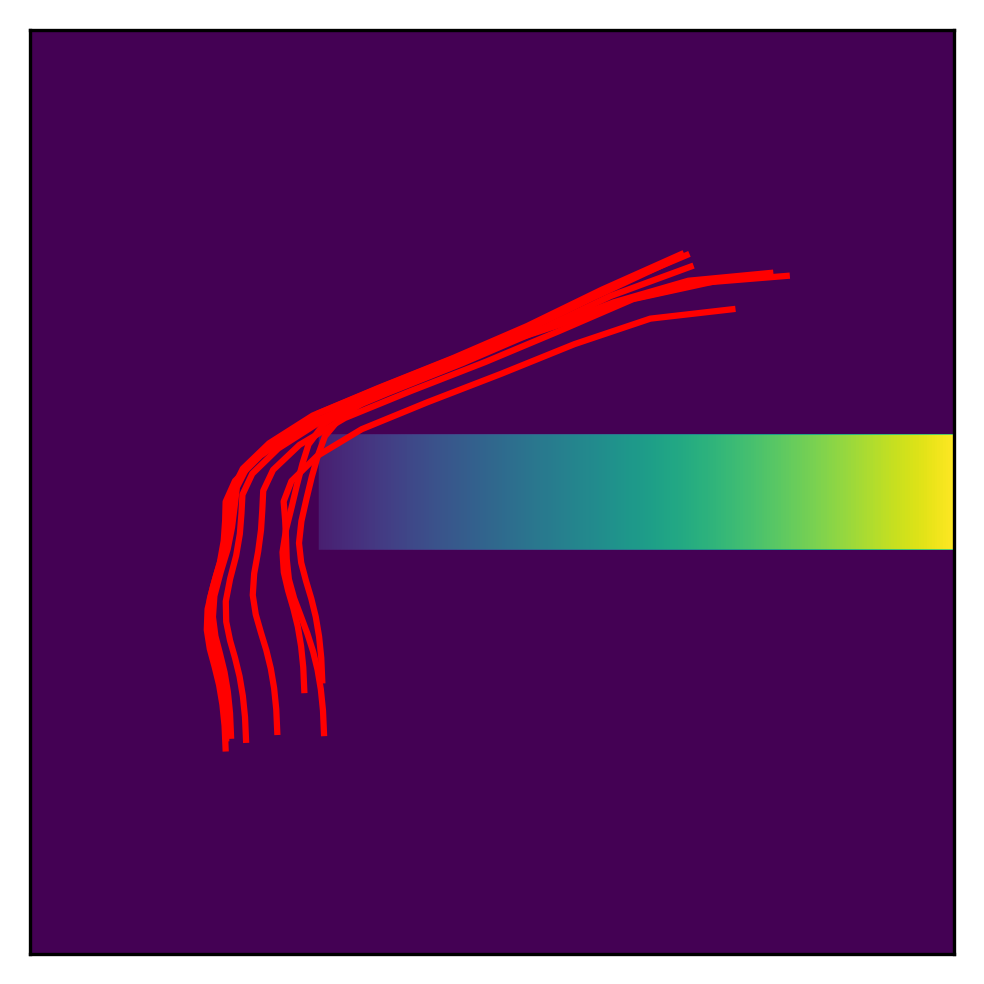}}
\caption{Illustration of crowd motion problems at $d=2$.}
\label{crowdmotion}
\end{center}
\end{figure}

The obstacles in the mazes are represented by an indicator function, with the regions defined as regular rectangles. 
To ensure that obstacles are successfully avoided, we initially apply a method similar to a Gaussian function to blur the rectangular obstacles, as illustrated in the last column (g) of Figure \ref{crowdmotion}.
The density evolution at times $t=0,0.25,0.5,0.75$ and $1$ is shown in (a)-(e). Furthermore, the red lines in (g) represent the characteristics (i.e. the learned trajectories) starting from randomly sampled points based on the initial density $\rho_0$. 
Besides, the push-forward of $\rho_0$ in (e) closely resembles the target density $\rho_1$ in (f), both in terms of position and colorbar values.
These results show that the mass can effectively avoid the obstacles, which demonstrates the effectiveness of our UOT algorithm in handing the crowd motion problem with varying total masses.

\section{Conclusion}\label{Conclusion}
In this paper, we propose a novel neural network-based method for solving high-dimensional dynamic UOT problems. 
To tackle the curse of dimensionality inherent in the traditional mesh-based discretization methods, we employ Lagrangian discretization to solve the continuity equation with the source term along the characteristics, while utilizing the Monte Carlo method to approximate the high-dimensional integrals. 
Additionally, we carefully designed neural networks to parametrize both the velocity field and the source function.
Our method is easy to implement through back-propagation and yields highly accurate results with a constant number of samples in the numerical experiments. 
In the future, we aim to explore the applications of the proposed framework in inverse problems and machine learning.

\section*{Acknowledgments}
This work is supported by the National Natural Science Foundation of China (grants 12301538 and 92370125) and the National Key R\&D Program of China (grant 2021YFA1001300).

\bibliographystyle{siam}
\bibliography{bibtex}

\appendix
\section{Gaussian examples for dimensions $d=10,30,60,80$ on V100}
\label{appA}


\begin{figure}[H]
\begin{center}
\includegraphics[width=2.4cm]{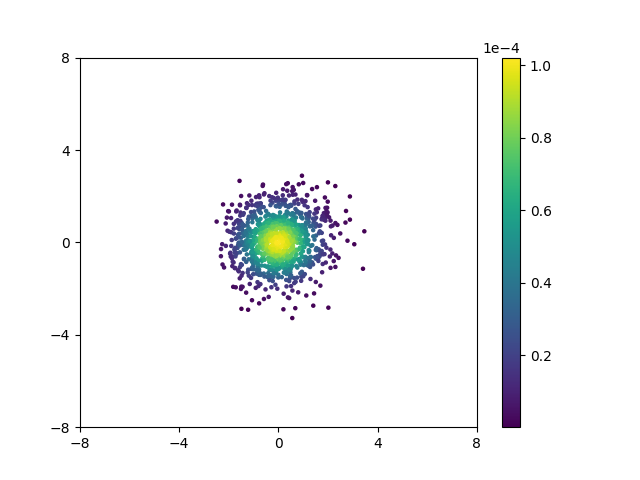}
\hspace{-5mm}
\includegraphics[width=2.4cm]{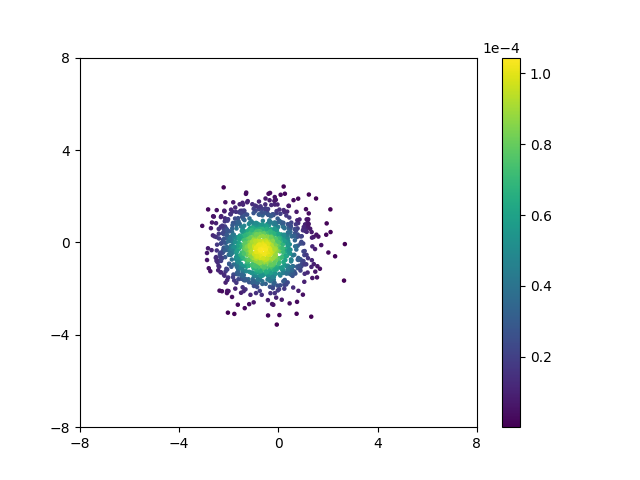}
\hspace{-5mm}
\includegraphics[width=2.4cm]{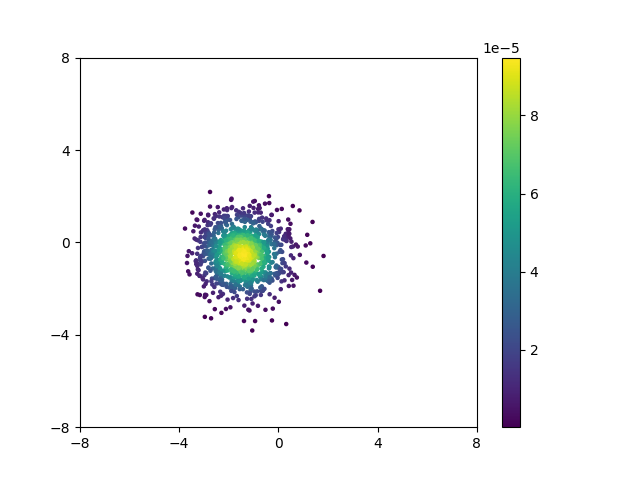}
\hspace{-5mm}
\includegraphics[width=2.4cm]{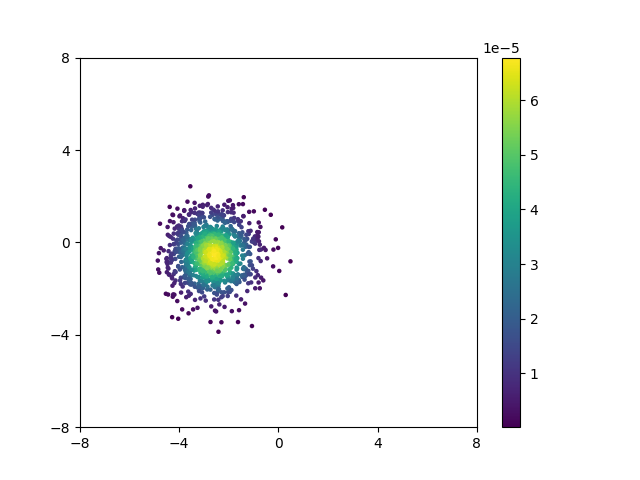}
\hspace{-5mm}
\includegraphics[width=2.4cm]{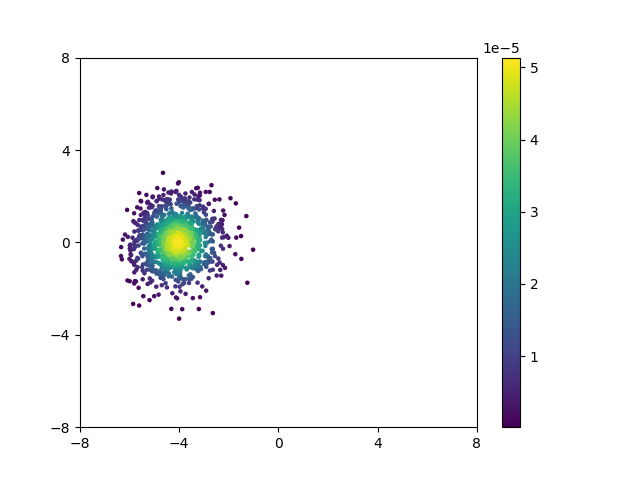}
\hspace{-5mm}
\includegraphics[width=2.4cm]{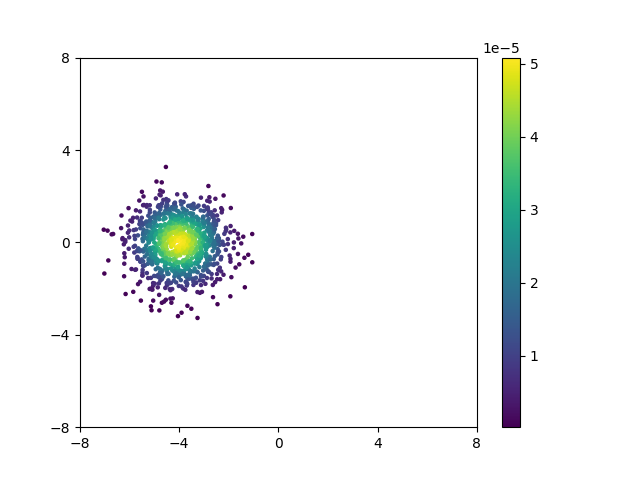}\\
\vspace{5pt}

\includegraphics[width=2.4cm]{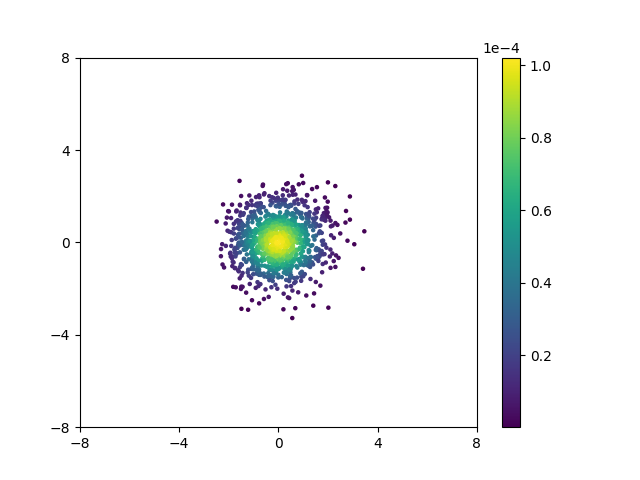}
\hspace{-5mm}
\includegraphics[width=2.4cm]{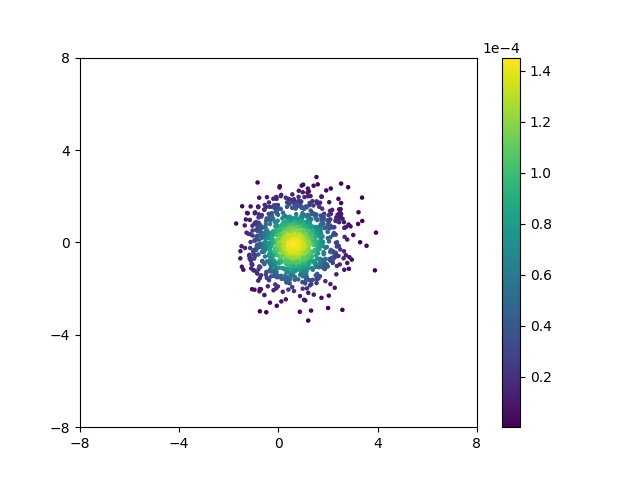}
\hspace{-5mm}
\includegraphics[width=2.4cm]{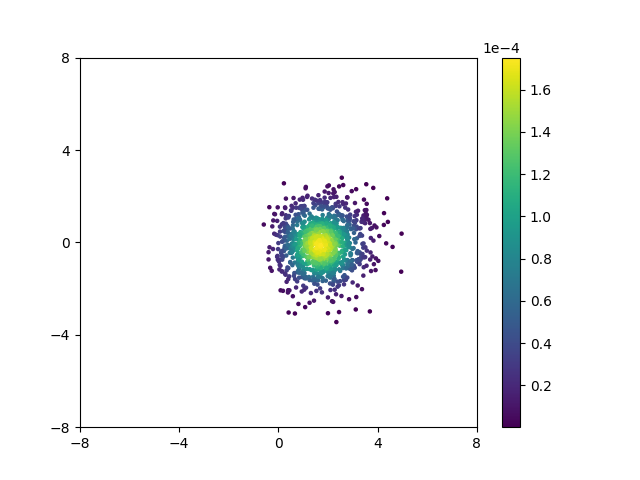}
\hspace{-5mm}
\includegraphics[width=2.4cm]{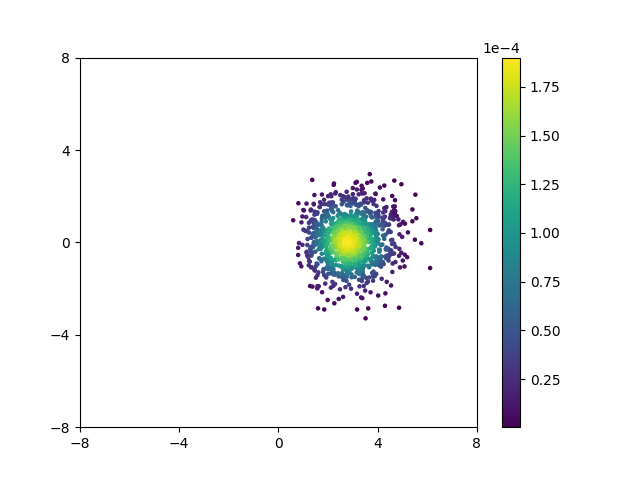}
\hspace{-5mm}
\includegraphics[width=2.4cm]{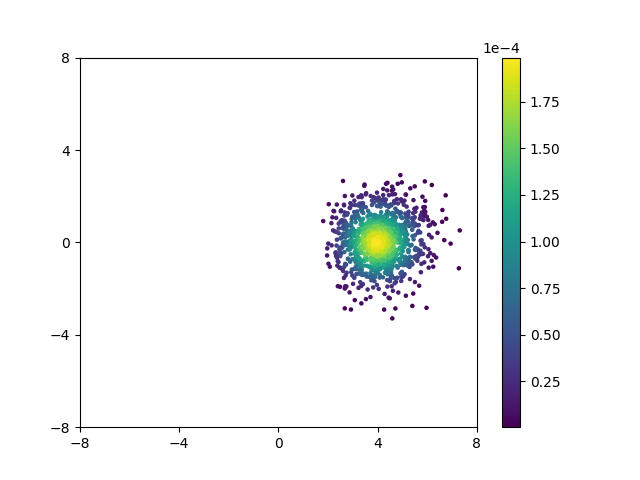}
\hspace{-5mm}
\includegraphics[width=2.4cm]{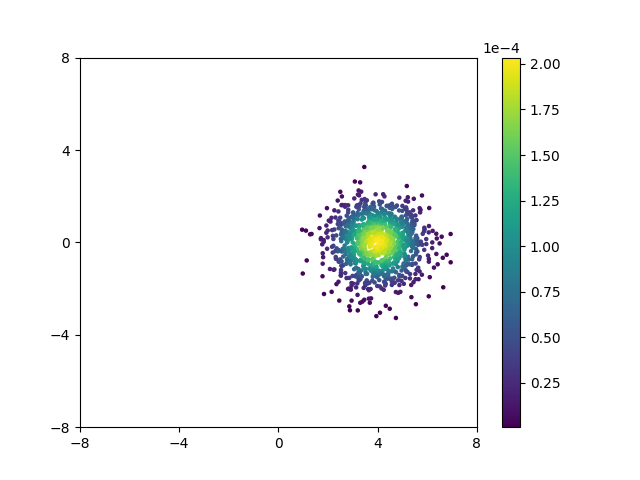}\\
\vspace{5pt}

\includegraphics[width=2.4cm]{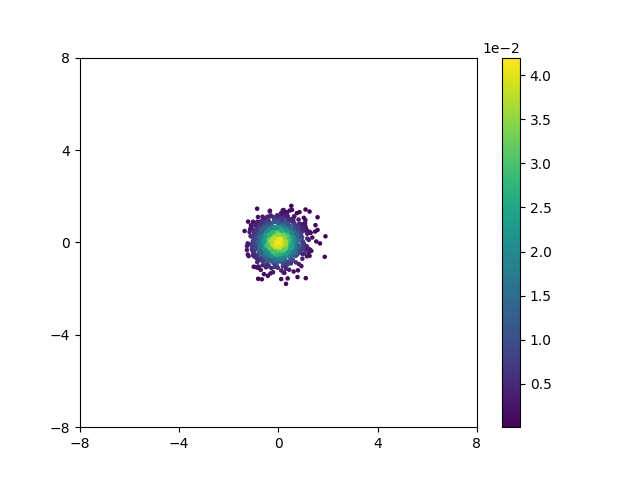}
\hspace{-5mm}
\includegraphics[width=2.4cm]{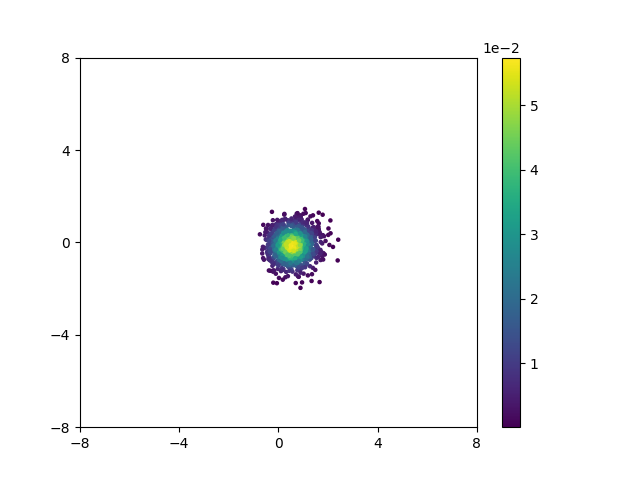}
\hspace{-5mm}
\includegraphics[width=2.4cm]{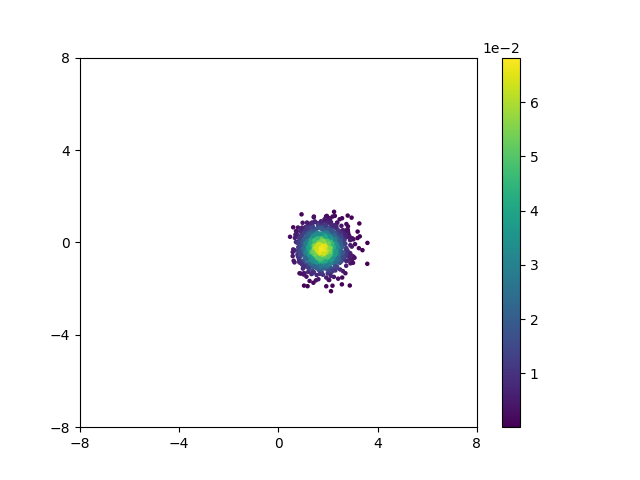}
\hspace{-5mm}
\includegraphics[width=2.4cm]{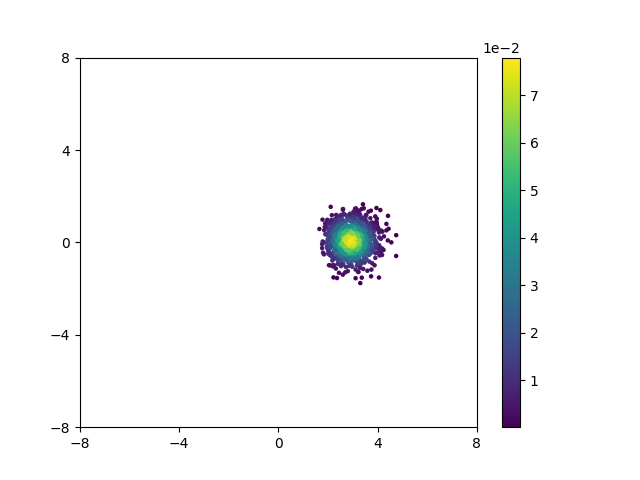}
\hspace{-5mm}
\includegraphics[width=2.4cm]{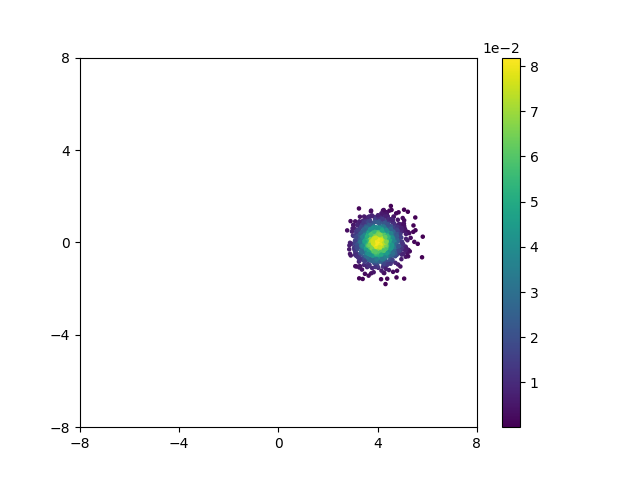}
\hspace{-5mm}
\includegraphics[width=2.4cm]{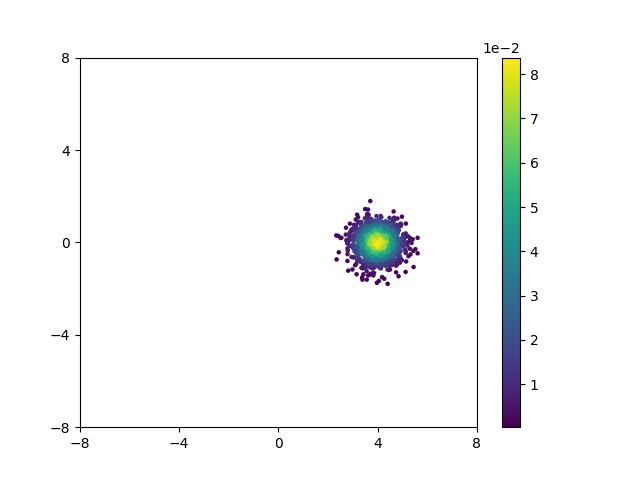}\\
\vspace{5pt}

\subfigure[$\rho_0(\boldsymbol{\boldsymbol{x}}_0)$]{\includegraphics[width=2.4cm]{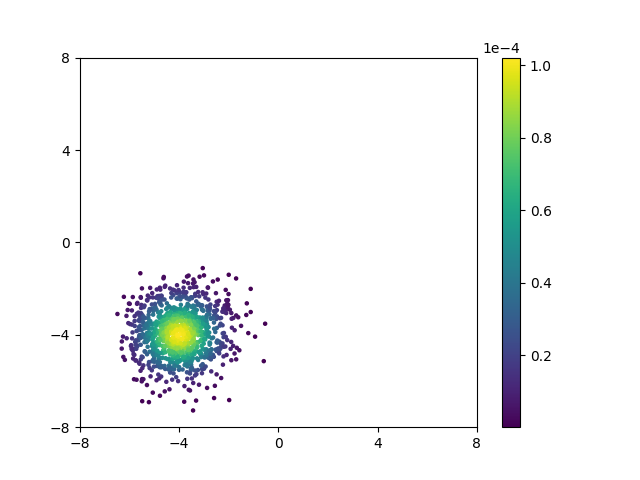}}
\hspace{-5mm}
\subfigure[$\tilde{\rho}_{1/4}(\tilde{\boldsymbol{\boldsymbol{x}}}_{1/4})$]{\includegraphics[width=2.4cm]{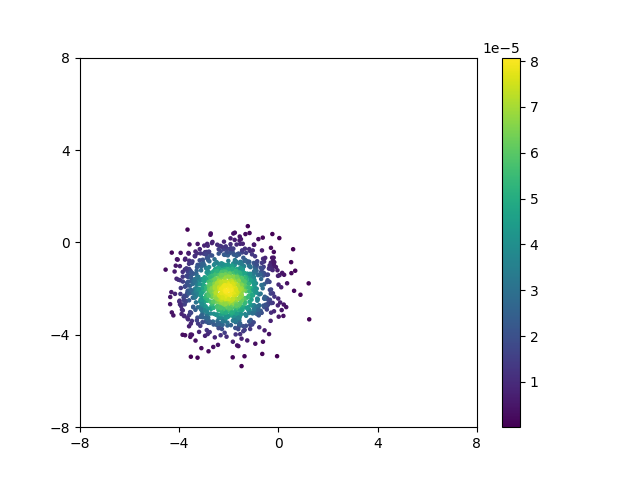}}
\hspace{-5mm}
\subfigure[$\tilde{\rho}_{1/2}(\tilde{\boldsymbol{\boldsymbol{x}}}_{1/2})$]{\includegraphics[width=2.4cm]{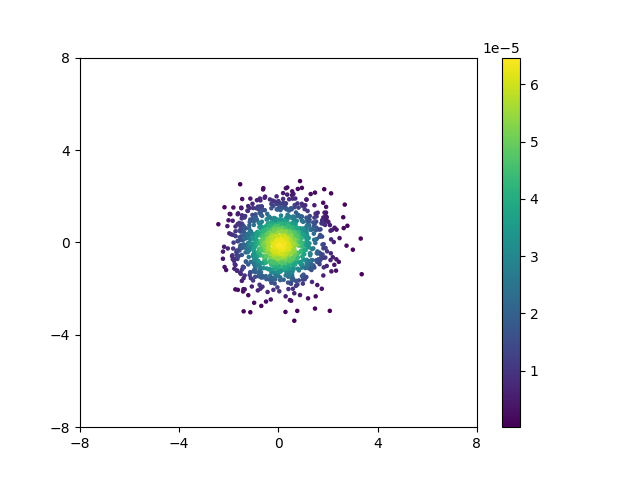}}
\hspace{-5mm}
\subfigure[$\tilde{\rho}_{3/4}(\tilde{\boldsymbol{\boldsymbol{x}}}_{3/4})$]{\includegraphics[width=2.4cm]{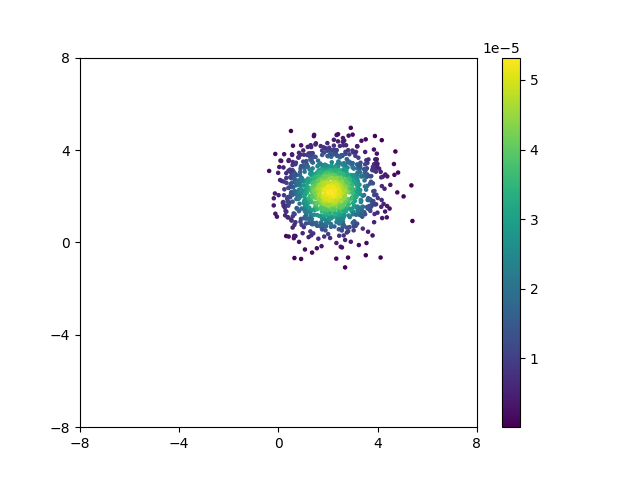}}
\hspace{-5mm}
\subfigure[$\tilde{\rho}_{1}(\tilde{\boldsymbol{\boldsymbol{x}}}_{1})$]{\includegraphics[width=2.4cm]{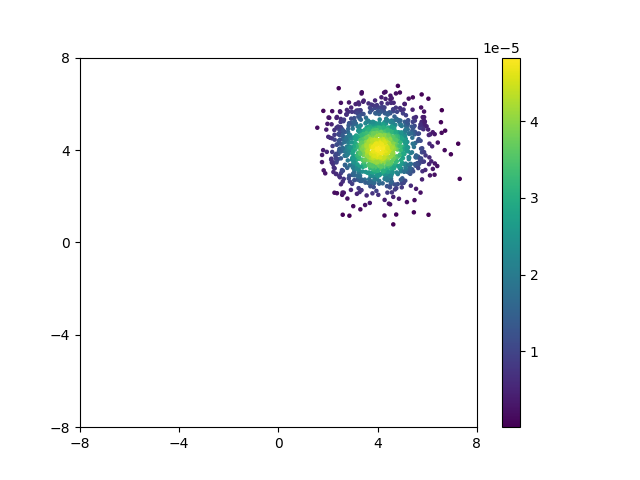}}
\hspace{-5mm}
\subfigure[$\rho_{1}(\boldsymbol{\boldsymbol{x}}_{1})$]{\includegraphics[width=2.4cm]{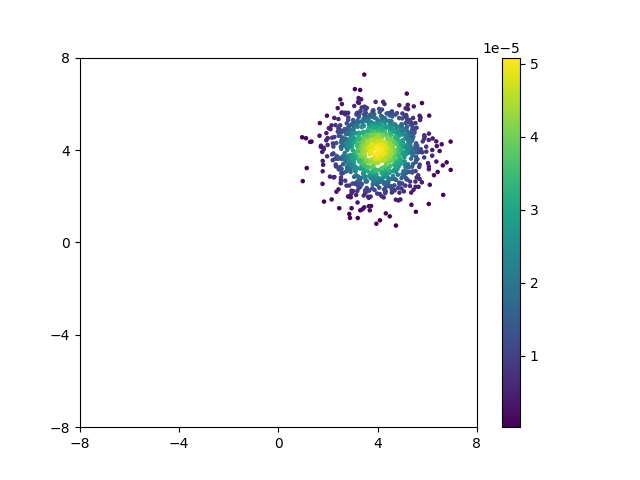}}

\caption{Illustration of these Gaussian problems at $d=10$.}
\label{compare_d10}
\end{center}
\end{figure}


\begin{figure}[H]
\begin{center}
\includegraphics[width=2.4cm]{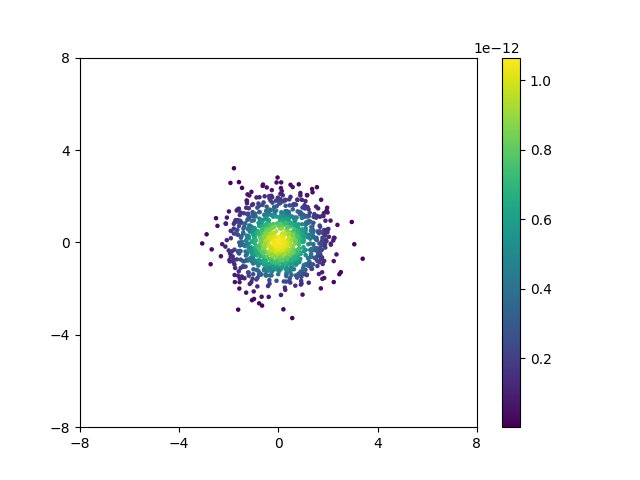}
\hspace{-5mm}
\includegraphics[width=2.4cm]{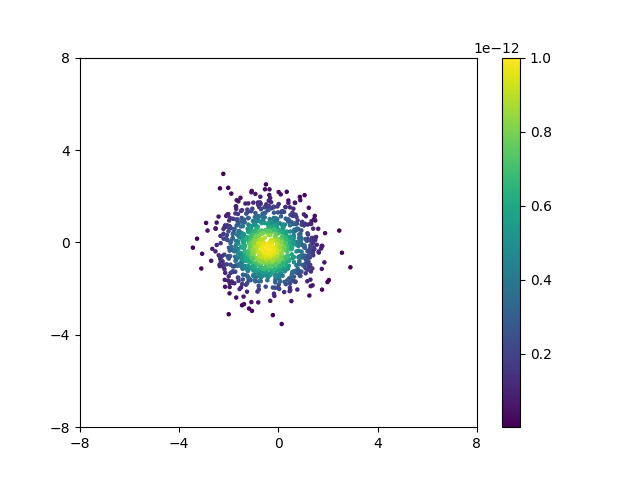}
\hspace{-5mm}
\includegraphics[width=2.4cm]{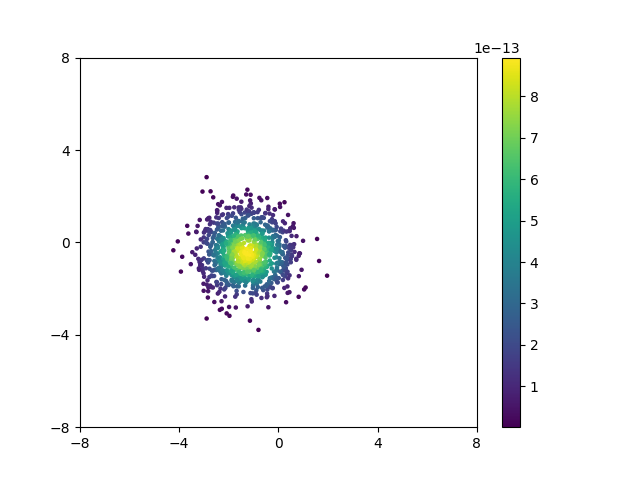}
\hspace{-5mm}
\includegraphics[width=2.4cm]{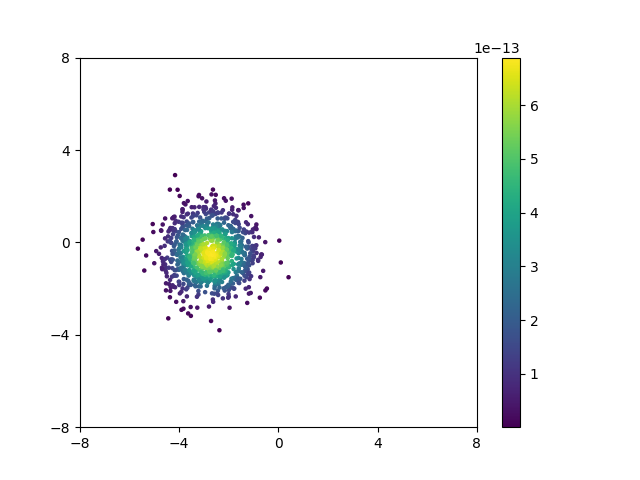}
\hspace{-5mm}
\includegraphics[width=2.4cm]{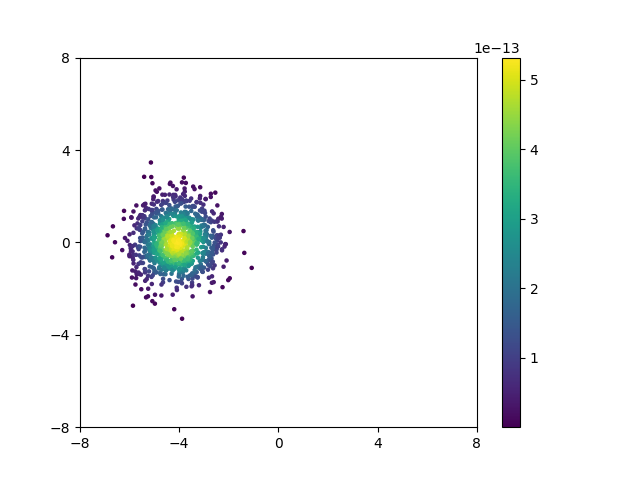}
\hspace{-5mm}
\includegraphics[width=2.4cm]{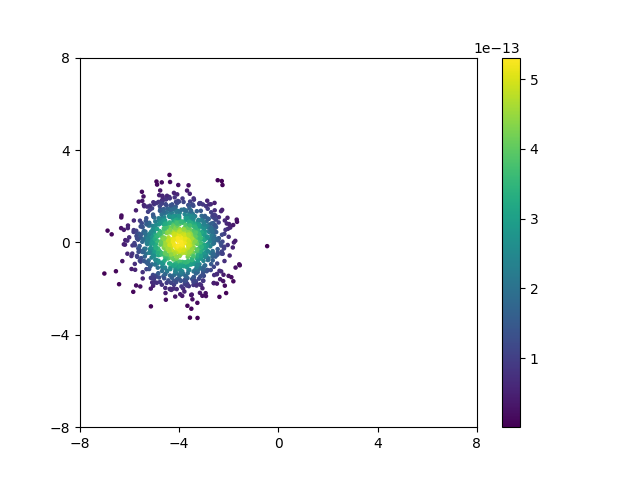}\\
\vspace{5pt}

\includegraphics[width=2.4cm]{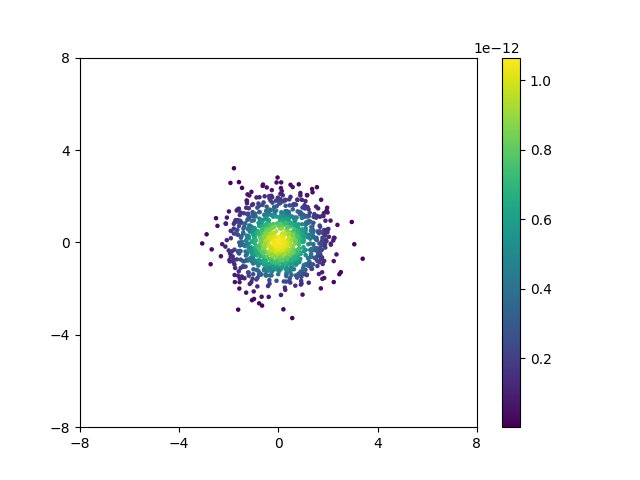}
\hspace{-5mm}
\includegraphics[width=2.4cm]{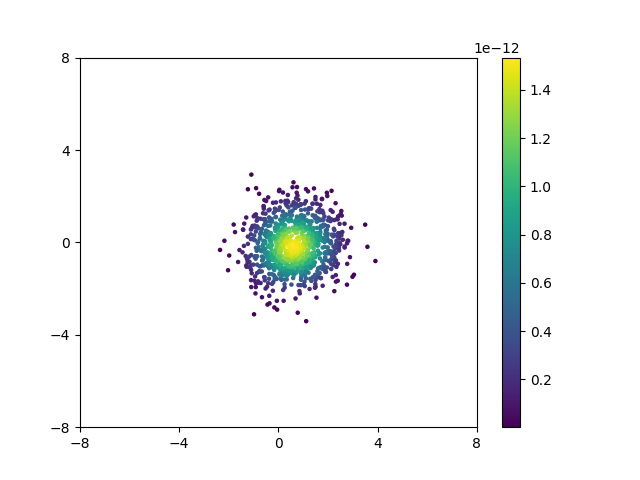}
\hspace{-5mm}
\includegraphics[width=2.4cm]{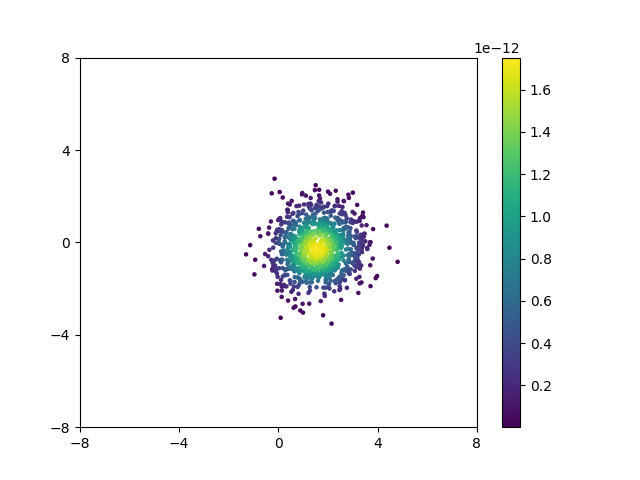}
\hspace{-5mm}
\includegraphics[width=2.4cm]{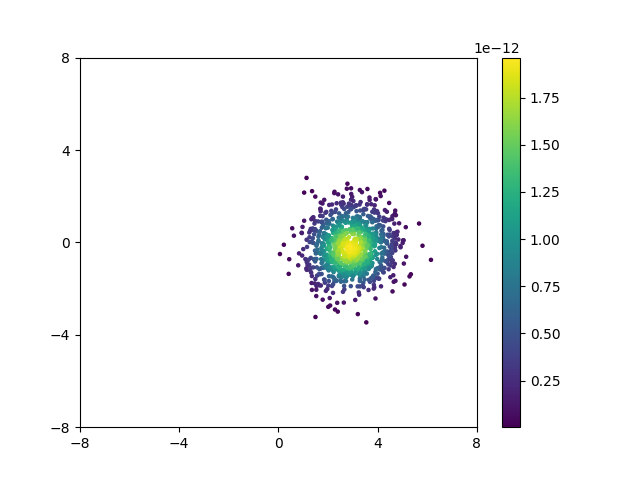}
\hspace{-5mm}
\includegraphics[width=2.4cm]{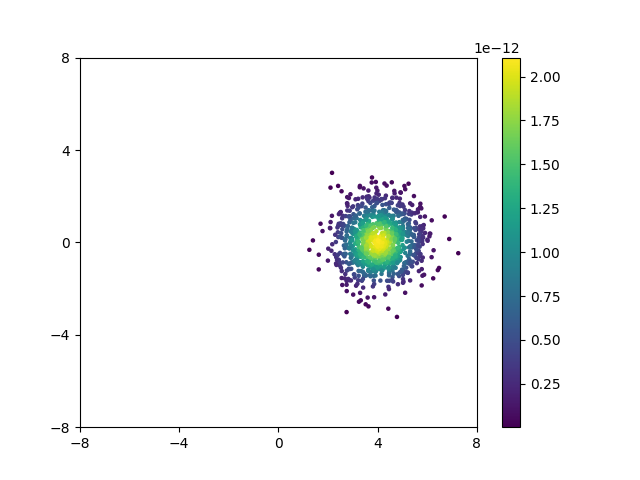}
\hspace{-5mm}
\includegraphics[width=2.4cm]{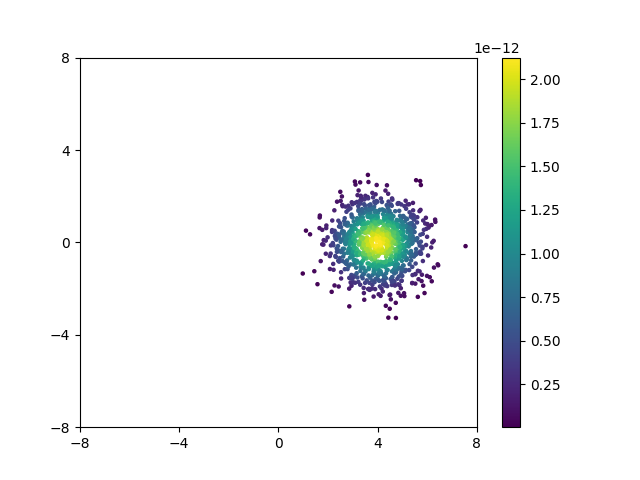}\\
\vspace{5pt}

\includegraphics[width=2.4cm]{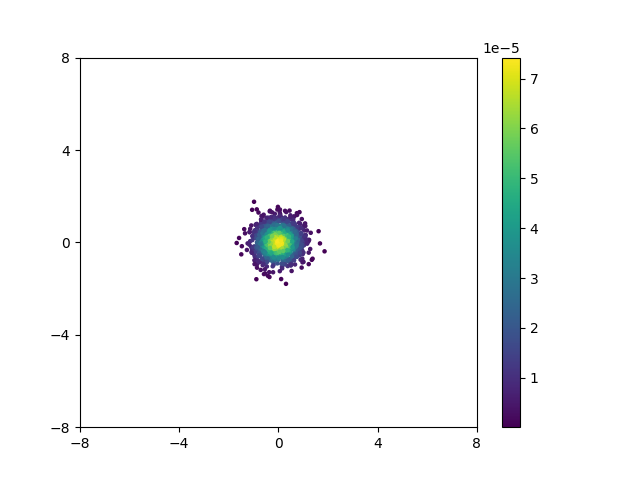}
\hspace{-5mm}
\includegraphics[width=2.4cm]{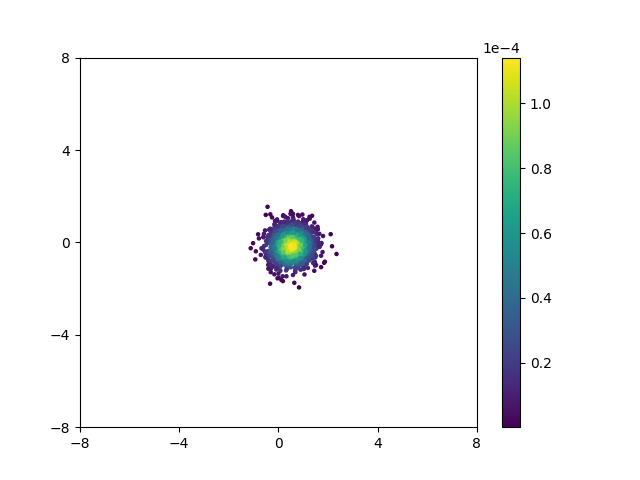}
\hspace{-5mm}
\includegraphics[width=2.4cm]{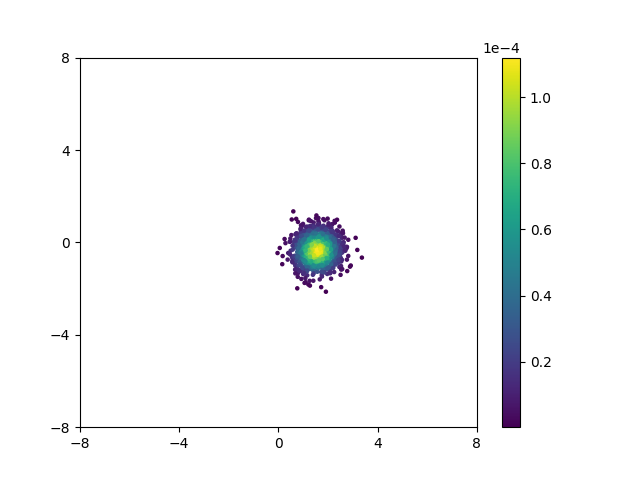}
\hspace{-5mm}
\includegraphics[width=2.4cm]{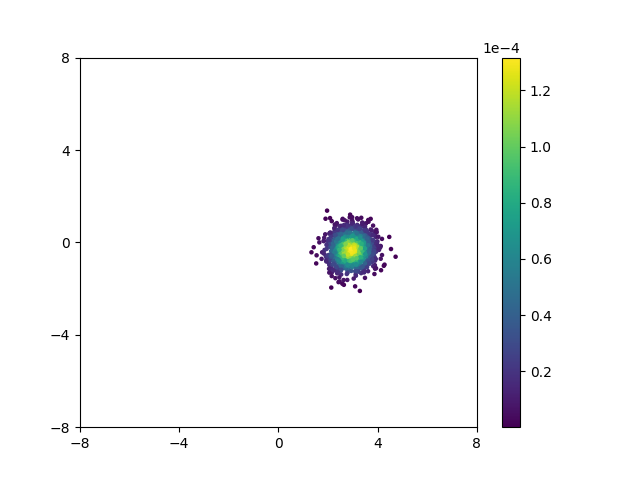}
\hspace{-5mm}
\includegraphics[width=2.4cm]{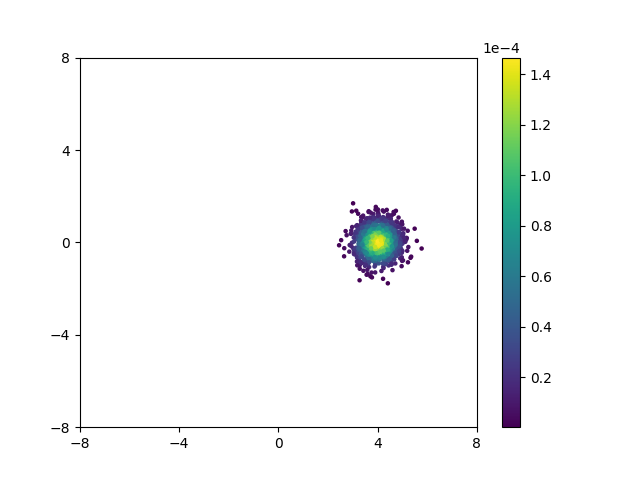}
\hspace{-5mm}
\includegraphics[width=2.4cm]{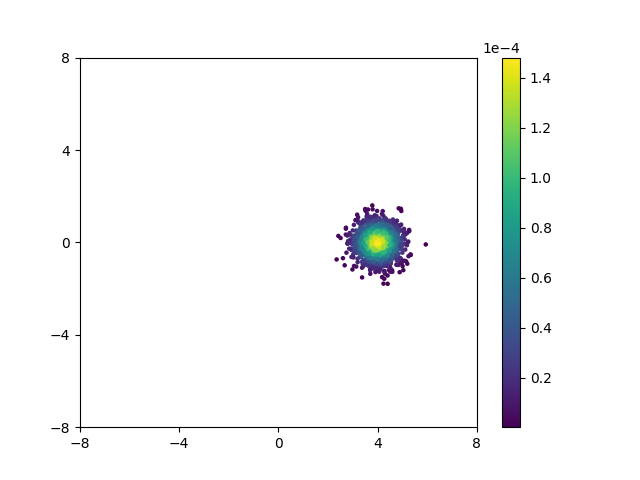}\\
\vspace{5pt}

\subfigure[$\rho_0(\boldsymbol{\boldsymbol{x}}_0)$]{\includegraphics[width=2.4cm]{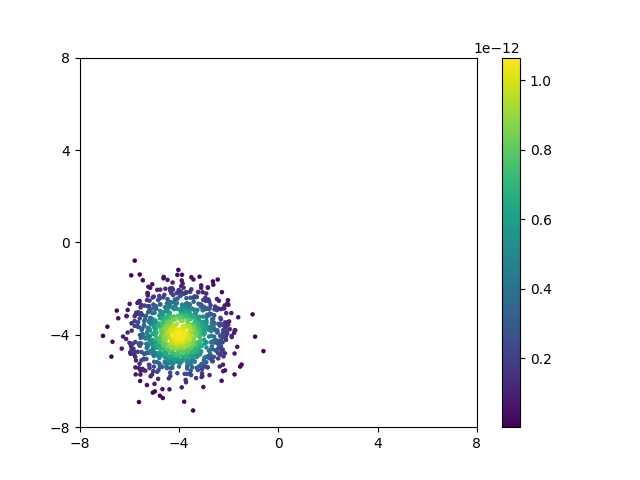}}
\hspace{-5mm}
\subfigure[$\tilde{\rho}_{1/4}(\tilde{\boldsymbol{\boldsymbol{x}}}_{1/4})$]{\includegraphics[width=2.4cm]{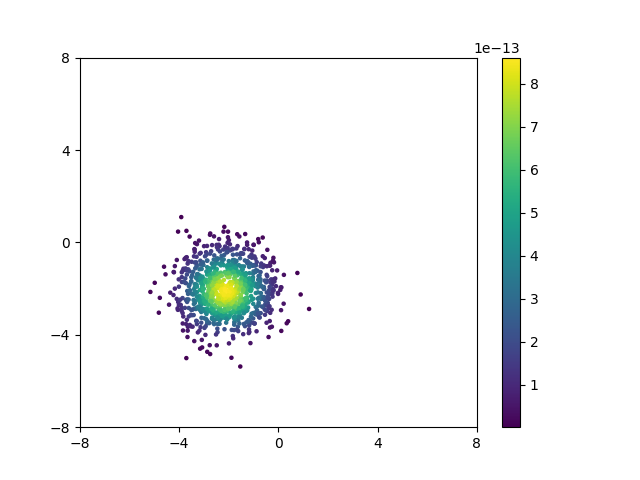}}
\hspace{-5mm}
\subfigure[$\tilde{\rho}_{1/2}(\tilde{\boldsymbol{\boldsymbol{x}}}_{1/2})$]{\includegraphics[width=2.4cm]{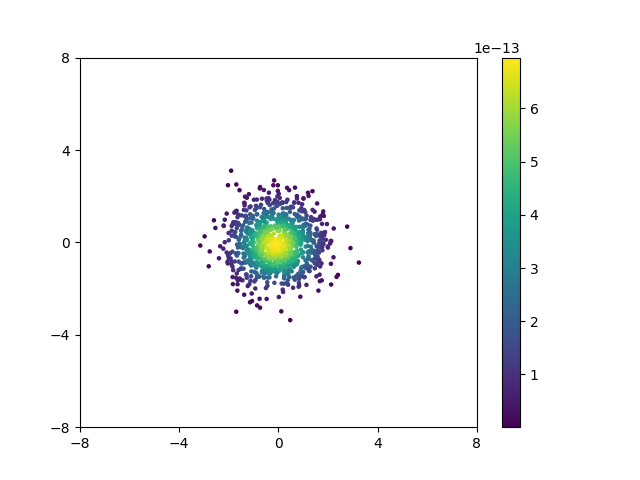}}
\hspace{-5mm}
\subfigure[$\tilde{\rho}_{3/4}(\tilde{\boldsymbol{\boldsymbol{x}}}_{3/4})$]{\includegraphics[width=2.4cm]{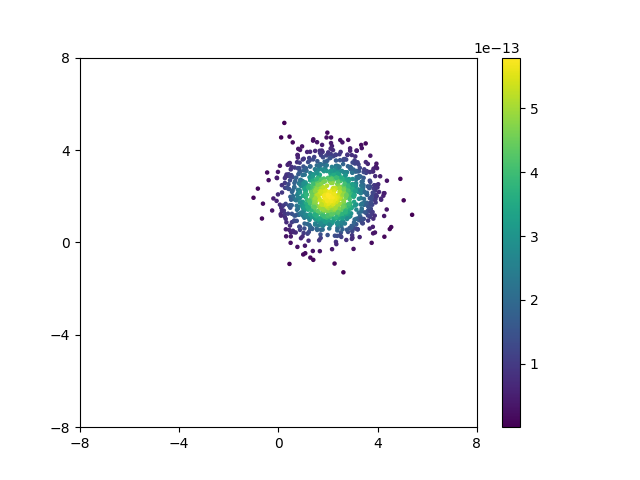}}
\hspace{-5mm}
\subfigure[$\tilde{\rho}_{1}(\tilde{\boldsymbol{\boldsymbol{x}}}_{1})$]{\includegraphics[width=2.4cm]{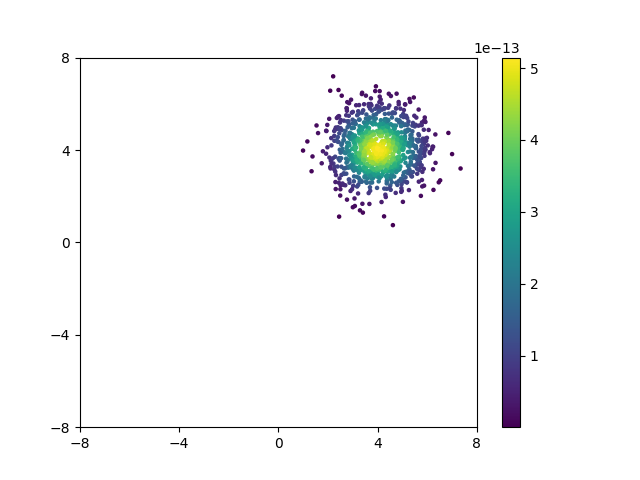}}
\hspace{-5mm}
\subfigure[$\rho_{1}(\boldsymbol{\boldsymbol{x}}_{1})$]{\includegraphics[width=2.4cm]{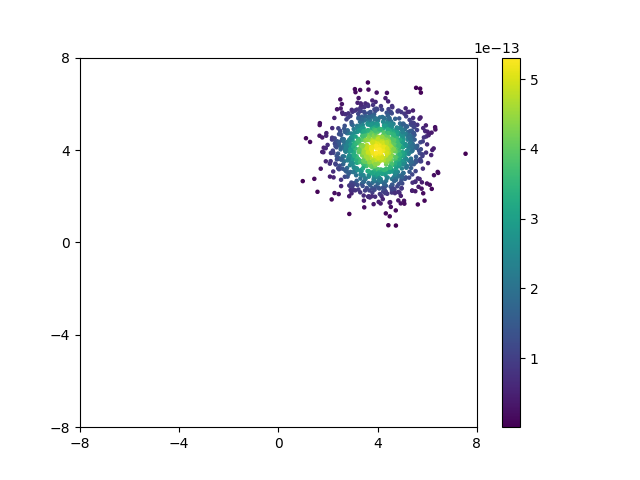}}

\caption{Illustration of these Gaussian problems at $d=30$.}
\label{compare_d30}
\end{center}
\end{figure}

\begin{figure}[H]
\begin{center}
\includegraphics[width=2.4cm]{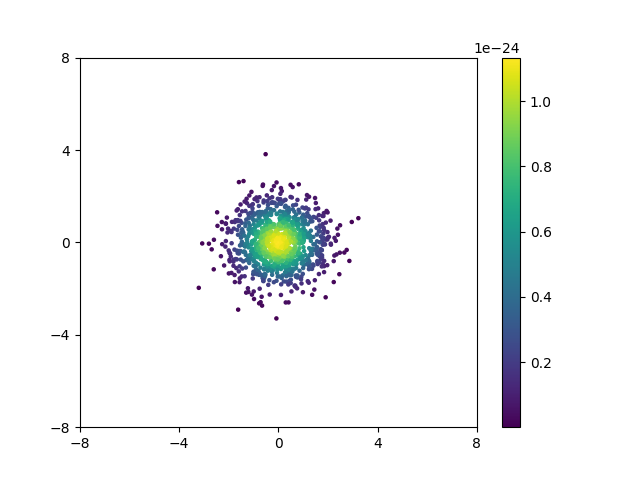}
\hspace{-5mm}
\includegraphics[width=2.4cm]{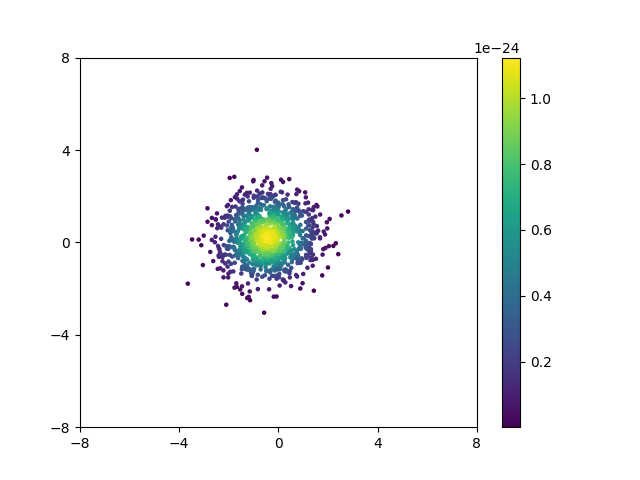}
\hspace{-5mm}
\includegraphics[width=2.4cm]{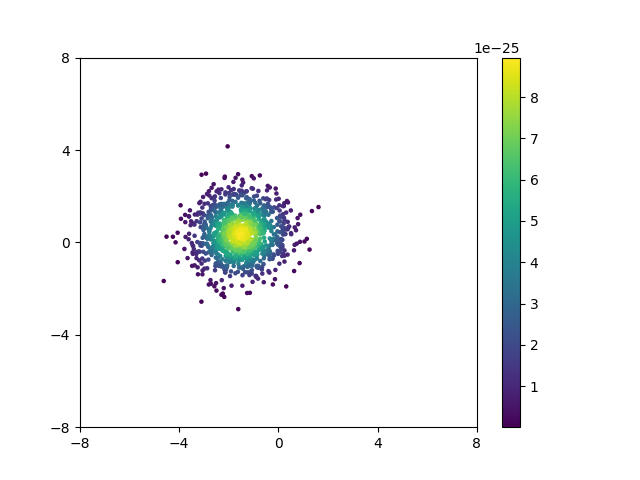}
\hspace{-5mm}
\includegraphics[width=2.4cm]{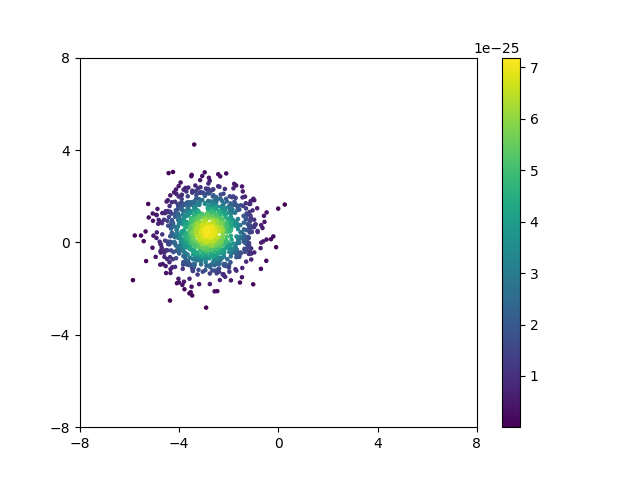}
\hspace{-5mm}
\includegraphics[width=2.4cm]{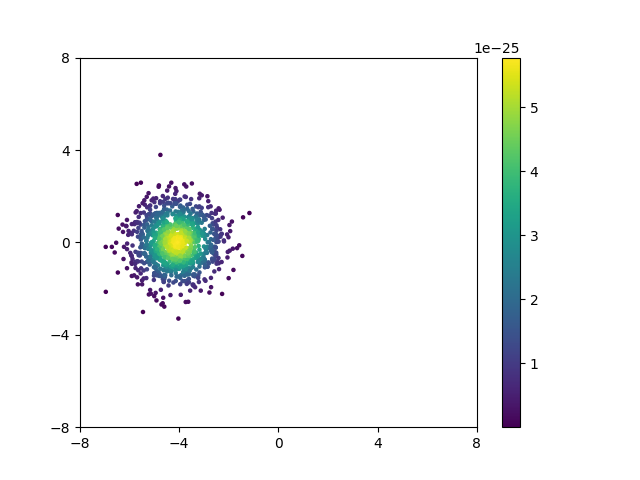}
\hspace{-5mm}
\includegraphics[width=2.4cm]{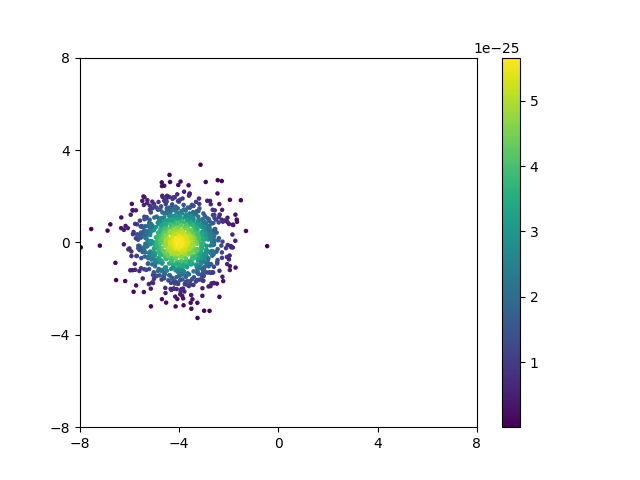}\\
\vspace{5pt}

\includegraphics[width=2.4cm]{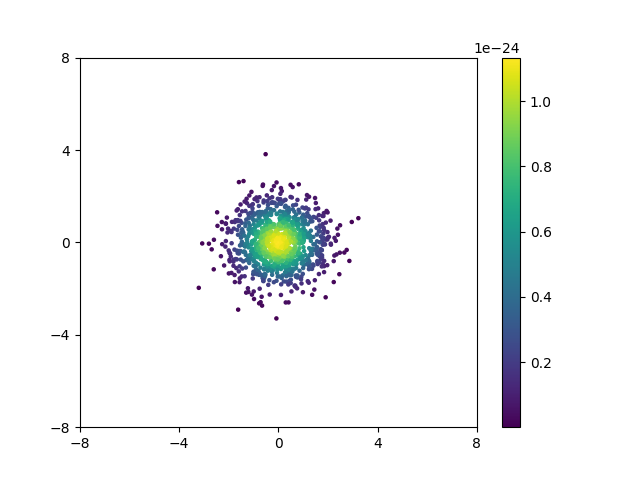}
\hspace{-5mm}
\includegraphics[width=2.4cm]{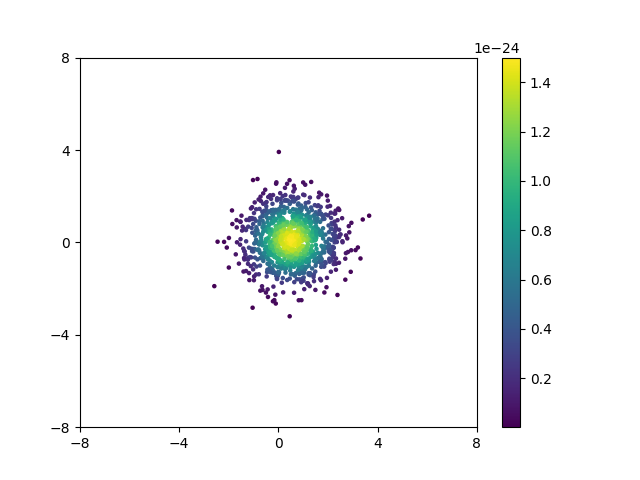}
\hspace{-5mm}
\includegraphics[width=2.4cm]{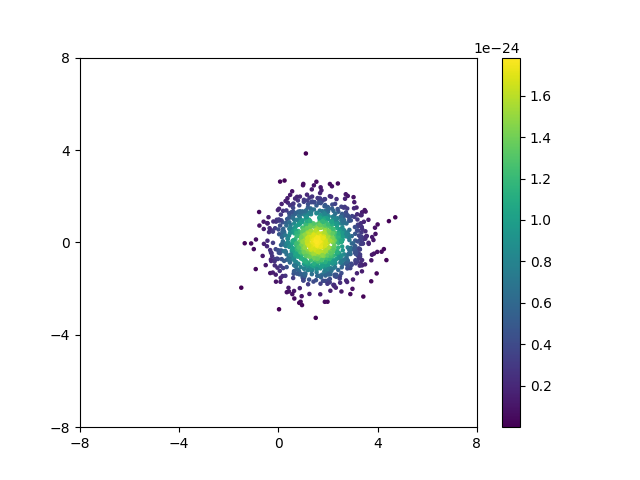}
\hspace{-5mm}
\includegraphics[width=2.4cm]{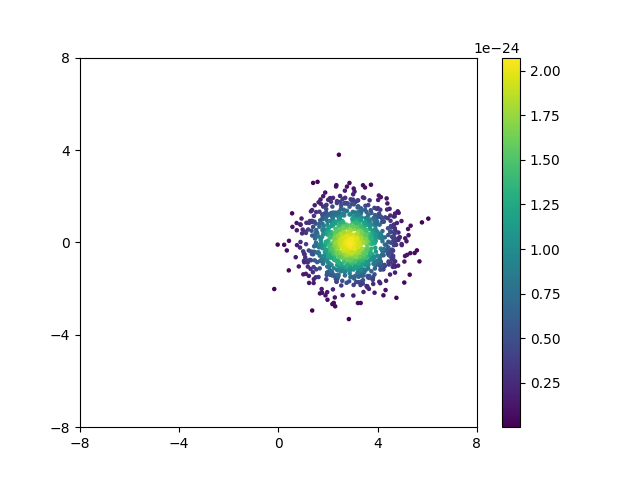}
\hspace{-5mm}
\includegraphics[width=2.4cm]{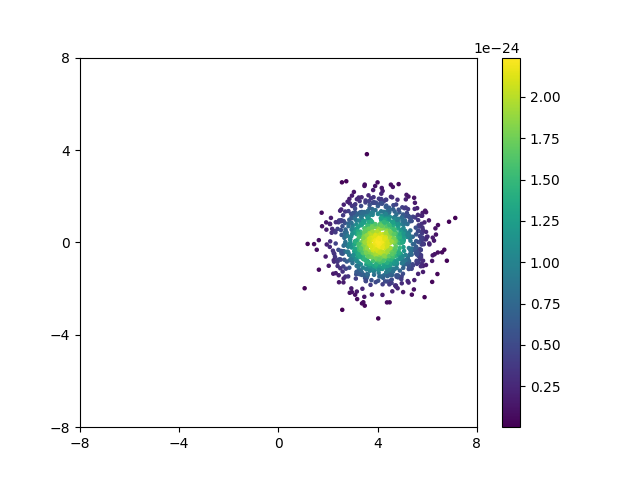}
\hspace{-5mm}
\includegraphics[width=2.4cm]{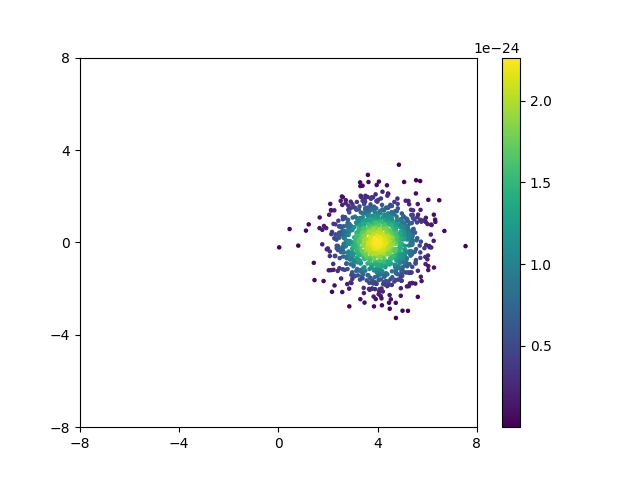}\\
\vspace{5pt}

\includegraphics[width=2.4cm]{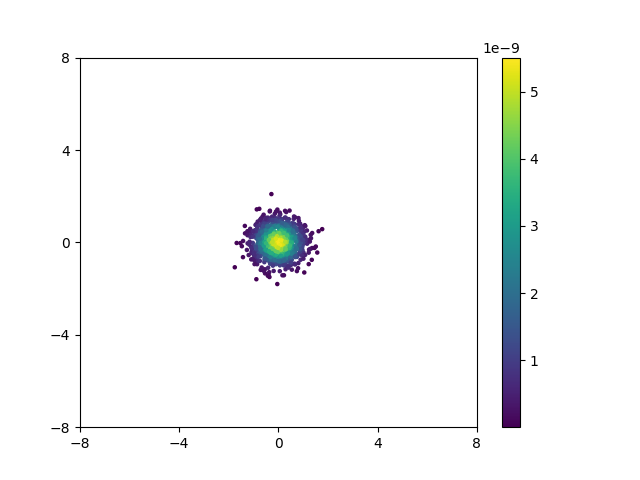}
\hspace{-5mm}
\includegraphics[width=2.4cm]{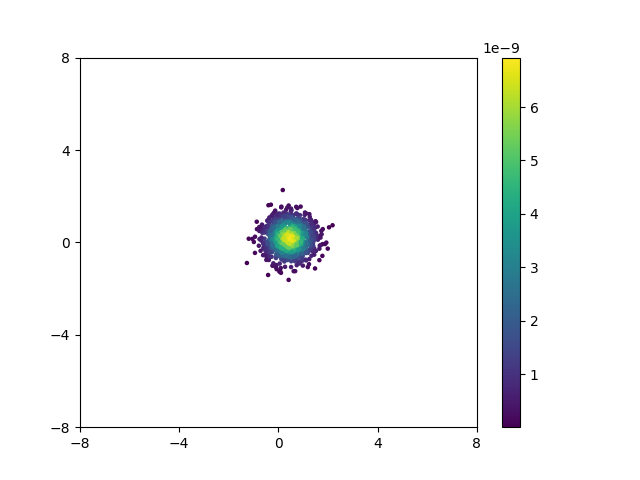}
\hspace{-5mm}
\includegraphics[width=2.4cm]{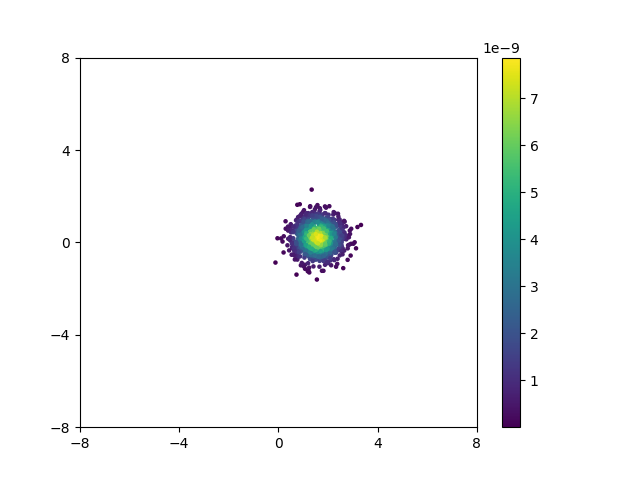}
\hspace{-5mm}
\includegraphics[width=2.4cm]{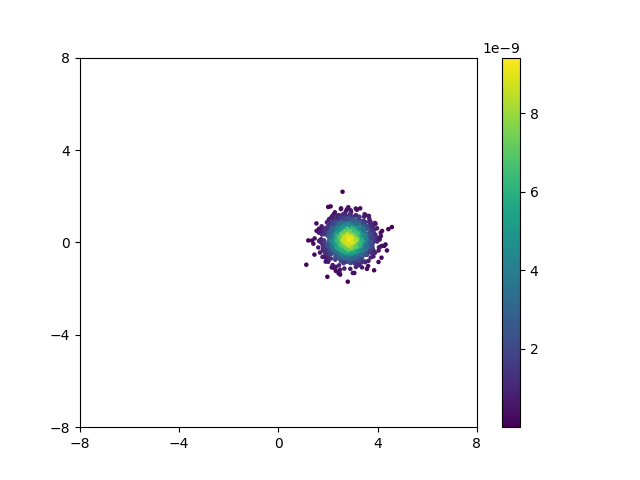}
\hspace{-5mm}
\includegraphics[width=2.4cm]{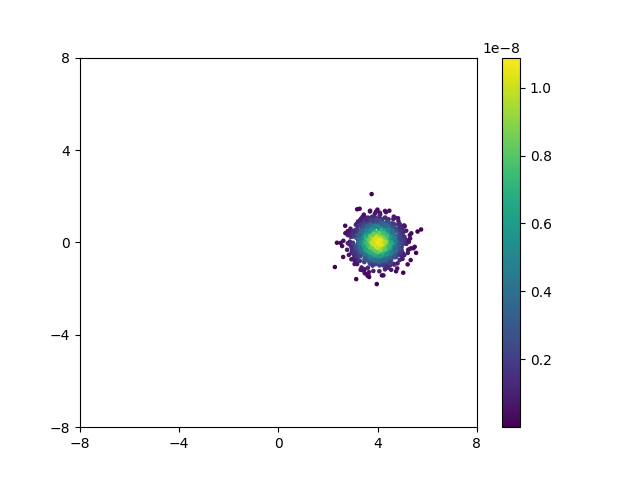}
\hspace{-5mm}
\includegraphics[width=2.4cm]{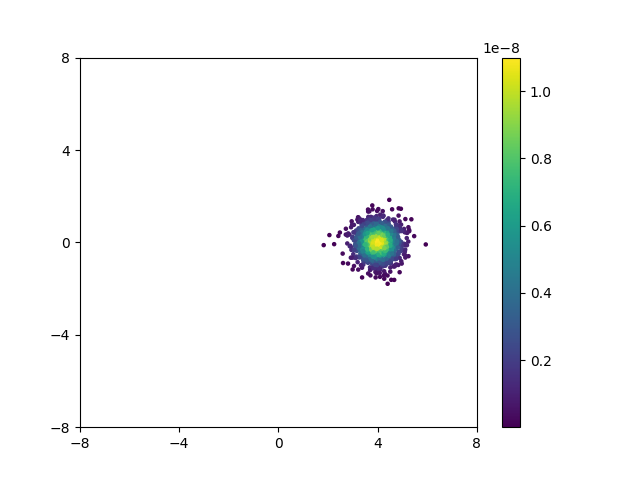}\\
\vspace{5pt}

\subfigure[$\rho_0(\boldsymbol{\boldsymbol{x}}_0)$]{\includegraphics[width=2.4cm]{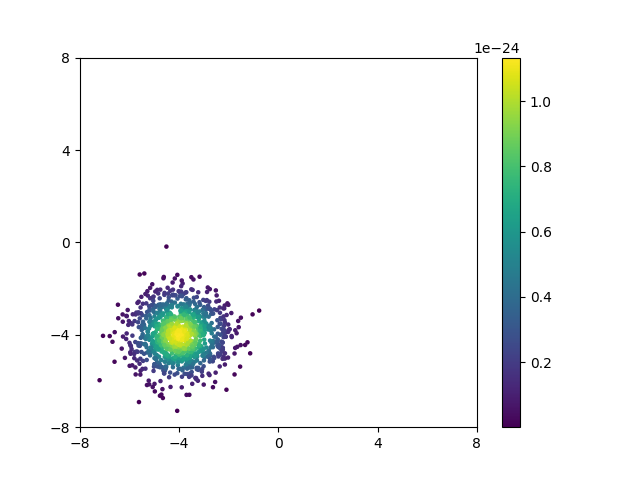}}
\hspace{-5mm}
\subfigure[$\tilde{\rho}_{1/4}(\tilde{\boldsymbol{\boldsymbol{x}}}_{1/4})$]{\includegraphics[width=2.4cm]{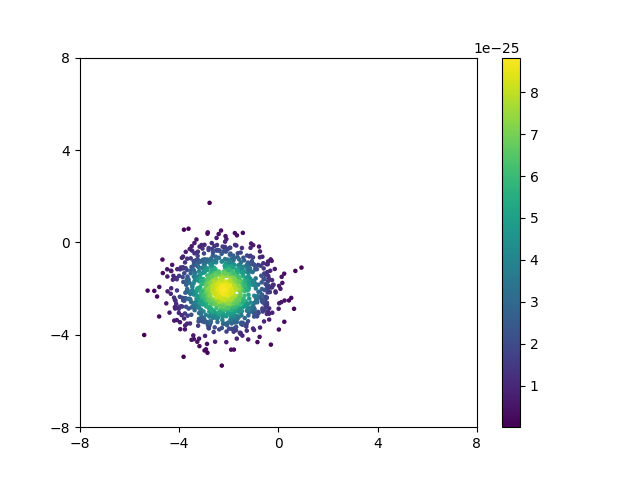}}
\hspace{-5mm}
\subfigure[$\tilde{\rho}_{1/2}(\tilde{\boldsymbol{\boldsymbol{x}}}_{1/2})$]{\includegraphics[width=2.4cm]{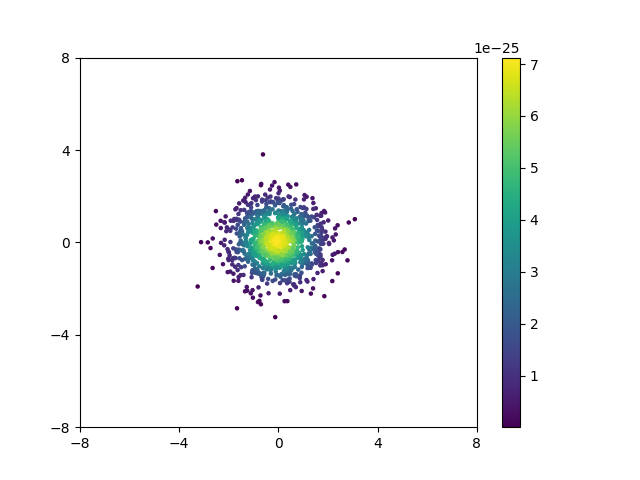}}
\hspace{-5mm}
\subfigure[$\tilde{\rho}_{3/4}(\tilde{\boldsymbol{\boldsymbol{x}}}_{3/4})$]{\includegraphics[width=2.4cm]{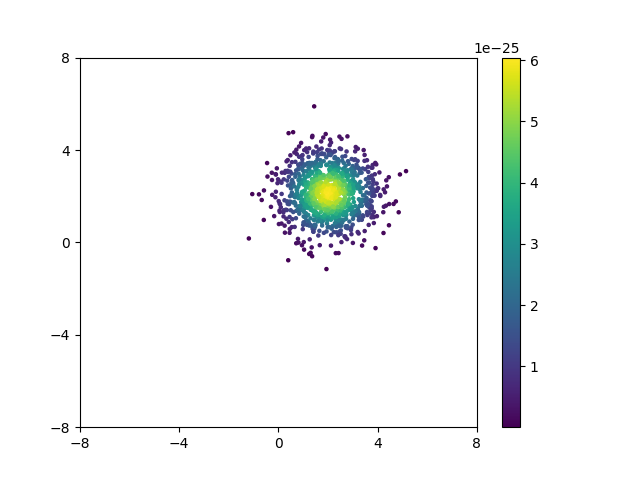}}
\hspace{-5mm}
\subfigure[$\tilde{\rho}_{1}(\tilde{\boldsymbol{\boldsymbol{x}}}_{1})$]{\includegraphics[width=2.4cm]{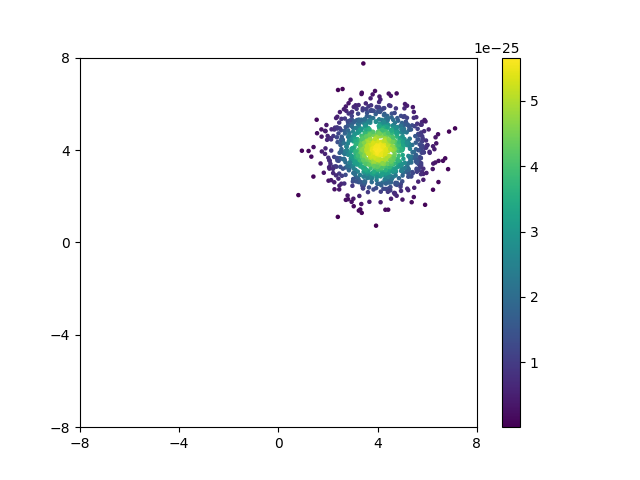}}
\hspace{-5mm}
\subfigure[$\rho_{1}(\boldsymbol{\boldsymbol{x}}_{1})$]{\includegraphics[width=2.4cm]{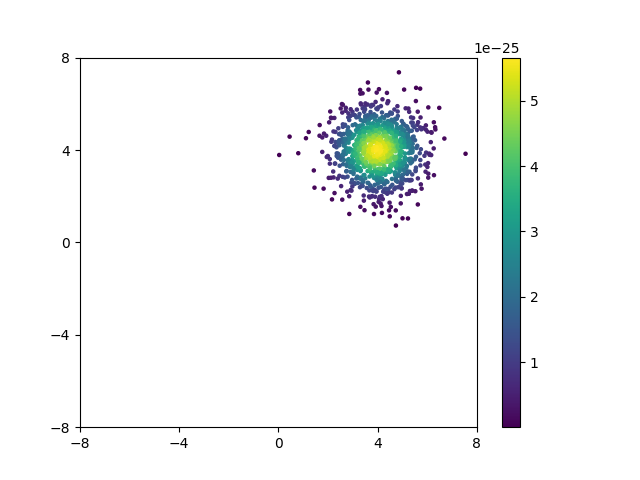}}

\caption{Illustration of these Gaussian problems at $d=60$.}
\label{compare_d60}
\end{center}
\end{figure}


\begin{figure}[H]
\begin{center}
\includegraphics[width=2.4cm]{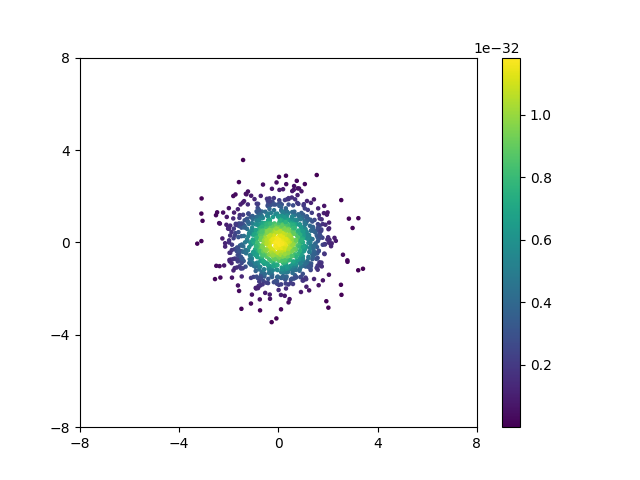}
\hspace{-5mm}
\includegraphics[width=2.4cm]{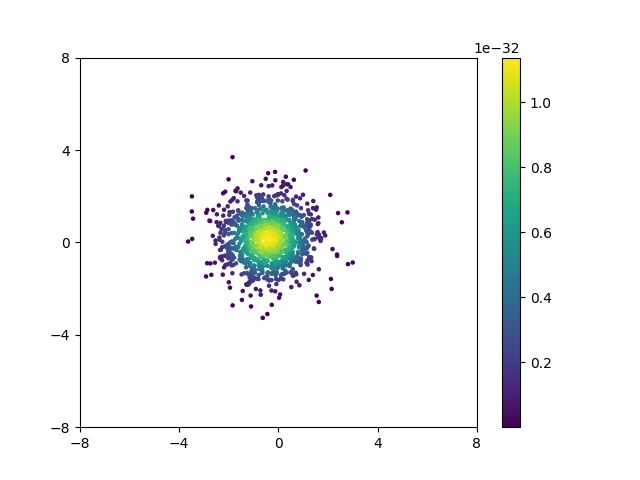}
\hspace{-5mm}
\includegraphics[width=2.4cm]{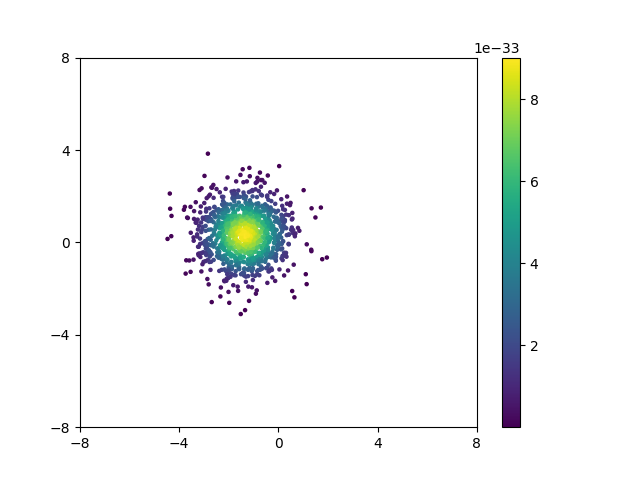}
\hspace{-5mm}
\includegraphics[width=2.4cm]{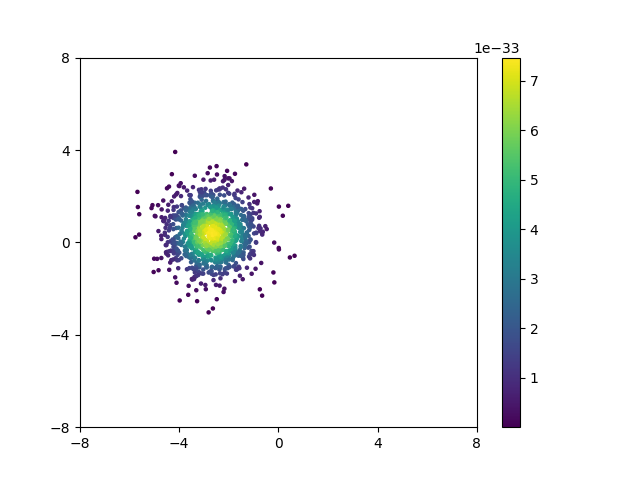}
\hspace{-5mm}
\includegraphics[width=2.4cm]{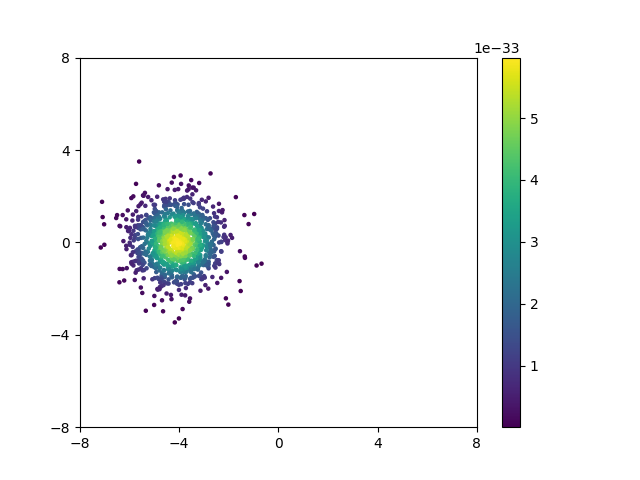}
\hspace{-5mm}
\includegraphics[width=2.4cm]{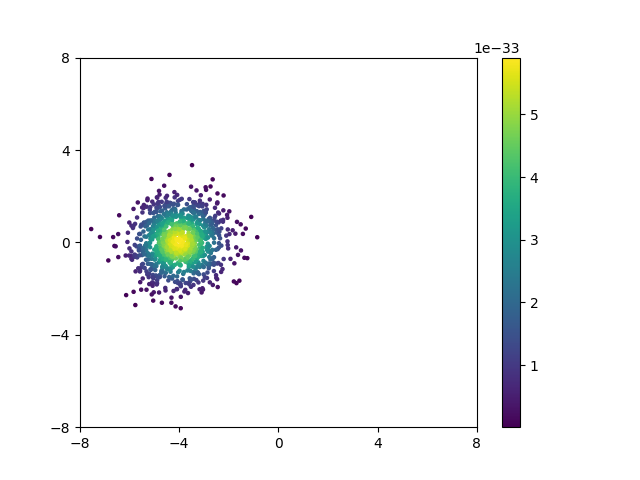}\\
\vspace{5pt}

\includegraphics[width=2.4cm]{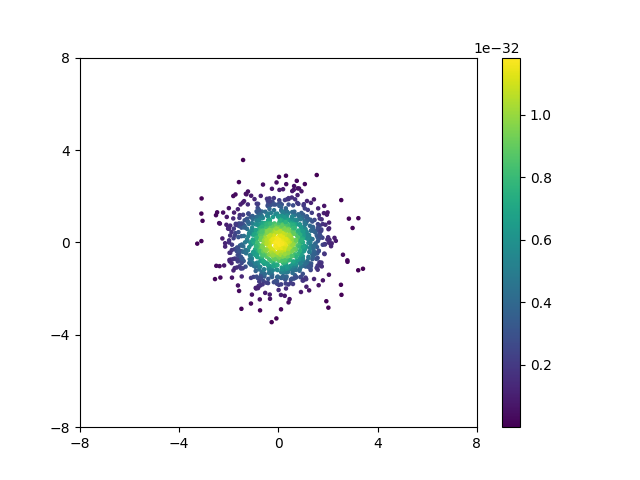}
\hspace{-5mm}
\includegraphics[width=2.4cm]{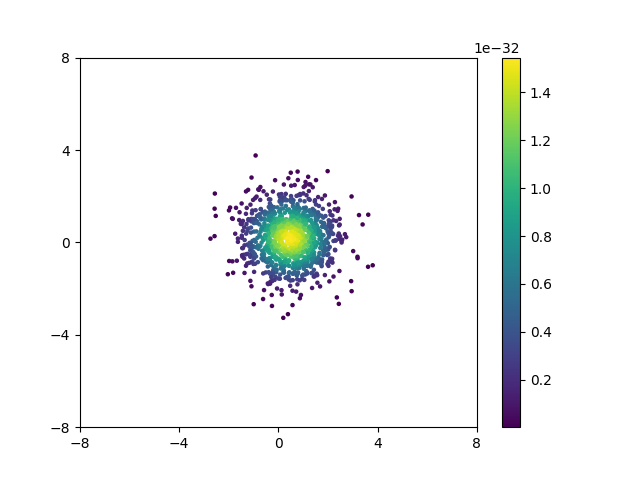}
\hspace{-5mm}
\includegraphics[width=2.4cm]{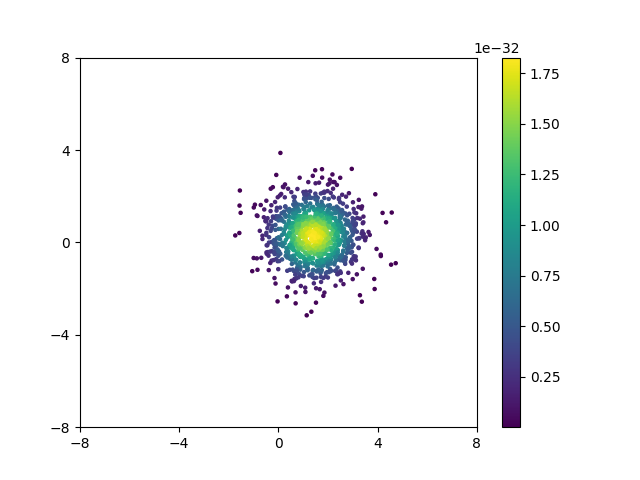}
\hspace{-5mm}
\includegraphics[width=2.4cm]{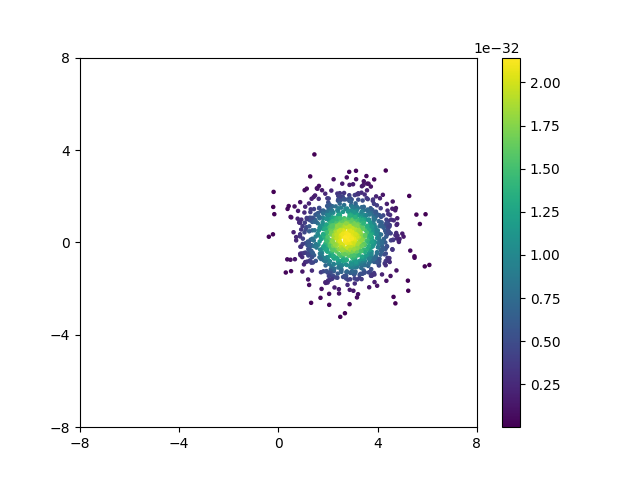}
\hspace{-5mm}
\includegraphics[width=2.4cm]{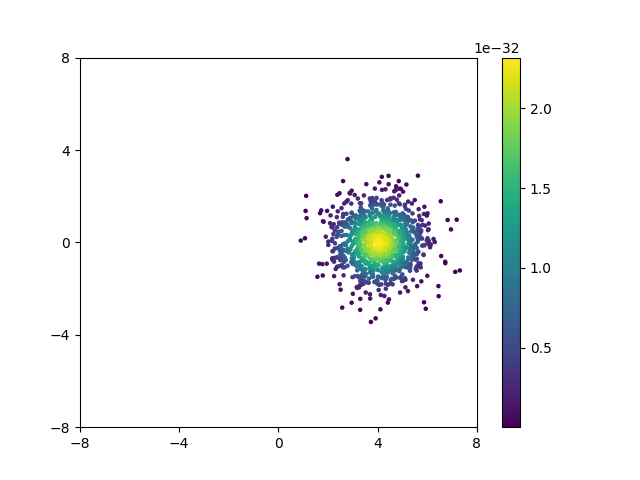}
\hspace{-5mm}
\includegraphics[width=2.4cm]{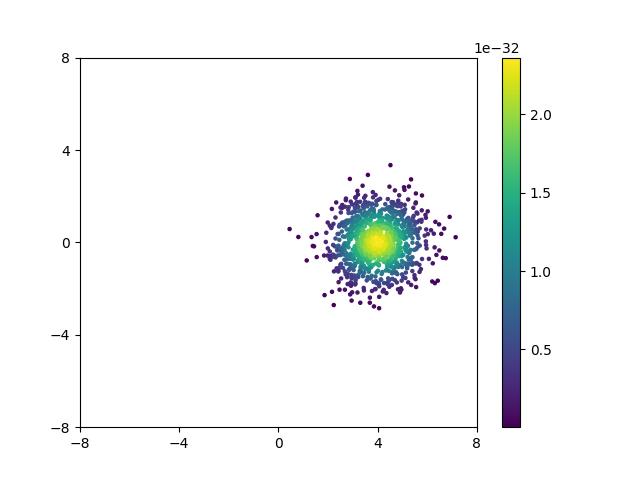}\\
\vspace{5pt}

\includegraphics[width=2.4cm]{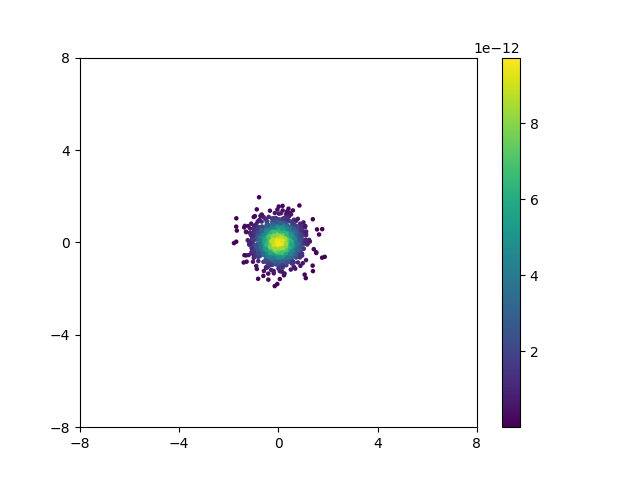}
\hspace{-5mm}
\includegraphics[width=2.4cm]{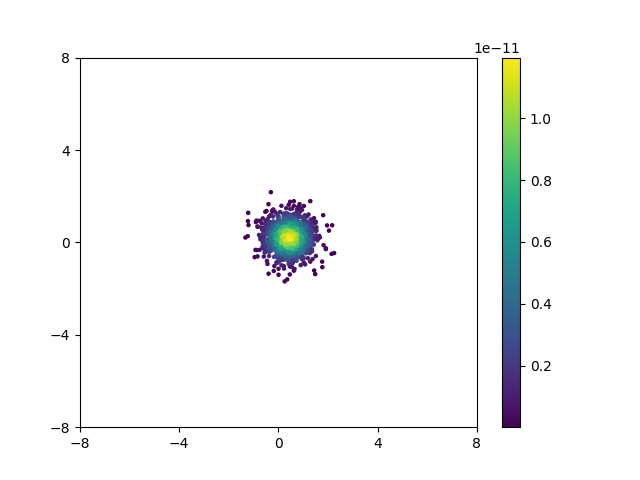}
\hspace{-5mm}
\includegraphics[width=2.4cm]{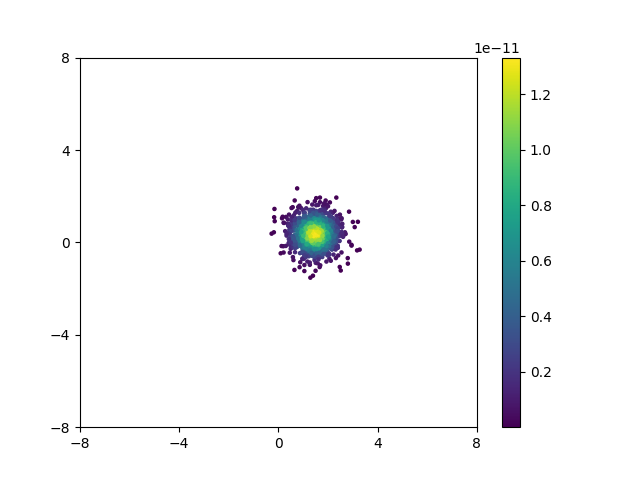}
\hspace{-5mm}
\includegraphics[width=2.4cm]{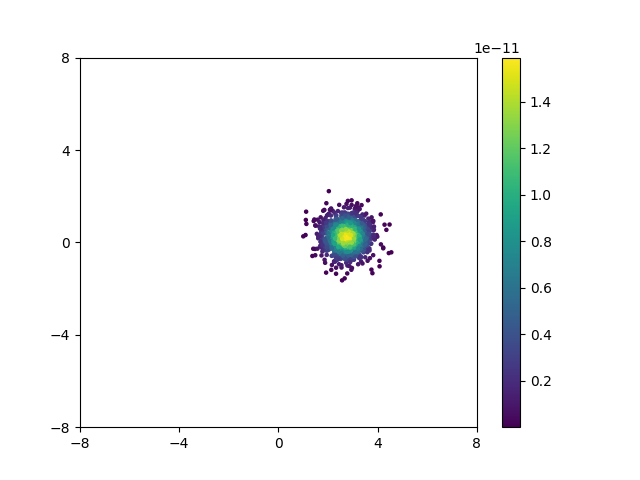}
\hspace{-5mm}
\includegraphics[width=2.4cm]{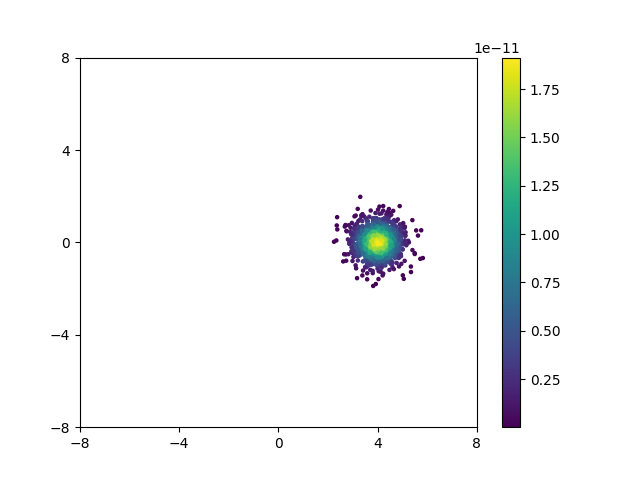}
\hspace{-5mm}
\includegraphics[width=2.4cm]{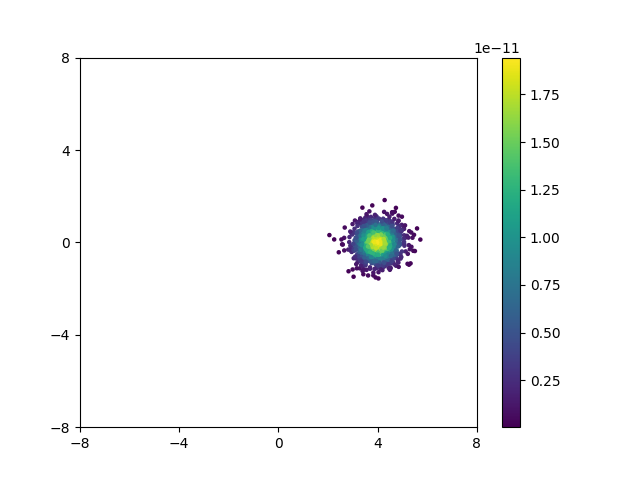}\\
\vspace{5pt}

\subfigure[$\rho_0(\boldsymbol{\boldsymbol{x}}_0)$]{\includegraphics[width=2.4cm]{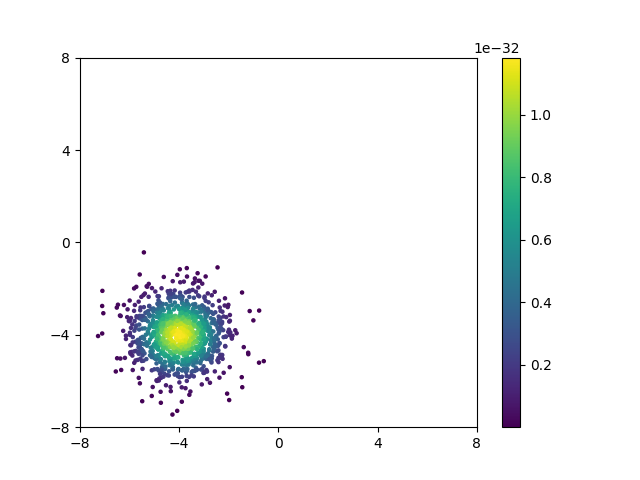}}
\hspace{-5mm}
\subfigure[$\tilde{\rho}_{1/4}(\tilde{\boldsymbol{\boldsymbol{x}}}_{1/4})$]{\includegraphics[width=2.4cm]{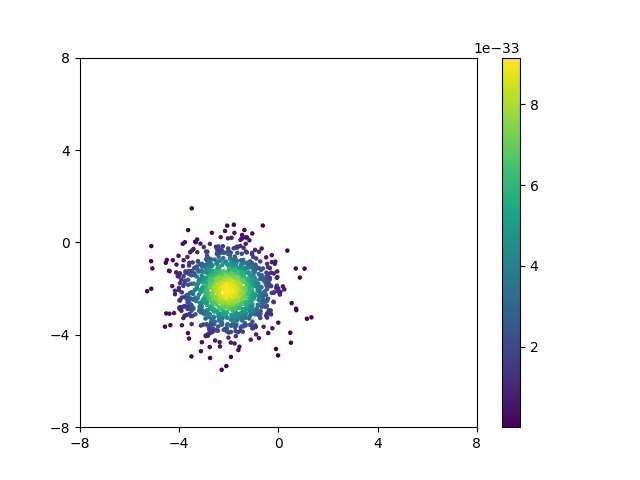}}
\hspace{-5mm}
\subfigure[$\tilde{\rho}_{1/2}(\tilde{\boldsymbol{\boldsymbol{x}}}_{1/2})$]{\includegraphics[width=2.4cm]{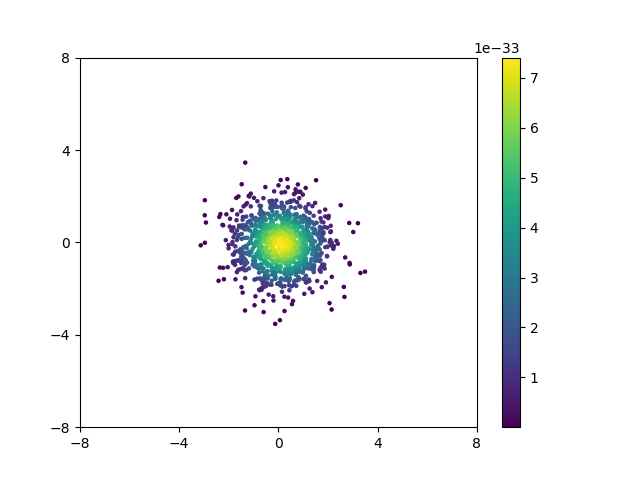}}
\hspace{-5mm}
\subfigure[$\tilde{\rho}_{3/4}(\tilde{\boldsymbol{\boldsymbol{x}}}_{3/4})$]{\includegraphics[width=2.4cm]{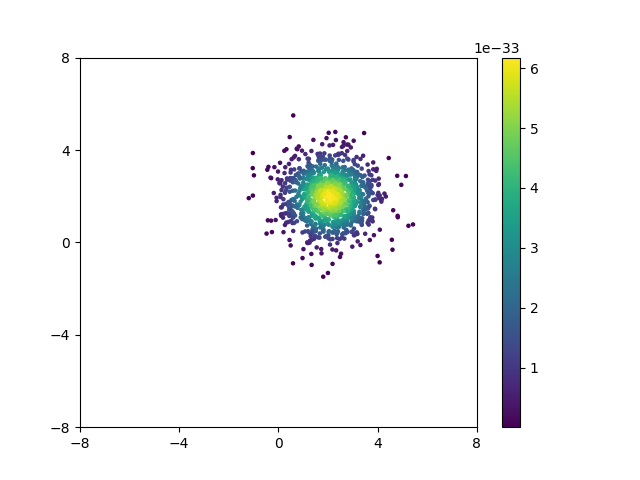}}
\hspace{-5mm}
\subfigure[$\tilde{\rho}_{1}(\tilde{\boldsymbol{\boldsymbol{x}}}_{1})$]{\includegraphics[width=2.4cm]{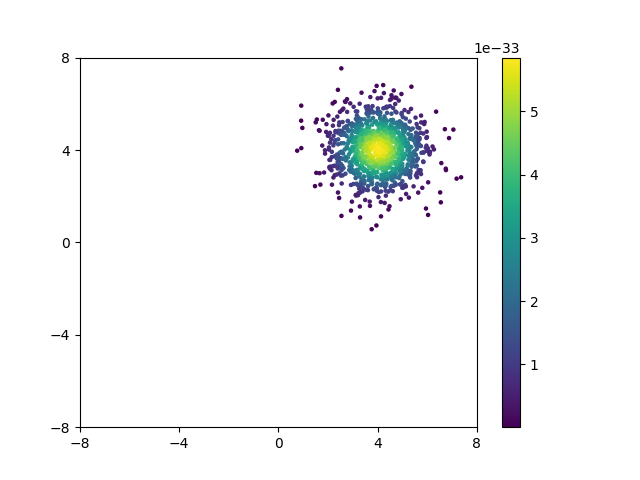}}
\hspace{-5mm}
\subfigure[$\rho_{1}(\boldsymbol{\boldsymbol{x}}_{1})$]{\includegraphics[width=2.4cm]{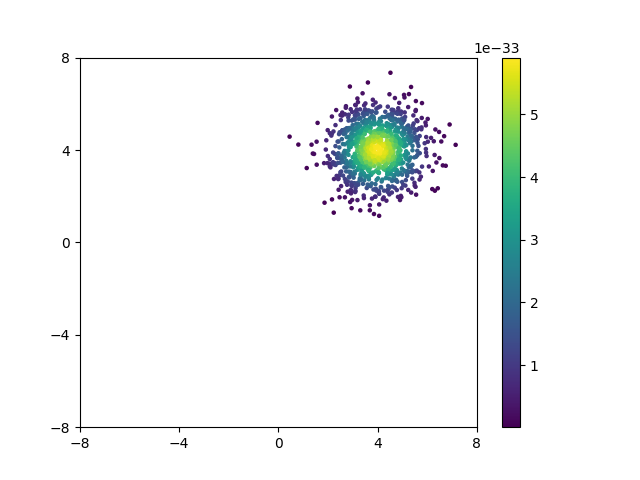}}

\caption{Illustration of these Gaussian problems at $d=80$.}
\label{compare_d80}
\end{center}
\end{figure}

\section{Gaussian examples for dimensions $d=300, 400$ on H800}
\label{appB}

\begin{figure}[H]
\begin{center}
\includegraphics[width=2.4cm]{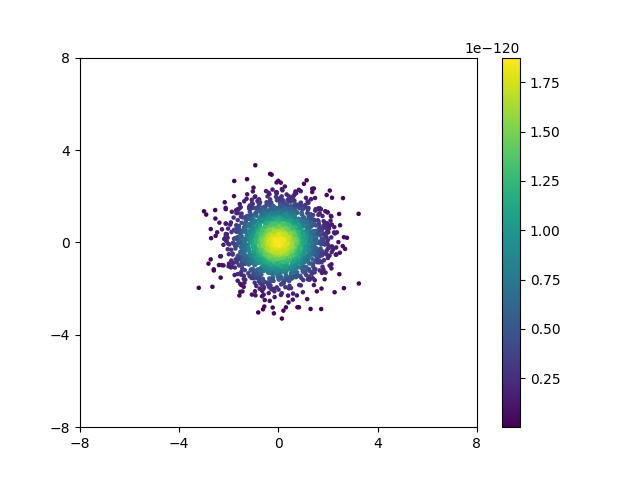}
\hspace{-5mm}
\includegraphics[width=2.4cm]{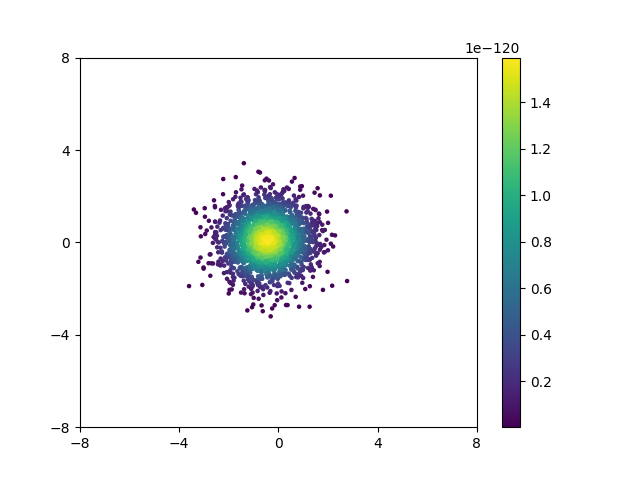}
\hspace{-5mm}
\includegraphics[width=2.4cm]{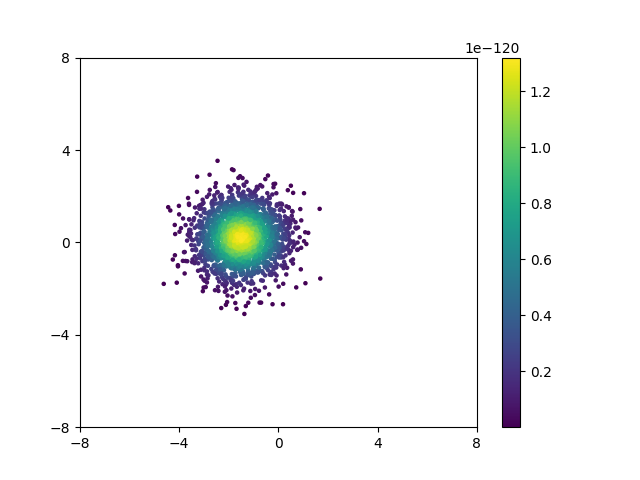}
\hspace{-5mm}
\includegraphics[width=2.4cm]{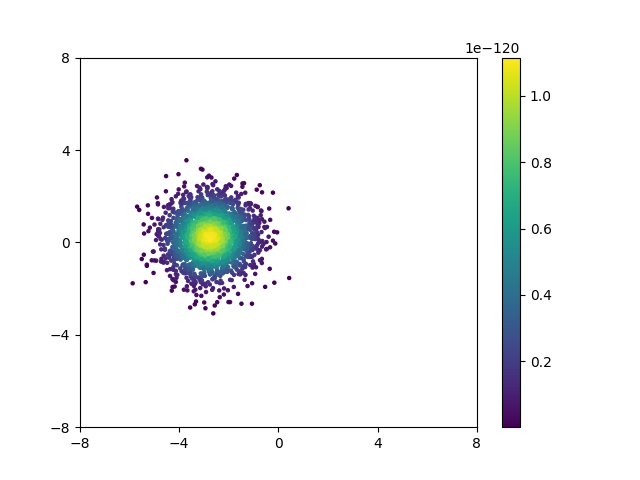}
\hspace{-5mm}
\includegraphics[width=2.4cm]{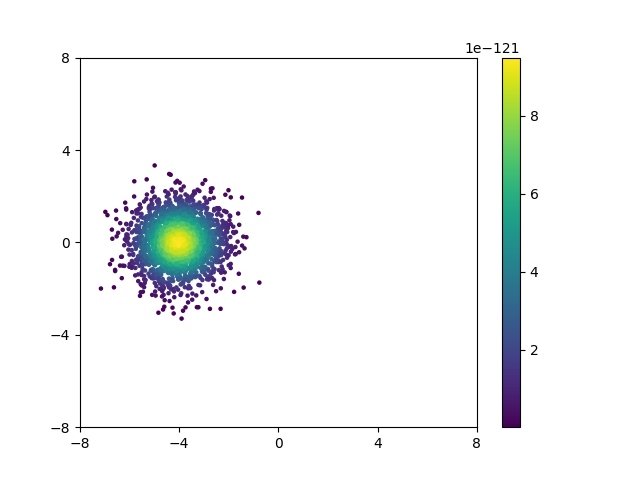}
\hspace{-5mm}
\includegraphics[width=2.4cm]{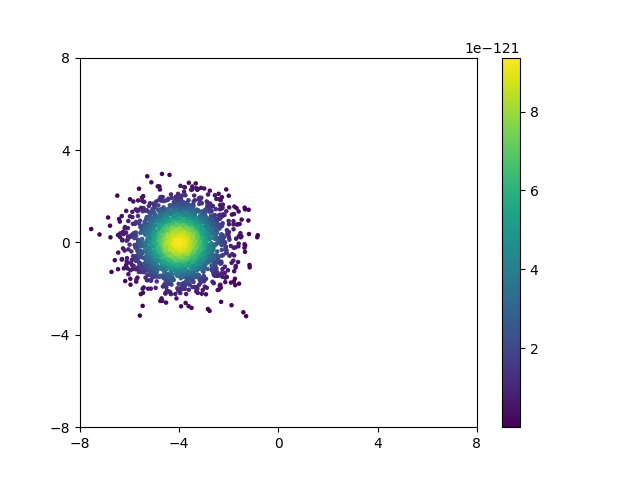}\\
\vspace{5pt}

\includegraphics[width=2.4cm]{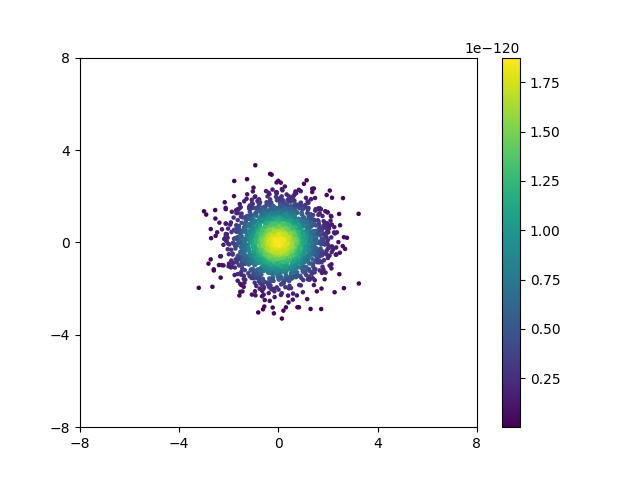}
\hspace{-5mm}
\includegraphics[width=2.4cm]{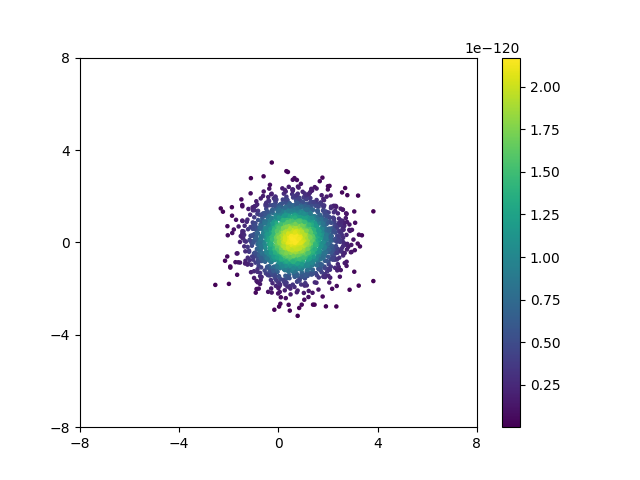}
\hspace{-5mm}
\includegraphics[width=2.4cm]{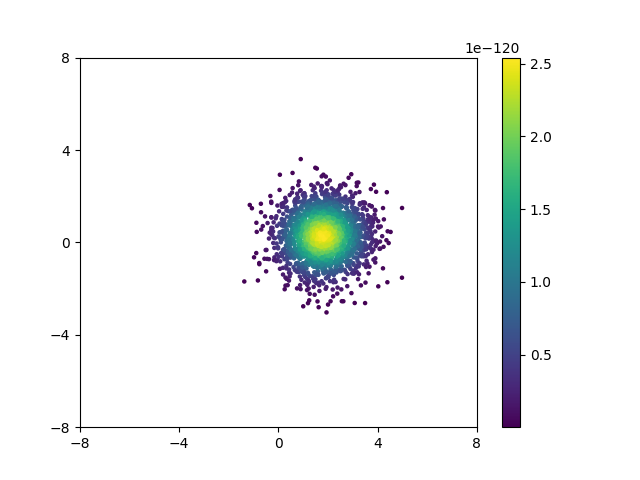}
\hspace{-5mm}
\includegraphics[width=2.4cm]{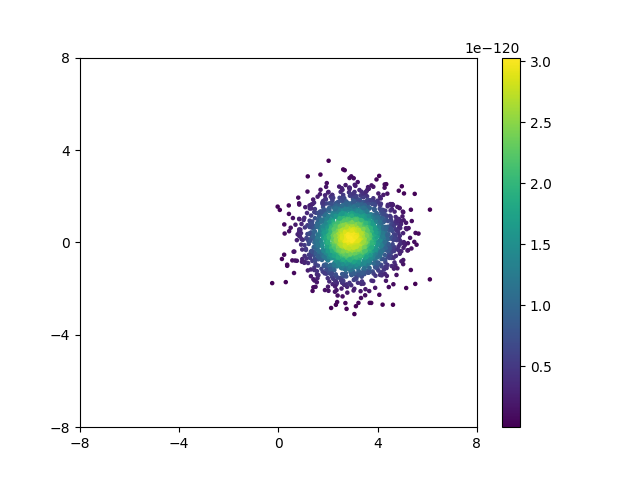}
\hspace{-5mm}
\includegraphics[width=2.4cm]{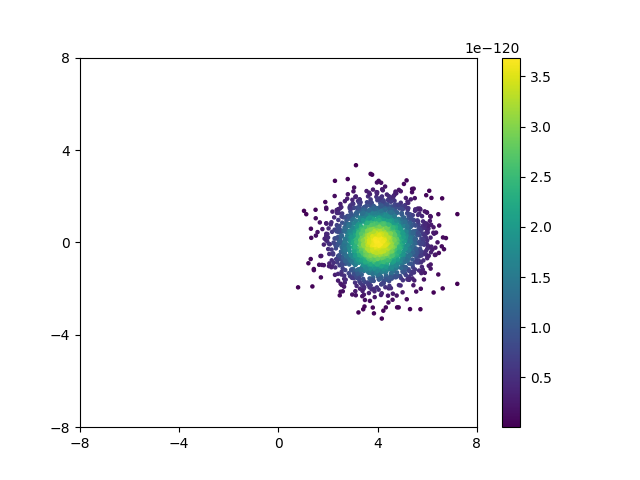}
\hspace{-5mm}
\includegraphics[width=2.4cm]{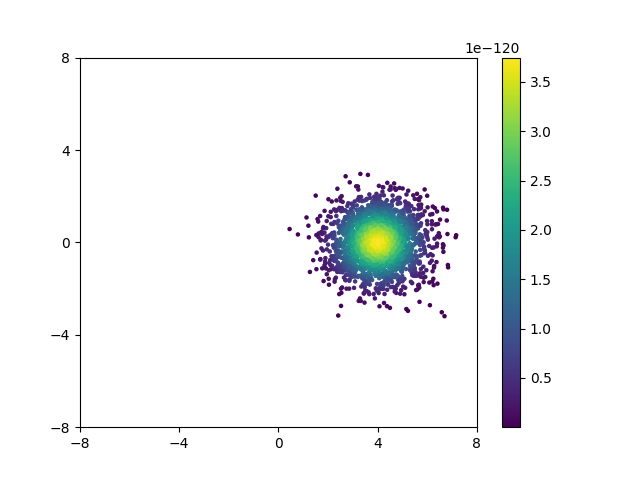}\\
\vspace{5pt}

\includegraphics[width=2.4cm]{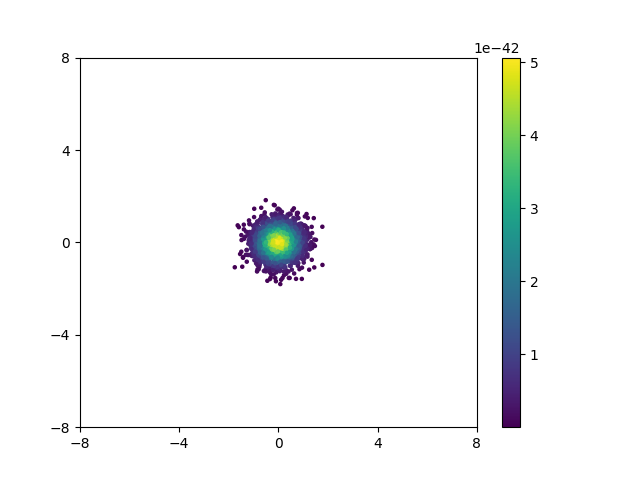}
\hspace{-5mm}
\includegraphics[width=2.4cm]{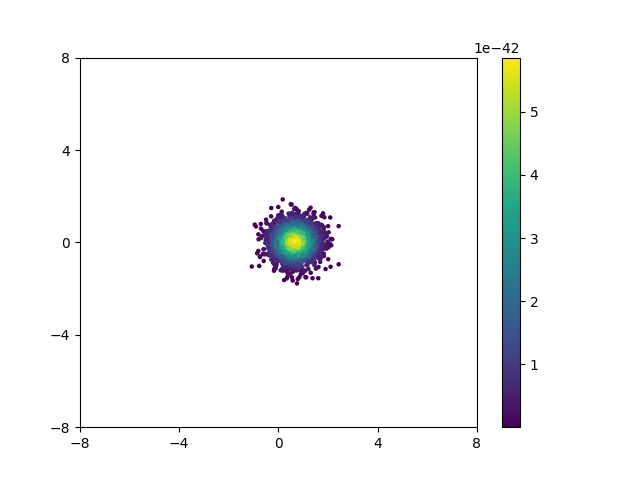}
\hspace{-5mm}
\includegraphics[width=2.4cm]{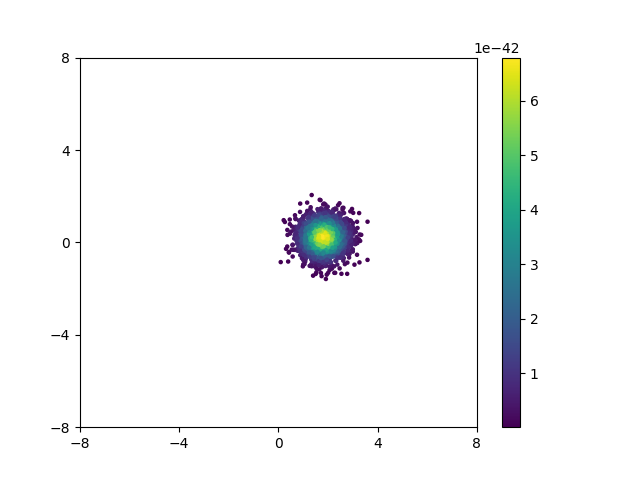}
\hspace{-5mm}
\includegraphics[width=2.4cm]{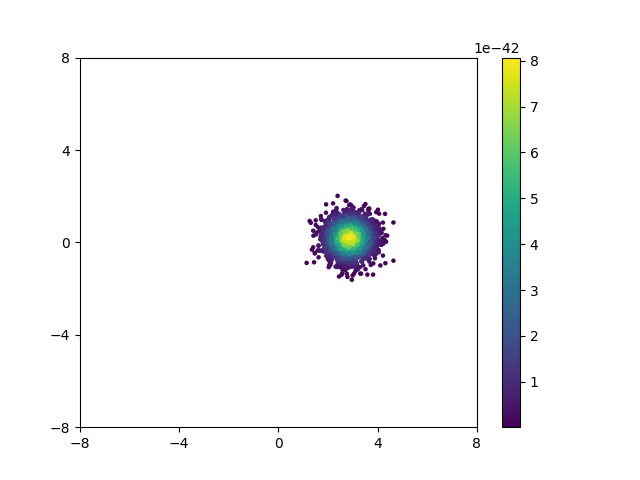}
\hspace{-5mm}
\includegraphics[width=2.4cm]{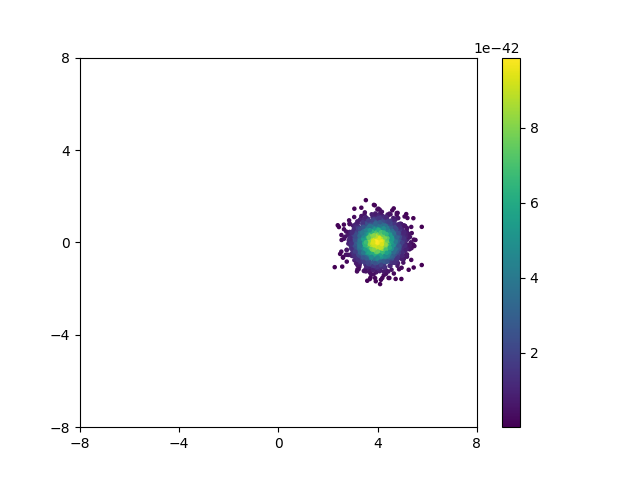}
\hspace{-5mm}
\includegraphics[width=2.4cm]{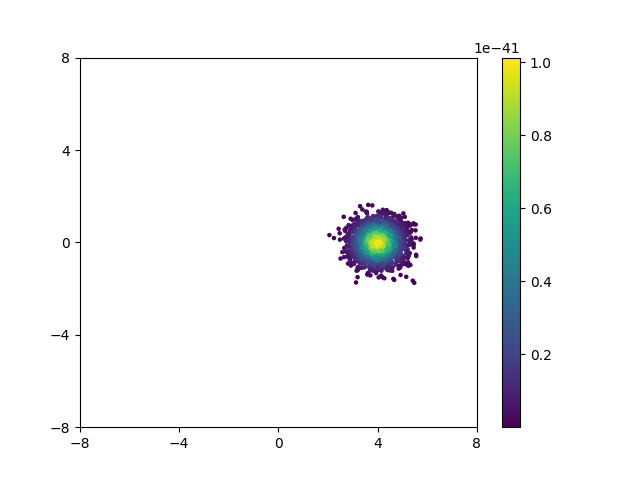}\\
\vspace{5pt}

\subfigure[$\rho_0(\boldsymbol{\boldsymbol{x}}_0)$]{\includegraphics[width=2.4cm]{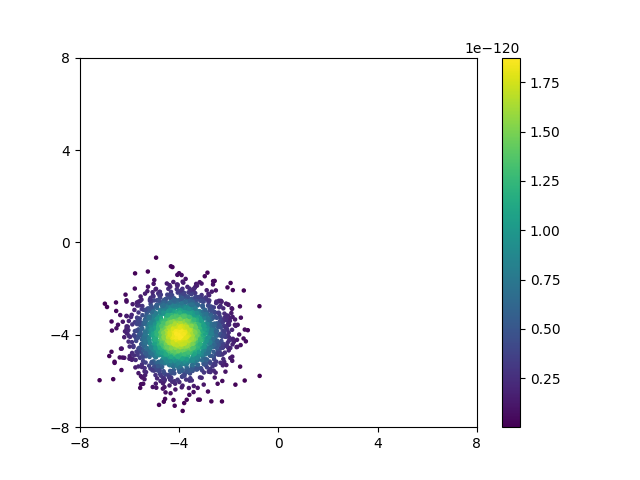}}
\hspace{-5mm}
\subfigure[$\tilde{\rho}_{1/4}(\tilde{\boldsymbol{\boldsymbol{x}}}_{1/4})$]{\includegraphics[width=2.4cm]{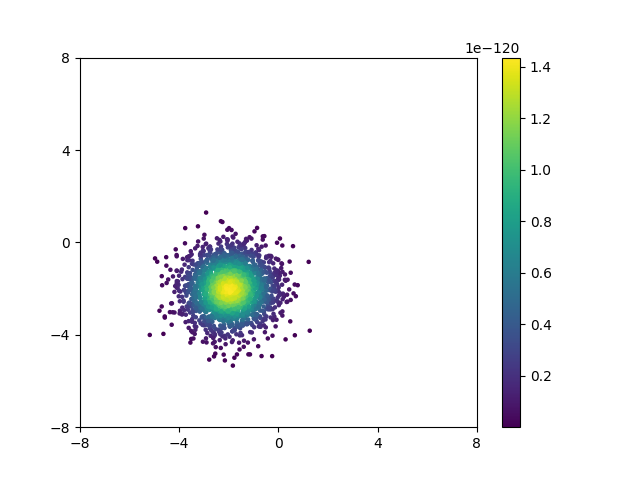}}
\hspace{-5mm}
\subfigure[$\tilde{\rho}_{1/2}(\tilde{\boldsymbol{\boldsymbol{x}}}_{1/2})$]{\includegraphics[width=2.4cm]{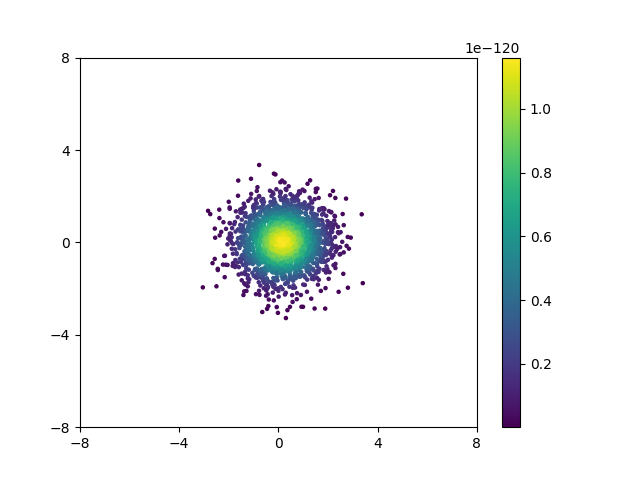}}
\hspace{-5mm}
\subfigure[$\tilde{\rho}_{3/4}(\tilde{\boldsymbol{\boldsymbol{x}}}_{3/4})$]{\includegraphics[width=2.4cm]{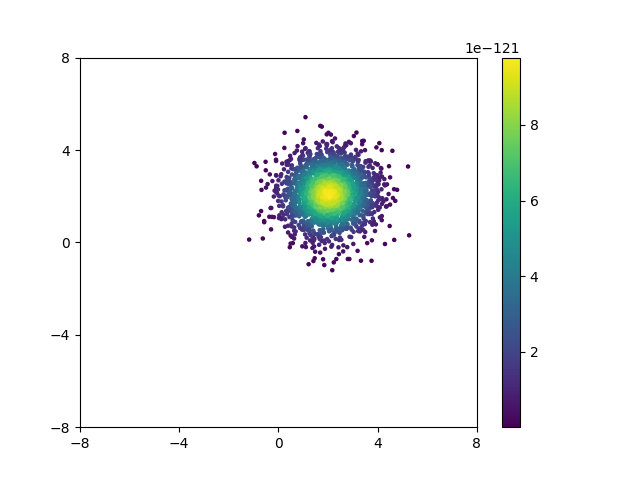}}
\hspace{-5mm}
\subfigure[$\tilde{\rho}_{1}(\tilde{\boldsymbol{\boldsymbol{x}}}_{1})$]{\includegraphics[width=2.4cm]{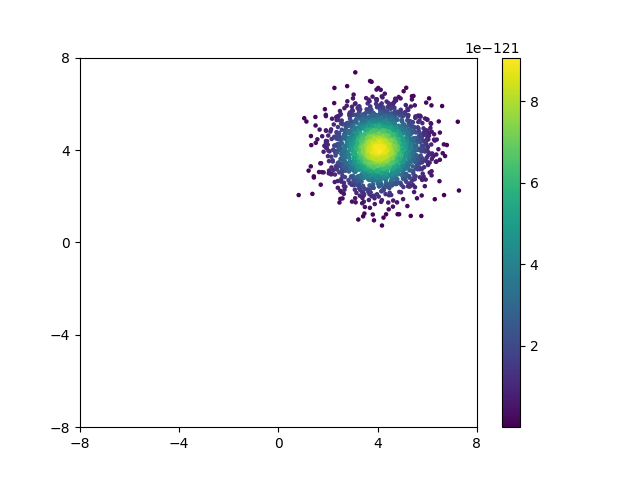}}
\hspace{-5mm}
\subfigure[$\rho_{1}(\boldsymbol{\boldsymbol{x}}_{1})$]{\includegraphics[width=2.4cm]{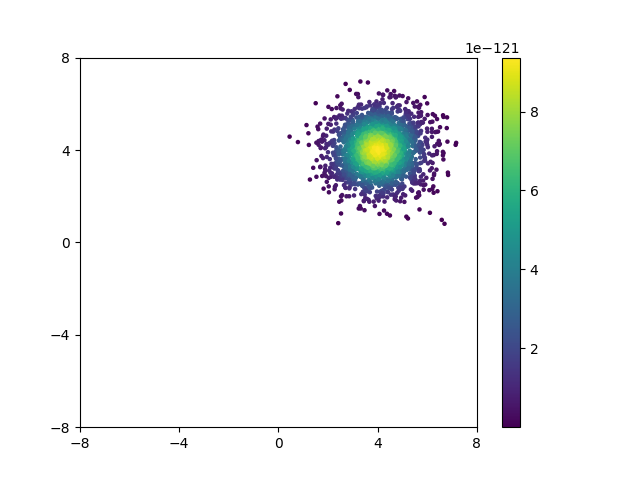}}

\caption{Illustration of these Gaussian problems at $d=300$.}
\label{compare_d300}
\end{center}
\end{figure}


\begin{figure}[H]
\begin{center}
\includegraphics[width=2.4cm]{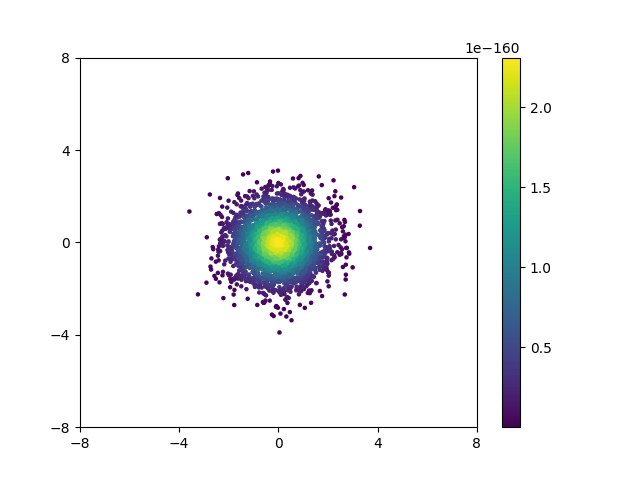}
\hspace{-5mm}
\includegraphics[width=2.4cm]{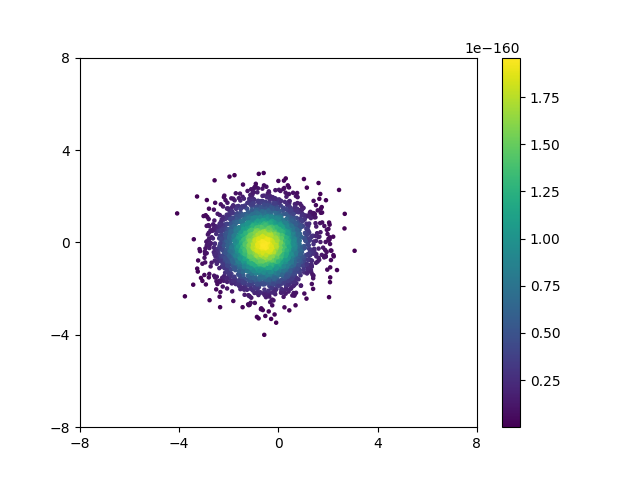}
\hspace{-5mm}
\includegraphics[width=2.4cm]{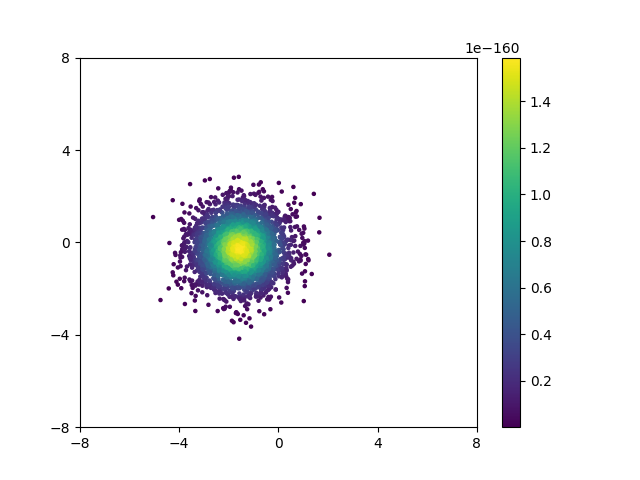}
\hspace{-5mm}
\includegraphics[width=2.4cm]{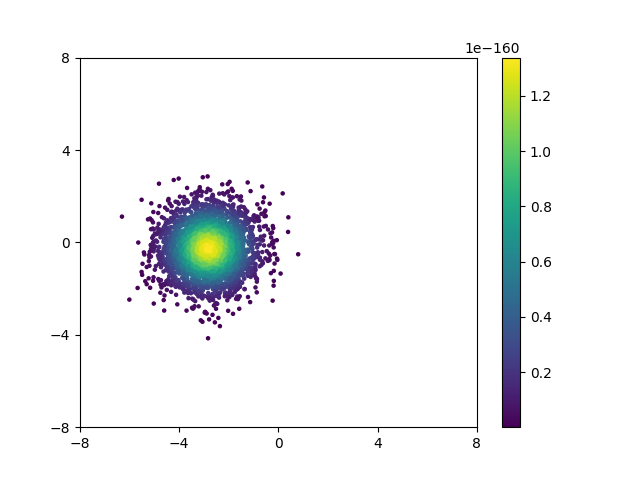}
\hspace{-5mm}
\includegraphics[width=2.4cm]{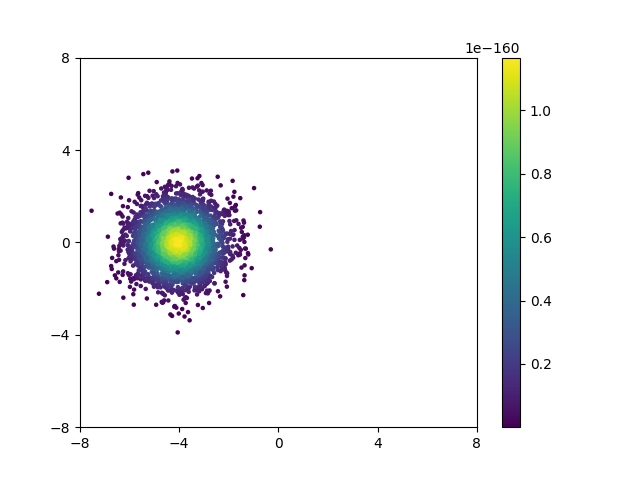}
\hspace{-5mm}
\includegraphics[width=2.4cm]{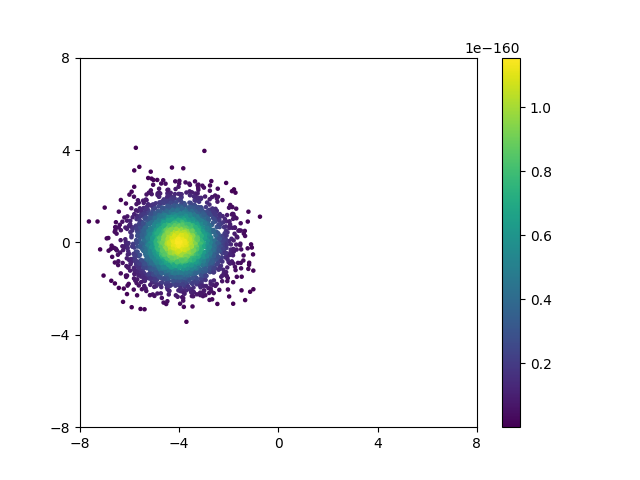}\\
\vspace{5pt}

\includegraphics[width=2.4cm]{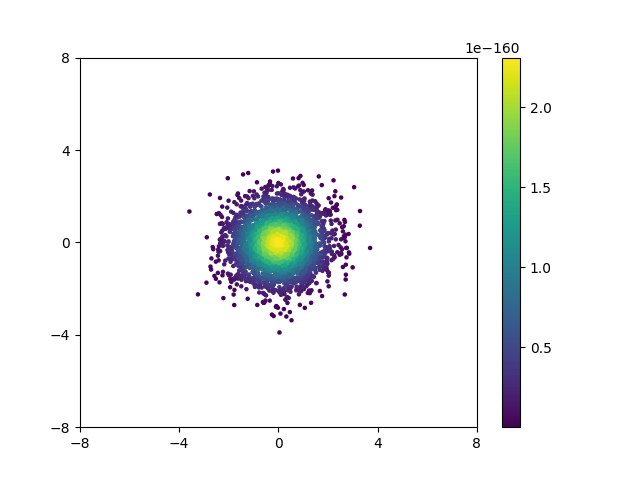}
\hspace{-5mm}
\includegraphics[width=2.4cm]{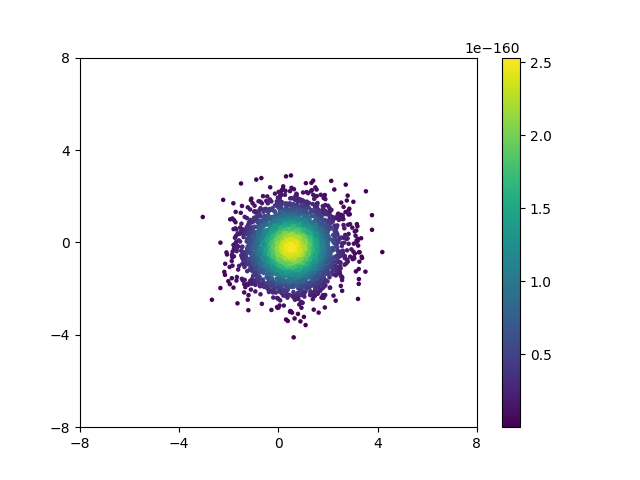}
\hspace{-5mm}
\includegraphics[width=2.4cm]{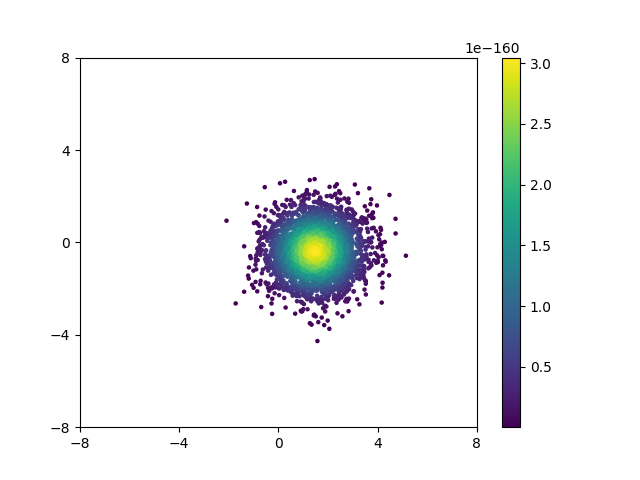}
\hspace{-5mm}
\includegraphics[width=2.4cm]{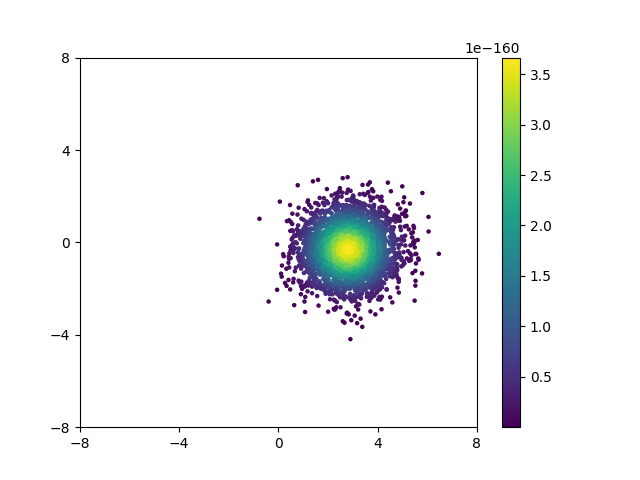}
\hspace{-5mm}
\includegraphics[width=2.4cm]{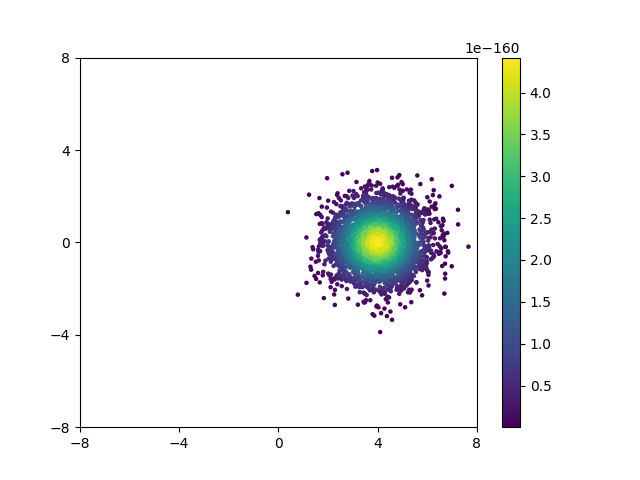}
\hspace{-5mm}
\includegraphics[width=2.4cm]{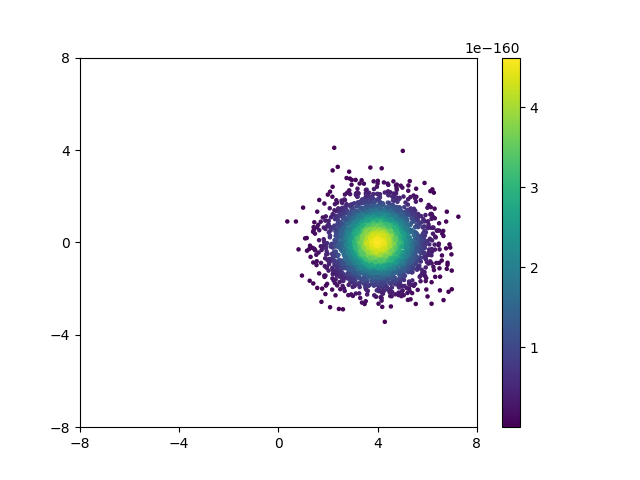}\\
\vspace{5pt}

\includegraphics[width=2.4cm]{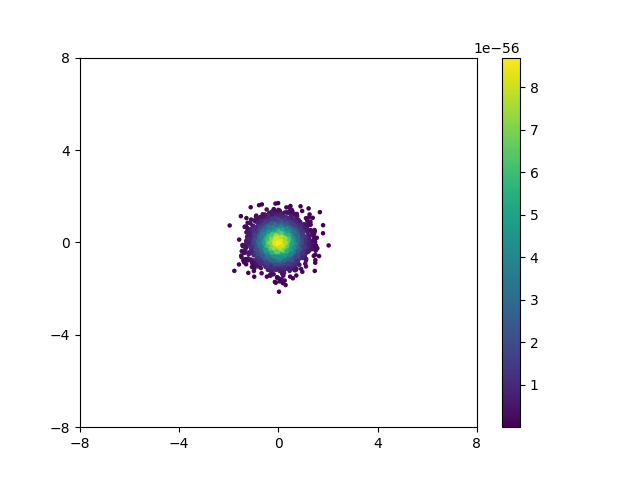}
\hspace{-5mm}
\includegraphics[width=2.4cm]{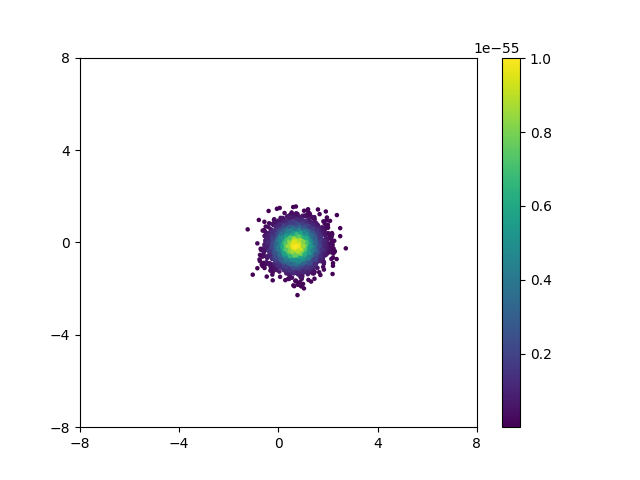}
\hspace{-5mm}
\includegraphics[width=2.4cm]{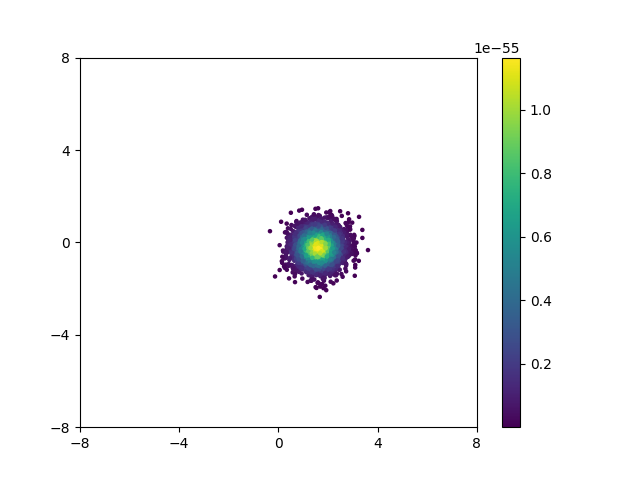}
\hspace{-5mm}
\includegraphics[width=2.4cm]{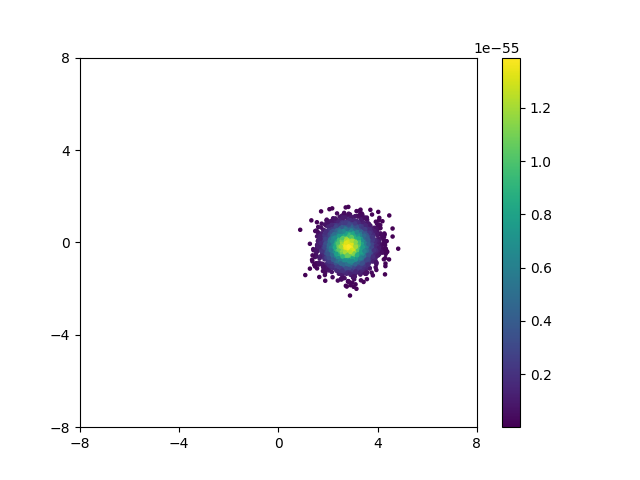}
\hspace{-5mm}
\includegraphics[width=2.4cm]{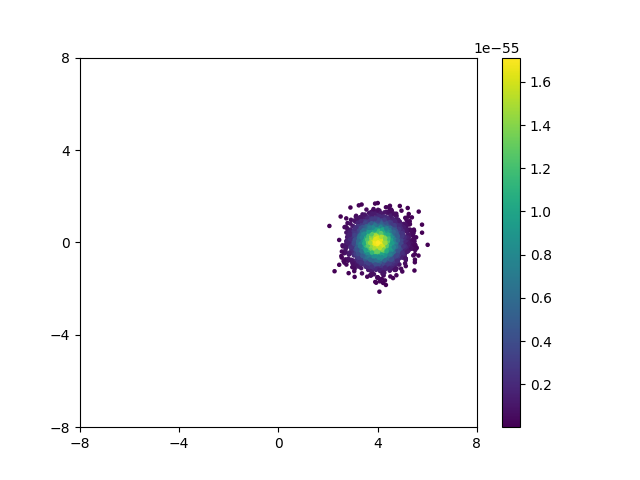}
\hspace{-5mm}
\includegraphics[width=2.4cm]{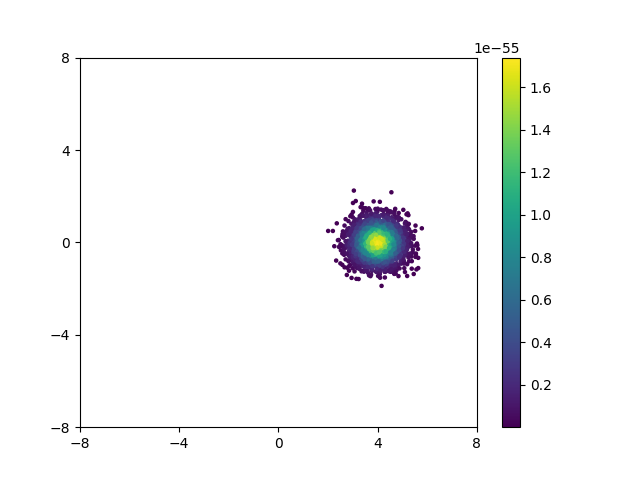}\\
\vspace{5pt}

\subfigure[$\rho_0(\boldsymbol{\boldsymbol{x}}_0)$]{\includegraphics[width=2.4cm]{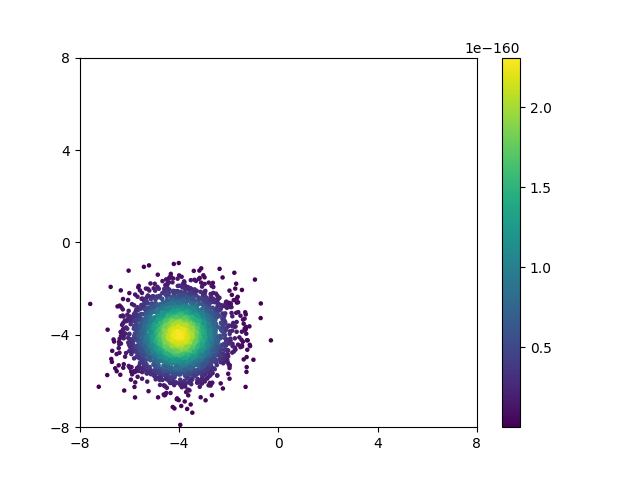}}
\hspace{-5mm}
\subfigure[$\tilde{\rho}_{1/4}(\tilde{\boldsymbol{\boldsymbol{x}}}_{1/4})$]{\includegraphics[width=2.4cm]{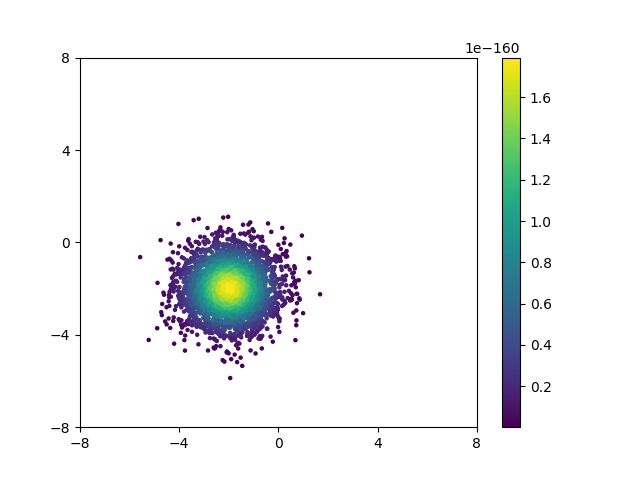}}
\hspace{-5mm}
\subfigure[$\tilde{\rho}_{1/2}(\tilde{\boldsymbol{\boldsymbol{x}}}_{1/2})$]{\includegraphics[width=2.4cm]{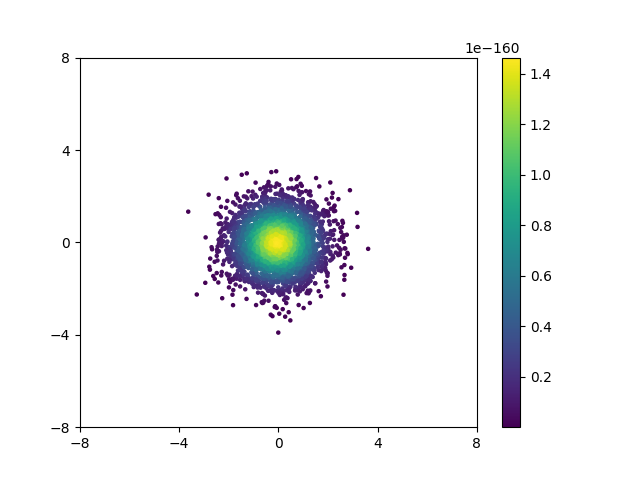}}
\hspace{-5mm}
\subfigure[$\tilde{\rho}_{3/4}(\tilde{\boldsymbol{\boldsymbol{x}}}_{3/4})$]{\includegraphics[width=2.4cm]{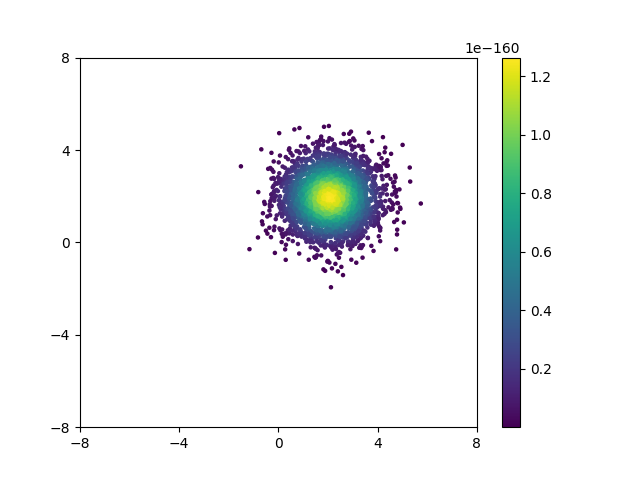}}
\hspace{-5mm}
\subfigure[$\tilde{\rho}_{1}(\tilde{\boldsymbol{\boldsymbol{x}}}_{1})$]{\includegraphics[width=2.4cm]{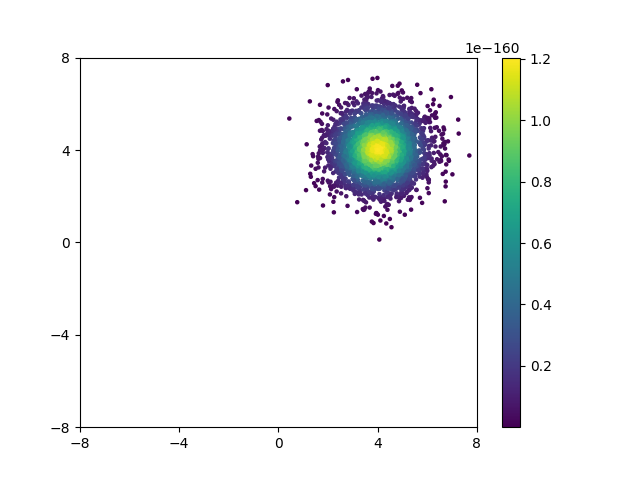}}
\hspace{-5mm}
\subfigure[$\rho_{1}(\boldsymbol{\boldsymbol{x}}_{1})$]{\includegraphics[width=2.4cm]{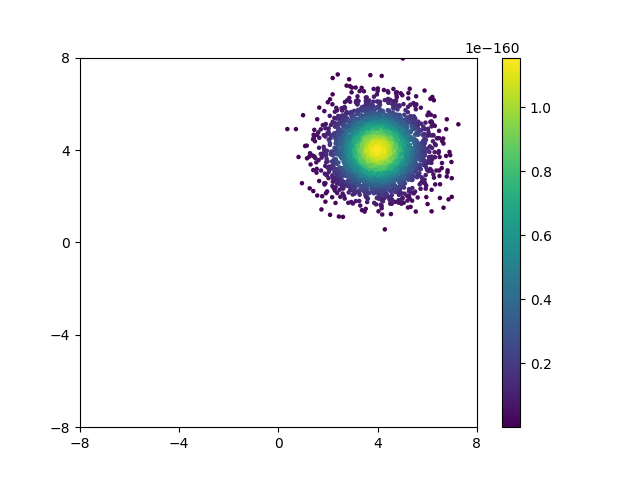}}

\caption{Illustration of these Gaussian problems at $d=400$.}
\label{compare_d400}
\end{center}
\end{figure}

\end{document}